\documentclass{amsart}
\usepackage{amssymb}

\numberwithin{equation}{section}

\newtheorem{thm}{Theorem}[section]
\newtheorem{defn}[thm]{Definition}
\newtheorem{prop}[thm]{Proposition}
\newtheorem{cor}[thm]{Corollary}
\newtheorem{lem}[thm]{Lemma}
\newtheorem{rema}[thm]{Remark}
\newcommand{\Z}{\mathbb{Z}_+}

\begin{document}

\title[$N=2$ vertex superalgebras]{On axiomatic aspects of $N=2$ vertex superalgebras with odd formal variables, and deformations of $N=1$ vertex superalgebras}

\author{Katrina Barron}
\address{Department of Mathematics, University of Notre Dame,
Notre Dame, IN 46556}
\email{kbarron@nd.edu}

\subjclass[2000]{17B68, 17B69, 17B81, 81R10, 81T40, 81T60}

\date{October 30, 2007}

\keywords{Superconformal field theory, Neveu-Schwarz Lie superalgebra, vertex operator superalgebra}

\begin{abstract}
The notion of ``$N = 2$ vertex superalgebra with two odd formal variables" is presented, the main axiom being a Jacobi identity with odd formal variables in which an $N=2$ superconformal shift is incorporated into the usual Jacobi identity for a vertex superalgebra.  It is shown that as a consequence of these axioms, the $N=2$ vertex superalgebra is naturally a representation of the Lie algebra isomorphic to the three-dimensional algebra of superderivations with basis consisting of the usual conformal operator and the two $N=2$ superconformal operators.  In addition, this superconformal shift in the Jacobi identity dictates the form of the odd formal variable components of the vertex operators, and allows one to easily derive the useful formulas in the theory.  The notion of $N=2$ Neveu-Schwarz vertex operator superalgebra with two odd formal variables is introduced, and consequences of this notion are derived.  Various other formulations of the notion of $N=2$ (Neveu-Schwarz) vertex (operator) superalgebra appearing in the mathematics and physics literature are discussed, and several mistakes in the literature are noted and corrected.  The notion of  ``$N=2$ (Neveu-Schwarz) vertex (operator) superalgebra with one odd formal variable" is formulated.  It is shown that this formulation naturally arises {}from alternate notions of $N=1$ superconformality and the continuous deformation of an $N=1$ (Neveu-Schwarz) vertex (operator) superalgebra with one odd formal variable.  This notion is formulated to reflect the underlying $N=1$ superanalytic geometry, and it is shown that the equivalence of the notions of $N=2$ (Neveu-Schwarz) vertex (operator) superalgebra with one and with two odd formal variables reflects the equivalence of $N=2$ superconformal and $N=1$ superanalytic geometry.   In particular we prove that the group of formal $N=2$ superconformal functions vanishing at zero and invertible in a neighborhood of zero is isomorphic to a certain subgroup of $N=1$ superanalytic functions vanishing at zero and invertible in a neighborhood of zero.
\end{abstract}

\maketitle

\tableofcontents

\section{Introduction and preliminaries}

\subsection{Introduction}
In this paper, we study several equivalent formulations of the notion of ``$N = 2$ (Neveu-Schwarz) vertex (operator) superalgebra" motivated principally by the two different but equivalent supergeometric settings underlying genus-zero, holomorphic $N=2$ superconformal field theory, and continuously deformed genus-zero, holomorphic $N=1$ superconformal field theory, respectively.

(Super)conformal field theory is an attempt at developing a physical theory that combines all fundamental interactions of particles (cf. \cite{BPZ}, \cite{Fd}, \cite{FS}, \cite{V}, \cite{S}).    In supersymmetric theories, symmetries between the two fundamental types of particles -- bosons and fermions -- are assumed.  $N=n$ superconformal field theory for $n$ a nonnegative integer explores the theory of conformal fields if $n$ boson-fermion symmetries are present  (cf. \cite{DPZ}, \cite{FMS}, \cite{Wa}, \cite{Ge}).  Algebraically the interactions of superparticles are described by vertex operators or more precisely by vertex (operator) superalgebras, also called ``chiral algebras" in the physics literature \cite{Bo}, \cite{FLM}, \cite{KT},  \cite{DPZ}, \cite{YZ}, \cite{LVW}, \cite{MSS}, \cite{FFR}, \cite{DL}.  If the theory is $N=1$ (resp., $N=2$) supersymmetric then in particular, the vertex operator superalgebra is a representation of the Lie superalgebra of infinitesimal superconformal symmetries which is called the $N=1$ (resp., $N=2$) Neveu-Schwarz algebra.  The natural geometric setting for such theories is supergeometry in which one (resp., two) odd variables are included, \cite{Fd}, \cite{D}, \cite{B-thesis}, \cite{B-memoir}, \cite{B-moduli}.

There are many notions of ``$N=2$ vertex (operator) superalgebra"  in the literature (cf.  \cite{Kac1997}, \cite{Adamovic1999}, \cite{Bo1}, \cite{HM}, \cite{Bo2}, \cite{HK}, \cite{He}) inspired by the notion in physics of the algebra of chiral $N=2$ superfields.  When odd variables are not included in the definition, these notions are all more or less equivalent.   Our formulations of the notions of ``$N=2$ (Neveu-Schwarz) vertex (operator) superalgebra with two odd formal variables" is equivalent to or similar in spirit to, some of the notions found in the literature such as in \cite{HK} and \cite{HM}.   However, there are typos in these works that make some of the formulations mathematically inconsistent and confusing to a novice reader.    

The first part of this paper is devoted to presenting what we feel to be the most natural formulation of the notion of $N=2$ vertex superalgebra with two odd formal variables.  We then show how our notion relates to others, and show how our formulation gives as a consequence the properties that have been troublesome to formulate in some of the literature.   Our notion of $N=2$ vertex superalgebra with two odd formal variables involves the usual notion of vertex superalgebra formulated using the Jacobi identity as in \cite{LL}, but extended as follows:  first odd variable components are added to the vertex operators; and second, we require that these operators satisfy the usual Jacobi identity for vertex superalgebras extended so that in place of the usual shift in the $\delta$-functions appearing in the Jacobi identity, we substitute  an $N=2$ superconformal shift.   That is instead of, for instance $x_1-x_2$ appearing, we substitute $x_1 - x_2 - \varphi_1^+ \varphi_2^- - \varphi_1^- \varphi_2^+$ where $x_1$ and $x_2$ are ``even" (or commuting) formal variables, and $\varphi^\pm_1$ and $\varphi_2^\pm$  are  ``odd"  (or anticommuting) formal variables.   In particular, we show that as a {\it consequence} of this simple two-step extension of the usual notion of vertex algebra, then the underlying vector space is naturally endowed with the structure of a module for a certain subalgebra of the orthogonal-symplectic Lie superalgebra, $\mathfrak{osp}(2|2)$.  We show that the three representative elements for the usual basis elements of this Lie superalgebra, when bracketed with a vertex operator, act as the conformal derivation and the two $N=2$ superconformal superderivations, respectively.   Thus we get as a result of the definition, the presence of endomorphisms on the underlying vector space that correspond to the usual conformal operator and the two $N=2$ superconformal operator.  This is in contrast to other formulations of the notion of $N=2$ vertex superalgebra that axiomatize the existence of such endomorphisms and axiomatize certain properties the vertex operators must satisfy in relation to these endomorphisms (cf. \cite{HK}).  In our formulation, all we are assuming is that the notion of a shift from one point in superconformal space to another, should be of the form of an $N=2$ superconformal shift so as to reflect some underlying $N=2$ superconformal geometry.   Furthermore, we show the form of the odd variable components of the vertex operators are completely dictated by this one simple modification of the shift appearing in the usual Jacobi identity for vertex superalgebras, and we show that concise formulas for commutators and iterates of vertex operators with odd formal variables, skew supersymmetry, and commutativity and associativity properties all follow easily {}from the Jacobi identity formulation.   

In addition, we formulate the notion of $N=2$ Neveu-Schwarz vertex operator superalgebra with two odd formal variables.  This is an $N=2$ vertex superalgebra that includes, among other things, a representation of the full $N=2$ Neveu-Schwarz algebra (which is an infinite extension of $\mathfrak{osp}(2|2)$) and an associated half-integer grading.  This formulation is in the spirit of \cite{HM} and previous work by the author in the $N=1$ setting \cite{B-announce}, \cite{B-thesis}, \cite{B-vosas}. In fact, the Jacobi identity for vertex operators involving two odd formal variables where the $N=2$ superconformal shift is incorporated in to the $\delta$-functions, was first given in \cite{HM} as a consequence of their notion of ``$N=2$ superconformal vertex operator superalgebra", which is (with a few minor typos corrected) equivalent to our notion of $N=2$ Neveu-Schwarz vertex operator superalgebra.  We give further consequences of this notion, including certain change of variables formulas related to the grading, as well as commutativity and associativity properties.  

Since there are principally two different basis that are widely used for the algebra of $N=2$ superconformal symmetries, i.e., the $N=2$ Neveu-Schwarz algebra, we formulate many of the notions presented here and the consequences of those notions in both basis so as to compare the present work to that in the literature.  We continue to use the terminology we introduced in \cite{B-moduli} for the two different basis -- ``homogeneous" and ``nonhomogeneous" -- and explain in Remarks \ref{homo-infinitesimals-remark} and \ref{last-nonhomo-terminology-remark} some of our reasons for using this terminology.   

We summarize the relationships between the notions of $N=2$ (Neveu-Schwarz) vertex (operator) superalgebra we present here and other notions found in the literature as follows: Our notion of ``$N=2$ Neveu-Schwarz vertex operator superalgebra with two odd formal variables in the homogeneous coordinates" is equivalent to that given in \cite{HM} in the homogeneous coordinates if one interprets their formulas in the homogeneous basis correctly; however most of the relevant formulas in \cite{HM} are given in the nonhomogeneous basis and there are a couple of typos involving incorrect signs in these formulas, as in for instance the formulation of  the vertex operators with odd formal variables, and the vertex operator corresponding to the $N=2$ superconformal element; see Remarks \ref{huang-milas-remark1} and  \ref{huang-milas-remark2}.   The notion of ``$N=2$ superconformal vertex algebra" given in \cite{Kac1997} which includes odd formal variables, is not equivalent to our notion of ``$N=2$ vertex superalgebra with odd formal variables" as equations (5.9.7) in \cite{Kac1997} do not endow the odd variable components of the vertex operators with the correct structure; see Remark \ref{nonhomo-equivalence-remark}.   The notion of ``$N=2$ vertex superalgebra with two odd formal variables in the nonhomogeneous basis" we present in Section \ref{nonhomo-section} is equivalent to the notion of ``$N_K = 2$ SUSY vertex superalgebra" as presented in \cite{HK} and \cite{He}, and some of the basic properties of such algebras that we present here appear in \cite{HK}.  However, we point out an inconsistency in the formulations given in \cite{HK} and \cite{He} involving the vertex operator corresponding to the $N=2$ superconformal element; see Remark \ref{HK-remark}.    In the early physics literature, the odd variables were often included in the fields and the vertex operator with odd formal variables corresponding to the $N=2$ superconformal element as well as the operator product expansion for this vertex operator with odd formal variables is given in, for instance \cite{P} and \cite{YZ}; see Remark \ref{mu-operator-physics-remark}.  As mentioned before, in this work, we show that the odd variable components for a vertex operator in general and for the vertex operator corresponding to the $N=2$ superconformal element in particular are completely determined by the form of the superconformal shift that appears in the Jacobi identity when two odd formal variables are included in the vertex operators.  In future work,  using some of the results presented here, we will give a geometric basis for the form that the vertex operator corresponding to the $N=2$ superconformal element must take.      

One of the main differences between our work in presenting the notion of $N=2$ vertex superalgebra with two odd formal variables and deriving the main properties of these algebras, in comparison to the work done in \cite{HK} is that we stress here (in the spirit of \cite{LL}) the Jacobi identity formulation of the notions of $N=1$ and $N=2$ vertex superalgebra and then show that virtually all the basic properties for such algebras follow quite easily and naturally {}from this one identity if the odd formal variables are included appropriately.   The generalization of the notion of vertex operator superalgebra to include odd formal variables was first done in \cite{B-thesis}; see also \cite{B-announce}, and for a more extensive treatment \cite{B-vosas}.  The notion of $N=1$ Neveu-Schwarz vertex operator superalgebra with odd formal variables presented in \cite{B-announce}--\cite{B-vosas} was the first instance of a Jacobi identity for vertex operators which included odd formal variables.   As mentioned before, the Jacobi identity was first formulated for vertex operators with two odd formal variables  in \cite{HM}.  Some of the properties that naturally follow {}from the Jacobi identity for $N=1$ Neveu-Schwarz vertex operator superalgebras with odd formal variables that were given in \cite{B-announce}--\cite{B-vosas} were generalized to other $N=n$ settings for $n>1$ in \cite{HK} and \cite{He}.  However the Jacobi identity with odd formal variables is not formulated in \cite{HK} or \cite{He}, and thus more work must be done to derive many of the useful formulas, and in addition, as we mentioned before, there are inconsistencies between or typos in some of the formulas in \cite{HM}, \cite{HK} and  \cite{He}.

After a review and clarification of the notion of $N=2$ (Neveu-Schwarz) vertex (operator) superalgebra, we present the main new results in this paper. These include a realization of the $N=2$ Neveu-Schwarz algebra as the algebra of infinitesimal $N=1$ superanalytic (as opposed to superconformal) coordinate transformations and an interpretation of this as arising {}from alternate notions of $N=1$ superconformality and the continuous deformation of $N=1$ superconformal functions by a one-parameter family corresponding to these alternate notions.    We then show that an $N=2$ (Neveu-Schwarz) vertex (operator) superalgebra is naturally a continuous deformation of an $N=1$ (Neveu-Schwarz) vertex (operator) superalgebra and is related to alternate notions of $N=1$ superconformality.  In particular we prove that the group of formal $N=2$ superconformal functions vanishing at zero and invertible in a neighborhood of zero is isomorphic to a certain subgroup of $N=1$ superanalytic functions vanishing at zero and invertible in a neighborhood of zero.   We also prove a similar statement about a group arising from $N=2$ superconformal functions vanishing at infinity and invertible in a neighborhood of infinity and the corresponding group arising from certain $N=1$ superanalytic functions vanishing at infinity and invertible in a neighborhood of infinity.  Thus one can view an $N=2$ (Neveu-Schwarz) vertex (operator) superalgebra as the algebra of correlations functions for $N=1$ superanalytic (as opposed to superconformal) fields.  It is important to note, however, that we do {\it not} use the transformation used in for instance \cite{DRS} to achieve a reduction {}from $N=2$ superconformal functions to $N=1$ superanalytic functions.  The transformation used in \cite{DRS} and, for instance \cite{He}, depends on a transformation of the even variable by a shift involving the product of the odd variables, and results in a basis for the $N=1$ superanalytic algebra which has zero central charge for the conformal elements.  Instead, we achieve a correspondence between the $N=2$ superconformal and $N=1$ superanalytic settings by continuously deforming the odd variable in the $N=1$ superconformal setting, and show that the $N=2$ Neveu-Schwarz algebra in its original basis is realized as the infinitesimal $N=1$ superanalytic transformations; see Remark \ref{compare-to-DRS-remark}.    We expect to use this result in making rigorous the connection between the algebraic and geometric aspects of $N=2$ holomorphic genus-zero superconformal field theory following the program of Huang in \cite{H-book} in the nonsuper case and the author in \cite{B-announce}--\cite{B-iso} in the $N=1$ super case as discussed below.  In addition, we expect this result to give insights into the construction of $N=2$ Neveu-Schwarz vertex operator superalgebras from $N=1$ Neveu-Schwarz vertex operator superalgebras.

The main reason for including the odd formal variables in the definition of $N=2$ (Neveu-Schwarz) vertex (operator) superalgebra is that this makes the correspondence with the underlying geometry  more transparent, giving a context to many identities that might otherwise seem arbitrary.   In fact, after each formulation of the notion of $N=2$ (Neveu-Schwarz) vertex (operator) superalgebra, we give consequences of the notion, some of which will be used in subsequent work to establish the rigorous correspondence between the algebraic settings presented here and the differential supergeometric settings which we allude to in this paper as motivation.   That is, we plan to use the results of this paper to extend the work of the author in \cite{B-announce}--\cite{B-iso} on a rigorous correspondence between the geometric and algebraic foundations of $N=1$ holomorphic genus-zero superconformal field theory to the $N=2$ case using the $N=2$ supergeometric results of \cite{B-moduli}.  This includes plans to establish the notion of ``$N=2$ supergeometric vertex operator superalgebra" and to define an isomorphism between the category of such objects and the category of $N=2$ Neveu-Schwarz vertex operator superalgebras.  The isomorphism of categories is much easier to formulate and prove when one includes the odd formal variables in the algebraic notion of $N=2$ Neveu-Schwarz vertex operator superalgebra.  This rigorous correspondence between the differential geometric and algebraic aspects of the theory of $N=2$ Neveu-Schwarz vertex operator superalgebras is necessary for the construction of $N=2$ superconformal field theory in the sense of Segal \cite{S} following the program of Huang \cite{H-thesis}--\cite{H2005-2}.  

In addition to being useful for establishing the correspondence between the geometric and algebraic aspects of $N=2$ genus-zero holomorphic superconformal field theory, including the odd variables makes it easier to formulate many of the identities arising as a consequences of the notion of $N=2$ (Neveu-Schwarz) vertex (operator) superalgebra; many vertex operator calculations become easier; and many problems such as determining change of variables formulas for $N=2$ Neveu-Schwarz vertex operator superalgebras become easier to formulate and solve (cf. \cite{B-change}).  We take this opportunity to point out that in \cite{HK}, Heluani and Kac formulate an axiomatic notion of what they call ``strongly conformal $N_K = N$ SUSY vertex superalgebras" and in  \cite{He}, Heluani  gives certain specialized change of variables formulas for these algebras.  For $N=1$ these algebras correspond to ``$N=1$ Neveu-Schwarz vertex operator superalgebras" as formulated and studied previously by the author in  \cite{B-announce}--\cite{B-vosas}, and change of variable formulas for such algebras had  been formulated and proved by the author in \cite{B-change}.  In fact in \cite{B-change}, the author not only proves general change of variables formulas for {\it any} $N=1$ superconformal change of coordinates (not just certain specialized ones near zero as given in \cite{He}) but also near infinity and more generally, on an annulus.  In addition in \cite{B-change}, the author concludes when and how one obtains isomorphic vertex operator superalgebras under general changes of variables, and proves certain convergence properties for the resulting correlation functions.  The subtleties entering into the formulation and the proof of the change of variable formulas were explained in detail in \cite{B-change}.   It was actually the rigorous treatment of the correspondence between certain geometric and algebraic aspects of $N=1$ superconformal field theory as developed in \cite{B-thesis} and  \cite{B-memoir} that led the author to formulate and study the notion of $N=1$ Neveu-Schwarz vertex operator superalgebra with odd formal variables and to introduce (in \cite{B-announce} and \cite{B-thesis}; see also \cite{B-vosas}) such properties as the $G(-1/2)$-derivative property, the Jacobi identity and the commutator formula with odd variables, supercommutativity and associativity with odd variables, and skew supersymmetry.    In fact, only using this correspondence with the $N=1$ supergeometry was the author able to, in \cite{B-change}, prove the convergence properties of the correlation functions arising {}from the change of variables formulas.  The present paper is the second in a series of papers (the first being \cite{B-moduli}) extending the author's entire algebraic and geometric program in the $N=1$ case \cite{B-announce}--\cite{B-iso}  to the $N=2$ case.

The results presented here concerning the Jacobi identity formulations of the notions of $N=1$ and $N=2$ vertex superalgebras with odd formal variables and how these formulations naturally gives rise to the desired conformal operator and the various superconformal operators were first derived in conjunction with a graduate course on vertex (super)algebras the author taught Spring 2006 at the University of Notre Dame using \cite{LL} as one of the texts.  As such, this paper, along with \cite{B-vosas} would be appropriate as part of an introductory graduate course on vertex algebras using \cite{LL} if the instructor desired a foray into the theory of vertex superalgebras, especially if one were to tie into the geometric foundations of (super)conformal field theory as presented in \cite{H-book} and \cite{B-memoir}--\cite{B-iso}, \cite{B-moduli}.  Most of the results presented here on axiomatic aspects of $N=2$ Neveu-Schwarz vertex operator superalgebras with one odd formal variable and their relationship to deformations of $N=1$ Neveu-Schwarz vertex operator superalgebras were first presented by the author in a talk at the Fields Institute in October 2000.

\subsection{Grassmann algebras and the $N=2$ Neveu-Schwarz Lie superalgebra}\label{prelim-section}

Let $\mathbb{F}$ be a field of characteristic zero,  let $\mathbb{Z}$ denote the integers, and let $\mathbb{Z}_2$ denote the integers modulo two.  For a $\mathbb{Z}_2$-graded vector space over $\mathbb{F}$, $V = V^0 \oplus V^1$, define the {\it sign function} $\eta$ on the homogeneous subspaces of $V$ by $\eta(v) = j$, for $v \in V^j$ and  $j = \mathbb{Z}_2$.  If $\eta(v) = 0$, we say that $v$ is {\it even}, and if $\eta(v) = 1$, we say that $v$ is {\it odd}. 

A {\it superalgebra} is an (associative) algebra $A$ (with identity $1 \in A$), such that: (i) $A$ is a $\mathbb{Z}_2$-graded algebra; (ii) $ab = (-1)^{\eta(a)\eta(b)} ba$ for $a,b$ homogeneous in $A$.
Note that when working over a field of characteristic zero or of characteristic greater than two, property (ii), supercommutativity, implies that the square of any odd element is zero.

The exterior algebra over a vector space $U$, denoted $\bigwedge(U)$, has the structure of a superalgebra.  Let $\mathbb{N}$ denote the natural numbers. Fix $U_L$ to be an $L$-dimensional vector space over $\mathbb{C}$ for $L \in \mathbb{N}$ such that $U_L \subset U_{L+1}$.  We denote $\bigwedge(U_L)$ by $\bigwedge_L$ and call this the {\it Grassmann algebra on $L$ generators}.  Note that $\bigwedge_L \subset \bigwedge_{L+1}$, and taking the direct limit as $L \rightarrow \infty$, we have the {\it infinite Grassmann algebra} denoted by $\bigwedge_\infty$.  Since this paper is motivated by the geometry of $N=2$ superconformal field theory  \cite{FMS}, \cite{DRS}, \cite{D}, \cite{B-moduli}, in which one considers superanalytic structures, for the purposes of this paper we have defined $\bigwedge_*$ over $\mathbb{C}$.  However, we could just as well have formulated the results that follow for Grassmann algebras over any field of characteristic zero.   

Throughout the rest of this paper we assume $\bigwedge_L$ for $L\in \mathbb{N}$ and $\bigwedge_\infty$ to be fixed.  (That is we assume $U_1 \subset U_2 \subset \cdots$ to be fixed.)  We use the notation $\bigwedge_*$ to denote the Grassmann algebra $\bigwedge_L$ for some $L \in \mathbb{N}$ or the infinite Grassmann algebra.  Note that if $L=0$, then $\bigwedge_0 = \mathbb{C}$.  We call this subalgebra of $\bigwedge_*$ the {\it body} of $\bigwedge_*$, denoted $(\bigwedge_*)_B$.  We call the orthogonal complement of the body in $\bigwedge_*$, the {\it soul} and denote it $(\bigwedge_*)_S$.  Let $\pi_B$ be the projection {}from $\bigwedge_*$ onto the body $(\bigwedge_*)_B$ (cf. \cite{D}, \cite{B-thesis}, \cite{B-memoir}).  For $a \in \bigwedge_*$, denote $\pi_B (a) = a_B$.

Let $(\bigwedge_*)^\times$ denote the set of invertible elements in $\bigwedge_*$.  Then
\[(\mbox{$\bigwedge_*$})^\times = \{a \in \mbox{$\bigwedge_*$} \; | \; a_B \neq 0 \} \]
since 
\[\frac{1}{a} = \frac{1}{a_B + a_S} = \sum_{n \in \mathbb{N}} \frac{(-1)^n a_S^n}{a_B^{n + 1}} \]
is well defined if and only if $a_B \neq 0$.

A $\mathbb{Z}_2$-graded vector space $\mathfrak{g} = \mathfrak{g}^0 \oplus \mathfrak{g}^1$ is said to be a {\it Lie superalgebra} if it has a bilinear operation $[\cdot,\cdot]$ on $\mathfrak{g}$ such that for $u,v$ homogeneous in $\mathfrak{g}$: (i) $\; [u,v] \in {\mathfrak g}^{(\eta(u) + \eta(v))\mathrm{mod} \; 2}$;  
(ii) skew symmetry holds $[u,v] = -(-1)^{\eta(u)\eta(v)}[v,u]$; (iii) the following Jacobi identity holds 
\begin{equation} \label{Jac-Lie}
(-1)^{\eta(u)\eta(w)}[[u,v],w] + (-1)^{\eta(v)\eta(u)}[[v,w],u]+ \; (-1)^{\eta(w)\eta(v)}[[w,u],v] = 0. 
\end{equation}

The {\it $N=2$ Neveu-Schwarz algebra} is the Lie superalgebra with basis consisting of the central element $d$, even elements $L_n$ and $J_n$, and odd elements $G^\pm_{n + 1/2}$ for $n \in \mathbb{Z}$, and commutation relations  
\begin{eqnarray}
\left[L_m ,L_n \right] \! \! &=& \! \! (m - n)L_{m + n} + \frac{1}{12} (m^3 - m) \delta_{m + n 
, 0} \; d , \label{Virasoro-relation} \\
\left[ J_m, J_n \right] \! \! &=& \! \! \frac{1}{3} m \delta_{m+n,0} d, \qquad \qquad \qquad \qquad \quad \ \ \left [L_m, J_n \right] \  = \ -n J_{m+n}, \\
\left[ L_m, G^\pm_{n + \frac{1}{2}} \right] \! \! &=& \! \! \left(\frac{m}{2} - n - \frac{1}{2} \right) G^\pm_{m+n+\frac{1}{2}} ,\\
\left[ J_m, G^\pm_{n + \frac{1}{2}} \right] \! \! &=& \! \!  \pm G^\pm_{m+n+\frac{1}{2}} , \qquad \qquad \qquad \ \ \ \left[ G^\pm_{m + \frac{1}{2}} , G^\pm_{n + \frac{1}{2}} \right] \ =\  0 , \\
\qquad \left[ G^+_{m + \frac{1}{2}} , G^-_{n - \frac{1}{2}} \right]  \! \! &=& \! \! 2L_{m + n} + (m-n+1) J_{m+n}  + \frac{1}{3} (m^2 + m) \delta_{m + n , 0} \; d ,\label{Neveu-Schwarz-relation-last} 
\end{eqnarray}
for $m, n \in \mathbb{Z}$. 

\begin{rema}\label{n2-generators-remark}
{\em The $N=2$ Neveu-Schwarz algebra is generated by $G^\pm_{n+1/2}$ for $n \in \mathbb{Z}$.}
\end{rema}

\begin{rema}\label{N2-subalgebras-remark}
{\em The subalgebra $\mathrm{span}_{\bigwedge_*} \{ L_{-1}, G^+_{- 1/2}, G^-_{- 1/2}, L_0, J_0, G^+_{1/2}, G^-_{1/2}, L_1 \}$ of the $N=2$ Neveu-Schwarz algebra is the Lie superalgebra $\mathfrak{osp}_{\bigwedge_*} (2|2)$ (cf. \cite{Kac}).  It is the orthogonal-symplectic superalgebra over the Grassmann algebra $\bigwedge_*$ and is also the Lie superalgebra of ``infinitesimal $N=2$ superconformal transformations of the $N=2$ super-Riemann sphere" as shown in \cite{B-moduli}.   We will denote the subalgebra of $\mathfrak{osp}_{\bigwedge_*}(2|2)$ given by $\mathrm{span}_{\bigwedge_*} \{ L_{-1}, G^+_{- 1/2}, G^-_{- 1/2} \}$ by $\mathfrak{osp}_{\bigwedge_*}(2|2)_{<0}$.
}
\end{rema}

\begin{rema}\label{auto-remark}
{\em The group of automorphisms of the $N=2$ Neveu-Schwarz algebra which preserve the Lie subalgebra generated by $L_n$ and $J_n$ for $n \in \mathbb{Z}$ are given by:
\begin{equation}\label{auto1}
G^\pm_{n + \frac{1}{2}} \mapsto b^{\pm 1} G^\pm_{n+ \frac{1}{2}}, \quad J_n \mapsto J_n, \quad L_n \mapsto L_n,
\end{equation}
for $b \in \mathbb{C}^\times$, if we are taking the algebra over $\mathbb{C}$, or more generally, for $b$ an invertible even element of the underlying Grassmann algebra, i.e., $b \in ( \bigwedge_*^0)^\times$.   In addition, we have the Virasoro-preserving automorphism 
\begin{equation}\label{auto2}
G^\pm_{n + \frac{1}{2}} \mapsto G^\mp_{n+ \frac{1}{2}}, \quad J_n \mapsto -J_n, \quad L_n \mapsto L_n.
\end{equation}
The composition of these give the Virasoro-preserving automorphisms
\begin{equation}\label{auto3}
G^\pm_{n + \frac{1}{2}} \mapsto b^{\mp 1} G^\mp_{n+ \frac{1}{2}}, \quad J_n \mapsto -J_n, \quad L_n \mapsto L_n,
\end{equation}
for $b \in (\bigwedge_*^0)^\times$, which together with (\ref{auto1}) give the group $\mathbb{Z}_2 \times \mathbb{C}^\times$ of automorphisms of the $N=2$ Neveu-Schwarz algebra.  
}
\end{rema}

Now consider the substitutions
\begin{equation}
G^{(1)}_{n + \frac{1}{2}} = \frac{1}{\sqrt{2}} \left( G^+_{n + \frac{1}{2}} + G^-_{n + \frac{1}{2}} \right), \qquad
G^{(2)}_{n + \frac{1}{2}} = \frac{i}{\sqrt{2}} \left( G^+_{n + \frac{1}{2}} - G^-_{n + \frac{1}{2}} \right) . 
\end{equation}
Note then that $G^\pm_{n+1/2} = \frac{1}{\sqrt{2}}(G^{(1)}_{n + 1/2} \mp i G^{(2)}_{n + 1/2} )$.   Under these substitutions, the $N=2$ Neveu-Schwarz algebra has a basis consisting of the central element $d$ and $L_n$, $J_n$, and $G^{(j)}_{n + 1/2}$ for $n \in \mathbb{Z}$, $j = 1,2$ and commutation relations  
\begin{eqnarray}
\left[L_m ,L_n \right] \! \! &=& \! \! (m - n)L_{m + n} + \frac{1}{12} (m^3 - m) \delta_{m + n 
, 0} \; d , \label{Virasoro-relation2} \\
\left[ J_m, J_n \right] \! \! &=& \! \! \frac{1}{3} m \delta_{m+n,0} d, \qquad \qquad \qquad \qquad \ \  \left [L_m, J_n \right] \  = \   -n J_{m+n}, \\
\left[ L_m, G^{(j)}_{n + \frac{1}{2}} \right] \! \! &=& \! \! \left(\frac{m}{2} - n - \frac{1}{2} \right) G^{(j)}_{m+n+\frac{1}{2}} ,\\
\left[ J_m, G^{(1)}_{n + \frac{1}{2}} \right] \! \! &=& \! \! - i G^{(2)}_{m+n+\frac{1}{2}} , \qquad  \qquad \qquad  \ \ \ \left[ J_m, G^{(2)}_{n + \frac{1}{2}} \right]  \  = \   i G^{(1)}_{m+n+\frac{1}{2}} , \\
\left[ G^{(j)}_{m + \frac{1}{2}} , G^{(j)}_{n - \frac{1}{2}} \right] \! \! &=& \! \! 2L_{m + n}  + \frac{1}{3} (m^2 + m) \delta_{m + n , 0} \; d , \\
\qquad \ \left[ G^{(1)}_{m + \frac{1}{2}} , G^{(2)}_{n - \frac{1}{2}} \right] \! \! &=& \! \!  - i  (m-n+1) J_{m+n} , \label{transformed-Neveu-Schwarz-relation-last}
\end{eqnarray}
for $m, n \in \mathbb{Z}$, $j=1,2$.  We call this the {\it nonhomogeneous basis} for the $N=2$ Neveu-Schwarz algebra.

\begin{rema}\label{homo-infinitesimals-remark}
{\em  Comparing the commutation relations (\ref{Virasoro-relation})--(\ref{Neveu-Schwarz-relation-last}) for the $N=2$ Neveu-Schwarz algebra in the homogeneous basis with the commutation relations (\ref{Virasoro-relation2})--(\ref{transformed-Neveu-Schwarz-relation-last}) in the nonhomogeneous basis, we see a justification for our terminology ``homogeneous" and ``nonhomogeneous".   In the former case, the commutator of $J_m$ with $G^\pm_{n + 1/2}$ is homogeneous in $G^\pm_{n+1/2}$, respectively, whereas in the later case, the commutator of $J_m$ with $G^{(j)}_{n + 1/2}$, for $j = 1,2$, is nonhomogeneous.   Another justification for our terminology comes {}from the underlying $N=2$ supergeometry.  In \cite{B-moduli}, we prove that the $N=2$ Neveu-Schwarz algebra is the algebra of infinitesimal $N=2$ superconformal transformations; see Proposition \ref{moduli-prop} below.  In the homogeneous coordinate system (which corresponds to the homogeneous basis for the $N=2$ Neveu-Schwarz algebra of infinitesimals) an $N=2$ superconformal change of variables transforms certain superderivations homogeneously of degree one (see Remark \ref{superconformal-remark}); whereas in the nonhomogeneous coordinate system (which corresponds to the nonhomogeneous basis for the $N=2$ Neveu-Schwarz algebra of infinitesimals) an $N=2$ superconformal functions is one that transforms the analogous  superderivations in the nonhomogeneous coordinates to a nonhomogeneous linear combination of the superderivations (see Remark \ref{superconformal-remark-nonhomo} and Proposition \ref{moduli-prop-nonhomo} below).  }
\end{rema} 

\begin{rema}\label{N2-subalgebras-remark-nonhomo}
{\em Following Remark \ref{N2-subalgebras-remark}, we see that in the nonhomogeneous basis the orthogonal-symplectic superalgebra is given by $\mathfrak{osp}_{\bigwedge_*} (2|2)=\mathrm{span}_{\bigwedge_*} \{ L_{-1}, G^{(1)}_{- 1/2},$ $G^{(2)}_{- 1/2}, L_0, J_0, G^{(1)}_{1/2}, G^{(2)}_{1/2}, L_1 \}$, and $\mathfrak{osp}_{\bigwedge_*}(2|2)_{<0} = \mathrm{span}_{\bigwedge_*} \{ L_{-1}, G^{(1)}_{- 1/2}, G^{(2)}_{- 1/2} \}$.
}
\end{rema}

\begin{rema}\label{nonhomo-auto-remark}
{\em In the nonhomogeneous basis, the automorphisms of the $N=2$ Neveu-Schwarz algebra which preserve the Lie subalgebra generated by $L_n$ and $J_n$ for $n \in \mathbb{Z}$ and which correspond to the automorphisms (\ref{auto1}) are given by
\begin{eqnarray}\label{nonhomo-auto1}
G^{(1)}_{n + \frac{1}{2}} &\mapsto& \frac{1}{2} \left( (b+ b^{-1}) G^{(1)}_{n+ \frac{1}{2}} - i (b-b^{-1}) G^{(2)}_{n+\frac{1}{2}}\right), \\
G^{(2)}_{n + \frac{1}{2}} &\mapsto&  \frac{1}{2} \left( i(b-b^{-1}) G^{(1)}_{n+ \frac{1}{2}} + (b+b^{-1}) G^{(2)}_{n+\frac{1}{2}}\right), \nonumber\\
J_n &\mapsto& J_n, \quad \mbox{and} \quad L_n \ \ \mapsto \ \ L_n, \nonumber
\end{eqnarray}
for $b \in (\bigwedge_*^0)^\times$.  Or writing $b = e^\beta$ for $\beta \in \bigwedge_*^0$,  this becomes
\begin{eqnarray}\label{nonhomo-auto1-hyperbolic}
G^{(1)}_{n + \frac{1}{2}} &\mapsto& (\cosh \beta) G^{(1)}_{n+ \frac{1}{2}} - i (\sinh \beta) G^{(2)}_{n+\frac{1}{2}} , \\
G^{(2)}_{n + \frac{1}{2}} &\mapsto& i(\sinh \beta) G^{(1)}_{n+ \frac{1}{2}} + (\cosh \beta) G^{(2)}_{n+\frac{1}{2}} , \nonumber\\
J_n &\mapsto& J_n, \quad \mbox{and} \quad L_n \ \ \mapsto \ \ L_n. \nonumber
\end{eqnarray}
The Virasoro-preserving automorphism corresponding to the automorphism (\ref{auto2}) is
\begin{equation}\label{nonhomo-auto2}
G^{(1)}_{n + \frac{1}{2}} \mapsto G^{(1)}_{n + \frac{1}{2}}, \quad G^{(2)}_{n + \frac{1}{2}} \mapsto - G^{(2)}_{n + \frac{1}{2}}, \quad J_n \mapsto -J_n, \quad L_n \mapsto L_n,
\end{equation}
and the composition of these give the continuous family of Virasoro-preserving automorphisms
\begin{eqnarray}\label{nonhomo-auto3}
G^{(1)}_{n + \frac{1}{2}} &\mapsto& \frac{1}{2} \left( (b + b^{-1}) G^{(1)}_{n+ \frac{1}{2}} + i (b - b^{-1}) G^{(2)}_{n+\frac{1}{2}}\right), \\
G^{(2)}_{n + \frac{1}{2}} &\mapsto&  \frac{1}{2} \left( i(b-b^{-1}) G^{(1)}_{n+ \frac{1}{2}} - (b+b^{-1}) G^{(2)}_{n+\frac{1}{2}}\right), \nonumber\\
J_n &\mapsto& -J_n, \quad \mbox{and} \quad L_n \ \ \mapsto \ \ L_n, \nonumber
\end{eqnarray}
for $b \in (\bigwedge_*^0)^\times$, or writing $b = e^\beta$ for $\beta \in \bigwedge_*^0$
\begin{eqnarray}\label{nonhomo-auto3-hyperbolic}
G^{(1)}_{n + \frac{1}{2}} &\mapsto& (\cosh \beta) G^{(1)}_{n+ \frac{1}{2}} + i (\sinh \beta) G^{(2)}_{n+\frac{1}{2}} , \\
G^{(2)}_{n + \frac{1}{2}} &\mapsto& i(\sinh \beta) G^{(1)}_{n+ \frac{1}{2}} - (\cosh \beta) G^{(2)}_{n+\frac{1}{2}} , \nonumber\\
J_n &\mapsto& -J_n, \quad \mbox{and} \quad L_n \ \ \mapsto \ \ L_n. \nonumber
\end{eqnarray}
}
\end{rema}

\subsection{Delta functions and formal calculus with odd formal variables}\label{delta-section}

This section we follow and extend the formal calculus as presented in for instance \cite{FHL}, \cite{LL}, \cite{B-thesis}, and \cite{B-vosas}.

Let $x,x_0, x_1$ and $x_2$ be formal variables which commute with each other and with $\bigwedge_*$.  We call such variables {\it even formal variables}.  Let $\varphi^\pm, \varphi^\pm_1$ and $\varphi_2^\pm$ be formal variables which commute with $x, x_0, x_1$, $x_2$ and $\bigwedge_*^0$ and anticommute with each other and $\bigwedge_*^1$.  We call such variables {\it odd formal variables}.  Note that the square of any odd formal variable is zero.   Thus for any formal Laurent series $f(x) \in  \bigwedge_*[[x,x^{-1}]]$,  we define 
\begin{eqnarray}\label{expansion}
\lefteqn{f(x + \varphi_1^+ \varphi_2^- + \varphi_1^- \varphi_2^+) }\\
&=& f(x) + (\varphi_1^+ \varphi_2^- + \varphi_1^- \varphi_2^+) f'(x) + \frac{1}{2}(\varphi_1^+ \varphi_2^- + \varphi_1^- \varphi_2^+)^2 f''(x) \nonumber \\
&=& f(x) + (\varphi_1^+ \varphi_2^- + \varphi_1^- \varphi_2^+) f'(x) +  \varphi_1^+ \varphi_1^-\varphi_2^+ \varphi_2^- f''(x) \nonumber
\end{eqnarray}
which is in $\bigwedge_* [[x,x^{-1}]] [\varphi_1^+, \varphi_1^-, \varphi_2^+, \varphi_2^-]$.

Recall (cf. \cite{FLM}, \cite{LL}) the {\it formal $\delta$-function at $x=1$} given by $\delta(x) = \sum_{n \in \mathbb{Z}} x^n$.  In \cite{B-vosas}, we extended the $\delta$-function to be defined for certain expressions involving both even and odd formal variables.   Although the terms involving the odd formal variables we will use here are slightly different than those used in \cite{B-vosas},  the arguments for extending the delta function and related identities to these expressions involving odd formal variables carry through to the present situation.  Thus using (\ref{expansion}) we have the following $\delta$-function of expressions involving three even formal variables and four odd formal variables
\begin{eqnarray}\label{delta-with}
\lefteqn{\delta \biggl( \frac{x_1 - x_2 - \varphi_1^+ \varphi_2^- - \varphi_1^- \varphi_2^+}{x_0} \biggr)}\\
 &=&  \sum_{n \in \mathbb{Z}} (x_1 - x_2 - \varphi_1^+ \varphi_2^- - \varphi_1^- \varphi_2^+)^n x_0^{-n} \nonumber \\
&=& \delta \biggl( \frac{x_1 - x_2}{x_0} \biggr)  - (\varphi_1^+ \varphi_2^- +  \varphi_1^- \varphi_2^+)x_0^{-1} \delta' \biggl( \frac{x_1 - x_2}{x_0} \biggr) \nonumber\\
& & \quad + \, \varphi^+_1 \varphi^-_1 \varphi^+_2 \varphi^-_2 x_0^{-2} \delta '' \biggl( \frac{x_1 - x_2}{x_0} \biggr) \nonumber
\end{eqnarray} 
where $\delta^k(x) = d^k/dx^k  \delta (x)$,  for $k \in \mathbb{N}$, and we use the conventions that a function of even and odd variables should be expanded about the even variables and any expression in two even variables (such as $(x_1 - x_2)^n$, for $n \in \mathbb{Z}$) should be expanded in positive powers of the second variable, (in this case $x_2$).   Note that 
\begin{equation}
\frac{\partial^k}{\partial x_0^k} x_0^{-1} \delta \biggl( \frac{x_1 - x_2 }{x_0} \biggr) = (-1)^k x_0^{-k-1} \delta^k \biggl( \frac{x_1 - x_2 }{x_0} \biggr),
\end{equation}
and thus
\begin{multline}\label{delta-with2}
\delta \biggl( \frac{x_1 - x_2 - \varphi_1^+ \varphi_2^- - \varphi_1^- \varphi_2^+}{x_0} \biggr)\\
= \delta \biggl( \frac{x_1 - x_2}{x_0} \biggr) + (\varphi_1^+ \varphi_2^- + \varphi_1^- \varphi_2^+) x_0 \frac{\partial}{\partial x_0} x_0^{-1} \delta  \biggl( \frac{x_1 - x_2}{x_0} \biggr)  \\
+  \varphi^+_1 \varphi^-_1 \varphi^+_2 \varphi^-_2  x_0 \frac{\partial^2}{\partial x_0^2} x_0^{-1} \delta  \biggl( \frac{x_1 - x_2}{x_0} \biggr) .
\end{multline}

Following the treatment of formal calculus for even formal variables given in \cite{FHL}, we have  
\begin{equation}\label{delta1}
f(x) \delta (x) = f(1) \delta (x) \qquad \mbox{for $f(x) \in \mathbb{C}[x,x^{-1}]$}.
\end{equation}
Let $V$ be a vector space over $\mathbb{C}$.  For $X(x_1,x_2) \in (\mbox{End} \; V)[[x_1, x_1^{-1},x_2, x_2^{-1}]]$ such that
\[\lim_{x_1 \rightarrow x_2} X(x_1,x_2) \quad \mbox{exists, i.e.,} \quad X(x_1,x_2)|_{x_1 = x_2} \quad \mbox{exists} \]
(that is, when $X(x_1,x_2)$ is applied to any element of $V$, setting the variables equal leads to only finite sums in $V$), we have
\begin{equation}\label{delta-multiplication}
X(x_1,x_2) \delta \biggl( \frac{x_1}{x_2} \biggr) = X(x_2,x_2) \delta \biggl( \frac{x_1}{x_2} \biggr) .
\end{equation}
For the three-variable generating function
\[\delta \biggl( \frac{x_1 - x_2}{x_0} \biggr) = \sum_{n \in \mathbb{Z}} \frac{(x_1 - x_2)^n}{x_0^n} = \sum_{n \in \mathbb{Z}, m \in \mathbb{N}} (-1)^m \binom{n}{m} x_0^{-n} x_1^{n-m} x_2^m ,\]  
we have
\begin{equation}\label{delta-2-terms}
x_1^{-1} \delta \biggl( \frac{x_2 + x_0}{x_1} \biggr) = x_2^{-1} \delta \biggl( \frac{x_1 - x_0}{x_2} \biggr) 
\end{equation}
and
\begin{equation}\label{delta-3-terms}
x_0^{-1} \delta \biggl( \frac{x_1 - x_2}{x_0} \biggr) - x_0^{-1} \delta \biggl( \frac{x_2 - x_1}{-x_0} \biggr) = x_2^{-1} \delta \biggl( \frac{x_1 - x_0}{x_2} \biggr) .
\end{equation} 
Notice that in the spirit of the $\delta$-function multiplication principle (\ref{delta-multiplication}), the expressions on both sides of (\ref{delta-2-terms}) and the three terms occurring in (\ref{delta-3-terms}) all correspond to the same formal substitution $x_0 = x_1 - x_2$. 

The above facts (\ref{delta1})--(\ref{delta-multiplication}) extend immediately to the following: 
\[f(x,\varphi^+, \varphi^-)\delta(x) = f(1,\varphi^+,\varphi^-)\delta(x) \qquad \mbox{for  $f(x,\varphi^+, \varphi^-) \in \bigwedge_*[x,x^{-1},\varphi^+, \varphi^-]$}.\]
Let $V$ be a $\bigwedge_*$-module. If 
\[X(x_1,\varphi_1^+, \varphi_1^-,x_2,\varphi_2^+, \varphi_2^-) \in (\mbox{End} \; V)[[x_1,x_1^{-1},x_2, x_2^{-1}]][\varphi_1^+, \varphi_1^-, \varphi_2^+, \varphi_2^-] \] 
such that  
\[\lim_{x_1 \rightarrow x_2} X(x_1,\varphi_1^+. \varphi_1^-,x_2,\varphi_2^+, \varphi_2^-) \]
exists, i.e.,
\[X(x_1,\varphi_1^+, \varphi_1^-,x_2,\varphi_2^+, \varphi_2^-)|_{x_1 = x_2} \]
exists, then 
\begin{equation}\label{delta-multiplication-with-phis}
X(x_1,\varphi_1^+, \varphi_1^-,x_2,\varphi_2^+, \varphi_2^-) \delta \biggl( \frac{x_1}{x_2} \biggr) =
X(x_2,\varphi_1^+, \varphi_1^-,x_2,\varphi_2^+, \varphi_2^-) \delta \biggl(\frac{x_1}{x_2} \biggr) . 
\end{equation}

Repeatedly taking $\frac{\partial}{\partial x_0}$ of both sides of (\ref{delta-2-terms}), we obtain the following identity (cf. \cite{B-vosas}):
\begin{equation}\label{delta-derivative-2-terms}
x_1^{-k-1} \delta^k \biggl( \frac{x_2 + x_0}{x_1} \biggr) = (-1)^k x_2^{-k-1} \delta^k \biggl( \frac{x_1 - x_0}{x_2} \biggr) .
\end{equation}
Repeatedly taking $\frac{\partial}{\partial x_1}$ of both sides of (\ref{delta-3-terms}), we obtain (cf. \cite{B-vosas}):
\begin{equation}\label{delta-derivative-3-terms}
x_0^{-k - 1} \delta^k \biggl( \frac{x_1 - x_2}{x_0} \biggr) - x_0^{-k-1} \delta^k \biggl( \frac{x_2 - x_1}{-x_0} \biggr) = x_2^{-k-1} \delta^k \biggl( \frac{x_1 - x_0}{x_2} \biggr) .
\end{equation}
        
Thus {}from (\ref{delta-with}), (\ref{delta-2-terms}), (\ref{delta-3-terms}), (\ref{delta-derivative-2-terms}) and (\ref{delta-derivative-3-terms}), we have 
\begin{equation}\label{delta-2-terms-with-phis}
x_1^{-1} \delta \biggl( \frac{x_2 + x_0 + \varphi_1^+ \varphi_2^- + \varphi_1^- \varphi_2^+}{x_1}
\biggr) = x_2^{-1} \delta \biggl( \frac{x_1 - x_0 - \varphi_1^+ \varphi_2^- - \varphi_1^- \varphi_2^+ }{x_2} \biggr)  
\end{equation}
and
\begin{multline}\label{delta-3-terms-with-phis}
x_0^{-1} \delta \biggl( \frac{x_1 - x_2 - \varphi_1^+ \varphi_2^- - \varphi_1^- \varphi_2^+}{x_0}
\biggr) - x_0^{-1} \delta \biggl( \frac{x_2 - x_1 + \varphi_1^+ \varphi_2^- + \varphi_1^- \varphi_2^+}{-x_0} \biggr)  \\
=  x_2^{-1} \delta \biggl( \frac{x_1 - x_0 - \varphi_1^+ \varphi_2^- - \varphi_1^- \varphi_2^+}{x_2} \biggr) . 
\end{multline}
Notice that in the spirit of the $\delta$-function multiplication principle (\ref{delta-multiplication-with-phis}), the expressions on both sides of (\ref{delta-2-terms-with-phis}) and the three terms
occurring in (\ref{delta-3-terms-with-phis}) all correspond to the same formal substitution 
\begin{equation}
x_0 = x_1 - x_2  - \varphi_1^+  \varphi_2^- - \varphi_1^- \varphi_2^+.
\end{equation}  
And thus  for
\[X(x_0,x_1,\varphi_1^+, \varphi_1^-,x_2,\varphi_2^+, \varphi_2^-) \in (\mbox{End} \;V)
[[x_0,x_0^{-1},x_1,x_1^{-1},x_2,x_2^{-1}]] [\varphi_1^+, \varphi_1^-, \varphi_2^+, \varphi_2^-] ,\] 
a formal substitution corresponding to $x_0 = x_1 - x_2 - \varphi_1^+ \varphi_2^- - \varphi_1^-
\varphi_2^+$ can be made as long as the resulting expression is well defined, e.g., if 
\[X(x_1,\varphi_1^+, \varphi_1^-,x_2,\varphi_2^+, \varphi_2^-) \in (\mbox{End} \;V) [[x_1,x_1^{-1}]] ((x_2)) [\varphi_1^+, \varphi_1^-, \varphi_2^+, \varphi_2^-] ,\] 
then  
\begin{multline}\label{delta-substitute}
\delta \biggl( \frac{x_2 + x_0 + \varphi_1^+ \varphi_2^- + \varphi_1^- \varphi_2^+}{x_1} \biggr)
X(x_1,\varphi_1^+,\varphi_1^-,x_2,\varphi_2^+, \varphi_2^-) \\
= \delta \biggl( \frac{x_2 + x_0 + \varphi_1^+ \varphi_2^- + \varphi_1^- \varphi_2^+}{x_1} \biggr) X(x_2 + x_0 + \varphi_1^+ \varphi_2^- + \varphi_1^- \varphi_2^+, \varphi_1^+,  \\
\varphi_1^-, x_2,\varphi_2^+, \varphi_2^-) , 
\end{multline}
and (\ref{delta-substitute}) holds for more general $X(x_1, \varphi_1^+,\varphi_1^-, x_2, \varphi_2^+, \varphi_2^-)$ as long as the expression $X(x_2 + x_0 + \varphi_1^+  \varphi_2^- + \varphi_1^-
\varphi_2^+, \varphi_1^+, \varphi_1^-,x_2,\varphi_2^+, \varphi_2^-)$ is well defined.  

\begin{rema}\label{superconformal-remark}
{\em The substitution  $x_0 = x_1 - x_2 - \varphi_1^+ \varphi_2^- - \varphi_1^- \varphi_2^+$ can be 
thought of as the even part of a ``superconformal shift of  $(x_1, \varphi_1^+, \varphi_1^-)$ by $(x_2, \varphi_2^+, \varphi_2^-)$".  That is, if we let $f(x, \varphi^+, \varphi^-) = (\tilde{x}, \tilde{\varphi}^+, \tilde{\varphi}^-)$ with $\tilde{x} \in (\bigwedge_*[[x, x^{-1}]][\varphi^+, \varphi^-])^0$ and $\tilde{\varphi}^\pm \in (\bigwedge_*[[x, x^{-1}]][\varphi^{(1)}, \varphi^{(2)}])^1$, then $f$ is said to be
{\it $N=2$ superconformal} in $x$ and $\varphi^\pm$ if and only if $f$ satisfies 
\begin{equation}\label{nice-superconformal-condition}
D^- \tilde{\varphi}^+ = D^+ \tilde{\varphi}^- = D^+\tilde{x} - \tilde{\varphi}^- D^+ \tilde{\varphi}^+ = D^-\tilde{x} - \tilde{\varphi}^+ D^-\tilde{\varphi}^- =0 
\end{equation}     
for 
\begin{equation}\label{define-Dpm-derivative}
D^\pm = \frac{\partial}{\partial \varphi^\pm} + \varphi^\mp  \frac{\partial}{\partial x}
\end{equation} 
(see \cite{B-moduli}).  For a formal superanalytic vector-valued function $f(x, \varphi^+, \varphi^-) = (\tilde{x}, \tilde{\varphi}^+$, $\tilde{\varphi}^-)$, the conditions (\ref{nice-superconformal-condition}) are equivalent to requiring that $f$ transform both of the super-differential operators $D^\pm$ homogeneously of degree one so that $D^\pm = (D^\pm \tilde{\varphi}^\pm ) \tilde{D}^\pm$. We observe that $f(x_1,\varphi_1^+, \varphi_1^-) = (x_1 - x_2 - \varphi_1^+ \varphi_2^- - \varphi_1^- \varphi_2^+, \varphi_1^+ - \varphi_2^+, \varphi_1^- - \varphi_2^-)$ is formally $N=2$ superconformal in $x_1$ and $\varphi_1^\pm$ since it satisfies (\ref{nice-superconformal-condition}). }
\end{rema}

{}From \cite{LL} we have the following fact about $\delta$-functions and expansions of zero:
\begin{equation}\label{expansion-without}
(x_1 - x_2)^{-n-1} - (-x_2 + x_1)^{-n-1} = \frac{(-1)^n}{n!} \left( \frac{\partial}{\partial x_1} \right)^n x_2^{-1} \delta \left(\frac{x_1}{x_2} \right) ,
\end{equation}
for $n \in \mathbb{N}$.  Thus including formal odd variables we have the following proposition:

\begin{prop}\label{expansion-prop}
For $n \in \mathbb{N}$, 
\begin{multline}
(x_1 - x_2 - \varphi^+_1 \varphi^-_2 - \varphi^-_1 \varphi^+_2)^{-n-1} -  (-x_2 + x_1 - \varphi^+_1 \varphi^-_2 - \varphi^-_1 \varphi^+_2)^{-n-1} \\
 = \frac{(-1)^n}{n!} \left( \frac{\partial}{\partial x_1} \right)^n x_2^{-1} \delta \left(\frac{x_1 - \varphi^+_1 \varphi^-_2 - \varphi^-_1 \varphi^+_2}{x_2} \right).
\end{multline}
\end{prop}

\begin{proof}
{}From (\ref{expansion}) and (\ref{expansion-without}), we have
\begin{eqnarray*}
\lefteqn{(x_1 - x_2 - \varphi^+_1 \varphi^-_2 - \varphi^-_1 \varphi^+_2)^{-n-1} -  (-x_2 + x_1 - \varphi^+_1 \varphi^-_2 - \varphi^-_1 \varphi^+_2)^{-n-1} }\\
&=&  (x_1 - x_2)^{-n-1} - (-x_2 + x_1)^{-n-1} \\
& & \quad + \, (n+1)(\varphi^+_1 \varphi^-_2 + \varphi^-_1 \varphi^+_2) \bigl( (x_1 - x_2)^{-n-2} - (-x_2 + x_1)^{-n-2}\bigr) \\
& & \quad + \, (n+1)(n+2)\varphi^+_1 \varphi_2^- \varphi_1^- \varphi_2^+ \bigl( (x_1 - x_2)^{-n-3} - (-x_2 + x_1)^{-n-3} \bigr) \\
\end{eqnarray*}
\begin{eqnarray*}
&=&  \frac{(-1)^n}{n!} \left( \frac{\partial}{\partial x_1} \right)^n x_2^{-1} \delta \left(\frac{x_1}{x_2} \right) \\
& & \quad + \, (n+1)(\varphi^+_1 \varphi^-_2 + \varphi^-_1 \varphi^+_2) 
\frac{(-1)^{n+1}}{(n+1)!} \left( \frac{\partial}{\partial x_1} \right)^{n+1} x_2^{-1} \delta \left(\frac{x_1}{x_2} \right) 
 \\
& & \quad + \,  (n+1)(n+2)\varphi^+_1 \varphi_2^- \varphi_1^- \varphi_2^+ \frac{(-1)^{n+2}}{(n+2)!} \left( \frac{\partial}{\partial x_1} \right)^{n+2} x_2^{-1} \delta \left(\frac{x_1}{x_2} \right) \\
&=& \frac{(-1)^n}{n!} \left( \frac{\partial}{\partial x_1} \right)^n  \biggl( 1- (\varphi^+_1 \varphi^-_2 + \varphi^-_1 \varphi^+_2)  \frac{\partial}{\partial x_1}   \\
& & \quad + \, \varphi^+_1 \varphi_2^- \varphi_1^- \varphi_2^+  \left( \frac{\partial}{\partial x_1} \right)^2 \biggr) x_2^{-1} \delta \left(\frac{x_1}{x_2} \right)  \\
&=&  \frac{(-1)^n}{n!} \left( \frac{\partial}{\partial x_1} \right)^n x_2^{-1} \delta \left(\frac{x_1 - \varphi^+_1 \varphi^-_2 - \varphi^-_1 \varphi^+_2}{x_2} \right).
\end{eqnarray*}
\end{proof}

We will often use a formal residue operator $\mathrm{Res}_x$ {}from $V[[x,x^{-1}]][\varphi^+, \varphi^-]$ to $V[\varphi^+, \varphi^-]$ which we define by
\begin{equation}
\mathrm{Res}_x \sum_{n \in \mathbb{Z}} a_n (\varphi^+, \varphi^-) x^n   = a_{-1}(\varphi^+, \varphi^-)
\end{equation}
where $a_n(\varphi^+, \varphi^-) \in V[\varphi^+, \varphi^-]$, for $n \in \mathbb{Z}$.   

We have the following extension of the the formal Taylor Theorem (cf. \cite{FLM}, \cite{LL}):
\begin{prop}\label{Taylor}
Let $v(x,\varphi^+, \varphi^-) \in V[[x,x^{-1}]][\varphi^+, \varphi^-]$. Then 
\begin{multline}
e^{x_0 \frac{\partial}{\partial x} + \varphi^+_0 \left(
\frac{\partial}{\partial \varphi^+} + \varphi^- \frac{\partial}{\partial x} \right)+ \varphi^-_0 \left(
\frac{\partial}{\partial \varphi^-} + \varphi^+ \frac{\partial}{\partial x} \right)} v(x, \varphi^+, \varphi^-) \\
= v(x + x_0 + \varphi^+_0 \varphi^- + \varphi_0^- \varphi^+, \varphi_0^+ + \varphi^+, \varphi_0^- + \varphi^-) \\
 \end{multline}
(and these expressions exist).
\end{prop}

\begin{proof}
By direct expansion, we see that 
\begin{multline}
e^{x_0 \frac{\partial}{\partial x} + \varphi^+_0 \left( \frac{\partial}{\partial \varphi^+} + \varphi^- \frac{\partial}{\partial x} \right)+ \varphi^-_0 \left( \frac{\partial}{\partial \varphi^-} + \varphi^+ \frac{\partial}{\partial x} \right)} (x, \varphi^+, \varphi^-) \\
= (x + x_0 + \varphi^+_0 \varphi^- + \varphi_0^- \varphi^+, \varphi_0^+ + \varphi^+, \varphi_0^- + \varphi^-) .
\end{multline}
By Proposition 6.20 in \cite{B-moduli}, we have that 
\begin{multline}
e^{x_0 \frac{\partial}{\partial x} + \varphi^+_0 \left( \frac{\partial}{\partial \varphi^+} + \varphi^- \frac{\partial}{\partial x} \right)+ \varphi^-_0 \left( \frac{\partial}{\partial \varphi^-} + \varphi^+ \frac{\partial}{\partial x} \right)} v(x, \varphi^+, \varphi^-) \\
= v\left(e^{x_0 \frac{\partial}{\partial x} + \varphi^+_0 \left(
\frac{\partial}{\partial \varphi^+} + \varphi^- \frac{\partial}{\partial x} \right)+ \varphi^-_0 \left(
\frac{\partial}{\partial \varphi^-} + \varphi^+ \frac{\partial}{\partial x} \right)} (x, \varphi^+, \varphi^-)\right).
\end{multline}
The result follows.
\end{proof}

\section{$N=2$ (Neveu-Schwarz) vertex (operator) superalgebras with two odd formal variables}\label{N2-with-section}

\subsection{The notion of $N=2$ vertex superalgebra with two odd formal variables}\label{vertex-superalgebra-section}

\begin{defn}\label{vertex-superalgebra-definition}
{\em An} $N=2$ vertex superalgebra over $\bigwedge_*$ and with two odd formal variables {\em consists of a $\mathbb{Z}_2$-graded $\bigwedge_*$-module (graded by} sign $\eta${\em) 
\begin{equation}\label{first-n2-superalgebra-with}
V = V^0 \oplus V^1
\end{equation} 
equipped, first, with a linear map 
\begin{eqnarray}\label{operator-with}
V &\longrightarrow&  (\mbox{End} \; V)[[x,x^{-1}]][\varphi^+, \varphi^-] \\
v  &\mapsto&  Y(v,(x,\varphi^+, \varphi^-)) \nonumber
\end{eqnarray}
with 
\begin{equation}
Y(v,(x,\varphi^+, \varphi^-)) = \sum_{n \in \mathbb{Z}} \Bigl( v_n  + \varphi^+  v_{n - \frac{1}{2}}^+ + \varphi^-  v_{n - \frac{1}{2}}^- + \varphi^+ \varphi^- v_{n-1}^{+-} \Bigr)  x^{-n-1},
\end{equation}
where for the $\mathbb{Z}_2$-grading of $\mathrm{End} \; V$ induced {}from that of $V$, we have
\begin{equation}
v_n, v_n^{+-} \in (\mbox{End} \; V)^{\eta(v)}, \quad \mathrm{and} \quad v_{n - \frac{1}{2}}^\pm \in  (\mbox{End} \; V)^{(\eta(v) + 1) \mbox{\begin{footnotesize} mod 
\end{footnotesize}} 2} 
\end{equation}
for $v$ of homogeneous sign in $V$, $x$ is an even formal variable, and $\varphi^\pm$ are odd formal variables, and where $Y(v,(x,\varphi^+, \varphi^-))$ denotes the} vertex operator associated with 
$v$.  {\em We also have a distinguished element $\mathbf{1}$ in $V$ (the} vacuum vector{\em ).  The following conditions are assumed for $u,v \in V$:  the} truncation conditions: {\em
\begin{equation}\label{truncation}
u_n v = u_{n}^{+-} v = u_{n-\frac{1}{2}}^\pm v  = 0 \qquad \mbox{for $n \in \mathbb{Z}$ sufficiently large,} 
\end{equation}
that is
\begin{equation}
Y(u, (x, \varphi^+, \varphi^-))v \in V((x))[\varphi^+, \varphi^-];
\end{equation}
next, the following} vacuum property{\em :
\begin{equation}\label{vacuum-identity}
Y(\mathbf{1}, (x, \varphi^+, \varphi^-)) = \mathrm{id}_V ;
\end{equation}
the} creation property {\em holds:
\begin{eqnarray}\
Y(v,(x,\varphi^+,\varphi^-)) \mathbf{1} &\in& V[[x]][\varphi^+,\varphi^-] \label{creation-property1} \\
\lim_{(x,\varphi^+,\varphi^-) \rightarrow 0} Y(v,(x, \varphi^+,\varphi^-)) \mathbf{1} &=& v ;  \label{creation-property2}
\end{eqnarray} 
and finally the} Jacobi identity {\em holds:  

\begin{eqnarray}\label{Jacobi-identity}
\lefteqn{\quad x_0^{-1} \delta \biggl( \frac{x_1 - x_2 - \varphi^+_1 \varphi^-_2 - \varphi^-_1 \varphi^+_2}{x_0} \biggr)Y(u,(x_1, \varphi^+_1, \varphi^-_1)) Y(v,(x_2, \varphi^+_2,\varphi^-_2))} \\
& & \quad - (-1)^{\eta(u)\eta(v)} x_0^{-1} \delta \biggl( \frac{x_2 - x_1 + \varphi^+_1 \varphi^-_2 + \varphi^-_1 \varphi^+_2}{-x_0} \biggr) Y(v,(x_2, \varphi^+_2,\varphi^-_2)) \nonumber\\
& & \quad \cdot Y(u,(x_1,\varphi^+_1,\varphi^-_1)) \nonumber \\
\quad &=& \! \! \! x_2^{-1} \delta \biggl( \frac{x_1 - x_0 - \varphi^+_1 \varphi^-_2 - \varphi^-_1 \varphi^+_2}{x_2} \biggr) Y(Y(u,(x_0,  \varphi^+_1 - \varphi^+_2,\varphi^-_1 - \varphi^-_2))v, \nonumber\\
& & \quad (x_2, \varphi^+_2,\varphi^-_2)) , \nonumber
\end{eqnarray}
for $u,v$ of homogeneous sign in $V$.
}
\end{defn}

The $N=2$ vertex superalgebra just defined is denoted by 
\[(V,Y(\cdot,(x, \varphi^+, \varphi^-)), \mathbf{1}),\] 
or for simplicity by $V$.

Since $\bigwedge_L \subset \bigwedge_\infty$ for $L\in \mathbb{N}$, any $N=2$ vertex superalgebra $V$ over $\bigwedge_L$ can be extended to an $N=2$ vertex superalgebra over $\bigwedge_\infty$ by taking the $\bigwedge_\infty$-module induced by the $\bigwedge_L$-module $V$.

Let $(V_1, Y_1(\cdot,(x,\varphi^+, \varphi^-)),\mbox{\bf 1}_1)$ and $(V_2, Y_2(\cdot,(x,\varphi^+, \varphi^-)),\mbox{\bf 1}_2)$ be two $N=2$ vertex superalgebras.  A {\it homomorphism of $N=2$ vertex superalgebras with odd formal variables} is a $\mathbb{Z}_2$-graded $\bigwedge_\infty$-module homomorphism $\gamma : V_1 \longrightarrow V_2 \;$ satisfying
\begin{equation}
\gamma (Y_1(u,(x,\varphi^+,\varphi^-))v) = Y_2(\gamma(u),(x,\varphi^+, \varphi^-))\gamma(v)
\quad \mbox{for} \quad u,v \in V_1 ,
\end{equation}
and $\gamma(\mbox{\bf 1}_1) = \mbox{\bf 1}_2$.  

Note that an $N=2$ vertex superalgebra over $\bigwedge_*$ is isomorphic to the same $N=2$ Neveu-Schwarz vertex operator superalgebra over a different subalgebra of the underlying Grassmann algebra $\bigwedge_\infty$.

We have the following consequences of the definition of $N=2$ vertex superalgebra with two odd formal variables:
By the creation property, we see that $\mathbf{1} \in V^0$, and 
\begin{equation}
v_{-1} \cdot \mathbf{1} = v.
\end{equation}

Taking $\mathrm{Res}_{x_0}$ of the Jacobi identity and using the $\delta$-function identity (\ref{delta-2-terms-with-phis}), we obtain the following {\it supercommutator formula} with odd formal variables
\begin{multline}\label{homo-bracket}
[ Y(u, (x_1,\varphi^+_1, \varphi_1^-)), Y(v,(x_2,\varphi_2^+, \varphi_2^-))]  \\
=  \mbox{Res}_{x_0} x_2^{-1} \delta \biggl( \frac{x_1 - x_0 - \varphi_1^+ \varphi_2^- - \varphi_1^- \varphi_2^+}{x_2} \biggr) \\
\cdot Y(Y(u,(x_0, \varphi_1^+ - \varphi_2^+, \varphi_1^- - \varphi_2^-))v,(x_2, \varphi_2^+, \varphi_2^-)) . 
\end{multline}
Taking $\mathrm{Res}_{x_1}$ of the Jacobi identity and using the $\delta$-function identity (\ref{delta-2-terms-with-phis}), we obtain the following {\it iterate formula} for $N=2$ vertex operators with two odd formal variables
\begin{multline}\label{iterate}
Y(Y(u,(x_0, \varphi_1^+ - \varphi_2^+, \varphi_1^- - \varphi_2^-))v,(x_2, \varphi_2^+, \varphi_2^-))\\
= \mathrm{Res}_{x_1} \biggl( x_0^{-1} \delta \biggl( \frac{x_1 - x_2 - \varphi^+_1 \varphi^-_2 - \varphi^-_1 \varphi^+_2}{x_0} \biggr)Y(u,(x_1, \varphi^+_1, \varphi^-_1)) Y(v,(x_2, \varphi^+_2,\varphi^-_2)) \\
- \, (-1)^{\eta(u)\eta(v)} x_0^{-1} \delta \biggl( \frac{x_2 - x_1 + \varphi^+_1 \varphi^-_2 + \varphi^-_1 \varphi^+_2}{-x_0} \biggr) Y(v,(x_2, \varphi^+_2,\varphi^-_2))\\
\cdot Y(u,(x_1,\varphi^+_1,\varphi^-_1)) \biggr).
\end{multline}

We will sometimes use the following notation for the vertex operator associated with $v \in V$:
\begin{eqnarray}
Y(v,(x,\varphi^+, \varphi^-)) &=& \sum_{n \in \mathbb{Z}} \Bigl( v_n  + \varphi^+  v_{n - \frac{1}{2}}^+ + \varphi^-  v_{n - \frac{1}{2}}^- + \varphi^+ \varphi^- v_{n-1}^{+-} \Bigr)  x^{-n-1} \nonumber \\
&=& v(x) + \varphi^+ v^+ (x) + \varphi^- v^-(x) + \varphi^+ \varphi^- v^{+-}(x). \label{v-notation}
\end{eqnarray}

\begin{prop}\label{Dpm-derivative-prop}
Let $V$ be an $N=2$ vertex superalgebra with two odd formal variables and let $\mathcal{D}^\pm$ be the odd endomorphisms of $V$ defined by
\begin{equation}\label{define-Dpm}
\mathcal{D}^\pm (v) = v^\pm_{-\frac{3}{2}} \mathbf{1} \qquad \mbox{for $v \in V$.}
\end{equation}
Then 
\begin{equation}\label{Dpm-derivative-property}
Y(\mathcal{D}^\pm v,(x,\varphi^+, \varphi^-)) = \biggl( \frac{\partial}{\partial \varphi^\pm} + \varphi^\mp \frac{\partial}{\partial x} \biggr) Y(v,(x,\varphi^+, \varphi^-))  . 
\end{equation} 
\end{prop}

\begin{proof}
Using the iterate formula (\ref{iterate}), the vacuum property (\ref{vacuum-identity}), Proposition \ref{expansion-prop}, (\ref{delta-2-terms-with-phis}) with $x_0=0$ (which is well-defined), (\ref{delta-substitute}), and the notation (\ref{v-notation}), we have
\begin{eqnarray*}
\lefteqn{Y(\mathcal{D}^\pm v,(x_2,\varphi^+_2, \varphi^-_2)) }\\
&=& Y(v^\pm_{-\frac{3}{2}} \mathbf{1},(x_2,\varphi^+_2, \varphi^-_2)) \\
&=& \left( \frac{\partial}{\partial \varphi_1^\pm} \pm \varphi^\mp_2 \frac{\partial}{\partial \varphi_1^-} \frac{\partial}{\partial \varphi_1^+} \right) \mathrm{Res}_{x_0} x_0^{-1} Y(Y(v,(x_0, \varphi_1^+ - \varphi_2^+, \varphi_1^- - \varphi_2^-)) \mathbf{1}, \\
& & \quad (x_2, \varphi_2^+, \varphi_2^-)) \biggr|_{\varphi_1^\mp = 0}\\
&=& \left( \frac{\partial}{\partial \varphi_1^\pm} \pm \varphi^\mp_2 \frac{\partial}{\partial \varphi_1^-} \frac{\partial}{\partial \varphi_1^+} \right) \mathrm{Res}_{x_0} x_0^{-1}  \mathrm{Res}_{x_1} \biggl( x_0^{-1} \delta \biggl( \frac{x_1 - x_2 - \varphi^+_1 \varphi^-_2 - \varphi^-_1 \varphi^+_2}{x_0} \biggr)\\
\end{eqnarray*}
\begin{eqnarray*}
& & \quad \cdot Y(v,(x_1, \varphi^+_1, \varphi^-_1))   -  x_0^{-1} \delta \biggl( \frac{x_2 - x_1 + \varphi^+_1 \varphi^-_2 + \varphi^-_1 \varphi^+_2}{-x_0} \biggr)\\
& & \quad  \cdot Y(v,(x_1,\varphi^+_1,\varphi^-_1)) \biggr) \biggr|_{\varphi_1^\mp = 0}\\
&=& \left( \frac{\partial}{\partial \varphi_1^\pm} \pm \varphi^\mp_2 \frac{\partial}{\partial \varphi_1^-} \frac{\partial}{\partial \varphi_1^+} \right)  \mathrm{Res}_{x_1} \biggl( \Bigl( (x_1 - x_2 - \varphi^+_1 \varphi^-_2 - \varphi^-_1 \varphi^+_2)^{-1}\\
& & \quad  -  \, (-x_2 + x_1 - \varphi^+_1 \varphi^-_2 - \varphi^-_1 \varphi^+_2)^{-1} \Bigr)  Y(v,(x_1,\varphi^+_1,\varphi^-_1)) \biggr) \biggr|_{\varphi_1^\mp = 0}\\
&=& \left( \frac{\partial}{\partial \varphi_1^\pm} \pm \varphi^\mp_2 \frac{\partial}{\partial \varphi_1^-} \frac{\partial}{\partial \varphi_1^+} \right)  \mathrm{Res}_{x_1} \biggl( x_2^{-1} \delta \left(\frac{x_1 - \varphi^+_1 \varphi^-_2 - \varphi^-_1 \varphi^+_2}{x_2} \right) \\
& & \quad \cdot Y(v,(x_1,\varphi^+_1,\varphi^-_1)) \biggr) \biggr|_{\varphi_1^\mp = 0}\\
&=& \left( \frac{\partial}{\partial \varphi_1^\pm} \pm \varphi^\mp_2 \frac{\partial}{\partial \varphi_1^-} \frac{\partial}{\partial \varphi_1^+} \right)  \mathrm{Res}_{x_1} \biggl( x_1^{-1} \delta \left(\frac{x_2 + \varphi^+_1 \varphi^-_2 + \varphi^-_1 \varphi^+_2}{x_1} \right) \\
& & \quad \cdot Y(v,(x_1,\varphi^+_1,\varphi^-_1)) \biggr) \biggr|_{\varphi_1^\mp = 0}\\
&=& \left( \frac{\partial}{\partial \varphi_1^\pm} \pm \varphi^\mp_2 \frac{\partial}{\partial \varphi_1^-} \frac{\partial}{\partial \varphi_1^+} \right)  \mathrm{Res}_{x_1} \biggl( x_1^{-1} \delta \left(\frac{x_2 + \varphi^+_1 \varphi^-_2 + \varphi^-_1 \varphi^+_2}{x_1} \right) \\
& & \quad \cdot Y(v,(x_2 + \varphi^+_1 \varphi^-_2 + \varphi^-_1 \varphi^+_2 ,\varphi^+_1,\varphi^-_1)) \biggr) \biggr|_{\varphi_1^\mp = 0}\\
&=& \left( \frac{\partial}{\partial \varphi_1^\pm} \pm \varphi^\mp_2 \frac{\partial}{\partial \varphi_1^-} \frac{\partial}{\partial \varphi_1^+} \right)  Y(v,(x_2 + \varphi^+_1 \varphi^-_2 + \varphi^-_1 \varphi^+_2 ,\varphi^+_1,\varphi^-_1)) \biggr|_{\varphi_1^\mp = 0}\\
&=& \left( \frac{\partial}{\partial \varphi_1^\pm} \pm \varphi^\mp_2 \frac{\partial}{\partial \varphi_1^-} \frac{\partial}{\partial \varphi_1^+} \right)  \biggl( Y(v,(x_2 ,\varphi^+_1,\varphi^-_1)) + ( \varphi^+_1 \varphi^-_2 + \varphi^-_1 \varphi^+_2) \\
& & \quad \frac{\partial}{\partial x_2}  Y(v,(x_2 ,\varphi^+_1,\varphi^-_1)) +  \varphi^+_1 \varphi^-_2 \varphi^-_1 \varphi^+_2 \left(\frac{\partial}{\partial x_2}\right)^2  Y(v,(x_2 ,\varphi^+_1,\varphi^-_1))\biggr) \biggr|_{\varphi_1^\mp = 0}\\
&=& \biggl( \frac{\partial}{\partial \varphi_1^\pm}  Y(v,(x_2 ,\varphi^+_1,\varphi^-_1)) 
+   \varphi^\mp_2 \frac{\partial}{\partial x_2} Y(v,(x_2 ,\varphi^+_1,\varphi^-_1))  \pm \varphi^\mp_2 \frac{\partial}{\partial \varphi_1^-} \frac{\partial}{\partial \varphi_1^+}    \\
& & \quad \cdot Y(v,(x_2 ,\varphi^+_1,\varphi^-_1)) \mp \varphi^+_2 \varphi^-_2 \frac{\partial}{\partial \varphi_1^\pm} \frac{\partial}{\partial x_2} Y(v,(x_2 ,\varphi^+_1,\varphi^-_1)) \biggr) \biggr|_{\varphi_1^\mp = 0}\\
&=& v^\pm (x_2)  +   \varphi^\mp_2 \frac{\partial}{\partial x_2}v(x_2)  \pm \varphi^\mp_2 v^{+-}(x_2) \mp \varphi^+_2 \varphi^-_2 \frac{\partial}{\partial x_2} v^\pm(x_2)\\
&=& \left( \frac{\partial}{\partial \varphi_2^\pm} + \varphi^\mp_2 \frac{\partial}{\partial x_2} \right)  Y(v,(x_2,\varphi^+_2, \varphi^-_2 )).
\end{eqnarray*}
\end{proof}

\begin{rema}\label{geometry-correspondence-remark}
{\em The $\mathcal{D}^\pm$-derivative properties (\ref{Dpm-derivative-property}) show that the operators $\mathcal{D}^\pm \in (\mathrm{End} \; V)^1$ correspond to the $N=2$ superconformal operators $D^\pm$ as defined in (\ref{define-Dpm-derivative});  see Remark \ref{superconformal-remark}. In fact this is one of the main motivations for adding the odd formal variables to the notion of vertex superalgebra.  The inclusion of odd formal variables in the notion of vertex operator superalgebra was first done in the $N=1$ case in \cite{B-announce} (see also \cite{B-thesis}--\cite{B-iso}), motivated by the supergeometry of genus-zero, holomorphic, $N=1$ superconformal field theory.  The rigorous correspondence between this geometry and $N=1$ Neveu-Schwarz vertex operator superalgebras was developed in \cite{B-announce}--\cite{B-iso}.   In Section \ref{iso-section}, we will present the notion of $N=2$ vertex superalgebra without odd formal variables and show that the category of such objects is isomorphic to the category of $N=2$ vertex superalgebras with odd formal variables.  Thus although the notion above has an equivalent counterpart without the odd formal variables, in adding the formal variables the correspondence with the underlying superconformal theory is emphasized.  We will discuss this correspondence with the geometry further in Section \ref{vosa-section}, and then in Sections \ref{nonhomo-section} and \ref{N=1-section} we show that there are other superconformal settings which will motivate other ways of formulating an equivalent notion of $N=2$ vertex superalgebra which incorporates odd formal variables.}
\end{rema}

As a consequences of Proposition \ref{Dpm-derivative-prop}, and the fact that
\begin{equation}
\left[D^+, D^-\right] = \left[ \frac{\partial}{\partial \varphi^+} + \varphi^- \frac{\partial}{\partial x}, \frac{\partial}{\partial \varphi^-} + \varphi^+ \frac{\partial}{\partial x} \right] = 2 \frac{\partial }{\partial x},
\end{equation}
we have that 
\begin{equation}\label{for-D-derivative}
Y([\mathcal{D}^+, \mathcal{D}^-]v,(x, \varphi^+, \varphi^-)) = 2 \frac{\partial}{\partial x} Y(v, (x, \varphi^+, \varphi^-)).
\end{equation}
Thus we have the following corollary to Proposition \ref{Dpm-derivative-prop}:

\begin{cor}\label{D-derivative-cor}
Let $V$ be an $N=2$ vertex superalgebra with two odd formal variables.  Let $\mathcal{D}^\pm$ be defined by (\ref{define-Dpm}), and let $\mathcal{D}$ be the even endomorphisms of $V$ defined by
\begin{equation}\label{define-D}
\mathcal{D} = \frac{1}{2} [\mathcal{D}^+, \mathcal{D}^-].
\end{equation}
Then 
\begin{equation}\label{D-derivative-property}
Y(\mathcal{D}v,(x, \varphi^+, \varphi^-)) =  \frac{\partial}{\partial x} Y(v, (x, \varphi^+, \varphi^-)).
\end{equation}
Furthermore, we have
\begin{equation}\label{D-is-what-it-should-be}
\mathcal{D}(v) = v_{-2}\mathbf{1}.
\end{equation}
\end{cor}

\begin{proof}
Since (\ref{D-derivative-property}) follows {}from (\ref{for-D-derivative}), we need only prove equation (\ref{D-is-what-it-should-be}).   {}From the $\mathcal{D}^\pm$-derivative properties (\ref{Dpm-derivative-property}), we have
\begin{equation}
(\mathcal{D}^\pm v)^\mp_{n-\frac{1}{2}} = \pm v_{n-1}^{+-} - n v_{n-1}.
\end{equation}
Thus
\begin{eqnarray*}
\mathcal{D}(v) &=&  \frac{1}{2} (\mathcal{D}^+ \mathcal{D}^- + \mathcal{D}^- \mathcal{D}^+)v  \  = \ \frac{1}{2} \Bigl( (\mathcal{D}^- v)^+_{-\frac{3}{2}} \mathbf{1} + (\mathcal{D}^+ v)^-_{-\frac{3}{2}} \mathbf{1} \Bigr)\\
&=& \frac{1}{2} \Bigl( (-v_{-2}^{+-} + v_{-2}) \mathbf{1} + (v_{-2}^{+-} + v_{-2}) \mathbf{1} \Bigr) \ = \ v_{-2} \mathbf{1}.
\end{eqnarray*}
\end{proof}

\begin{rema}\label{N2-subalgebras-representation-remark}
{\em Proposition \ref{Dpm-derivative-prop} and Corollary \ref{D-derivative-cor} show that an $N=2$ vertex superalgebra with two odd formal variables is naturally a representation of $\mathfrak{osp}_{\bigwedge_*}(2|2)_{<0}$ under the isomorphism $\mathcal{D}^\pm \mapsto G^\pm _{-1/2}$ and $\mathcal{D} \mapsto L_{-1}$; see Remark \ref{N2-subalgebras-remark}.
}
\end{rema}

Note that as a consequence of Proposition \ref{Dpm-derivative-prop} and Corollary \ref{D-derivative-cor}, we have
\begin{equation}
\mathcal{D}^\pm (\mathbf{1}) = \mathcal{D} (\mathbf{1}) = 0.
\end{equation}

{}From Proposition \ref{Taylor} and the $\mathcal{D}^\pm$- and $\mathcal{D}$-derivative properties (\ref{Dpm-derivative-property}) and (\ref{D-derivative-property}), we have 
\begin{eqnarray}
\lefteqn{Y(e^{x_0 \mathcal{D} + \varphi_0^+ \mathcal{D}^+ + \varphi_0^- \mathcal{D}^-}v, (x,\varphi^+, \varphi^-)) } \label{exponential-L(-1)-G(-1/2)} \\
&=& e^{x_0 \frac{\partial}{\partial x} + \varphi^+_0 \left(
\frac{\partial}{\partial \varphi^+} + \varphi^- \frac{\partial}{\partial x} \right)+ \varphi^-_0 \left(
\frac{\partial}{\partial \varphi^-} + \varphi^+ \frac{\partial}{\partial x} \right)} Y(v,(x,\varphi^+, \varphi^-)) \nonumber \\  
&=& Y(v,(x + x_0 + \varphi_0^+ \varphi^- + \varphi_0^- \varphi^+, \varphi_0^+ + \varphi^+, \varphi_0^- + \varphi^-)) . \nonumber
\end{eqnarray}

{}From the creation property and (\ref{exponential-L(-1)-G(-1/2)}), we have
\begin{equation}\label{for-skew-symmetry}
e^{x \mathcal{D} + \varphi^+ \mathcal{D}^+ +  \varphi^- \mathcal{D}^-} v \; = \; Y(v,(x,\varphi^+, \varphi^-))
\mathbf{1} .
\end{equation}

\begin{rema}\label{Jacobi-symmetry-remark}
{\em The left-hand side of the Jacobi identity (\ref{Jacobi-identity}) is invariant under  the transformation 
\begin{equation}
(u,v,x_0,x_1,x_2,\varphi_1^+,\varphi_1^-,\varphi_2^+, \varphi_2^-) \longleftrightarrow ((-1)^{\eta(u)
\eta(v)} v,u,- x_0,x_2,x_1,\varphi_2^+, \varphi_2^-,\varphi_1^+, \varphi_1^-) .
\end{equation}
Thus the right-hand side of the Jacobi identity must be symmetric with respect to this also.  }
\end{rema}

\begin{prop} {\bf (skew supersymmetry)}
Let $V$ be an $N=2$ vertex superalgebra with two odd formal variables.  Recall the operators $\mathcal{D}^\pm$ defined in Proposition \ref{Dpm-derivative-prop}, and let $\mathcal{D} = \frac{1}{2} [\mathcal{D}^+, \mathcal{D}^-]$ as in (\ref{define-D}).  Then
\begin{equation}\label{skew-supersymmetry}
Y(u,(x,\varphi^+, \varphi^-))v =(-1)^{\eta(u) \eta(v)}  e^{x \mathcal{D} + \varphi^+ \mathcal{D}^+ + \varphi^- \mathcal{D}^-} Y(v,(-x,-\varphi^+, - \varphi^-))u 
\end{equation}
for $u,v$ of homogeneous sign in $V$.
\end{prop}

\begin{proof}
Using the symmetry of the left-hand side of the Jacobi identity as noted in Remark \ref{Jacobi-symmetry-remark}, property (\ref{delta-substitute}), and (\ref{exponential-L(-1)-G(-1/2)}), we have
\begin{eqnarray*}
\lefteqn{x_2^{-1} \delta \biggl( \frac{x_1 - x_0 - \varphi^+_1 \varphi^-_2 - \varphi^-_1 \varphi^+_2}{x_2} \biggr) Y(Y(u,(x_0,  \varphi^+_1 - \varphi^+_2,\varphi^-_1 - \varphi^-_2))v,(x_2, \varphi^+_2,\varphi^-_2))}\\
&=& (-1)^{\eta(u) \eta(v)}  x_1^{-1} \delta \biggl( \frac{x_2 + x_0 - \varphi^+_2 \varphi^-_1 - \varphi^-_2 \varphi^+_1}{x_1} \biggr) Y(Y(v,(-x_0,  \varphi^+_2 - \varphi^+_1,\varphi^-_2 \\
& & \quad - \, \varphi^-_1))u,(x_1, \varphi^+_1,\varphi^-_1)) \\
&=& (-1)^{\eta(u) \eta(v)}  x_1^{-1} \delta \biggl( \frac{x_2 + x_0 - \varphi^+_2 \varphi^-_1 - \varphi^-_2 \varphi^+_1}{x_1} \biggr) Y(Y(v,(-x_0,  \varphi^+_2 - \varphi^+_1,\varphi^-_2 \\
& & \quad - \,  \varphi^-_1))u,(x_2 + x_0 - \varphi^+_2 \varphi^-_1 - \varphi^-_2 \varphi^+_1, \varphi^+_1,\varphi^-_1)) \\
&=& (-1)^{\eta(u) \eta(v)}  x_1^{-1} \delta \biggl( \frac{x_2 + x_0 - \varphi^+_2 \varphi^-_1 - \varphi^-_2 \varphi^+_1}{x_1} \biggr) Y(e^{x_0 \mathcal{D} + (\varphi^+_1 - \varphi_2^+) \mathcal{D}^+ +  (\varphi^-_1 - \varphi_2^-) \mathcal{D}^-} \\
& & \quad Y(v,(-x_0,  \varphi^+_2 - \varphi^+_1,\varphi^-_2  - \varphi^-_1))u,(x_2 , \varphi^+_2,\varphi^-_2)) .
\end{eqnarray*}  
Taking $\mathrm{Res}_{x_1}$ and using  (\ref{delta-2-terms-with-phis}), we have
\begin{multline*}
Y(Y(u,(x_0,  \varphi^+_1 - \varphi^+_2,\varphi^-_1 - \varphi^-_2))v,(x_2, \varphi^+_2,\varphi^-_2)) \\
= 
 (-1)^{\eta(u) \eta(v)} Y(e^{x_0 \mathcal{D} + (\varphi^+_1 - \varphi_2^+) \mathcal{D}^+ +  (\varphi^-_1 - \varphi_2^-) \mathcal{D}^-}  Y(v,(-x_0,  \varphi^+_2 \\
- \, \varphi^+_1,\varphi^-_2  - \varphi^-_1))u,
(x_2 , \varphi^+_2,\varphi^-_2)).
\end{multline*}
Then acting on $\mathbf{1}$, taking the limit as $(x_2, \varphi_2^+, \varphi_2^-)$ goes to zero, and using the creation property (\ref{creation-property2}), the result follows.
\end{proof}

We also have the following bracket formulas involving $\mathcal{D}^\pm$ and $\mathcal{D}$:
\begin{prop}\label{bracket-prop}
Let $V$ be an $N=2$ vertex superalgebra with two odd formal variables, and let $\mathcal{D}^\pm$ and $\mathcal{D}$ be defined by (\ref{define-Dpm}) and (\ref{define-D}), respectively.   Then for $v \in V$ the following $\mathcal{D}^\pm$- and $\mathcal{D}$-bracket-derivative properties hold
\begin{eqnarray}
\left[ \mathcal{D}^\pm, Y(v, (x, \varphi^+, \varphi^-)) \right] &=& \left( \frac{\partial}{\partial \varphi^\pm} - \varphi^\mp \frac{\partial}{\partial x} \right) Y(v, (x, \varphi^+, \varphi^-)) \label{Dpm-bracket-derivative}\\
\left[ \mathcal{D}, Y(v, (x, \varphi^+, \varphi^-)) \right] &=& \frac{\partial}{\partial x} Y(v, (x,\varphi^+, \varphi^-)) , \label{D-bracket-derivative}
\end{eqnarray}
and the following $\mathcal{D}^\pm$- and $\mathcal{D}$-bracket properties hold
\begin{eqnarray}
\left[ \mathcal{D}^\pm, Y(v, (x, \varphi^+, \varphi^-)) \right] &=& Y(\mathcal{D}^\pm v, (x, \varphi^+, \varphi^-)) \label{Dpm-bracket}\\
& & \quad - \, 2 \varphi^\mp Y(\mathcal{D}v, (x, \varphi^+, \varphi^-)) \nonumber\\
\left[ \mathcal{D}, Y(v, (x, \varphi^+, \varphi^-)) \right] &=& Y(\mathcal{D}v, (x, \varphi^+, \varphi^-)) \label{D-bracket}.
\end{eqnarray}
\end{prop}

\begin{proof}
Let  $T = x \mathcal{D} + \varphi^+ \mathcal{D}^+ + \varphi^- \mathcal{D}^-$.  Note that $[\mathcal{D}^\pm, \mathcal{D}] = [\mathcal{D}^\pm, \mathcal{D}^\pm] = 0$.  Thus since $(\partial/\partial x) T = \mathcal{D}$, we have
\begin{equation}\label{for-proof-1}
\frac{\partial}{\partial x} e^T = \mathcal{D} e^T + e^T \frac{\partial}{\partial x} .
\end{equation}
However since $(\partial/\partial \varphi^\pm)T = \mathcal{D}^\pm$, we have that for $n \in \mathbb{N}$,
\begin{eqnarray*}
\frac{\partial}{\partial \varphi^\pm} T^n &=& \mathcal{D}^\pm T^{n-1} + \sum_{k=1}^{n-1} \left( [T^k, \mathcal{D}^\pm] T^{n-1-k} + \mathcal{D}^\pm  T^{n-1}\right) + T^n \frac{\partial}{\partial \varphi^\pm} \\
&=&  \sum_{k=0}^{n-1} \left( 2k \varphi^\mp \mathcal{D} T^{k-1} T^{n-1-k} + \mathcal{D}^\pm  T^{n-1}\right) + T^n \frac{\partial}{\partial \varphi^\pm} \\
&=& n(n-1) \varphi^\mp \mathcal{D} T^{n-2} + n\mathcal{D}^\pm  T^{n-1} + T^n \frac{\partial}{\partial \varphi^\pm} ,
\end{eqnarray*}
and thus
\begin{equation}\label{for-proof-2}
\frac{\partial}{\partial \varphi^\pm} e^T = (\varphi^\mp \mathcal{D} + \mathcal{D}^\pm ) e^T + e^T \frac{\partial}{\partial \varphi^\pm} .
\end{equation}
Therefore for $u,v \in V$ of homogeneous sign in $V$, using skew supersymmetry (\ref{skew-supersymmetry}), equations (\ref{for-proof-1}) and (\ref{for-proof-2}), and the $\mathcal{D}^\pm$-derivative properties (\ref{Dpm-derivative-property}), we obtain
\begin{eqnarray*}
\lefteqn{\left( \frac{\partial}{\partial \varphi^\pm} - \varphi^\mp \frac{\partial}{\partial x} \right) Y(u, (x, \varphi^+, \varphi^-))v}\\
&=& (-1)^{\eta(v) \eta(u)}   \left( \frac{\partial}{\partial \varphi^\pm} - \varphi^\mp \frac{\partial}{\partial x} \right) e^{x \mathcal{D} + \varphi^+ \mathcal{D}^+ + \varphi^- \mathcal{D}^-} Y(v,(-x,-\varphi^+, - \varphi^-))u \\
&=& (-1)^{\eta(v) \eta(u)}   \biggl(\Bigl(\varphi^\mp \mathcal{D} + \mathcal{D}^\pm - \varphi^\mp \mathcal{D}  \Bigr) e^{x \mathcal{D} + \varphi^+ \mathcal{D}^+ + \varphi^- \mathcal{D}^-} Y(v,(-x,-\varphi^+, - \varphi^-))u \\
& & \quad + \, e^{x \mathcal{D} + \varphi^+ \mathcal{D}^+ + \varphi^- \mathcal{D}^-} \left( \frac{\partial}{\partial \varphi^\pm}  - \varphi^\mp \frac{\partial}{\partial x}  \right) Y(v,(-x,-\varphi^+, - \varphi^-))u  \biggr) \\
&=& (-1)^{\eta(v) \eta(u)}   \biggl( \mathcal{D}^\pm  e^{x \mathcal{D} + \varphi^+ \mathcal{D}^+ + \varphi^- \mathcal{D}^-} Y(v,(-x,-\varphi^+, - \varphi^-))u 
\end{eqnarray*}
\begin{eqnarray*}
& & \quad - \, e^{x \mathcal{D} + \varphi^+ \mathcal{D}^+ + \varphi^- \mathcal{D}^-} \left( \frac{\partial}{\partial (-\varphi^\pm)}  + (-\varphi^\mp) \frac{\partial}{\partial (-x)}  \right) Y(v,(-x,-\varphi^+, - \varphi^-))u  \biggr) \\
&=& (-1)^{\eta(v) \eta(u)}   \biggl(\mathcal{D}^\pm e^{x \mathcal{D} + \varphi^+ \mathcal{D}^+ + \varphi^- \mathcal{D}^-} Y(v,(-x,-\varphi^+, - \varphi^-))u \\
& & \quad - \, e^{x \mathcal{D} + \varphi^+ \mathcal{D}^+ + \varphi^- \mathcal{D}^-}Y(\mathcal{D}^\pm v,(-x,-\varphi^+, - \varphi^-))u  \biggr) \\
&=&  \mathcal{D}^\pm Y(u,(x,\varphi^+,  \varphi^-))v  -  (-1)^{\eta(u)}   Y(u,(x,\varphi^+,  \varphi^-)) \mathcal{D}^\pm v \\
&=& \left[ \mathcal{D}^\pm , Y(u,(x,\varphi^+,  \varphi^-)) \right]v .
\end{eqnarray*}
Using these $\mathcal{D}^\pm$-bracket properties and the $\mathcal{D}$-derivative property (\ref{D-derivative-property}), we have
\begin{eqnarray*}
\lefteqn{\frac{\partial}{\partial x} Y(u, (x, \varphi^+, \varphi^-))}\\
&=& - \frac{1}{2} \left[ \frac{\partial}{\partial \varphi^-} - \varphi^+ \frac{\partial}{\partial x} ,  \frac{\partial}{\partial \varphi^+} - \varphi^- \frac{\partial}{\partial x} \right] Y(u, (x, \varphi^+, \varphi^-))\\
&=& \frac{1}{2} \left([\mathcal{D}^+, [\mathcal{D}^-, Y(u, (x, \varphi^+, \varphi^-)) ]] + [\mathcal{D}^-, [\mathcal{D}^+, Y(u, (x, \varphi^+, \varphi^-)) ]] \right)\\
&=& [\frac{1}{2} [\mathcal{D}^+, \mathcal{D}^-], Y(u, (x, \varphi^+, \varphi^-)) ]\\
&=& [\mathcal{D}, Y(u, (x, \varphi^+, \varphi^-)) ].
\end{eqnarray*}

Finally (\ref{Dpm-bracket}) and (\ref{D-bracket}) follow {}from the $\mathcal{D}^\pm$- and $\mathcal{D}$-derivative properties (\ref{Dpm-derivative-property}) and (\ref{D-derivative-property}), respectively. 
\end{proof}

Making repeated use of the $\mathcal{D}^\pm$ -and $\mathcal{D}$-bracket properties  (\ref{Dpm-bracket}) and (\ref{D-bracket}), and using (\ref{exponential-L(-1)-G(-1/2)}), we have
\begin{eqnarray}
e^{x_0\mathcal{D}} Y(v,(x,\varphi^+, \varphi^-)) e^{- x_0 \mathcal{D}} &=& Y(e^{x_0 \mathcal{D}}v,(x,\varphi^+, \varphi^-)) \label{conjugation1} \\
&=& Y(v,(x +x_0, \varphi^+, \varphi^-)) \nonumber \\
e^{\varphi^+_0  \mathcal{D}^+} Y(v,(x,\varphi^+, \varphi^-)) e^{- \varphi_0^+ \mathcal{D}^+ } &=& Y(e^{ \varphi^+_0  \mathcal{D}^+ + 2\varphi^- \varphi_0^+  \mathcal{D}}v,(x,\varphi^+, \varphi^-)) \\
&=& Y(v,(x + \varphi^- \varphi_0^+, \varphi^+ + \varphi_0^+, \varphi^-)) \nonumber\\ 
\quad e^{\varphi^-_0  \mathcal{D}^-} Y(v,(x,\varphi^+, \varphi^-)) e^{- \varphi_0^-  \mathcal{D}^- } &=& Y(e^{ \varphi^-_0  \mathcal{D}^-  + 2\varphi^+ \varphi_0^-  \mathcal{D}}v,(x,\varphi^+, \varphi^-))\label{conjugation3} \\
&=& Y(v,(x + \varphi^+ \varphi_0^-,\varphi^+, \varphi^- + \varphi_0^-)). \nonumber
\end{eqnarray} 
By the Campbell-Baker-Hausdorff formula (cf. \cite{BHL}), we have
\begin{equation}\label{CBH}
e^{x_0  \mathcal{D}+ \varphi^+_0  \mathcal{D}^+ + \varphi^-_0  \mathcal{D}^- }
= e^{(x_0 + \varphi_0^+ \varphi_0^-)  \mathcal{D}} e^{\varphi^+_0  \mathcal{D}^+} e^{ \varphi^-_0  \mathcal{D}^-} ,
\end{equation}
and thus we have the following corollary to Proposition \ref{bracket-prop}:

\begin{cor}
Let $V$ be an $N=2$ vertex superalgebra with two odd formal variables, and let $\mathcal{D}^\pm$ and $\mathcal{D}$ be defined by 
(\ref{define-Dpm}) and (\ref{define-D}), respectively.   Then for $v \in V$
\begin{eqnarray}\label{conjugate-shift}
\lefteqn{ e^{x_0  \mathcal{D} + \varphi^+_0  \mathcal{D}^+ +  \varphi^-_0  \mathcal{D}^-} Y(v,(x,\varphi^+, \varphi^-)) e^{- x_0  \mathcal{D} - \varphi_0^+  \mathcal{D}^+ -  \varphi^-_0  \mathcal{D}^-}}  \\
&=& Y(e^{(x_0 + 2 \varphi^+ \varphi_0^- + 2 \varphi^- \varphi_0^+)  \mathcal{D} + \varphi^+_0  \mathcal{D}^+ + \varphi^-_0  \mathcal{D}^-}v,(x,\varphi^+, \varphi^-)) \nonumber \\
&=& Y(v,(x +x_0  + \varphi^+ \varphi_0^- + \varphi^- \varphi_0^+, \varphi^+ + \varphi_0^+, \varphi^- + \varphi_0^-)).  \nonumber
\end{eqnarray} 
\end{cor}

Formulas (\ref{D-bracket}) and (\ref{Dpm-bracket}) imply that if we let $Y(v,x) = Y(v,(x, \varphi^+, \varphi^-)) |_{\varphi^+ = \varphi^- = 0}$, then
\begin{equation} \label{Dpm-bracket-without}
\left[ \mathcal{D}, Y(v, x) \right] = Y (\mathcal{D}v, x), \quad \mathrm{and} \quad \left[ \mathcal{D}^\pm , Y(v, x) \right] = Y(\mathcal{D}^\pm v, x).
\end{equation}
Continuing to use this notation, we have {}from the $\mathcal{D}^\pm$-bracket formulas (\ref{Dpm-bracket-derivative})--(\ref{D-bracket}), that 
\begin{multline}\label{homo-odd-relationship1}
Y(v,(x,\varphi^+, \varphi^-)) = \sum_{n \in \mathbb{Z}} \Bigl( v_n  +
\varphi^+  [\mathcal{D}^+ , v_n ]  +  \varphi^-  [\mathcal{D}^-, v_n ]  \\
+  \frac{1}{2} \varphi^+ \varphi^- \left( [\mathcal{D}^-, [ \mathcal{D}^+, v_n ] ] - [\mathcal{D}^+, [ \mathcal{D}^-, v_n ] ] \right) \Bigr) x^{-n-1},
\end{multline}
i.e., 
\begin{eqnarray}
v^\pm_{n - 1/2} &=&  [\mathcal{D}^\pm, v_n] \label{Dvpm-bracket-condition}\\ 
v^{+-}_{n-1} &=& \frac{1}{2}\left( [\mathcal{D}^-, [ \mathcal{D}^+, v_n ] ] - [\mathcal{D}^+, [ \mathcal{D}^-, v_n ] ] \right) \label{Dv+--bracket-condition}\\
&=& [\mathcal{D}^-, [ \mathcal{D}^+, v_n ] ] - [\mathcal{D}, v_n]. \nonumber
\end{eqnarray}  
Using (\ref{Dpm-bracket-without}), this is equivalent to
\begin{equation}\label{Dvpm-condition} 
v^\pm_{n - 1/2} =  (\mathcal{D}^\pm v)_n,  \quad \mathrm{and} \quad v^{+-}_{n-1} =   \frac{1}{2} \left( (\mathcal{D}^- \mathcal{D}^+ v)_n  - (\mathcal{D}^+ \mathcal{D}^- v)_n \right),\end{equation}
i.e.,
\begin{multline}
Y(v,(x,\varphi^+, \varphi^-)) = Y(v,x) + \varphi^+ Y(\mathcal{D}^+ v,x) + \varphi^- Y(\mathcal{D}^- v,x) \\
+\frac{1}{2} \varphi^+ \varphi^- Y((\mathcal{D}^- \mathcal{D}^+ - \mathcal{D}^+ \mathcal{D}^-  )v,x) .
\end{multline}

\begin{rema}
{\em  In Section \ref{nonhomo-section} we will show that the notion of $N=2$ vertex superalgebra is, for instance, equivalent to the notion of  ``$N_K=2$ SUSY vertex algebra" in \cite{HK} and \cite{He}, see Remark \ref{nonhomo-equivalence-remark}, but is not equivalent to the notion of ``$N=2$ superconformal vertex algebra" given in \cite{Kac1997}. }
\end{rema}

Finally, we note that taking $\mbox{Res}_{x_1}$ of the supercommutator formula (\ref{homo-bracket}), we find that for $u,v \in V$ and $n \in \mathbb{Z}$, we have for instance
\begin{eqnarray}
\left[u_{0}, v_n \right] &=& (u_{0}v)_n \\
\left[u^\pm_{-\frac{1}{2}}, v^\pm_{n-\frac{1}{2}} \right] &=& (-1)^{\eta(u)+1}(u^\pm_{-\frac{1}{2}} v)^\pm_{n-\frac{1}{2}} \\
\left[u^{+-}_{-1}, v_n^{+-}\right] &=& (u^{+-}_{-1}v)_n^{+-}
\end{eqnarray}
and in particular,
\begin{eqnarray}
\left[u_{0}, v_0 \right] &=& (u_{0}v)_0 \label{form-Lie1} \\
\left[u^\pm_{-\frac{1}{2}}, v^\pm_{-\frac{1}{2}} \right] &=& (-1)^{\eta(u)+1} (u^\pm_{-\frac{1}{2}} v)^\pm_{-\frac{1}{2}} \label{form-Lie2} \\
\left[u^{+-}_{-1}, v_{-1}^{+-}\right] &=& (u^{+-}_{-1}v)_{-1}^{+-}. \label{form-Lie3}
\end{eqnarray}
Thus the operators $u_0$ form a Lie superalgebra, as do the operators $u^+_{-1/2}$, the operators $u^-_{-1/2}$, and the operators $u^{+-}_{-1}$, respectively.

\subsection{The notion of $N=2$ Neveu-Schwarz vertex operator superalgebra with two odd formal variables}\label{vosa-section} 

\begin{defn}\label{VOSA-definition}
{\em An} $N = 2$ Neveu-Schwarz vertex operator superalgebra over $\bigwedge_*$ and with two odd variables {\em is a $\frac{1}{2} \mathbb{Z}$-graded $\bigwedge_*$-module (graded by} weights{\em)
\begin{equation}\label{vosa1}
V = \coprod_{n \in \frac{1}{2} \mathbb{Z}} V_{(n)} 
\end{equation}  
such that 
\begin{equation}\label{vosa2}
\dim V_{(n)} < \infty \qquad \mbox{for  $n \in \frac{1}{2} \mathbb{Z}$,} 
\end{equation}
\begin{equation}\label{positive-energy}
V_{(n)} = 0 \qquad \mbox{for $n$ sufficiently negative},
\end{equation}
equipped with an $N=2$ vertex superalgebra structure $(V, Y(\cdot, (x, \varphi^+, \varphi^-)), \mathbf{1})$, and a distinguished homogeneous vector $\mu \in V_{(1)}^0$ (the {\em $N=2$ Neveu-Schwarz element} or {\em $N=2$ superconformal element})), satisfying the following conditions:
the $N=2$ Neveu-Schwarz algebra relations hold:
\begin{eqnarray}
\left[L(m) ,L(n) \right] \! \! \! \! &=& \! \! \! \! (m - n)L(m + n) + \frac{1}{12} (m^3 - m) \delta_{m + n, 0} \; c_V   \label{first-n2}\\
\left[J(m) , J(n) \right] \! \! \!  \! &=& \! \! \!  \! \frac{1}{3} m \delta_{m + n , 0} \;c_V \\ 
\left[L(m) , J(n) \right] \! \! \! \! &=&\!  \! \! \! -nJ(m+n) ,\\
\bigl[L(m),G^\pm(n +1/2)\bigr] \! \! \! \! &=& \! \! \!  \! \Bigl(\frac{m}{2} - n - \frac{1}{2} \Bigr) G^\pm (m + n +1/2)  \\
\bigl[ J(m) , G^\pm (n + 1/2)\bigr] \! \! \! \! &=& \! \! \! \!  \pm G^\pm (m + n + 1/2) \\
\qquad  \  \bigl[ G^\pm (m + 1/2) , G^\pm (n + 1/2) \bigr] \! \! \! \! &=& \! \! \! \! 0 \\
\bigl[ G^+ (m + 1/2) , G^- (n - 1/2) \bigr] \! \! \! \! &=& \! \! \! \! 2L(m + n) + (m-n+1) J (m+n) \label{last-n2} \\
& & \quad + \, \frac{1}{3}(m^2 + m) \delta_{m + n , 0} \;c_V, \nonumber
\end{eqnarray}

for $m,n \in \mathbb{Z}$, where 
\begin{eqnarray}
J(n) = \mu_n, \quad  \mp G^\pm(n - 1/2) = \mu_{n-\frac{1}{2}}^\pm, \quad \mbox{and} \quad -2L(n) = \mu_n^{+-} , 
\end{eqnarray}
i.e., 
\begin{multline}\label{homo-mu}
Y(\mu,(x,\varphi^+,\varphi^-)) = \sum_{n \in \mathbb{Z}} \Bigl( J(n) x^{- n - 1} -  \varphi^+ G^+ (n+1/2) x^{-n-2}  \\
+ \varphi^- G^- (n+1/2) x^{-n-2}  - 2 \varphi^+ \varphi^- L(n) x^{- n - 2} \Bigr)
\end{multline}
and $c_V \in \mathbb{C}$ (the} central charge{\em);
for $n \in \frac{1}{2} \mathbb{Z}$ and $v \in V_{(n)}$
\begin{equation}\label{grading-for-vosa-with}
L(0)v = nv
\end{equation}
and in addition, $V_{(n)}$ is the direct sum of eigenspaces for $J(0)$ such that if $v \in V_{(n)}$ is also an eigenvector for $J(0)$ with eigenvalue $k$, i.e., if  
\begin{equation}\label{J(0)-grading}
J(0)v = kv, \quad \mbox{then $k \equiv 2n \; \mathrm{mod} \; 2$};
\end{equation} 
and finally, the} $G^\pm(-1/2)$-derivative properties {\em hold:
\begin{equation}\label{G(-1/2)-derivative}
\biggl( \frac{\partial}{\partial \varphi^\pm} + \varphi^\mp \frac{\partial}{\partial x} \biggr) Y(v,(x,\varphi^+, \varphi^-)) =  Y(G^\pm(-  1/2)v,(x,\varphi^+, \varphi^-)) . 
\end{equation} }
\end{defn}

The $N=2$ Neveu-Schwarz vertex operator superalgebra just defined is denoted
by 
\[(V,Y(\cdot,(x,\varphi^+, \varphi^-)),\mathbf{1},\mu),\] 
or for simplicity by $V$.

\begin{rema}\label{mu-operator-physics-remark}
{\em  In the physics literature, when a system of fields is said to be ``$N=2$ supersymmetric", this is denoted by specifying the ``energy momentum stress tensor" which is the vertex operator corresponding to $\mu$ (\ref{homo-mu}).  However, often a slightly different basis is used for the $N=2$ Neveu-Schwarz algebra and a different normalization for the variables then we have used above.  Thus in for instance \cite{P} and \cite{YZ}, we see that letting $T_n = (1/2)J(n)$, for $n \in \mathbb{Z}$, $G_r = G^+(r)$ and $\bar{G}_r = G^-(r)$, for $r \in \frac{1}{2} + \mathbb{Z}$, $C_2 = c_V/3$, $\theta = i \sqrt{2} \varphi^-$, and $\bar{\theta} = - i \sqrt{2} \varphi^+$, then the vertex operator corresponding to the $N=2$ superconformal element $\mu$ given by equation (5) in \cite{P} and equation (2.13) in \cite{YZ} is exactly (\ref{homo-mu}).  Moreover, using these substitutions, the ``operator product expansion" which we give below in (\ref{OPE}) for the vertex operator corresponding to $\mu$ is exactly that given in \cite{YZ}.  This operator product expansion is a way of encoding the $N=2$ Neveu-Schwarz algebra relations (\ref{first-n2})--(\ref{last-n2}). In addition to the commonly used notation mentioned above, one often encounters the switching of the notation for the variables corresponding to $\varphi^\pm \mapsto -\varphi^\mp$ as in \cite{MSS}.  Such changes in notational convention do not correspond to an $N=2$ superconformal change of variables in the spirit of \cite{B-change} and do not arise {}from the automorphisms of the $N=2$ Neveu-Schwarz algebra as discussed in Remark \ref{auto-remark2} below.   Our choice of notation has more to do with the aesthetics arising {}from the mathematics rather than with the mathematics itself.  
}
\end{rema}

\begin{rema}\label{huang-milas-remark1}  
{\em Our formulation of the notion of $N=2$ Neveu-Schwarz vertex operator superalgebra with two odd formal variables is equivalent to that given in \cite{HM} under the relabeling of coordinates such that  their $\varphi^\pm$ is our $i\varphi^\pm$ and  their $G^\pm(n-1/2)$ is our $\pm iG^\pm(n-1/2)$, and as long as one judiciously interprets  the transformation of coordinates {}from nonhomogeneous to homogeneous (as is discussed below in Section \ref{nonhomo-section}).  However, in \cite{HM} the equation for how the vertex operators with odd variables are constructed {}from vertex operators without odd variables is only given in the nonhomogeneous coordinates and has a sign typo.  The vertex operator corresponding to $\mu$ has a couple of sign typos as well; see Remark \ref{huang-milas-remark2}.  Finally we note that in the coordinates used in \cite{HM}, the $G^\pm(-1/2)$-derivative properties are then given by $Y(G^\pm (-1/2)v, (x, \varphi^+, \varphi^-)) = \mp \Bigl( \frac{\partial}{\partial \varphi^\pm} - \varphi^\mp \frac{\partial}{\partial x}  \Bigr)Y(v, (x, \varphi^+, \varphi^-))$, cf. equation (\ref{G(-1/2)-derivative}).}
\end{rema}

For $V$ an $N=2$ Neveu-Schwarz vertex operator superalgebra over $\bigwedge_*$ for $\bigwedge_* = \bigwedge_0 = \mathbb{C}$, the $\mathbb{Z}_2$-grading is given by
\[ V^0 = \coprod_{n\in \mathbb{Z}} V_{(n)}, \qquad \mathrm{and}  \qquad V^1 =  \coprod_{n\in \mathbb{Z} + \frac{1}{2}} V_{(n)} .\]
Extending $V$ to a vertex operator superalgebra over a general $\bigwedge_*$ by $\bigwedge_* \otimes V$ changes the  $\mathbb{Z}_2$-grading via $(\bigwedge_* \otimes V)^0 = \bigwedge_*^0 
\otimes V^0 + \bigwedge_*^1 \otimes V^1$ and  $(\bigwedge_*  \otimes V)^1 = \bigwedge_*^1 \otimes V^0 + \bigwedge_*^0 \otimes V^1$.

We have the following consequences of the definition of $N=2$ Neveu-Schwarz vertex operator superalgebra with two odd formal variables:
By the creation property, (\ref{homo-mu}) and the fact that $G^\pm (n + \frac{1}{2})$ for $n \geq -1$ generate $L(n)$ for $n \geq -1$ and $J(n)$ for $n \geq 0$, we have 
\begin{equation}\label{kill-vacuum}
L(n) \mathbf{1} = J(n+1) \mathbf{1} =  G^\pm(n + 1/2) \mathbf{1} = 0, \quad \mbox{for
$n \geq -1$}.
\end{equation}
By the creation property, we have
\begin{equation}\label{mu-vacuum}
J(-1) \mbox{\bf 1} = \mu_{-1} \mathbf{1} = \mu .
\end{equation}

{}From (\ref{homo-mu}) and the $G^\pm(-1/2)$-derivative properties (\ref{G(-1/2)-derivative}), we have that there exist two vectors 
\begin{equation}
\tau^{(\pm)} = G^\pm(-3/2) \mathbf{1} = \mp G^\pm (-1/2) \mu \in V_{(3/2)}^1
\end{equation} 
such that
\begin{multline}\label{tau-pm-with}
Y(\tau^{(\pm)},(x,\varphi^+,\varphi^-)) =  \sum_{n \in \mathbb{Z}} \Bigl( G^\pm (n+ 1/2) x^{- n - 2} + 2\varphi^\mp L(n) x^{-n-2}   \\ 
\pm (n+1)\varphi^\mp J(n) x^{-n-2}    \pm (n+2)  \varphi^+ \varphi^- G^\pm(n+ 1/2) x^{- n - 3} \Bigr)
\end{multline}
as well as a third vector 
\begin{equation}
\omega = L(-2) \mathbf{1} = \frac{1}{4}(G^-(-1/2)\tau^{(+)} + G^+(-1/2) \tau^{(-)}) \in V_{(2)}^0
\end{equation} 
with 
\begin{multline}
Y(\omega,(x,\varphi^+,\varphi^-)) =  \sum_{n \in \mathbb{Z}} \Bigl( L(n) x^{- n - 2} - \frac{1}{2}(n+2)  \varphi^+ G^+(n+ 1/2) x^{-n-3}  \\ 
- \frac{1}{2} (n+2) \varphi^- G^-(n+ 1/2 ) x^{-n-3}    - \frac{1}{2}(n+1)(n+2) \varphi^+ \varphi^- J(n) x^{- n - 3} \Bigr).
\end{multline}

Using the $N=2$ Neveu-Schwarz algebra supercommutation relations (\ref{first-n2})--(\ref{last-n2})  as well as (\ref{kill-vacuum}) and (\ref{mu-vacuum}),  we see that  
\begin{eqnarray}
\lefteqn{\quad Y(\mu, (x,\varphi^+, \varphi^-))\mu } \label{mu-singular}\\
&=& \! \! \! \! J(1) \mu x^{-2} - \varphi^+ G^+(- 1/2) \mu  x^{-1} + \varphi^- G^-(- 1/2) \mu x^{-1}  - 2\varphi^+ \varphi^- L(0) \mu x^{-2} \nonumber\\
& & \quad - \, 2\varphi^+ \varphi^-L(-1) \mu x^{-1} + y (x, \varphi^+, \varphi^-) \nonumber\\
\qquad &=& \! \! \! \! \frac{1}{3}c_V \mathbf{1} x^{-2}  + \varphi^+ G^+(- 3/2) \mathbf{1} x^{-1} + \varphi^- G^-(-3/2) \mathbf{1} x^{-1}  - 2\varphi^+ \varphi^- J(-1) \mathbf{1} x^{-2} \nonumber \\
& & \quad - \, 2\varphi^+ \varphi^- J(-2) \mathbf{1} x^{-1} + y(x, \varphi^+, \varphi^-) \nonumber \\
&=& \! \! \! \! \frac{1}{3}c_V \mathbf{1} x^{-2} + \varphi^+ \tau^{(+)} x^{-1} + \varphi^- \tau^{(-)} x^{-1}  - 2\varphi^+ \varphi^- \mu x^{-2} \nonumber\\
& & \quad - \, 2\varphi^+ \varphi^-L(-1) \mu x^{-1} + y (x, \varphi^+, \varphi^-) \nonumber
\end{eqnarray}
where $y(x, \varphi^+, \varphi^-) \in V[[x]][\varphi^+, \varphi^-]$.  

Using the notation and results of Section \ref{vertex-superalgebra-section}, we observe that the $G^\pm(-1/2)$-derivative properties (\ref{G(-1/2)-derivative}), the creation property, Proposition \ref{Dpm-derivative-prop}, and Corollary \ref{D-derivative-cor} imply that 
\begin{equation}\label{D-correspondence}
G^\pm(-1/2) = \mathcal{D}^\pm, \qquad \mbox{and} \qquad L(-1) = \mathcal{D},
\end{equation}
and thus the identities for $\mathcal{D}^\pm$ and $\mathcal{D}$ of Section \ref{vertex-superalgebra-section} apply using $G^\pm(-1/2)$ and $L(-1)$, respectively.  That is, we now have the $L(-1)$-derivative property (\ref{D-derivative-property}), skew supersymmetry (\ref{skew-supersymmetry}),  the  $G^\pm(-1/2)$- and $L(-1)$-bracket-derivative properties (\ref{Dpm-bracket-derivative}) and (\ref{D-bracket-derivative}), and the $G^\pm(-1/2)$- and $L(-1)$-bracket properties (\ref{Dpm-bracket}) and (\ref{D-bracket}) with $\mathcal{D}^\pm= G^\pm(-1/2)$ and $\mathcal{D} = L(-1)$.   

Using (\ref{kill-vacuum}), we have the following special cases of identities using the supercommutator formula (\ref{homo-bracket}) with $u = \omega$, $u = \mu$, or $u = \tau^{(\pm)}$, which we will find useful later: 
\begin{eqnarray} 
\left[ L(0), Y(v, (x, \varphi^+, \varphi^-)) \right]  \! \! \!  \! &= &  \! \! \!  Y\! \biggl( \! \Bigl( L(0)+ \frac{1}{2} \varphi^+ G^+(-1/2)  + \frac{1}{2} \varphi^- G^-(-1/2) \label{L0pm-bracket} \\
& & \quad +   xL(-1)\Bigr) v, (x, \varphi^+, \varphi^-) \! \biggr) \nonumber \\
\left[ J(0), Y(v, (x, \varphi^+, \varphi^-)) \right]  \! \! \!  \! &=&  \! \! \!  Y\left(\left(J(0)+ \varphi^+G^+(-1/2) -  \varphi^-G^-(-1/2)  \right. \right.  \label{J0pm-bracket}  \\
& & \quad \left. \left. - 2 \varphi^+ \varphi^- L(-1)\right)v,(x, \varphi^+, \varphi^- )\right) \nonumber\\
\qquad  \ \  \left[ G^\pm(1/2), Y(v, (x, \varphi^+, \varphi^-)) \right]  \! \! \! \!  &=&  \! \! \!  Y\left(\left(G^\pm(1/2)- \varphi^\mp (2 L(0) \pm J(0))  \right. \right. \label{G-1/2-bracket}\\
& & \quad \pm \varphi^+ \varphi^- G^\pm (-1/2) +  x G^\pm (-1/2) \nonumber \\
& & \quad  \left. \left.  -  2 x \varphi^\mp L(-1) \right) v, (x, \varphi^+, \varphi^-)\right).\nonumber
\end{eqnarray}
And of course any other brackets involving $N=2$ Neveu-Schwarz algebra elements can be computed similarly using the supercommutator formula.

\begin{rema}\label{Lie-algebras-remark}
{\em Recall that by (\ref{form-Lie1}), the operators $u_0$ form a Lie superalgebra; denote this Lie superalgebra by $\mathcal{L}_0$.  By (\ref{form-Lie2}), the operators $u^+_{-1/2}$ and  $u^-_{-1/2}$ also form two Lie superalgebras; denote these by $\mathcal{L}^+_{-1/2}$ and $\mathcal{L}^-_{-1/2}$, respectively.  Finally, by (\ref{form-Lie3}), the operators $u^{+-}_{-1}$ form a Lie superalgebra; denote this Lie superalgebra by $\mathcal{L}^{+-}_{-1}$.   Since $\tau^{(\pm)}_0 = G^\pm(-1/2) = \mathcal{D}^\pm$, by (\ref{Dvpm-bracket-condition}) and (\ref{Dv+--bracket-condition}), we have that $\mathcal{L}^\pm_{-1/2}$ and $\mathcal{L}_{-1}^{+-}$ are subalgebras of  $\mathcal{L}_0$.  }
\end{rema}

We will use the notation $\mathrm{wt} \ v = n$ for $v\in V$ satisfying $L(0)v = nv$, and $\mathrm{wt}^J \ v = n$ for $v \in V$ satisfying $J(0)v = nv$.   Note that by (\ref{kill-vacuum}), we have
\begin{equation}\label{weight-of-vacuum}
\mathrm{wt} \ \mathbf{1} = \mathrm{wt}^J \ \mathbf{1} = 0,
\end{equation}
and using the $N=2$ Neveu-Schwarz relations along with (\ref{mu-vacuum}), we have
\begin{equation}\label{weight-of-mu}
\mathrm{wt}^J \ \mu=  0.
\end{equation}

{}From (\ref{L0pm-bracket}), and using the $L(-1)$- and $G^\pm(-1/2)$-derivative properties we obtain the following:
\begin{eqnarray}
\mbox{wt} \; v_n  =  \mbox{wt} \; v_n^{+-} = \mbox{wt} \; v - n - 1 , \label{L(0)-grading} \\
\mbox{wt} \; v_{n- \frac{1}{2}}^\pm = \mbox{wt} \; v - \Bigl(n - \frac{1}{2}\Bigr) - 1 , 
\end{eqnarray}
for $n \in \mathbb{Z}$ and for $v \in V$ of homogeneous weight.  {}From (\ref{J0pm-bracket}), and using the $L(-1)$- and $G^\pm(-1/2)$-derivative properties we obtain the following:
\begin{eqnarray}
\mathrm{wt}^J \ v_n = \mathrm{wt}^J \ v^{+-}_n = \mathrm{wt}^J \ v  \label{J(0)-grading1} \\
\mathrm{wt}^J \ v_{n-\frac{1}{2}}^\pm = \mathrm{wt}^J \ v \pm 1 \label{J(0)-grading2}
\end{eqnarray}
for $n \in \mathbb{Z}$ and for $v \in V$ of homogeneous weight with respect to the $J(0)$ grading.

For $v$ of homogeneous $L(0)$-weight in $V$, let $x_0^{2L(0)} v = x_0^{2\mathrm{wt} \, v}$, and for $v$ of homogeneous $J(0)$-weight, let $x_0^{J(0)} v = x_0^{\mathrm{wt}^J \, v}$.  As a consequence of the above $L(0)$- and $J(0)$-grading properties (\ref{L(0)-grading}) - (\ref{J(0)-grading2}), we have the following  
\begin{eqnarray} 
x_0^{2L(0)} Y(v, (x, \varphi^+, \varphi^-)) x_0^{-2 L(0)} &=& Y(x_0^{2L(0)}v, (x_0^2 x,
x_0 \varphi^+, x_0 \varphi^-)) \label{conjugate-by-L(0)} \\ 
x_0^{J(0)} Y(v, (x, \varphi^+, \varphi^-)) x_0^{-J(0)} &=& Y(x_0^{J(0)}v, (x,
x_0 \varphi^+, x_0^{-1} \varphi^-)) \label{conjugate-by-J(0)}.
\end{eqnarray}

Let $(V_1, Y_1(\cdot,(x,\varphi^+, \varphi^-)),\mbox{\bf 1}_1,\mu_1)$ and $(V_2, Y_2(\cdot,(x,\varphi^+, \varphi^-)),\mbox{\bf 1}_2,\mu_2)$ be two $N=2$ Neveu-Schwarz vertex operator superalgebras with two odd formal variables.  A {\it homomorphism of $N=2$ Neveu-Schwarz vertex operator superalgebras with two odd formal variables} is an $N=2$ vertex superalgebra homomorphism $\gamma : V_1 \longrightarrow V_2 \;$ such that $\gamma(\mu_1) = \mu_2$.  

We note that this last condition implies that $\gamma$ is grade-preserving for the three gradings: sign, $L(0)$-weight,  and $J(0)$-weight.  That is $\gamma : (V_1)_{(n)}^j \longrightarrow  (V_2)_{(n)}^j$ for $n \in \frac{1}{2} \mathbb{Z}$, and $j \in \mathbb{Z}_2$ and if $J(0)v = kv$ for $v \in V_1$, and $k \in \mathbb{Z}$, then $J(0) \gamma(v) = k\gamma(v)$.   In addition, $V_1$ and $V_2$ must have the same central charge. 

For $a,b \in (\bigwedge_*^0)^\times$ let  
\begin{eqnarray}\label{define-gamma}
\gamma_{(a,b)} : V &\longrightarrow& V \\
v & \mapsto& a^{2L(0)} b^{J(0)} v.\nonumber
\end{eqnarray}
Then $\gamma_{(a,b)}$ is a bijection, with $\gamma_{(a,b)}^{-1} = \gamma_{(a^{-1}, b^{-1})}$.  Note that by (\ref{weight-of-vacuum}) and (\ref{weight-of-mu}), we have that $\gamma_{(a,b)} \mathbf{1} = \mathbf{1}$ and $\gamma_{(a,b)} \mu = a^2 \mu$. {}From (\ref{conjugate-by-L(0)}) and (\ref{conjugate-by-J(0)}), we have that 
\begin{equation}\label{gamma-iso}
\gamma_{(a,b)} \circ Y(v, (x, \varphi^+, \varphi^-)) \circ \gamma_{(a,b)}^{-1} = Y(\gamma_{(a,b)}(v), (a^2x, ab \varphi^+, ab^{-1} \varphi^-)).
\end{equation}
Thus we have the following lemma:

\begin{lem}
Let $(V, Y(\cdot, (x, \varphi^+, \varphi^-)), \mathbf{1}, \mu)$ be an $N=2$ Neveu-Schwarz vertex operator superalgebra.  Then $(V, Y(\cdot, (a^2x, ab \varphi^+, ab^{-1} \varphi^-)), \mathbf{1}, a^2\mu)$ is isomorphic to  $(V, Y(\cdot, (x, \varphi^+, \varphi^-)), \mathbf{1}, \mu)$, for $a,b \in (\bigwedge_*^0)^\times$ .
\end{lem}

\begin{rema}\label{auto-remark2} 
{\em Given an $N=2$ Neveu-Schwarz vertex operator superalgebra $(V, Y(\cdot,$ $(x, \varphi^+, \varphi^-)), \mathbf{1}, \mu)$, the automorphisms of the $N=2$ Neveu-Schwarz algebra (\ref{auto1})--(\ref{auto3}), give rise to the following $N=2$ Neveu-Schwarz vertex operator superalgebras which are isomorphic to $(V, Y(\cdot, (x, \varphi^+, \varphi^-)), \mathbf{1}, \mu)$:
We have the following family of $N=2$ Neveu-Schwarz vertex operator superalgebras which have the same $N=2$ Neveu-Schwarz algebra element $\mu$
\begin{equation}
(V, Y(\cdot,(x, \varphi^+, \varphi^-)), \mathbf{1}, \mu) \cong (V, Y(\cdot, (x, b\varphi^+, b^{-1} \varphi^-)), \mathbf{1}, \mu) 
\end{equation}
for $b \in (\bigwedge_*^0)^\times$ and correspond to the isomorphisms $\gamma_{(1,b)}$ defined in (\ref{define-gamma}).  In addition, we have the isomorphism
\begin{equation}
(V, Y(\cdot,(x, \varphi^+, \varphi^-)), \mathbf{1}, \mu) \cong (V, Y(\cdot, (x, \varphi^-, \varphi^+)), \mathbf{1}, -\mu)
\end{equation}
which does not preserve the $N=2$ Neveu-Schwarz algebra element but does preserve the ``Virasoro element" $\omega = L(-2) \mathbf{1}$.  And of course, the composition of these two isomorphisms gives the following continuous family of isomorphic $N=2$ Neveu-Schwarz vertex operator superalgebras which preserve $\omega$ but not $\mu$
\begin{equation}
(V, Y(\cdot, (x, \varphi^+, \varphi^-)), \mathbf{1}, \mu) \cong (V, Y(\cdot, (x, b^{-1}\varphi^-, b\varphi^+)), \mathbf{1}, -\mu)
\end{equation}
for $b \in (\bigwedge_*^0)^\times$.  Cf. Remark 3.2 of \cite{B-vosas}. }
\end{rema}

We have formulated the notion of $N=2$ Neveu-Schwarz vertex operator superalgebra with two odd formal variables so as to emphasis the relationship to the underlying supergeometry of propagating worldsheets in $N=2$ superconformal field theory (cf.  Remark \ref{geometry-correspondence-remark}).  In \cite{B-moduli}, we define the moduli space of $N=2$ super-Riemann spheres with tubes.   In analogy to the author's work on the geometric interpretation of $N=1$ Neveu-Schwarz vertex operator superalgebras \cite{B-announce}--\cite{B-iso}, we plan in subsequent work to show that the there is a sewing operator on the moduli space of $N=2$ super-Riemann spheres with tubes and that this along with the action of the symmetric groups defined in \cite{B-moduli} gives the moduli space the structure of a partial-pseudo operad \cite{operads}, in analogy to the non-super case \cite{HL2}, \cite{HL3} and the $N=1$ super case  \cite{B-memoir}, \cite{B-iso}.  Then we will show that the category of $N=2$ Neveu-Schwarz vertex operator superalgebras is isomorphic to the category of algebras over this partial-pseudo operad.   One can already see some aspects of the correspondence with the geometry of $N=2$ super-Riemann spheres with tubes in that such objects are topological superspaces whose transition functions are $N=2$ superconformal, meaning they transform the superdifferential operators $D^\pm$ defined by (\ref{define-Dpm-derivative}) homogeneously of degree one; see Remark \ref{superconformal-remark}.  As noted in Remark \ref{geometry-correspondence-remark} the $\mathcal{D}^\pm = G^\pm (-1/2)$-derivative properties in the notion of $N=2$ (Neveu-Schwarz) vertex (operator) superalgebra with two odd formal variables show some of the correspondence with the notion of $N=2$ superconformality.  In addition, here we note that in \cite{B-moduli} we prove the following proposition:

\begin{prop}\label{moduli-prop}(\cite{B-moduli})
Let $\mathbb{Z}_+$ denote the positive integers.  There is a bijection between formal $N=2$ superconformal functions vanishing at zero and invertible in a neighborhood of zero and expressions of the form
\begin{multline}\label{change-of-variables}
\exp \Biggl(-  \sum_{n \in \Z} \Bigl( A^{+}_n L_n(x,\varphi^+, \varphi^-) + A^{-}_n J_n(x,\varphi^+, \varphi^-) +  M^{+}_{n - \frac{1}{2}} G^+_{n -\frac{1}{2}} (x, \varphi^+, \varphi^-) \\
+  M^{-}_{n - \frac{1}{2}} G^-_{n -\frac{1}{2}} (x, \varphi^+, \varphi^-) \Bigr) \Biggr)  \cdot  (a_0^{+})^{-2L_0(x,\varphi^+, \varphi^-)} \cdot  (a_0^{-})^{-J_0(x,\varphi^+, \varphi^-)} \cdot (x, \varphi^+, \varphi^-) 
\end{multline}
for $(a_0^+, a_0^-) \in ((\bigwedge_*^0)^\times)^2/\langle (-1,-1) \rangle$, and $A_n^\pm \in \bigwedge_*^0$ and $M^\pm_{n- 1/2} \in \bigwedge_*^1$, for $n \in \Z$, 
where 
\begin{eqnarray}
L_n(x,\varphi^+,\varphi^-) &=& - \biggl( x^{n + 1} \frac{\partial}{\partial x} + \Bigl(\frac{n + 1}{2}\Bigr) x^n \Bigl( \varphi^+ \frac{\partial}{\partial \varphi^+} + \varphi^- \frac{\partial}{\partial \varphi^-}\Bigr) \biggr) \label{L-notation}\\
J_n(x,\varphi^+,\varphi^-) &=& - x^n\Bigl(\varphi^+\frac{\partial}{\partial \varphi^+} - \varphi^- \frac{\partial}{\partial \varphi^-}\Bigr)  \label{J-notation}\\
\quad G^\pm_{n -\frac{1}{2}} (x,\varphi^+,\varphi^-) &=& - \biggl( x^n \Bigl( \frac{\partial}{\partial \varphi^\pm} - \varphi^\mp \frac{\partial}{\partial x}\Bigr) \pm nx^{n-1} \varphi^+ \varphi^- \frac{\partial}{\partial \varphi^\pm} \biggr) \label{G-notation}
\end{eqnarray}
are superderivations in  $\mbox{Der} (\bigwedge_*[[x,x^{-1}]] [\varphi^+,\varphi^-])$, for $n \in \mathbb{Z}$, which give a representation of the $N=2$ Neveu-Schwarz algebra with central charge zero; that is  (\ref{L-notation})--(\ref{G-notation}) satisfy (\ref{Virasoro-relation})--(\ref{Neveu-Schwarz-relation-last}) with $d = 0$.  Similarly, there is a bijection between formal $N=2$ superconformal functions vanishing at $(\infty, 0,0)$ and invertible in a neighborhood of $(\infty, 0,0)$ and expressions of the form
\begin{multline}\label{change-of-variables-infty}
\exp \Biggl(\sum_{n \in \Z} \Bigl( A^{+}_n L_{-n}(x,\varphi^+, \varphi^-) - A^{-}_n J_{-n}(x,\varphi^+, \varphi^-) +  iM^{+}_{n - \frac{1}{2}} G^+_{-n +\frac{1}{2}} (x, \varphi^+, \varphi^-) \\
+ \, iM^{-}_{n - \frac{1}{2}} G^-_{-n +\frac{1}{2}} (x, \varphi^+, \varphi^-) \Bigr) \Biggr)  \cdot  (a_0^{+})^{2L_0(x,\varphi^+, \varphi^-)} \cdot  (a_0^{-})^{-J_0(x,\varphi^+, \varphi^-)} \cdot \left(\frac{1}{x}, \frac{i\varphi^+}{x},\frac{i \varphi^-}{x} \right) .
\end{multline}
\end{prop}

In light of this proposition, properties such as (\ref{gamma-iso}) can be rewritten as 
\begin{multline}\label{LJ-change-of-variables}
a^{2L(0)} b^{J(0)}  Y(v, (x, \varphi^+, \varphi^-)) a^{-2L(0)} b^{-J(0)}  \\
= Y(a^{2L(0)} b^{J(0)} v, a^{-2L_0(x, \varphi^+, \varphi^-)} b^{-J_0(x, \varphi^+, \varphi^-)} (x,  \varphi^+,  \varphi^-)),
\end{multline}
where using (\ref{L-notation}) and (\ref{J-notation}), we have the following $N=2$ superconformal change of variables vanishing at zero and invertible in a neighborhood of zero
\begin{eqnarray}
\lefteqn{a^{-2L_0(x, \varphi^+, \varphi^-)} b^{-J_0(x, \varphi^+, \varphi^-)} (x,  \varphi^+,  \varphi^-)}\\
&=&  a^{\left(2x\frac{\partial}{\partial x} + \varphi^+ \frac{\partial}{\partial \varphi^+} + \frac{1}{2} \varphi^- \frac{\partial}{\partial \varphi^-}\right)} b^{\left( \varphi^+ \frac{\partial}{\partial \varphi^+} - \varphi^- \frac{\partial}{\partial \varphi^-}\right)} (x,  \varphi^+, \varphi^-) \nonumber\\
&=& (a^2 x, ab \varphi^+, a b^{-1} \varphi^-). \nonumber
\end{eqnarray}
Another example of change of variables interpretation is that (\ref{conjugate-shift}) can be expressed as
\begin{eqnarray}\label{shift-change-of-variables}
\lefteqn{ \qquad e^{x_0 L(-1) + \varphi^+_0  G^+(-1/2) +  \varphi^-_0  G^-(-1/2)} Y(v,(x,\varphi^+, \varphi^-)) }\\
& & \quad \cdot e^{- x_0 L(-1)- \varphi_0^+  G^+(-1/2) -  \varphi^-_0  G^-(-1/2)}\nonumber  \\
&=& Y(v,e^{-x_0 L_{-1}(x, \varphi^+, \varphi^-)  - \varphi^+_0  G^+_{-1/2} (x, \varphi^+, \varphi^-)  -  \varphi^-_0  G^-_{-1/2}(x, \varphi^+, \varphi^-)}  (x , \varphi^+, \varphi^-)). \nonumber
\end{eqnarray} 

\begin{rema}\label{change-of-variables-remark}{\em 
We see {}from the discussion above that one can think of (\ref{LJ-change-of-variables}) and (\ref{shift-change-of-variables}) as a change of variables formulas related to certain changes of variables of the form (\ref{change-of-variables}) or (\ref{change-of-variables-infty}), respectively.  This leads to the question of what are the change of variables formulas for a general change of variables of the form (\ref{change-of-variables}), or a change of variables that vanishes at infinity instead of zero (cf. \cite{B-moduli}).  In the nonsuper case  in \cite{H-book}, change of variables formulas for a general conformal change of variables for a vertex operator algebra were given and certain convergences were also proved.  These formulas were derived {}from the sewing together of spheres with tubes arising {}from the worldsheet geometry underlying conformal field theory and the correspondence between the geometry of spheres with tubes and a sewing operation and the notion of vertex operator algebra.  In \cite{B-change}, the change of variables formulas for a general $N=1$ superconformal change of variables for an $N=1$ Neveu-Schwarz vertex operator superalgebra were given and the convergence of the correlation functions that arise was proved given the convergence of the change of variables.  This work used the sewing together of $N=1$ superspheres with tubes arising {}from the worldsheet geometry underlying $N=1$ superconformal field theory and the correspondence between the geometry of $N=1$ superspheres with tubes and a sewing operation and the notion of $N=1$ Neveu-Schwarz vertex operator superalgebra as developed in \cite{B-announce}--\cite{B-iso}.   In \cite{He}, Heluani gives change of variables formulas for a certain subset of  superconformal change of variables vanishing at zero for $N=1,2$ assuming the version of Proposition \ref{moduli-prop} in the nonhomogeneous coordinate system, Proposition \ref{moduli-prop-nonhomo} below.   As mentioned above, the change of variables formulas had already been proven by the author for the $N=1$ case in \cite{B-change} along with more general change of variables formulas and in addition questions of convergence of the resulting correlations functions were addressed.   This was achieved not only by using the underlying supergeometry as developed in \cite{B-announce}--\cite{B-iso} as motivation, but by using the actual results of this work. 
}
\end{rema}

Recalling (\ref{mu-singular}) we have the following proposition which is useful in proving whether an $N=2$ vertex superalgebra is an $N=2$ Neveu-Schwarz vertex operator superalgebra.   

\begin{prop}\label{mu-prop}
Let $(V, Y(\cdot, (x, \varphi^+, \varphi^-)), \mathbf{1})$ be an $N=2$ vertex superalgebra and let 
$\mathcal{D}^\pm$ and $\mathcal{D}$ be defined by (\ref{define-Dpm}) and (\ref{define-D}), respectively.  If $\mu \in V$, satisfies
\begin{multline}\label{for-mu-prop}
Y(\mu, (x,\varphi^+, \varphi^-))\mu =  \frac{1}{3}c \mathbf{1} x^{-2}  - \varphi^+ \mathcal{D}^+ \mu x^{-1} + \varphi^- \mathcal{D}^-  \mu x^{-1} \\
 - 2\varphi^+ \varphi^- \mu x^{-2} - 2\varphi^+ \varphi^-\mathcal{D} \mu x^{-1} + y(x, \varphi^+, \varphi^-) 
\end{multline}
for some $c \in \mathbb{C}$ and $y \in V[[x]][\varphi^+, \varphi^-]$, that is if
\begin{eqnarray}
\mu_1 \mu &=& \frac{1}{3} c \mathbf{1} \\
\mu_{-\frac{1}{2}}^\pm \mu &=& \mp \mathcal{D}^\pm  \mu \ = \ \mp \mu_{-\frac{3}{2}}^\pm \mathbf{1} \\
\mu_0^{+-} \mu &=& -2 \mu \\
\mu_{-1}^{+-} \mu  &=& -2 \mathcal{D} \mu  \ = \ -2\mu_{-2}\mathbf{1}
\end{eqnarray}
and
\begin{equation}
\mu_{n-1} \mu = \mu_{n-\frac{1}{2}}^\pm \mu = \mu_n^{+-} \mu =  0 \qquad \mbox{for $n\geq 1$},
\end{equation}
then defining $L(n), J(n), G^\pm(n - 1/2) \in \mathrm{End} \, V$, for $n \in \mathbb{Z}$, by $\mu_n = J(n)$, $\mu_{n-1/2}^\pm = \mp G^\pm(n - 1/2)$, and $\mu_n^{+-} = -2L(n)$, we have that $L(n), J(n)$, and $G^\pm(n-1/2)$ satisfy the relations for the $N=2$ Neveu-Schwarz algebra (\ref{first-n2})--(\ref{last-n2}) with central charge $c$. 
\end{prop}
 
\begin{proof}
By Proposition \ref{Dpm-derivative-prop} and Corollary \ref{D-derivative-cor}, we have that the $\mathcal{D}^\pm$- and $\mathcal{D}$-derivative properties (\ref{Dpm-derivative-property}) and (\ref{D-derivative-property}) hold.   Using these derivative properties, the supercommutator formula (\ref{homo-bracket}), and the $\delta$-function identity (\ref{delta-2-terms-with-phis}), we have that if $\mu \in V$ satisfies (\ref{for-mu-prop}), then for $k \in \mathbb{Z}$,
\begin{eqnarray}
\lefteqn{\qquad \mathrm{Res}_{x_1} x_1^k  \left[Y(\mu, (x_1, \varphi_1^+, \varphi_1^-)), Y(\mu, (x_2, \varphi_2^+, \varphi_2^-)) \right] }  \label{for-brackets}\\
&=&  \mathrm{Res}_{x_1} x_1^k  \mathrm{Res}_{x_0} x_2^{-1} \delta \biggl( \frac{x_1 - x_0 - \varphi_1^+ \varphi_2^- - \varphi_1^- \varphi_2^+}{x_2} \biggr)  \nonumber\\
& & \quad \cdot Y(Y(\mu,(x_0, \varphi_1^+ - \varphi_2^+, \varphi_1^- - \varphi_2^-))\mu,(x_2, \varphi_2^+, \varphi_2^-))\nonumber \\
&=& \mathrm{Res}_{x_1} x_1^k  \mathrm{Res}_{x_0}  x_1^{-1} \delta \biggl( \frac{x_2 + x_0 + \varphi_1^+ \varphi_2^- + \varphi_1^- \varphi_2^+}{x_1} \biggr) Y\Bigl(\frac{1}{3}c \mathbf{1} x_0^{-2}\nonumber \\
& & \quad  - \,  (\varphi^+_1 - \varphi_2^+) \mathcal{D}^+ \mu x_0^{-1} + (\varphi^-_1 - \varphi_2^-) \mathcal{D}^-  \mu x_0^{-1} - 2(\varphi^+_1 - \varphi_2^+) (\varphi^-_1 \nonumber\\
& & \quad - \, \varphi_2^-) \mu x_0^{-2}  - 2(\varphi^+_1 - \varphi_2^+) (\varphi^-_1 - \varphi_2^-)\mathcal{D} \mu x_0^{-1} + y(x_0, \varphi^+_1 - \varphi_2^+, \varphi^-_1 \nonumber\\
& & \quad - \, \varphi_2^-) ,(x_2, \varphi_2^+, \varphi_2^-)\Bigr)\nonumber  \\
&=& \mathrm{Res}_{x_0}  (x_2 + x_0 + \varphi_1^+ \varphi_2^- + \varphi_1^- \varphi_2^+)^k \biggl( \frac{1}{3}c  x_0^{-2} + \biggl( -(\varphi^+_1 - \varphi_2^+) x_0^{-1} \biggl( \frac{\partial}{\partial \varphi^+_2} \nonumber\\
& & \quad + \, \varphi^-_2 \frac{\partial}{\partial x_2} \biggr)  + (\varphi^-_1 - \varphi_2^-) x_0^{-1}  \biggl( \frac{\partial}{\partial \varphi^-_2} + \varphi^+_2 \frac{\partial}{\partial x_2} \biggr)    - 2(\varphi^+_1 - \varphi_2^+) (\varphi^-_1 \nonumber \\
& & \quad  - \, \varphi_2^-) x_0^{-2} - 2(\varphi^+_1 - \varphi_2^+) (\varphi^-_1 - \varphi_2^-) x_0^{ -1}\frac{\partial}{\partial x_2} \biggr) Y( \mu , (x_2, \varphi_2^+, \varphi_2^-)) \biggr) \nonumber \\
&=& ( kx_2^{k-1}  + (\varphi_1^+ \varphi_2^- + \varphi_1^- \varphi_2^+)k(k-1) x_2^{k-2} + \varphi_1^+ \varphi_1^- \varphi_2^+\varphi_2^-k(k-1)(k  \nonumber \\
& & \quad -\, 2) x_2^{k-3} ) \frac{1}{3}c + \biggl( -( x_2^k + (\varphi_1^+ \varphi_2^- + \varphi_1^- \varphi_2^+)k x_2^{k-1} ) (\varphi^+_1 - \varphi_2^+) \biggl( \frac{\partial}{\partial \varphi^+_2}   \nonumber 
\end{eqnarray}
\begin{eqnarray}
& & \quad + \, \varphi^-_2 \frac{\partial}{\partial x_2} \biggr)+  ( x_2^k + (\varphi_1^+ \varphi_2^- + \varphi_1^- \varphi_2^+)k x_2^{k-1} ) (\varphi^-_1 - \varphi_2^-)  \biggl( \frac{\partial}{\partial \varphi^-_2} \nonumber  \\
& &  \quad  +\,  \varphi^+_2 \frac{\partial}{\partial x_2} \biggr) - 2( kx_2^{k-1}  + (\varphi_1^+ \varphi_2^- + \varphi_1^- \varphi_2^+)k(k-1) x_2^{k-2} )(\varphi^+_1 \nonumber\\
& & \quad - \, \varphi_2^+) (\varphi^-_1 - \varphi_2^-) - 2 ( x_2^k + (\varphi_1^+ \varphi_2^- + \varphi_1^- \varphi_2^+)k x_2^{k-1} ) (\varphi^+_1 - \varphi_2^+) (\varphi^-_1  \nonumber \\
& & \quad - \, \varphi_2^-) \frac{\partial}{\partial x_2} \biggr) Y( \mu , (x_2, \varphi_2^+, \varphi_2^-)) \nonumber
.
\end{eqnarray}

Thus for $m,n \in \mathbb{Z}$, we have 
\begin{eqnarray*}
\lefteqn{\left[ L (m), L(n) \right] = \frac{1}{4} \left[\mu^{+-}_m, \mu^{+-}_n\right] } \\
&=& \frac{1}{4} \mathrm{Res}_{x_1} x_1^{m+1} \frac{\partial}{\partial \varphi_1^-} \frac{\partial}{\partial \varphi_1^+}  \mathrm{Res}_{x_2} x_2^{n+1} \frac{\partial}{\partial \varphi_2^-} \frac{\partial}{\partial \varphi_2^+}\left[Y(\mu, (x_1, \varphi_1^+, \varphi_1^-)), \right.\\
& & \quad \left. Y(\mu, (x_2, \varphi_2^+, \varphi_2^-)) \right]\\
&=&  \frac{1}{4} \frac{\partial}{\partial \varphi_2^-} \frac{\partial}{\partial \varphi_2^+} \mathrm{Res}_{x_2} x_2^{n+1}  \biggl(  \varphi_2^+\varphi_2^-(m+1)m(m-1) x_2^{m-2}) \frac{1}{3}c  \\
& & \quad + \biggl( -\varphi_2^+(m+1) x_2^m \biggl( \frac{\partial}{\partial \varphi^+_2}  + \varphi^-_2 \frac{\partial}{\partial x_2} \biggr)  - \varphi_2^- (m+1) x_2^m   \biggl( \frac{\partial}{\partial \varphi^-_2} + \varphi^+_2 \frac{\partial}{\partial x_2} \biggr) \\
& &  \quad   - \, 2(m+1)  x_2^m   - 2  x_2^{m+1}  \frac{\partial}{\partial x_2} \biggr)  Y( \mu , (x_2, \varphi_2^+, \varphi_2^-)) \biggr) \\
&=&  \frac{1}{12}(m+1)m(m-1)c  \delta_{m+n,0}  + \frac{1}{2} (m+1) L(m+n)  \\
& & \quad - \, \frac{1}{4} (m+1) (-m-n-1) J(m+n)   + \frac{1}{2}  (m+1) L(m+n) \\
& & \quad  + \,  \frac{1}{4} (m+1)(-m-n-1)J(m+n)  + (m+1)  L(m+n) \\
& &  \quad   + \, (-m-n-2)L(m+n)\\
&=&  \frac{1}{12} (m^3-m)  c  \delta_{m+n,0}   + 2(m+1) L(m+n) + (-m-n-2)L(m+n),
\end{eqnarray*}
\begin{eqnarray*}
\lefteqn{\left[J (m), J(n) \right] = \left[\mu_m, \mu_n\right] } \\
&=&\left. \mathrm{Res}_{x_1} x_1^m \mathrm{Res}_{x_2} x_2^n \left[Y(\mu, (x_1, \varphi_1^+, \varphi_1^-)),  Y(\mu, (x_2, \varphi_2^+, \varphi_2^-)) \right] \right|_{\varphi_1^+ =\varphi_1^-= \varphi_2^+ = \varphi_2^-= 0}\\
&=& \mathrm{Res}_{x_2} x_2^n  mx_2^{m-1}  \frac{1}{3}c \ = \ \frac{1}{3}mc \delta_{m+n,0},
\end{eqnarray*}
and
\begin{eqnarray*}
\lefteqn{\left[ G^+ (m + 1/2), G^- (n - 1/2) \right] = \left[ -\mu^+_{m+ \frac{1}{2}},\mu^-_{n- \frac{1}{2}} \right]}\\
&=& \! \! \! \! \! - \frac{\partial}{\partial \varphi_1^+} \mathrm{Res}_{x_1} x_1^{m+1} \frac{\partial}{\partial \varphi_2^-} \mathrm{Res}_{x_2} x_2^n  \!  \left[Y(\mu, (x_1, \varphi_1^+, \varphi_1^-)), 
Y(\mu, (x_2, \varphi_2^+, \varphi_2^-)) \right] \! \biggr|_{\varphi_1^- = \varphi_2^+ = 0} \\
&=&  \! \! \! \ \frac{\partial}{\partial \varphi_2^-} \mathrm{Res}_{x_2} x_2^n \biggl(  \varphi_2^-(m+1)m x_2^{m-1}\frac{1}{3}c + \biggl(  - x_2^{m+1}\biggl( \frac{\partial}{\partial \varphi^+_2}  + \varphi^-_2 \frac{\partial}{\partial x_2} \biggr)  \\
& & \quad  + \, 2(m+1) \varphi_2^- x_2^m  + \, 2  \varphi_2^- x_2^{m+1}  \frac{\partial}{\partial x_2} \biggr) Y(\mu, (x_2, \varphi_2^+, \varphi_2^-) \biggr) \biggr|_{\varphi_2^+ = 0} \\
&=&  \! \! \!  \frac{1}{3}(m^2+m) \delta_{m+n,0} c + 2L(m+n) - (-m-n-1)J(m+n) 
\end{eqnarray*}
\begin{eqnarray*}
& & \quad + \, 2(m+1) J(m+n)  + 2(-m-n-1)J(m+n) \\
&=& \! \! \! 2L(m+n) + (m- n+ 1)J(m+n)  + \frac{1}{3}(m^2+m) \delta_{m+n,0} c .
\end{eqnarray*}
The other $N=1$ Neveu-Schwarz relations are proved similarly; that is, by taking the appropriate term in (\ref{for-brackets}) and simplifying.
\end{proof}

\begin{rema}
{\em Proposition \ref{mu-prop} implies that if $(V, Y(\cdot, (x, \varphi^+, \varphi^-)), \mathbf{1})$ is an $N=2$ vertex superalgebra such that there exists an element $\mu \in V$ satisfying (\ref{for-mu-prop}), and if in addition $V$ satisfies the $L(0)$- and $J(0)$-grading axioms (\ref{grading-for-vosa-with}) and (\ref{J(0)-grading}) and the $G^\pm(-1/2)$-derivative properties (\ref{G(-1/2)-derivative}), then $(V, Y(\cdot, (x, \varphi^+, \varphi^-)), \mathbf{1}, \mu)$ is an $N=2$ Neveu-Schwarz vertex operator superalgebra with central charge $c$.}
\end{rema}

\begin{rema}\label{OPE-remark}
{\em Typically in the physics literature, instead of encoding the $N=2$ Neveu-Schwarz algebra relations as we did in equation (\ref{for-mu-prop}), they are encoded using the ``operator product expansion", which in our formulation with two odd formal variables is given in Section \ref{duality-section}, equation (\ref{OPE}). } 
\end{rema}

\section{Supercommutativity and associativity properties}\label{commut-assoc-section}

\subsection{Weak supercommutativity and weak associativity for $N=2$ vertex superalgebras with two odd formal variables}\label{weak-section}

In this section we show that  an $N=2$ vertex superalgebra with two odd formal variables satisfies properties called ``weak supercommutativity" and ``weak associativity".  Such properties for vertex superalgebras without odd formal variables and without superconformal elements were first formulated and studied for instance in \cite{DL}  and \cite{Li}.  In \cite{B-thesis} and \cite{B-vosas} we gave the first instance of these properties formulated to include odd formal variables, namely we formulated the properties of weak supercommutativity and weak associativity with odd formal variables for an $N=1$ Neveu-Schwarz vertex operator superalgebra with odd formal variables.  These properties are referred to as ``weak'' supercommutativity and ``weak'' associativity because in Section \ref{duality-section}, we will prove slightly stronger statements about the nature of certain rational functions associated with products and iterates of vertex operators.

In this section we are following and extending the corresponding results and arguments of \cite{Li} and \cite{B-thesis}; see also \cite{B-vosas} and \cite{LL}.

\begin{prop}\label{weak-supercommutativity}
{\bf (weak supercommutativity)} Let $(V,Y(\cdot,(x, \varphi^+, \varphi^-)), \mathbf{1})$ be an $N=2$ vertex  superalgebra with odd formal variables.  There exists $k \in \Z$ such that 
\begin{equation}\label{weak-supercommutativity-formula}
(x_1 - x_2 - \varphi_1^+ \varphi_2^- - \varphi_1^- \varphi_2^+)^k \Bigl[Y(u,(x_1,\varphi_1^+, \varphi_1^-))
,Y(v,(x_2,\varphi_2^+, \varphi_2^-)) \Bigr] = 0 .
\end{equation}
Furthermore this weak supercommutativity follows {}from the truncation condition (\ref{truncation}) and the Jacobi identity (\ref{Jacobi-identity}). 
\end{prop}

\begin{proof}
Let $k \in \Z$.  Taking $\mbox{Res}_{x_0} x_0^k$ of the Jacobi identity, we have
\begin{eqnarray*}
\lefteqn{(x_1 - x_2 - \varphi_1^+ \varphi_2^- - \varphi_1^- \varphi_2^+)^k \left[ Y(u,(x_1,\varphi_1^+, \varphi_1^-)), Y(v,(x_2,\varphi_2^+, \varphi_2^-)) \right] } \\
&=& \mbox{Res}_{x_0} x_2^{-1} \delta \biggl( \frac{x_1 - x_0 -
\varphi_1^+ \varphi_2^- - \varphi_1^- \varphi_2^+}{x_2}\biggr) x_0^k Y(Y(u,(x_0,\varphi_1^+ -
\varphi_2^+, \varphi_1^-\\
& & \quad  - \,  \varphi_2^-)) v, (x_2, \varphi_2^+, \varphi_2^-))\\ 
&=& \sum_{n \in \mathbb{Z}} \mbox{Res}_{x_0} x_2^{-1} \delta \biggl( \frac{x_1 - x_0 -
\varphi_1^+ \varphi_2^- - \varphi_1^- \varphi_2^+}{x_2}\biggr) x_0^k x_0^{-n-1} Y((u_n + (\varphi_1^+-
\varphi_2^+)  
\end{eqnarray*}
\begin{eqnarray*}
& & \quad \cdot u^+_{n - \frac{1}{2}} + (\varphi_1^- - \varphi_2^-)u^-_{n - \frac{1}{2}} + (\varphi_1^+ - \varphi_2^+)(\varphi_1^- - \varphi_2^-) u_{n-1}^{+-}) v, (x_2, \varphi_2^+, \varphi_2^-))  . 
\end{eqnarray*}
Now using the truncation condition (\ref{truncation}), choose $k$ such that $u_n v = u^\pm_{n-\frac{1}{2}} v = u_n^{+-} v = 0$ for all $n \in 
\Z$, $n \geq k $.   
\end{proof}

\begin{prop}\label{weak-associativity}
{\bf (weak associativity)} Let $(V,Y(\cdot,(x, \varphi^+,\varphi^-)), \mathbf{1})$ be an $N=2$ vertex  superalgebra with odd formal variables.  Then for $u, w \in V$, there exists $k \in \Z$  such that for any $v \in V$ 
\begin{multline}\label{weak-associativity-formula}
(x_0 + x_2 + \varphi_1^+ \varphi_2^- + \varphi_1^- \varphi_2^+)^k Y(Y(u,(x_0,\varphi_1^+ -
\varphi_2^+, \varphi_1^- - \varphi_2^-)v, (x_2,\varphi_2^+, \varphi_2^-))w \\
= (x_0 + x_2 + \varphi_1^+ \varphi_2^- + \varphi_1^- \varphi_2^+)^k Y(u,(x_0 + x_2 +
\varphi_1^+ \varphi_2^- + \varphi_1^- \varphi_2^+ ,\varphi_1^+, \varphi_1^-)) \\
\cdot Y(v,(x_2,\varphi_2^+, \varphi_2^-))w  .
\end{multline}
Furthermore this weak associativity follows {}from the truncation condition (\ref{truncation}) and the Jacobi identity (\ref{Jacobi-identity}). 
\end{prop}  

\begin{proof}
Using the Jacobi identity, the $\delta$-function identities (\ref{delta-2-terms-with-phis}) and (\ref{delta-3-terms-with-phis}), and (\ref{delta-substitute}),  we obtain the following iterate
\begin{eqnarray*}
\lefteqn{Y(Y(u,(x_0,\varphi_1^+ - \varphi_2^+, \varphi_1^- - \varphi_2^-)v, (x_2,\varphi_2^+, \varphi_2^-))}\\
&=& \mbox{Res}_{x_1} x_1^{-1} \delta \biggl( \frac{x_2 + x_0 + \varphi_1^+
\varphi_2^- + \varphi_1^- \varphi_2^+ }{x_1}\biggr) Y(Y(u,(x_0,\varphi_1^+ - \varphi_2^+, \varphi_1^- - \varphi_2^-)v,\\
& & \quad (x_2,\varphi_2^+, \varphi_2^-)) \\
&=& \mbox{Res}_{x_1} x_2^{-1} \delta \biggl( \frac{x_1 - x_0 - \varphi_1^+ \varphi_2^- - \varphi_1^- \varphi_2^+ }{x_2}\biggr) Y(Y(u,(x_0,\varphi_1^+ - \varphi_2^+, \varphi_1^- - \varphi_2^-)v,\\
& & \quad  (x_2,\varphi_2^+, \varphi_2^-)) \\
&=& \mbox{Res}_{x_1} \Biggl( 
x_0^{-1} \delta \biggl( \frac{x_1 - x_2 - \varphi^+_1 \varphi^-_2 - \varphi^-_1 \varphi^+_2}{x_0} \biggr)Y(u,(x_1, \varphi^+_1, \varphi^-_1)) Y(v,(x_2, \varphi^+_2,\varphi^-_2)) \\
& & \quad -  (-1)^{\eta(u)\eta(v)} x_0^{-1} \delta \biggl( \frac{x_2 - x_1 + \varphi^+_1 \varphi^-_2 + \varphi^-_1 \varphi^+_2}{-x_0} \biggr) Y(v,(x_2, \varphi^+_2,\varphi^-_2))\\
& & \quad \cdot Y(u,(x_1,\varphi^+_1,\varphi^-_1)) \Biggr) \\
&=& \mbox{Res}_{x_1} \left( x_1^{-1} \delta \biggl( \frac{x_0 + x_2 +
\varphi_1^+ \varphi_2^- + \varphi_1^- \varphi_2^+}{x_1} \biggr) Y(u,(x_1, \varphi_1^+, \varphi_1^-))Y(v,(x_2, \varphi_2^+, \varphi_2^-)) \right. \\
& & \quad - (-1)^{\eta(u)\eta(v)} Y(v,(x_2, \varphi_2^+, \varphi_2^-)) \left(x_0^{-1} \delta \biggl(\frac{x_1 - x_2 - \varphi_1^+ \varphi_2^- - \varphi_1^- \varphi_2^+}{x_0} \biggr) \right. \\
& & \quad \left. \left. - \, x_2^{-1} \delta \biggl(\frac{x_1 - x_0 - \varphi_1^+ \varphi_2^- - \varphi_1^- \varphi_2^+}{x_2} \biggr) \right) Y(u,(x_1, \varphi_1^+, \varphi_1^-))
\right) \\ 
&=& \mbox{Res}_{x_1} \left( x_1^{-1} \delta \biggl( \frac{x_0 + x_2 +
\varphi_1^+ \varphi_2^- + \varphi_1^- \varphi_2^+}{x_1} \biggr) Y(u,(x_1, \varphi_1^+, \varphi_1^-))Y(v,(x_2, \varphi_2^+, \varphi_2^-)) \right. \\
& & \quad - (-1)^{\eta(u)\eta(v)} Y(v,(x_2, \varphi_2^+, \varphi_2^-)) \left(x_1^{-1} \delta \biggl(\frac{x_0 + x_2 + \varphi_1^+ \varphi_2^- + \varphi_1^- \varphi_2^+}{x_1} \biggr) \right. \\
& & \quad \left. \left. - \, x_1^{-1} \delta \biggl(\frac{x_2  +x_0 + \varphi_1^+ \varphi_2^- + \varphi_1^- \varphi_2^+}{x_1} \biggr) \right) Y(u,(x_1, \varphi_1^+, \varphi_1^-))
\right) 
\end{eqnarray*}
\begin{eqnarray*}
&=& Y(u,(x_0 + x_2 + \varphi_1^+ \varphi_2^- + \varphi_1^- \varphi_2^+ ,\varphi_1^+, \varphi_1^-))
Y(v,(x_2,\varphi_2^+, \varphi_2^-)) \\
& & \quad - (-1)^{\eta(u)\eta(v)} Y(v,(x_2, \varphi_2^+, \varphi_2^-))\Bigl(Y(u,(x_0 + x_2 + \varphi_1^+ \varphi_2^- + \varphi_1^- \varphi_2^+ ,\varphi_1^+, \varphi_1^-))  \\
& & \quad - Y(u,(x_2 + x_0 + \varphi_1^+ \varphi_2^- + \varphi_1^- \varphi_2^+
,\varphi_1^+, \varphi_1^-)) \Bigr) .\\  
\end{eqnarray*}
Let $k \in \Z$ be such that $x^{k-2} Y(u,(x,\varphi^+, \varphi^-))w$ involves only positive powers of $x$.  
Then 
\begin{multline*}
(x_0 + x_2 + \varphi_1^+ \varphi_2^- + \varphi_1^- \varphi_2^+)^k \Bigl(Y(u,(x_0 + x_2 + \varphi_1^+
\varphi_2^- + \varphi_1^- \varphi_2^+ ,\varphi_1^+, \varphi_1^-))  \\ 
- Y(u,(x_2 + x_0 + \varphi_1^+ \varphi_2^- + \varphi_1^- \varphi_2^+,\varphi_1^+, \varphi_1^-)) \Bigr)w   = 0 
\end{multline*}
and weak associativity follows. 
\end{proof}

\begin{prop}\label{equals-Jacobi1}  
In the presence of the other axioms in the definition of $N=2$ vertex superalgebra with odd formal variables, the Jacobi identity is equivalent to weak supercommutativity and weak associativity. 
\end{prop}

\begin{proof}
Propositions \ref{weak-supercommutativity} and \ref{weak-associativity} show that in the presence of the
other axioms for an $N=2$ vertex superalgebra with two odd formal variables, the Jacobi identity implies weak supercommutativity and weak associativity. 

Assume weak supercommutativity and weak associativity hold.  Let $u,v,w, \in V$, and choose $k \in \Z$ such that weak supercommutativity (\ref{weak-supercommutativity-formula}) and weak associativity (\ref{weak-associativity-formula}) hold.  Then using the $\delta$-function identity (\ref{delta-3-terms-with-phis}), we have
\begin{eqnarray}
\lefteqn{\qquad x_0^k x_1^k \left(x_0^{-1} \delta \biggl(\frac{x_1 - x_2 - \varphi_1^+ \varphi_2^- - \varphi_1^-
\varphi_2^+}{x_0} \biggr)  Y(u,(x_1, \varphi_1^+, \varphi_1^-)) 
\right.} \label{assoc1} \\
& & \quad   \cdot Y(v,(x_2, \varphi_2^+, \varphi_2^-))w - (-1)^{\eta(u)\eta(v)} x_0^{-1} \delta \biggl(\frac{x_2 - x_1 + \varphi_1^+ \varphi_2^- + \varphi_1^- \varphi_2^+}{-x_0} \biggr)  \nonumber\\  
& & \quad \left. \cdot Y(v,(x_2, \varphi_2^+, \varphi_2^-)) Y(u,(x_1, \varphi_1^+, \varphi_1^-))w \right) \nonumber\\
&=& x_0^{-1} \delta \biggl(\frac{x_1 - x_2 - \varphi_1^+ \varphi_2^- - \varphi_1^- \varphi_2^+}{x_0}
\biggr) x_1^k (x_1 - x_2 - \varphi_1^+ \varphi_2^- - \varphi_1^- \varphi_2^+)^k \nonumber \\
& &  \quad \cdot Y(u,(x_1, \varphi_1^+, \varphi_1^-))Y(v,(x_2, \varphi_2^+, \varphi_2^-))w   - (-1)^{\eta(u)\eta(v)} x_0^{-1} \nonumber 
\end{eqnarray}
\begin{eqnarray}
& & \quad \cdot \delta \biggl(\frac{x_2 - x_1 + \varphi_1^+ \varphi_2^- + \varphi_1^- \varphi_2^+}{-x_0} \biggr) x_1^k (- x_2 + x_1 - \varphi_1^+ \varphi_2^- - \varphi_1^- \varphi_2^+)^k \nonumber\\
& & \quad \cdot Y(v,(x_2, \varphi_2^+, \varphi_2^-)) Y(u,(x_1, \varphi_1^+, \varphi_1^-))w \nonumber \\
&=& x_2^{-1} \delta \biggl(\frac{x_1 - x_0 - \varphi_1^+ \varphi_2^-  - \varphi_1^- \varphi_2^+}{x_2}
\biggr) x_1^k (x_1 - x_2 - \varphi_1^+ \varphi_2^- - \varphi_1^- \varphi_2^+)^k  \nonumber\\
& & \quad \cdot Y(u,(x_1,\varphi_1^+, \varphi_1^-))Y(v,(x_2, \varphi_2^+, \varphi_2^-))w .\nonumber
\end{eqnarray}
By weak supercommutativity, 
\[x_1^k (x_1 - x_2 - \varphi_1^+ \varphi_2^- - \varphi_1^- \varphi_2^+)^k
Y(u,(x_1, \varphi_1^+, \varphi_1^-))Y(v,(x_2, \varphi_2^+, \varphi_2^-))w\] 
involves only nonnegative powers of $x_1$.  Thus in this case, we can replace $x_1$ by $x_2 + x_0 + \varphi_1^+  \varphi_2^- + \varphi_1^- \varphi_2^+$ or $x_0 +x_2 + \varphi_1^+ \varphi_2^- + \varphi_1^- \varphi_2^+$.  Therefore since weak associativity holds for our choice of $k$, (\ref{assoc1}) is equal to
\begin{eqnarray*}
\lefteqn{x_0^k (x_0 + x_2 + \varphi_1^+  \varphi_2^- + \varphi_1^- \varphi_2^+)^k x_2^{-1} \delta \biggl(\frac{x_1 - x_0 - \varphi_1^+ \varphi_2^- - \varphi_1^- \varphi_2^+}{x_2}
\biggr) }  \\
& & \quad \cdot Y(u,(x_0 + x_2 + \varphi_1^+ \varphi_2^- + \varphi_1^- \varphi_2^+, \varphi_1^+, \varphi_1^-))Y(v,(x_2, \varphi_2^+, \varphi_2^-))w \\ 
&=& x_2^{-1} \delta \biggl(\frac{x_1 - x_0 - \varphi_1^+ \varphi_2^- - \varphi_1^- \varphi_2^+}{x_2}
\biggr) x_0^k (x_0 + x_2 + \varphi_1^+ \varphi_2^- + \varphi_1^- \varphi_2^+)^k  \\
& & \quad \cdot Y(Y(u,(x_0, \varphi_1^+ - \varphi_2^+, \varphi_1^- - \varphi_2^-))v,(x_2, \varphi_2^+, \varphi_2^-))w \\
&=& x_2^{-1} \delta \biggl(\frac{x_1 - x_0 - \varphi_1^+ \varphi_2^- - \varphi_1^- \varphi_2^+}{x_2}
\biggr) x_0^k x_1^k Y(Y(u,(x_0, \varphi_1^+ - \varphi_2^+, \varphi_1^- - \varphi_2^-))v,\\
& & \quad (x_2, \varphi_2^+, \varphi_2^-))w \\
&=& x_0^k x_1^k x_2^{-1} \delta \biggl(\frac{x_1 - x_0 - \varphi_1^+ \varphi_2^- - \varphi_1^- \varphi_2^+}{x_2}\biggr) Y(Y(u,(x_0, \varphi_1^+ - \varphi_2^+, \varphi_1^- - \varphi_2^-))v,\\
& & \quad (x_2, \varphi_2^+, \varphi_2^-))w ,
\end{eqnarray*}
which implies the Jacobi identity. 
\end{proof}

\begin{thm}\label{get-Jacobi2} 
The Jacobi identity for an $N=2$ vertex superalgebra with odd formal variables follows {}from weak supercommutativity (\ref{weak-supercommutativity-formula}) in the presence of the other axioms together with the $\mathcal{D}^\pm$-bracket-derivative properties (\ref{Dpm-bracket-derivative}).  In particular, in the definition of the notion of $N=2$ vertex superalgebra with odd formal variables, the Jacobi identity can be replaced by these properties.
\end{thm}

\begin{proof}  We generalize the proof of Theorem 3.5.1 in \cite{LL} to this setting.    Letting $\mathcal{D} = \frac{1}{2}\left[ \mathcal{D}^+, \mathcal{D}^- \right]$, the $\mathcal{D}^\pm$-bracket-derivative properties (\ref{Dpm-bracket-derivative}) imply the $\mathcal{D}$-bracket-derivative property since 
\begin{eqnarray*}
\lefteqn{\left[ \mathcal{D}, Y(v, (x,\varphi^+, \varphi^-))\right] }\\
&=& \left[\frac{1}{2}\left[\mathcal{D}^+, \mathcal{D}^- \right], Y(v, (x,\varphi^+, \varphi^-))\right] \\
&=& \frac{1}{2} \left( \left[ \mathcal{D}^+, \left[ \mathcal{D}^-, Y(v, (x,\varphi^+, \varphi^-))\right] \right] + \left[ \mathcal{D}^-, \left[ \mathcal{D}^+, Y(v, (x,\varphi^+, \varphi^-))\right] \right] \right)\\
\end{eqnarray*}
\begin{eqnarray*}
&=&\frac{1}{2} \biggl( \left( \frac{\partial}{\partial \varphi^-} - \varphi^+ \frac{\partial}{\partial x}\right)\left( \frac{\partial}{\partial \varphi^+} - \varphi^- \frac{\partial}{\partial x}\right) + \left( \frac{\partial}{\partial \varphi^+} - \varphi^- \frac{\partial}{\partial x}\right)\biggl( \frac{\partial}{\partial \varphi^-} \\
& & \quad - \, \varphi^+ \frac{\partial}{\partial x}\biggr) \biggr)  Y(v, (x, \varphi^+, \varphi^-))\\
&=& \frac{\partial}{\partial x} Y(v, (x, \varphi^+, \varphi^-)),
\end{eqnarray*}
for $v \in V$.  Repeatedly using the $\mathcal{D}^\pm$- and $\mathcal{D}$-bracket-derivative properties as well as (\ref{CBH}) and Proposition \ref{Taylor}, we have that
\begin{eqnarray}\label{conjugation}
\lefteqn{e^{x\mathcal{D} + \varphi^+\mathcal{D}^+ + \varphi^- \mathcal{D}^-} Y(v,(x_0, \varphi^+_0, \varphi^-_0))e^{-x\mathcal{D} - \varphi^+\mathcal{D}^+ - \varphi^- \mathcal{D}^-} }\\
&=& e^{x \frac{\partial}{\partial x_0} + \varphi^+ \left(
\frac{\partial}{\partial \varphi^+_0} + \varphi^-_0 \frac{\partial}{\partial x_0} \right)+ \varphi^- \left(
\frac{\partial}{\partial \varphi^-_0} + \varphi^+_0 \frac{\partial}{\partial x_0} \right)} Y(v,(x_0, \varphi^+_0, \varphi^-_0)) \nonumber \\
&=& Y(v,(x_0 + x + \varphi^+_0 \varphi^- + \varphi^-_0 \varphi^+, \varphi^+_0 + \varphi^+, \varphi^-_0 + \varphi^-)) .\nonumber
\end{eqnarray}
{}From the vacuum property (\ref{vacuum-identity}), and the definition of $\mathcal{D}^\pm$ (\ref{define-Dpm}), we have that $\mathcal{D}^\pm (\mathbf{1}) = \mathbf{1}^\pm_{-3/2} \mathbf{1} = 0$.   This and (\ref{conjugation}), imply that
\begin{eqnarray}\label{act-on-vacuum}
\lefteqn{Y(v,(x_0 + x + \varphi^+_0 \varphi^- + \varphi^-_0 \varphi^+, \varphi^+_0 + \varphi^+, \varphi^-_0 + \varphi^-)) \mathbf{1} }\\
&=& e^{x \mathcal{D} + \varphi^+ \mathcal{D}^+ + \varphi^- \mathcal{D}^-} Y(v,(x_0, \varphi^+_0, \varphi^-_0))e^{-x\mathcal{D} - \varphi^+\mathcal{D}^+ - \varphi^- \mathcal{D}^-} \mathbf{1}\nonumber \\
&=& e^{x\mathcal{D} + \varphi^+\mathcal{D}^+ + \varphi^- \mathcal{D}^-} Y(v,(x_0, \varphi^+_0, \varphi^-_0))\mathbf{1}. \nonumber
\end{eqnarray}
The creation property (\ref{creation-property1}), (\ref{creation-property2})  implies that we can set $x_0$, $\varphi^+_0$, and $\varphi^-_0$ equal to zero in (\ref{act-on-vacuum}), giving
\begin{equation}\label{for-skew-again}
Y(v,(x ,\varphi^+,\varphi^-)) \mathbf{1} = e^{x\mathcal{D} + \varphi^+ \mathcal{D}^+ + \varphi^- \mathcal{D}^-} v ,
\end{equation}
(cf. (\ref{for-skew-symmetry})).

Let $u,v \in V$.  By the truncation condition (\ref{truncation}), there exists $n \in \mathbb{N}$ be such that 
\begin{equation}\label{for-set-to-zero}
(x_1 - x_2 - \varphi_1^+ \varphi_2^- - \varphi_1^- \varphi_2^+)^n Y(v, (x_2- x_1 - \varphi_2^+ \varphi_1^- - \varphi_2^- \varphi_1^+, \varphi_2^+- \varphi_1^+, \varphi_2^-- \varphi_1^-)) u
\end{equation}
involves only nonnegative powers of $x_2 - x_1$.   Let $u,v \in V$ and let $k \in \mathbb{N}$ with $k \geq n$ be such that weak supercommutativity holds.  Then using weak supercommutativity, (\ref{for-skew-again}) and (\ref{conjugation}), we have
\begin{eqnarray*}
\lefteqn{(x_1 - x_2 - \varphi_1^+ \varphi_2^- - \varphi_1^- \varphi_2^+)^k Y(u,(x_1, \varphi^+_1, \varphi^-_1)) Y(v, (x_2, \varphi_2^+, \varphi_2^-)) \mathbf{1} }\\
&=& (-1)^{\eta(u)\eta(v)}(x_1 - x_2 - \varphi_1^+ \varphi_2^- - \varphi_1^- \varphi_2^+)^k Y(v, (x_2, \varphi_2^+, \varphi_2^-)) Y(u,(x_1, \varphi^+_1, \varphi^-_1)) \mathbf{1}\\
&=& (-1)^{\eta(u)\eta(v)}(x_1 - x_2 - \varphi_1^+ \varphi_2^- - \varphi_1^- \varphi_2^+)^k Y(v, (x_2, \varphi_2^+, \varphi_2^-)) e^{x_1\mathcal{D} + \varphi^+_1 \mathcal{D}^+ + \varphi^-_1 \mathcal{D}^-} u\\
&=& (-1)^{\eta(u)\eta(v)}(x_1 - x_2 - \varphi_1^+ \varphi_2^- - \varphi_1^- \varphi_2^+)^k e^{x_1\mathcal{D} + \varphi^+_1 \mathcal{D}^+ + \varphi^-_1 \mathcal{D}^-} Y(v, (x_2- x_1 \\
& & \quad - \, \varphi_2^+ \varphi_1^- - \varphi_2^- \varphi_1^+, \varphi_2^+- \varphi_1^+, \varphi_2^-- \varphi_1^-)) u.
\end{eqnarray*}
Since (\ref{for-set-to-zero}) involves only nonnegative powers of $x_2- x_1$, we may set $x_2 = 0$ and $\varphi_2^\pm = 0$.  Thus using the creation property, we have
\[x_1^k Y(u,(x_1, \varphi^+_1, \varphi^-_1)) v
= (-1)^{\eta(u)\eta(v)}x_1^k e^{x_1\mathcal{D} + \varphi^+_1 \mathcal{D}^+ + \varphi^-_1 \mathcal{D}^-} Y(v, (- x_1, - \varphi_1^+, - \varphi_1^-)) u.\]
Multiplying both sides by $x_1^{-k}$, we obtain skew-supersymmetry (\ref{skew-supersymmetry}).

For any $u,v,w \in V$ with $v$ and $w$ of homogeneous sign, let $k \in \mathbb{N}$ be such that weak supercommutativity (\ref{weak-supercommutativity-formula}) holds for $u$ and $w$.  Then using skew supersymmetry for $v$ and $w$, the conjugation formula (\ref{conjugation}), weak supercommutativity, and skew supersymmetry for $Y(u,(x_0, \varphi_1^+- \varphi_2^+, \varphi_1^- - \varphi_2^-))v$ and $w$, we have
\begin{eqnarray*}
\lefteqn{(x_0 + x_2 + \varphi_1^+ \varphi_2^- + \varphi_1^- \varphi_2^+)^k Y(u,(x_0 + x_2 +
\varphi_1^+ \varphi_2^- + \varphi_1^- \varphi_2^+ ,\varphi_1^+, \varphi_1^-))} \\
& & \quad \cdot Y(v,(x_2,\varphi_2^+, \varphi_2^-))w  \\
&=& \! \! \! (x_0 + x_2 + \varphi_1^+ \varphi_2^- + \varphi_1^- \varphi_2^+)^k Y(u,(x_0 + x_2 +
\varphi_1^+ \varphi_2^- + \varphi_1^- \varphi_2^+ ,\varphi_1^+, \varphi_1^-)) \\
& & \quad \cdot (-1)^{\eta(v) \eta(w)}  e^{x_2 \mathcal{D} + \varphi^+_2 \mathcal{D}^+ + \varphi^-_2 \mathcal{D}^-} Y(w,(-x_2,-\varphi^+_2, - \varphi^-_2))v    \\
&=&  \! \! \! (-1)^{\eta(v) \eta(w)}  e^{x_2 \mathcal{D} + \varphi^+_2 \mathcal{D}^+ + \varphi^-_2 \mathcal{D}^-} (x_0 + x_2 + \varphi_1^+ \varphi_2^- + \varphi_1^- \varphi_2^+)^k Y(u,(x_0 ,\varphi_1^+ - \varphi_2^+, \\
& & \quad \varphi_1^-- \varphi_2^-))  Y(w,(-x_2,-\varphi^+_2, - \varphi^-_2))v    \\
&=&  \! \! \!  (-1)^{\eta(v) \eta(w)}  e^{x_2 \mathcal{D} + \varphi^+_2 \mathcal{D}^+ + \varphi^-_2 \mathcal{D}^-} (-1)^{\eta(u)\eta(w)} (x_0 + x_2 + \varphi_1^+ \varphi_2^- + \varphi_1^- \varphi_2^+)^k Y(w,\\
& & \quad (-x_2,-\varphi^+_2, - \varphi^-_2)) Y(u,(x_0 ,\varphi_1^+ - \varphi_2^+, \varphi_1^-- \varphi_2^-))  v    \\
&=&  \! \! \!  (x_0 + x_2 + \varphi_1^+ \varphi_2^- + \varphi_1^- \varphi_2^+)^k Y(Y(u,(x_0,\varphi_1^+ -
\varphi_2^+, \varphi_1^- - \varphi_2^-)v, (x_2,\varphi_2^+, \varphi_2^-))w,
\end{eqnarray*}
proving weak associativity.  By Proposition \ref{equals-Jacobi1}, the result follows. 
\end{proof}

\begin{rema}
{\em
Theorem \ref{get-Jacobi2} gives a practical way of determining if a vector space $V$ equipped with a linear map $Y(\cdot, (x, \varphi^+, \varphi^-))$ and a distinguished vector $\mathbf{1}$ is in fact an $N=2$ vertex superalgebra.  All the axioms for an $N=2$ vertex superalgebra are relatively straightforward to check except for the Jacobi identity.  The weak commutativity and the $\mathcal{D}^\pm$-bracket-derivative properties are in practice usually easier to prove than the Jacobi identity, especially using the theory of local fields \cite{Li}. (See  Remark \ref{is-a-vertex-superalgebra-remark} and Theorem \ref{superalgebras} below for why the theory of local fields for vertex superalgebras as presented in \cite{Li} applies.)  However, the Jacobi identity formulation of $N=2$ vertex superalgebra allows one to immediately formulate such useful properties as the commutator and iterate formulas (\ref{homo-bracket}) and  (\ref{iterate}) (which themselves can be thought of as encoding an infinite number of formulas) and is in general much more useful for proving consequences and properties of vertex superalgebras as we did throughout Section \ref{N2-with-section}, for instance.   In fact the single Jacobi identity contains a wealth of information about the individual components $v_n, v_{n-1/2}^\pm, v_n^{+-} \in \mathrm{End} \, V$, and is easy to work with in its formal calculus form using the powerful techniques of formal calculus as presented in for instance \cite{FLM} and \cite{LL} (cf. Remark 3.1.2. \cite{LL}).  Furthermore, we see {}from Proposition \ref{Dpm-derivative-prop}, Corollary \ref{D-derivative-cor} and Remark \ref{N2-subalgebras-representation-remark}, that when the odd variables are included in the Jacobi identity, it naturally imposes a representation of $\mathfrak{osp}_{\bigwedge_*}(2|2)_{<0}$ on the vertex superalgebra reflecting the $N=2$ superconformal structure. 
}
\end{rema}

\subsection{Expansions of rational superfunctions}

In order to formulate the notions of associativity and supercommutativity for $N=2$ vertex operators, we will need to interpret correlation functions of vertex operators with odd formal variables as expansions of certain rational superfunctions.  In this section we follow and extend the treatment of rational functions as presented in \cite{FHL}, and the treatment of $N=1$ rational superfunctions as treated previously by the author in \cite{B-thesis} and \cite{B-vosas}.

Let $\bigwedge_*[x_1,x_2,\dots,x_n]_P$ be the ring of rational functions obtained by inverting (localizing with respect to) the set 
\[P = \biggl\{\sum_{i = 1}^{n} a_i x_i : a_i \in \mbox{$\bigwedge_*^0$}, \; \mbox{not all} \; (a_i)_B = 0\biggr\} . \] 
Recall the map $\iota_{i_1 \dots i_2} : \mathbb{F}[x_1,\dots,x_n]_P \longrightarrow \mathbb{F}[[x_1, x_1^{-1},\dots, x_n, x_n^{-1}]]$ defined in \cite{FLM} where in \cite{FLM} $P$ is denoted $S$ and the coefficients of elements in $S$ are restricted to the underlying field (in this case $\mathbb{C}$).  We extend this map $\iota_{i_1 \dots i_2}$ to 
\begin{multline}
\mbox{$\bigwedge_*$}[x_1,x_2,\dots,x_n]_P [\varphi_1^+, \varphi_1^-, \varphi_2^+, \varphi_2^-,\dots, \varphi_n^+, \varphi_n^-] \\
= \mbox{$\bigwedge_*$} [x_1, \varphi_1^+, \varphi_1^-, x_2, \varphi_2^+, \varphi_2^-, \dots, x_n, \varphi_n^+, \varphi_2^-]_P
\end{multline}
in the obvious way (cf. \cite{B-vosas}) obtaining
\begin{multline}
\iota_{i_1 \dots i_2} : \mbox{$\bigwedge_*$} [x_1,\varphi_1^+, \varphi_1^-,\dots,x_n,\varphi_n^+, \varphi_n^-]_P \longrightarrow \\
\mbox{$\bigwedge_*$}[[x_1, x_1^{-1},\dots, x_n, x_n^{-1}]] [\varphi_1^+, \varphi_1^-, \dots, \varphi_n^+, \varphi_n^-] .
\end{multline}
Let $\bigwedge_*[x_1, \varphi_1^+, \varphi_1^-, x_2, \varphi_2^+, \varphi_2^-, \dots, x_n,
\varphi_n^+, \varphi_n^-]_{P'}$ be the ring of rational functions obtained by
inverting the set   
\[P' = \biggl\{\sum_{\stackrel{i,j = 1}{i<j}}^{n} (a_i x_i + a^+_{ij} \varphi_i^+
\varphi_j^- + a^-_{ij} \varphi_i^-\varphi_j^+) : a_i, a^\pm_{ij} \in \mbox{$\bigwedge_*^0$}, \; \mbox{not all} \; (a_i)_B = 0\biggr\}. \]
Since we use the convention that a function of even and odd variables should be expanded about the even variables, we have 
\[\mbox{$\bigwedge_*$}[x_1,\varphi_1^+, \varphi_1^-, \dots, x_n, \varphi_n^+, \varphi_n^-]_{P'}
\subseteq \mbox{$\bigwedge_*$}[x_1,\varphi_1^+, \varphi_1^-, \dots, x_n, \varphi_n^+, \varphi_n^-]_P ,
\]   
and $\iota_{i_1 \dots i_2}$ is well defined on $\bigwedge_*[x_1, \varphi_1^+, \varphi_1^-, \dots, x_n, \varphi_n^+, \varphi_n^-]_{P'}$.

In the case $n = 2$, 
\[\iota_{12} : \mbox{$\bigwedge_*$} [x_1,\varphi_1^+, \varphi_1^-, x_2, \varphi_2^+, \varphi_2^-]_{P'} 
\longrightarrow \mbox{$\bigwedge_*$} [[x_1,x_2]] [\varphi_1^+, \varphi_1^-, \varphi_2^+, \varphi_2^-] \]   
is given by first expanding an element of $\bigwedge_* [x_1,\varphi_1^+, \varphi_1^-, x_2, \varphi_2^+, \varphi_2^-]_{P'}$ as a formal series in $\bigwedge_*[x_1,\varphi_1^+, \varphi_1^-, x_2, \varphi_2^+, \varphi_2^-]_P$ and then expanding each term as a series in $\bigwedge_* [[x_1,x_2]] [\varphi_1^+, \varphi_1^-,\varphi_2^+, \varphi_2^-]$ containing at most finitely many negative powers of $x_2$ (using binomial expansions for negative powers of linear polynomials involving both $x_1$ and $x_2$).

\subsection{Duality for $N=2$ Neveu-Schwarz vertex operator superalgebras with two odd formal variables}\label{duality-section}

In \cite{B-moduli}, we began to develop the geometric foundation necessary to formulate the notion of {\it $N=2$ supergeometric vertex operator superalgebra}.  In subsequent work we will show that any such object defines an $N=2$ Neveu-Schwarz vertex operator superalgebra with odd formal 
variables.  To show that the alleged vertex operator superalgebra satisfies the Jacobi identity, we will need the notions of associativity and supercommutativity for an $N=2$ Neveu-Schwarz vertex operator superalgebra with odd formal variables.  Together, these notions of associativity and  (super)commutativity are known as ``duality'', a term which arose {}from physics.   

In fact, it is supercommutativity and associativity which can easily be seen as coming {}from the geometry of the ``sewing" together of ``$N=2$ super-Riemann spheres with tubes".  These $N=2$ super-Riemann spheres with tubes correspond to the supersurfaces swept out by superstrings in the physical model of $N=2$ superconformal field theory.   It is this geometry underlying $N=2$ superconformal field theory in the genus zero case that is studied in \cite{B-moduli} and will be shown in a rigorous way to have deep connections to the algebras we study in this paper.

Throughout this section we follow and extend the treatment of duality as presented in \cite{FHL} for the non-super case and in \cite{B-vosas} for the $N=1$ super case.

For a $\frac{1}{2}\mathbb{Z}$-graded $\bigwedge_*$-module $V = \coprod_{\frac{1}{2}\mathbb{Z}} V_{(n)}$, let $V_{(n)}^*$ be the dual module of $V_{(n)}$ for $n \in \frac{1}{2} \mathbb{Z}$, i.e.,$V_{(n)}^* = \mbox{Hom}_{\bigwedge_*} (V, \bigwedge_*)$.  Let  
\begin{equation}
V' = \coprod_{n \in \frac{1}{2} \mathbb{Z}} V_{(n)}^* 
\end{equation}
be the graded dual space of $V$, 
\begin{equation}
\bar{V} = \prod _{n \in \frac{1}{2} \mathbb{Z}} V_{(n)} = V'^* 
\end{equation}
the algebraic completion of $V$, and $\langle \cdot , \cdot \rangle$ the natural pairing between $V'$ and $\bar{V}$.  We now formulate the weak supercommutativity and weak associativity properties of an $N=2$ Neveu-Schwarz vertex operator superalgebra with odd formal variables into slightly stronger
statements about ``matrix coefficients'' of products and iterates of vertex operators with odd formal variables.  

\begin{prop}\label{supercommutativity} Let $(V, Y(\cdot,(x,\varphi^+, \varphi^-)), \mathbf{1}, \mu)$ be an $N=2$ Neveu-Schwarz vertex operator superalgebra, and let $u, v, w \in V$, with $u$ and $v$
of homogeneous sign, and $v' \in V'$ be arbitrary.  

{\bf (a) (rationality of products)}  The formal series 
\begin{equation}
\langle v', Y(u,(x_1,\varphi_1^+, \varphi_1^-)) Y(v,(x_2,\varphi_2^+, \varphi_2^-)) w \rangle, 
\end{equation}
which involves only finitely many negative powers of $x_2$ and only finitely many positive powers of $x_1$, lies in the image of the map $\iota_{12}$: 
\begin{equation}
\langle v', Y(u,(x_1,\varphi_1^+, \varphi_1^-)) Y(v,(x_2,\varphi_2^+, \varphi_2^-)) w \rangle
= \iota_{12} f(x_1,\varphi_1^+, \varphi_1^-,x_2,\varphi_2^+, \varphi_2^-) , 
\end{equation}
where the (uniquely determined) element $f \in \bigwedge_* [x_1, \varphi_1^+, \varphi_1^-, x_2, \varphi_2^+, \varphi_2^-]_{P'}$ is of the form  
\begin{equation}
f(x_1,\varphi_1^+, \varphi_1^-,x_2,\varphi_2^+, \varphi_2^-) = \frac{g(x_1,\varphi_1^+, \varphi_1^-, x_2, \varphi_2^+, \varphi_2^-)}{x_1^r x_2^s (x_1 - x_2 - \varphi_1^+ \varphi_2^- - \varphi_1^- \varphi_2^+)^t} \end{equation}
for some $g \in \bigwedge_* [x_1,\varphi_1^+, \varphi_1^-,x_2,\varphi_2^+, \varphi_2^-]$ and $r, s, t
\in \mathbb{Z}$.  

{\bf (b) (supercommutativity)}  We also have
\begin{multline}
\langle v', Y(v,(x_2,\varphi_2^+, \varphi_2^-))Y(u,(x_1,\varphi_1^+, \varphi_1^-)) w \rangle \\
= (-1)^{\eta(u)\eta(v)} \iota_{21} f(x_1,\varphi_1^+,\varphi_1^-,x_2,\varphi_2^+, \varphi_2^-), 
\end{multline}
i.e,
\begin{multline}
\iota_{12}^{-1} \langle v', Y(u,(x_1,\varphi_1^+, \varphi_1^-))Y(v,(x_2,\varphi_2^+, \varphi_2^-)) w \rangle\\
= (-1)^{\eta(u)\eta(v)} \iota_{21}^{-1} \langle v', Y(v,(x_2,\varphi_2^+, \varphi_2^-)) Y(u,(x_1,\varphi_1^+, \varphi_1^-)) w \rangle .
\end{multline}  
\end{prop}

\begin{proof}
Part (a) follows {}from the positive energy axiom (\ref{positive-energy}) for an $N=2$ Neveu-Schwarz vertex operator superalgebra and the truncation condition (\ref{truncation}) for an $N=2$ vertex superalgebra.    For part (b), we note that by weak supercommutativity, there exists $k \in \Z$ such that
\begin{multline}\label{commute}
(x_1 - x_2 - \varphi_1^+ \varphi_2^- - \varphi_1^- \varphi_2^+)^k \langle v', Y(u,(x_1,\varphi_1^+, \varphi_1^-)) Y(v,(x_2,\varphi_2^+, \varphi_2^-)) w\rangle \\
= (-1)^{\eta(u)\eta(v)} (x_1 - x_2 - \varphi_1^+ \varphi_2^- - \varphi_1^- \varphi_2^+)^k \langle v', Y(v,(x_2,\varphi_2^+, \varphi_2^-)) Y(u,(x_1,\varphi_1^+, \varphi_1^-)) w\rangle 
\end{multline}
for all $w \in V$ and $v' \in V'$. {}From (a), we know the left-hand side of (\ref{commute}) involves only
finitely many negative powers of $x_2$ and only finitely many positive powers of $x_1$.  However, the right-hand side of (\ref{commute}) involves only finitely many negative powers of $x_1$ and only finitely
many positive powers of $x_2$.  Thus multiplying both sides of (\ref{commute}) by $(x_1 - x_2 - \varphi_1^+ \varphi_2^- - \varphi_1^- \varphi_2^+)^{-k}$ results in well-defined power series as long as on the left-hand side we expand  $(x_1 - x_2 - \varphi_1^+ \varphi_2^- - \varphi_1^- \varphi_2^+)^{-k}$ in positive powers of $x_2$ and on the right-hand side we expand $(x_1 - x_2 - \varphi_1^+ \varphi_2^- - \varphi_1^- \varphi_2^+)^{-k}$ in positive powers of $x_1$.  The result follows. 
\end{proof}

\begin{prop} \label{iterates}
Let $(V, Y(\cdot,(x,\varphi^+, \varphi^-)), \mathbf{1}, \mu)$ be an $N=2$ Neveu-Schwarz vertex operator superalgebra, and let $u, v, w \in V$, and $v' \in V'$ be arbitrary.  

{\bf (a) (rationality of iterates)} The formal series $\langle v', Y(Y(u,(x_0,\varphi_1^+ - \varphi_2^+, \varphi_1^- - \varphi_2^-))v,(x_2,\varphi_2^+, \varphi_2^-)) w \rangle$, which involves only finitely many negative powers of $x_0$ and only finitely many positive powers of $x_2$, lies in the image of the map $\iota_{20}$: 
\begin{multline}
\langle v', Y(Y(u,(x_0,\varphi_1^+ - \varphi_2^+, \varphi_1^- - \varphi_2^-))v,(x_2,\varphi_2^+, \varphi_2^-))w\rangle\\  
= \iota_{20} h(x_0,\varphi_1^+ - \varphi_2^+, \varphi_1^- - \varphi_2^-,x_2,\varphi_2^+, \varphi_2^-) , 
\end{multline}
where the (uniquely determined) element $h \in \bigwedge_*[x_0,\varphi_1^+, \varphi_1^-,x_2,\varphi_2^+, \varphi_2^-]_{P'}$ is of the form      
\begin{equation}
h(x_0,\varphi_1^+ - \varphi_2^+, \varphi_1^- - \varphi_2^- ,x_2,\varphi_2^+, \varphi_2^-) = \frac{k(x_0,\varphi_1^+ -\varphi_2^+, \varphi_1^- - \varphi_2^-,x_2,\varphi_2^+, \varphi_2^-)}{x_0^r x_2^s (x_0 + x_2 - \varphi_1^+ \varphi_2^- - \varphi_1^- \varphi_2^+)^t} 
\end{equation}
for some $k \in \bigwedge_* [x_0,\varphi_1^+, \varphi_1^-,x_2,\varphi_2^+, \varphi_2^-]$ and $r, s, t
\in \mathbb{Z}$.  

{\bf (b)}  The formal series 
\begin{equation}
\langle v', Y(u,(x_0 + x_2  + \varphi_1^+\varphi_2^- + \varphi_1^- \varphi_2^+, \varphi_1^+, \varphi_1^-)) Y(v,(x_2,\varphi_2^+, \varphi_2^-)) w \rangle,
\end{equation}
which involves only finitely many negative powers of $x_2$ and only finitely many positive powers of $x_0$, lies in the image of $\iota_{02}$, and in fact  
\begin{multline}
\langle v', Y(u,(x_0 + x_2  + \varphi_1^+ \varphi_2^- + \varphi_1^- \varphi_2^+, \varphi_1^+, \varphi_1^-))
Y(v,(x_2,\varphi_2^+, \varphi_2^-)) w \rangle \\
=  \iota_{02} h(x_0,\varphi_1^+ - \varphi_2^+, \varphi_1^- - \varphi_2^-, x_2, \varphi_2^+, \varphi_2^-) . 
\end{multline}    
\end{prop}

\begin{proof}
Part (a) follows {}from the positive energy axiom (\ref{positive-energy}) and truncation condition
(\ref{truncation}) for a vertex operator superalgebra.  For part (b), we note that {}from weak associativity, there exists $k \in \Z$ such that
\begin{multline}\label{assoc}
(x_0 + x_2 + \varphi_1^+ \varphi_2^- + \varphi_1^- \varphi_2^+)^k \langle v', Y(Y(u,(x_0,\varphi_1^+ -
\varphi_2^+, \varphi_1^- - \varphi_2^-))v,(x_2,\varphi_2^+, \varphi_2^-))w \rangle \\
= (x_0 + x_2 + \varphi_1^+ \varphi_2^- + \varphi_1^- \varphi_2^+)^k \langle v', Y(u, (x_0 + x_2  + \varphi_1^+ \varphi_2^- + \varphi_1^- \varphi_2^+, \varphi_1^+, \varphi_1^-)) \\
\cdot Y(v,(x_2, \varphi_2^+, \varphi_2^-)) w  \rangle 
\end{multline}
for all $v' \in V'$. {}From (a), we know the left-hand side of (\ref{assoc}) involves only finitely many negative powers of $x_0$ and only finitely many positive powers of $x_2$.  However, the right-hand
side of (\ref{assoc}) involves only finitely many negative powers of $x_2$ and only finitely many positive powers of $x_0$.  Thus multiplying both sides of (\ref{commute}) by $(x_0 + x_2 + \varphi_1^+
\varphi_2^- + \varphi_1^- \varphi_2^+)^{-k}$ results in a well-defined power series as long as on the left-hand side we expand $(x_0 + x_2 + \varphi_1^+ \varphi_2^- + \varphi_1^- \varphi_2^+)^{-k}$ in positive powers of $x_0$ and on the right-hand side we expand $(x_0 + x_2 + \varphi_1^+ \varphi_2^- + \varphi_1^- \varphi_2^+)^{-k}$ in positive powers of $x_2$.  The result follows.
\end{proof}  

\begin{prop}\label{associativity} 
{\bf (associativity)}  
Let $(V, Y(\cdot,(x,\varphi^+, \varphi^-)), \mathbf{1}, \mu)$ be an $N=2$ Neveu-Schwarz vertex operator superalgebra, and let $u, v, w \in V$, and $v' \in V'$ be arbitrary.  We have the following equality of rational functions: 
\begin{multline}
\iota^{-1}_{12} \langle v', Y(u,(x_1,\varphi_1^+, \varphi_1^-)) Y(v,(x_2,\varphi_2^+, \varphi_2^-)) w \rangle =   \Bigl( \iota^{-1}_{20} \langle v', Y(Y(u,\\
\left. (x_0,\varphi_1^+ - \varphi_2^+, \varphi_1^- - \varphi_2^- ))v,(x_2,\varphi_2^+, \varphi_2^-)) w \rangle \Bigr) \right|_{x_0 = x_1 - x_2 - \varphi_1^+ \varphi_2^- - \varphi_1^- \varphi_2^+} .
\end{multline}
\end{prop}

\begin{proof}
Let $f(x_1,\varphi_1^+, \varphi_1^-,x_2,\varphi_2^+, \varphi_2^-)$ be the rational function in Proposition \ref{supercommutativity}.  Then $f$ satisfies 
\begin{multline}
\iota_{02} f(x_0 + x_2 + \varphi_1^+ \varphi_2^- + \varphi_1^- \varphi_2^+,\varphi_1^+, \varphi_1^-,x_2,\varphi_2^+, \varphi_2^-) \\
= \left. \left( \iota_{12} f(x_1,\varphi_1^+, \varphi_1^-,x_2,\varphi_2^+, \varphi_2^-) \right) \right|_{x_1 = x_0 + x_2 +\varphi_1^+ \varphi_2^- + \varphi_1^- \varphi_2^+} . 
\end{multline}
Thus for $h(x_0,\varphi_1^+ -\varphi_2^+, \varphi_1^- - \varphi_2^-, x_2, \varphi_2^+, \varphi_2^-)$ {}from Proposition \ref{iterates}, we have 
\[h(x_0,\varphi_1^+ - \varphi_2^+, \varphi_1^- - \varphi_2^-, x_2, \varphi_2^+, \varphi_2^-) \! = \! f(x_0 + x_2 + \varphi_1^+ \varphi_2^- + \varphi_1^- \varphi_2^+,\varphi_1^+, \varphi_1^-,x_2,\varphi_2^+, \varphi_2^-).\]  
The result follows {}from Propositions \ref{supercommutativity} and \ref{iterates}. 
\end{proof} 

\begin{rema}\label{OPE-remark2}
{\em Proposition \ref{associativity} asserts that the expansion of
\[Y(Y(u,( x_1 - x_2 - \varphi_1^+ \varphi_2^- - \varphi_1^- \varphi_2^+,\varphi_1^+ - \varphi_2^+, \varphi_1^- - \varphi_2^- ))v,(x_2,\varphi_2^+, \varphi_2^-))   \]
in powers of $x_1 - x_2 - \varphi_1^+ \varphi_2^- - \varphi_1^- \varphi_2^+$ represents the ``operator product" 
\[Y(u,(x_1,\varphi_1^+, \varphi_1^-)) Y(v,(x_2,\varphi_2^+, \varphi_2^-)).\]
Thus using the fact that only the singular terms enter into the commutator formula, one can instead of Proposition \ref{mu-prop},  encode the $N=2$ Neveu-Schwarz algebra relations as  the ``operator product expansion", which in our formulation with two odd formal variables is given by 
\begin{eqnarray}\label{OPE}
\lefteqn{\qquad \ \ Y(\mu, (x_1, \varphi^+_1, \varphi^-_1)) Y(\mu, (x_2, \varphi^+_2, \varphi^-_2)) } \\
&=&\! \! \! \!   \frac{c_V/3}{(x_1 - x_2 - \varphi_1^+ \varphi_2^- - \varphi_1^- \varphi_2^+)^2}
+ \frac{(\varphi_1^+ - \varphi_2^+) Y(\tau^{(+)}, (x_2, \varphi^+_2, \varphi^-_2))}{x_1 - x_2 - \varphi_1^+ \varphi_2^- - \varphi_1^- \varphi_2^+}  \nonumber\\
& & \! \! \! \! \! \! \! \! + \frac{(\varphi_1^- - \varphi_2^-) Y(\tau^{(-)}, (x_2, \varphi^+_2, \varphi^-_2))}{x_1 - x_2 - \varphi_1^+ \varphi_2^- - \varphi_1^- \varphi_2^+}  - \frac{2(\varphi_1^+- \varphi_2^+)(\varphi_1^- - \varphi_2^-) Y(\mu,(x_2, \varphi^+_2, \varphi^-_2))}{(x_1 - x_2 - \varphi_1^+ \varphi_2^- - \varphi_1^- \varphi_2^+)^2} \nonumber\\
& & \! \! \! \! \! \! \! \! -  \frac{2(\varphi_1^+- \varphi_2^+)(\varphi_1^- - \varphi_2^-) \frac{\partial}{\partial x_2} Y(\mu,(x_2, \varphi^+_2, \varphi^-_2))}{x_1 - x_2 - \varphi_1^+ \varphi_2^- - \varphi_1^- \varphi_2^+} 
+ y(x_1, x_2, \varphi_1^+, \varphi_1^-, \varphi_2^+, \varphi_2^-) \nonumber
\end{eqnarray}
where $y(x_1, x_2, \varphi_1^+, \varphi_1^-, \varphi_2^+, \varphi_2^-) \in \mathbb{C}[x_1  - x_2][[x_2, x_2^{-1}]][\varphi_1^+, \varphi_1^-, \varphi_2^+, \varphi_2^-]$. That is $y$ contains only nonsingular terms in $x_1 - x_2 - \varphi_1^+ \varphi_2^- - \varphi_1^- \varphi_2^+$. Note that given an $N=2$ vertex superalgebra, $V$, if one wants to investigate whether a given element $\mu \in V$ is an $N=2$ Neveu-Schwarz vector, then checking that the vertex operators corresponding to $\mu$ satisfy the conditions of Proposition \ref{mu-prop} (which is equivalent to (\ref{OPE})) is much easier in comparison to checking the eleven operator product expansion terms usually presented in the literature as in for instance \cite{Greene}; see also \cite{Bo2}.  With the changes in notation given in Remark \ref{mu-operator-physics-remark}, we see that (\ref{OPE}) is exactly that given in the physics literature as in for instance \cite{YZ}.}
\end{rema}


Note that rationality of products and iterates and supercommutativity and associativity imply weak supercommutativity and weak associativity.   Thus by Proposition \ref{equals-Jacobi1}, we have: 

\begin{thm}\label{duality}
In the presence of the other axioms in the definition of $N=2$ Neveu-Schwarz vertex operator superalgebra with odd variables, the Jacobi identity follows {}from the rationality of products and iterates and supercommutativity and associativity.  In particular, the Jacobi identity may be replaced by these properties. 
\end{thm}

This theorem can also be proved by using the delta-function identity (\ref{delta-3-terms-with-phis}) and the substitution rule (\ref{delta-substitute}).  

\section{$N=2$ (Neveu-Schwarz) vertex (operator) superalgebras without odd formal variables}\label{N2-without-section}

\subsection{$N=2$ vertex superalgebras without  odd formal variables}\label{without-section}

In most of the literature on vertex superalgebras with superconformal properties, the vertex operators are not formulated with odd formal variables \cite{Scheit}, \cite{Adamovic1999}, \cite{Bo1}, \cite{Bo2}.  In fact, although most of the early physics literature contains odd variables in the expressions for fields (cf. \cite{DPZ}, \cite{YZ}), the first formal axiomatic treatment of vertex superalgebras with odd formal variables was in \cite{B-thesis}, \cite{B-vosas}.  In this section and the next, we see how our formulation of the notion of $N=2$ vertex superalgebra and $N=2$ Neveu-Schwarz vertex operator superalgebra relate to these notions without odd formal variables.     

\begin{defn}\label{vertex-superalgebra-definition-without}
{\em An} $N=2$ vertex superalgebra over $\bigwedge_*$  {\em consists of a $\mathbb{Z}_2$-graded (by} sign{\em) $\bigwedge_*$-module
\begin{equation}\label{first-n2-superalgebra-without}
V = V^0 \oplus V^1
\end{equation} 
equipped, first, with a linear map 
\begin{eqnarray}\label{operator-without}
V &\longrightarrow&  (\mbox{End} \; V)[[x,x^{-1}]] \\
v  &\mapsto&  Y(v,x) \nonumber
\end{eqnarray}
with 
\begin{equation}
Y(v,x) = \sum_{n \in \mathbb{Z}}  v_n  x^{-n-1},
\end{equation}
where for the $\mathbb{Z}_2$-grading of $\mathrm{End} \; V$ induced {}from that of $V$, we have
\begin{equation}
v_n \in  (\mbox{End} \; V)^{\eta(v)}
\end{equation}
for $v$ of homogeneous sign in $V$ and $x$ an even formal variable, and where $Y(v,x)$ denotes the} vertex operator associated with $v$.  {\em We also have a distinguished element $\mathbf{1}$ in $V$ (the} vacuum vector{\em ).  The following conditions are assumed for $u,v \in V$:  the} truncation condition: {\em
\begin{equation}\label{truncation-without}
u_n v  = 0 \qquad \mbox{for $n \in \mathbb{Z}$ sufficiently large,} 
\end{equation}
that is
\begin{equation}
Y(u, x)v \in V((x));
\end{equation}
next, the following} vacuum property{\em :
\begin{equation}\label{vacuum-identity-without}
Y(\mathbf{1}, x) = \mathrm{id}_V ;
\end{equation}
the} creation property {\em holds:
\begin{equation}
Y(v,x) \mathbf{1} \in V[[x]], \quad \mathrm{and} \quad \lim_{x\rightarrow 0} Y(v,x) \mathbf{1} = v ; 
\end{equation} 
the} Jacobi identity {\em holds:  
\begin{multline}\label{Jacobi-identity-without}
x_0^{-1} \delta \biggl( \frac{x_1 - x_2}{x_0} \biggr)Y(u,x_1) Y(v,x_2) -  (-1)^{\eta(u)\eta(v)} x_0^{-1} \delta \biggl( \frac{x_2 - x_1}{-x_0} \biggr) Y(v,x_2)Y(u,x_1) \\
= x_2^{-1} \delta \biggl( \frac{x_1 - x_0 }{x_2} \biggr) Y(Y(u,x_0)v,x_2) , 
\end{multline}
for $u,v$ of homogeneous sign in $V$; and finally, defining $\mathcal{D} \in (\mathrm{End} \, V)^0$ by
\begin{equation}\label{define-D-for vosa-without}
\mathcal{D}(v) = v_{-2} \mathbf{1},
\end{equation}
there exist $\mathcal{D}^\pm \in (\mathrm{End} \, V)^1$ satisfying 
\begin{equation}\label{need-D-relations}
\left[\mathcal{D}^+, \mathcal{D}^-\right] = 2 \mathcal{D}, \quad \mathrm{and} \quad \left[\mathcal{D}^\pm , \mathcal{D}^\pm \right] = 0,
\end{equation}
and such that the} $\mathcal{D}^\pm$-bracket relations {\em hold:
\begin{equation}\label{Dpm-bracket-without-definition}
\left[ \mathcal{D}^\pm, Y(v,x)\right] = Y(\mathcal{D}^\pm v,x).
\end{equation}
}
\end{defn}

The $N=2$ vertex superalgebra without odd formal variables just defined is denoted by \[(V,Y(\cdot,x),\mathbf{1}, \mathcal{D}^+, \mathcal{D}^-),\] 
or when no confusion will arise by $V$.

\begin{rema} \label{need-Dpm-if-without}
{\em In the definition above, we require that the two odd endomorphisms $\mathcal{D}^\pm$ be specified and that they satisfy (\ref{need-D-relations}), (\ref{Dpm-bracket-without-definition}) in contrast to the notion of $N=2$ vertex superalgebra with two odd formal variables in which we need not require that such endomorphisms be specified because the existence of $\mathcal{D}^\pm$ and their bracket properties naturally arise {}from the odd variable components of the vertex operators and the other axioms of the definition of $N=2$ vertex superalgebra with two odd formal variables.  } 
\end{rema}

\begin{rema}\label{is-a-vertex-superalgebra-remark}
{\em Given an $N=2$ vertex superalgebra without odd formal variables $(V,Y(\cdot,$ $x),\mathbf{1}, \mathcal{D}^+, \mathcal{D}^-)$ over $\mathbb{C}$, then $(V,Y(\cdot,$ $x),\mathbf{1}, \mathcal{D})$ is a vertex superalgebra in the sense of for instance \cite{Li} and \cite{T}.}
\end{rema}

Let $(V_1, Y_1(\cdot,x),\mbox{\bf 1}_1,\mathcal{D}^+_1,\mathcal{D}^-_1)$ and $(V_2, Y_2(\cdot,x),\mbox{\bf 1}_2,\mathcal{D}^+_2, \mathcal{D}^-_2)$ be two $N=2$ vertex superalgebras.  A {\it homomorphism of $N=2$ vertex superalgebras without odd formal variables} is a $\mathbb{Z}_2$-graded $\bigwedge_\infty$-module homomorphism $\gamma : V_1 \longrightarrow V_2$ satisfying
\[\gamma (Y_1(u,x)v) = Y_2(\gamma(u),x)\gamma(v) \quad \mbox{for}
\quad u,v \in V_1 ,\] 
$\gamma(\mbox{\bf 1}_1) = \mbox{\bf 1}_2$, and $\gamma(\mathcal{D}_1^\pm) = \mathcal{D}_2^\pm$.  

In the next section, we will show that the category of $N=2$ vertex superalgebras without odd formal variables is naturally isomorphic to the category of $N=2$ vertex superalgebras with two odd formal variables.  To prove this, we will need the following consequences of the definition of $N=2$ vertex superalgebra without  odd formal variables:  
Taking $\mathrm{Res}_{x_0}$ of the Jacobi identity and using the $\delta$-function identity (\ref{delta-2-terms}), we obtain the following supercommutator formula 
\begin{equation}\label{homo-bracket-without}
[ Y(u, x_1), Y(v,x_2)]  =  \mbox{Res}_{x_0} x_2^{-1} \delta \biggl( \frac{x_1 - x_0}{x_2} \biggr) Y(Y(u,x_0)v,x_2) . 
\end{equation}
{}From (\ref{need-D-relations}) we have that
\begin{equation}\label{from-D-relations}
(\mathcal{D}^\pm)^2 = 0, \quad \mathrm{and} \quad \left[\mathcal{D}, \mathcal{D}^\pm \right] = 0.
\end{equation}

Following the proof of Proposition 3.1.18 in \cite{LL}, which holds because the vacuum is an even vector in $V$, i.e., $\eta(\mathbf{1}) = 0$, we have that $V$ satisfies the $\mathcal{D}$-derivative property 
\begin{equation}\label{D-derivative-without}
\frac{\partial}{\partial x} Y(v,x) = Y(\mathcal{D}v, x).
\end{equation}

\subsection{$N=2$ Neveu-Schwarz vertex operator superalgebras without odd formal variables}

\begin{defn} {\em An} $N = 2$ Neveu-Schwarz vertex operator superalgebra over $\bigwedge_*$  and without odd variables {\em is a $\frac{1}{2} \mathbb{Z}$-graded $\bigwedge_*$-module (graded by} weights{\em)
\begin{equation} \label{first-n2-without}
V = \coprod_{n \in \frac{1}{2} \mathbb{Z}} V_{(n)} 
\end{equation}
such that
\begin{eqnarray}
& \dim V_{(n)} < \infty \quad \mbox{for} \quad n \in \frac{1}{2} \mathbb{Z} , & \\
& V_{(n)} = 0 \quad \mbox{for $n$ sufficiently negative} , & \label{third-n2-without}
\end{eqnarray}
equipped with an $N=2$ vertex superalgebra structure $(V,Y(\cdot, x), \mathbf{1}, \mathcal{D}^+, \mathcal{D}^-)$, and two distinguished homogeneous vectors $\tau^{(\pm)}  \in V_{(3/2)}^1$ (the {\em $\tau^{(\pm)}$ $N=2$ Neveu-Schwarz generators} or {\em $\tau^{(\pm)}$ $N=2$ superconformal generators}), satisfying the following conditions:  the $N=2$ Neveu-Schwarz algebra relations (\ref{first-n2})--(\ref{last-n2}) hold for $L(n), J(n), G^\pm(n - 1/2) \in \mathrm{End} \; V$, for $n \in \mathbb{Z}$, and $c_V \in \mathbb{C}$ (the} central charge{\em) where 
\begin{equation}\label{tau-no-variables}
G^\pm(n - 1/2) = \tau^{(\pm)}_n \quad \mbox{for} \; n \in \mathbb{Z}, \mbox{i.e.,} \quad Y(\tau^{(\pm)},x) = \sum_{n \in \mathbb{Z}} G^\pm (n + 1/2) x^{- n - \frac{1}{2} - \frac{3}{2}};
\end{equation} 
and
\begin{equation}\label{Dpm-compatibility}
G^\pm(-1/2) = \mathcal{D}^\pm;
\end{equation}
for $n \in \frac{1}{2} \mathbb{Z}$ and $v \in V_{(n)}$
\begin{equation}\label{L0-grading-without}
L(0)v = nv  
\end{equation}
and in addition, $V_{(n)}$ is the direct sum of eigenspaces for $J(0)$ such that if $v \in V_{(n)}$ is also an eigenvector for $J(0)$ with eigenvalue $k$, i.e., if 
\begin{equation}\label{J0-grading-without}
J(0)v = kv,  \qquad \mbox{then $k \equiv 2n \; \mathrm{mod} \; 2$}.
\end{equation}}
\end{defn}
The $N=2$ Neveu-Schwarz vertex operator superalgebra just defined is denoted by 
\[(V,Y(\cdot,x),\mbox{\bf 1},\tau^{(+)}, \tau^{(-)}) . \]

We have the following consequences of the definition:  Taking $\mathrm{Res}_{x_1} x_1$ of both sides of the commutator formula (\ref{homo-bracket-without})  with $u = \tau^{(\pm)}$, we have the following $G^\pm(1/2)$-bracket formula without odd formal variables
\begin{equation}\label{G-1/2-bracket-without}
\left[ G^\pm(1/2), Y(v, x) \right]  = Y(G^\pm (1/2)v, x) + x Y(G^\pm (-1/2)v,x),
\end{equation}
(cf. (\ref{G-1/2-bracket})).

There exist $\mu = \frac{1}{2} G^+(1/2) \tau^{(-)}\in V_{(3/2)}^1$ and $\omega  = L(-2) \mathbf{1} = \frac{1}{4} (G^-(-1/2) \tau^{(+)} + G^+(-1/2) \tau^{(-)}) \in V_{(2)}^0$.  Using (\ref{G-1/2-bracket-without}) and  the $\mathcal{D}^\pm = G^\pm (-1/2)$-bracket relations (\ref{Dpm-bracket-without-definition}), we see that 
\begin{equation}\label{get-mu-omega}
Y(\mu,x) =  \sum_{n \in \mathbb{Z}} J(n) x^{-n-1}, \quad \mathrm{and} \quad 
Y(\omega, x) = \sum_{n \in \mathbb{Z}} L(n) x^{-n-2}. 
\end{equation}

\begin{rema} \label{need-both-if-without}
{\em In the definition above, we require that the two elements $\tau^{(+)}$ and $\tau^{(-)}$ in $V_{(3/2)}$ be specified in contrast to the notion of $N=2$ Neveu-Schwarz vertex operator superalgebra with odd formal variables in which we specify just one vector giving rise to the $N=2$ Neveu-Schwarz algebra, that is $\mu \in V_{(1)}$.  This is because if we do not include odd formal variables we need these two vectors in order to give a full generating set for the $N=2$ Neveu-Schwarz algebra.  On the other hand, specifying the two vectors $\tau^{(\pm)}$ is enough since $G^\pm(n - 1/2)$, for $n \in \mathbb{Z}$, generate the $N=2$ Neveu-Schwarz algebra; see Remark \ref{n2-generators-remark} and (\ref{get-mu-omega}).} 
\end{rema}

\begin{rema}
{\em Given an $N=2$ Neveu-Schwarz vertex operator superalgebra without odd formal variables $(V,Y(\cdot,$ $x),\mathbf{1}, \tau^{(+)}, \tau^{(-)})$ over $\mathbb{C}$, and letting $ \omega = L(-2) \mathbf{1} = \frac{1}{4} (G^-(-1/2) \tau^{(+)} + G^+(-1/2) \tau^{(-)})$, then $(V,Y(\cdot,$ $x),\mathbf{1}, \omega)$ is a vertex operator superalgebra in the sense of \cite{Li} and \cite{T}.}
\end{rema}

\begin{rema}{\em  The notion of $N=2$ Neveu-Schwarz vertex operator superalgebra without odd formal variables is equivalent to the notion of  ``$N=2$ superconformal vertex operator superalgebra" without odd formal variables given in \cite{HM}.  However, we point out the minor simplification that in our notion of $N=2$ Neveu-Schwarz vertex operator superalgebra, we only specify the distinguished vectors $\tau^{(+)}$ and $\tau^{(-)}$ (and the vacuum), whereas in the definition in \cite{HM}, the vector $\mu$ is also specified.   The notion of ``$N=2$ SVOA" given in \cite{Adamovic1999} is equivalent to our notion of  $N=2$ Neveu-Schwarz vertex operator superalgebra without odd formal variables with the exception of the $J(0)$-grading condition (\ref{J0-grading-without}) which is not imposed in \cite{Adamovic1999}.  In addition, we note that in \cite{Adamovic1999}, the definition of $N=2$ SVOA includes four distinguished vectors corresponding to $\mu, \tau^{(\pm)}$, and $\omega$.  As mentioned above as a minor point, one only needs to specify the two vectors $\tau^{(\pm)}$. 
}
\end{rema}

Note that (\ref{Dpm-compatibility}) implies that 
\begin{equation}
L(-1) = \mathcal{D}.
\end{equation}

Let $(V_1, Y_1(\cdot,x),\mbox{\bf 1}_1,\tau^{(+)}_1,\tau^{(-)}_1)$ and $(V_2, Y_2(\cdot,x),\mbox{\bf 1}_2,\tau^{(+)}_2, \tau^{(-)}_2)$ be two $N=2$ Neveu-Schwarz vertex operator  superalgebras.  A {\it homomorphism of $N=2$ Neveu-Schwarz vertex operator superalgebras without odd formal variables} is an $N=2$ vertex superalgebra homomorphism $\gamma : V_1 \longrightarrow V_2$ such that 
$\gamma(\tau_1^{(\pm)}) = \tau_2^{(\pm)}$.

\subsection{Isomorphism between the categories of $N=2$ (Neveu-Schwarz) vertex (operator) superalgebras with two odd formal variables and without odd formal variables}\label{iso-section}

\begin{prop}\label{get-a-vosa-without}
Let $(V,Y(\cdot,(x,\varphi^+, \varphi^-)), \mathbf{1})$ be an $N=2$ vertex superalgebra with two odd formal variables, and let $\mathcal{D}^\pm \in \mathrm{End} \, V$ be defined by $\mathcal{D}^\pm (v)  = v^\pm_{-3/2} \mathbf{1}$.  Then $(V,Y(\cdot,(x,0,0)), \mathbf{1}, \mathcal{D}^+, \mathcal{D}^-)$ is an $N=2$ vertex superalgebra without odd formal variables. 

Let $(V,Y(\cdot,(x,\varphi^+, \varphi^-)), \mbox{\bf 1}, \mu)$ be an $N=2$ Neveu-Schwarz vertex
operator superalgebra with two odd formal variables, and let $\tau^{(\pm)} = \mp G^\pm(-1/2) \mu$.  Then $(V,Y(\cdot,(x,0,0)), \mbox{\bf 1}, \tau^{(+)}, \tau^{(-)})$ is an $N=2$ Neveu-Schwarz vertex operator superalgebra without odd formal variables. 
\end{prop}

\begin{proof}
Setting the odd variables $\varphi^\pm = 0$ in the axioms for $(V,Y(\cdot,(x,\varphi^+, \varphi^-)), \mathbf{1})$, it is clear that $(V,Y(\cdot,(x,0,0)), \mbox{\bf 1})$ satisfies all of the axioms for an $N=2$ vertex superalgebra without odd formal variables except for the axioms involving $\mathcal{D}$ and $\mathcal{D}^\pm$.  The remaining axioms (\ref{define-D-for vosa-without})--(\ref{Dpm-bracket-without-definition}) follow {}from Proposition \ref{Dpm-derivative-prop}, Corollary \ref{D-derivative-cor}, and  (\ref{Dpm-bracket-without}).

Similarly, it is trivial that $(V,Y(\cdot,(x,0,0)), \mbox{\bf 1}, \tau^{(+)}, \tau^{(-)})$ satisfies axioms (\ref{first-n2-without})--(\ref{third-n2-without}), as well as axioms (\ref{L0-grading-without}) and (\ref{J0-grading-without}).  By the argument above, $V$ is an $N=2$ vertex superalgebra without odd formal variables.   The axioms regarding the $N=2$ Neveu-Schwarz relations follow {}from consequences of the definition of an $N=2$ Neveu-Schwarz vertex operator superalgebra with two odd formal variables which show that the vertex operators for $\tau^{(\pm)} = \mp G^\pm(-1/2) \mu$ are given by (\ref{tau-pm-with}), and that $G^\pm(-1/2) = \mathcal{D}^\pm$ (\ref{D-correspondence}).  
\end{proof}

Proposition \ref{get-a-vosa-without} implies that all of the consequences of the definition of $N=2$ vertex superalgebra with two odd formal variables hold for an $N=2$ vertex superalgebra without odd formal variables if we set the odd formal variables equal to zero.  In particular, we have the following corollary to Proposition \ref{weak-supercommutativity} which we single out because we will use it in the proof of the next proposition.

\begin{cor}\label{weak-supercommutativity-without}
{\bf (weak supercommutativity)} Let $(V,Y(\cdot,x), \mathbf{1}, \mathcal{D}^+, \mathcal{D}^-)$ be an $N=2$ vertex  superalgebra without odd formal variables.  There exists $k \in \Z$ such that 
\begin{equation}\label{weak-supercommutativity-formula-without}
(x_1 - x_2)^k \left[Y(u,x_1),Y(v,x_2) \right] = 0 .
\end{equation}
Furthermore this weak supercommutativity follows {}from the truncation condition (\ref{truncation-without}) and the Jacobi identity (\ref{Jacobi-identity-without}). 
\end{cor} 

This corollary can also be thought of as following {}from Remark  \ref{is-a-vertex-superalgebra-remark} and the results on vertex superalgebras in \cite{Li}.

\begin{prop}\label{get-a-vosa-with}
Let $(V,Y(\cdot,x), \mathbf{1}, \mathcal{D}^+, \mathcal{D}^-)$ be an $N=2$ vertex superalgebra without odd formal variables.  Define 
\begin{multline}\label{Y-from-VOSA-without}
\tilde{Y}(v, (x,\varphi^+, \varphi^-)) = Y(v,x) + \varphi^+ Y(\mathcal{D}^+ v,x)  + \varphi^- Y(\mathcal{D}^- v,x) \\
+ \frac{1}{2} \varphi^+ \varphi^- Y((\mathcal{D}^- \mathcal{D}^+ - \mathcal{D}^+  \mathcal{D}^-) v, x).
\end{multline}
Then $(V,\tilde{Y}(\cdot,(x,\varphi^+, \varphi^-)), \mbox{\bf 1})$ is an $N=2$ vertex superalgebra with two odd formal variables. 

Furthermore, if $(V,Y(\cdot,x), \mathbf{1}, \tau^{(+)}, \tau^{(-)})$ is an $N=2$ Neveu-Schwarz vertex operator superalgebra without odd formal variables, then letting 
\begin{equation}
\mu  = \frac{1}{2} G^+(1/2) \tau^{(-)},
\end{equation}
gives $(V,\tilde{Y}(\cdot,(x,\varphi^+, \varphi^-)), \mbox{\bf 1}, \mu)$ the structure of an $N=2$ Neveu-Schwarz vertex operator superalgebra with two odd formal variables. 
\end{prop}

\begin{proof}
Assume that $(V,Y(\cdot,x), \mathbf{1}, \mathcal{D}^+, \mathcal{D}^-)$ is an $N=2$ vertex superalgebra without odd formal variables.  To show that $(V,\tilde{Y}(\cdot,(x, \varphi^+, \varphi^-)), \mathbf{1})$ is an $N=2$ vertex superalgebra with odd formal variables, we first note that axiom (\ref{first-n2-superalgebra-with}) follows trivially.  {}From the definition of $\tilde{Y}$, we have that  $v^\pm_{n-1/2} = (\mathcal{D}^\pm v)_n$, and $v_n^{+-} = 1/2((\mathcal{D}^- \mathcal{D}^+ - \mathcal{D}^+  \mathcal{D}^-) v)_n$.  Axioms (\ref{operator-with})--(\ref{truncation}) follow easily {}from this and the fact that $\mathcal{D}^\pm \in (\mathrm{End} \, V)^1$.

By the vacuum property (\ref{vacuum-identity-without}) for $(V,Y(\cdot,x), \mathbf{1},\mathcal{D}^+, \mathcal{D}^- )$, we have that $\mathbf{1}_{-3/2} = 0$ and thus $\mathcal{D}^\pm  \mathbf{1} = \mathbf{1}^\pm_{-3/2} \mathbf{1} = 0$. Therefore since 
\[\tilde{Y}(\mathbf{1}, (x.\varphi^+, \varphi^-)) = Y(\mathbf{1},x) + \varphi^+ Y(0, x) + \varphi^- Y(0, x) + \varphi^+ \varphi^- Y(0,x) = Y(\mathbf{1},x) = \mathrm{id}_V \]
the vacuum property (\ref{vacuum-identity}) holds. 

The creation property (\ref{creation-property1}) and (\ref{creation-property2}) holds since by the creation property for $(V,Y(\cdot,x), \mathbf{1}, \mathcal{D}^+, \mathcal{D}^-)$, we have that  $Y(v,x)\mathbf{1}$, $Y(\mathcal{D}^\pm v, x) \mathbf{1}$, and $Y((\mathcal{D}^- \mathcal{D}^+ - \mathcal{D}^+  \mathcal{D}^-)v,x) \mathbf{1}$ are all in $V[[x]]$, implying
\[\tilde{Y}(v,(x,\varphi^+, \varphi^-)) \mathbf{1}  \in V[[x]] [\varphi^+, \varphi^-] ,\]
and 
\[\lim_{(x,\varphi^+, \varphi^-) \rightarrow 0} \tilde{Y}(v,(x,\varphi^+, \varphi^-)) \mathbf{1}
= \lim_{x \rightarrow 0} Y(v,x)\mathbf{1} = v . \]

We have proved all the axioms for an $N=2$ vertex superalgebra with two odd formal variables except for the Jacobi identity.   To prove the Jacobi identity with odd formal variables, we will show that weak supercommutativity (\ref{weak-supercommutativity-formula}) and the $\mathcal{D}^\pm$-bracket-derivative properties (\ref{Dpm-bracket-derivative}) hold for $(V,\tilde{Y}(\cdot,(x, \varphi^+, \varphi^-)), \mathbf{1})$.  This will imply the Jacobi identity with odd formal variables by Theorem \ref{get-Jacobi2}.

To prove weak supercommutativity for $(V,\tilde{Y}(\cdot,(x,\varphi^+, \varphi^-)), \mbox{\bf 1})$, we let $S_1 = \{u, \mathcal{D}^+u, \mathcal{D}^-u, (\mathcal{D}^- \mathcal{D}^+- \mathcal{D}^+\mathcal{D}^-)u \}$ and $S_2 = \{v, \mathcal{D}^+v, \mathcal{D}^-v, (\mathcal{D}^- \mathcal{D}^+- \mathcal{D}^+\mathcal{D}^-)v \}$.  By Corollary \ref{weak-supercommutativity-without}, there exists $k(w_1, w_2) \in \Z$ such that 
\[(x_1-x_2)^k \left[Y(w_1,x_1), Y(w_2, x_2) \right] = 0\]
for any pair $(w_1, w_2) \in S_1 \times S_2$.  Let $k = \mathrm{max} \{ k(w_1, w_2) \; | \;  (w_1, w_2) \in S_1 \times S_2 \}$.  Now take $n = k+2$.  This implies that 
\begin{eqnarray*}
\lefteqn{(x_1 - x_2 - \varphi_1^+ \varphi_2^- - \varphi_1^- \varphi_2^+)^n \left[\tilde{Y}(u,(x_1,\varphi_1^+, \varphi_1^-)),\tilde{Y}(v,(x_2,\varphi_2^+, \varphi_2^-)) \right]}\\
&=& \Bigl( (x_1 - x_2)^n  - n(\varphi_1^+ \varphi_2^- +\varphi_1^- \varphi_2^+)(x_1- x_2)^{n-1} + n(n-1) \varphi_1^+ \varphi_1^- \varphi_2^+ \varphi_2^-(x_1\\
& & \quad -\, x_2 )^{n-2}\Bigr) \Bigl[ Y(u,x_1) + \varphi^+_1 Y(\mathcal{D}^+ u,x_1)  + \varphi^-_1Y(\mathcal{D}^- u,x_1) \\
& & \quad + \, \frac{1}{2} \varphi^+_1 \varphi^-_1 Y((\mathcal{D}^- \mathcal{D}^+ - \mathcal{D}^+  \mathcal{D}^-) u, x_1), Y(v,x_2) \\
& & \quad + \, \varphi^+_2 Y(\mathcal{D}^+ v,x_2)  + \varphi^-_2 Y(\mathcal{D}^- v,x_2) 
+ \frac{1}{2} \varphi^+_2 \varphi^-_2 Y((\mathcal{D}^- \mathcal{D}^+ - \mathcal{D}^+  \mathcal{D}^-) v, x_2) \Bigr]\\
&=& 0,
\end{eqnarray*}
proving weak supercommutativity for $(V,\tilde{Y}(\cdot,(x,\varphi^+, \varphi^-)), \mbox{\bf 1})$.

Now we prove the $\mathcal{D}^\pm$-bracket-derivative properties (\ref{Dpm-bracket-derivative}) by noting that (\ref{need-D-relations}), the $\mathcal{D}^\pm$-bracket relations (\ref{Dpm-bracket-without-definition}) and the $\mathcal{D}$-derivative property (\ref{D-derivative-without})  for $(V,Y(\cdot,x), \mathbf{1},$ $\mathcal{D}^+, \mathcal{D}^-)$ imply that
\begin{eqnarray*}
\lefteqn{ \left[ \mathcal{D}^\pm, \tilde{Y}(v,(x,\varphi^+, \varphi^-)) \right] }\\
&=&  Y(\mathcal{D}^\pm v,x) - \varphi^+ Y(\mathcal{D}^\pm\mathcal{D}^+ v,x)  - \varphi^- Y(\mathcal{D}^\pm \mathcal{D}^- v,x)  \\
& & \quad +\,  \frac{1}{2} \varphi^+ \varphi^- Y(\mathcal{D}^\pm (\mathcal{D}^- \mathcal{D}^+ - \mathcal{D}^+  \mathcal{D}^-) v, x) \\
&=&  Y(\mathcal{D}^\pm v,x) \pm \frac{1}{2} \varphi^\mp Y((\mathcal{D}^- \mathcal{D}^+ - \mathcal{D}^+ \mathcal{D}^-) v,x)  - \frac{1}{2} \varphi^\mp Y((\mathcal{D}^+ \mathcal{D}^- \\
& & \quad + \, \mathcal{D}^- \mathcal{D}^+) v,x)   \pm  \frac{1}{2} \varphi^+ \varphi^- Y(\mathcal{D}^\pm \mathcal{D}^\mp \mathcal{D}^\pm v, x) \\
&=&  Y(\mathcal{D}^\pm v,x) \pm \frac{1}{2} \varphi^\mp Y((\mathcal{D}^- \mathcal{D}^+ - \mathcal{D}^+ \mathcal{D}^-) v,x)  -  \varphi^\mp Y(\mathcal{D}v,x) \\
& & \quad  \pm \, \varphi^+ \varphi^- Y( \mathcal{D} \mathcal{D}^\pm v, x) \\
&=&  Y(\mathcal{D}^\pm v,x) \pm \frac{1}{2} \varphi^\mp Y((\mathcal{D}^- \mathcal{D}^+ - \mathcal{D}^+ \mathcal{D}^-) v,x)  -  \varphi^\mp \frac{\partial}{\partial x}Y(v,x) \\
& & \quad  \pm \, \varphi^+ \varphi^- \frac{\partial}{\partial x} Y(  \mathcal{D}^\pm v, x) \\
&=& \left( \frac{\partial}{\partial \varphi^\pm} - \varphi^\mp \frac{\partial}{\partial x} \right) \Bigl(Y(v,x) + \varphi^+ Y(\mathcal{D}^+ v,x)  + \varphi^- Y(\mathcal{D}^- v,x)  \\
& & \quad + \, \frac{1}{2} \varphi^+ \varphi^- Y((\mathcal{D}^- \mathcal{D}^+ - \mathcal{D}^+  
\mathcal{D}^-) v, x) \Bigr)
\end{eqnarray*}
\begin{eqnarray*}
&=& \left( \frac{\partial}{\partial \varphi^\pm} - \varphi^\mp \frac{\partial}{\partial x} \right) \tilde{Y}(v,(x, \varphi^+, \varphi^-)). 
\end{eqnarray*}
The Jacobi identity for $(V,\tilde{Y}(\cdot,(x,\varphi^+, \varphi^-)), \mbox{\bf 1})$ follows by Theorem \ref{get-Jacobi2}, proving that $(V,\tilde{Y}(\cdot,(x,\varphi^+, \varphi^-)), \mbox{\bf 1})$ is an $N=2$ vertex superalgebra with two odd formal variables.

Now assume $(V,Y(\cdot,x), \mathbf{1}, \tau^{(+)}, \tau^{(-)})$ is an $N=2$ Neveu-Schwarz vertex operator superalgebra without odd formal variables.  To show that $(V,\tilde{Y}(\cdot,(x, \varphi^+, \varphi^-)), \mathbf{1}, \mu)$ is an $N=2$ Neveu-Schwarz vertex operator superalgebra with two odd formal variables, we note that axioms (\ref{vosa1})--(\ref{positive-energy}), (\ref{grading-for-vosa-with}) and (\ref{J(0)-grading}) follow immediately.  For the $N=2$ Neveu-Schwarz element $\mu$, using the $\mathcal{D}^\pm= G^\pm(-1/2)$-bracket properties (\ref{Dpm-bracket-without-definition}), axiom (\ref{Dpm-compatibility}) and equation (\ref{get-mu-omega}), we have
\begin{eqnarray*}
\lefteqn{\tilde{Y} (\mu, (x,\varphi^+,\varphi^-)) }\\
&=& Y(\mu,x) + \varphi^+ Y(\mathcal{D}^+\mu,x) + \varphi^-Y(\mathcal{D}^-\mu,x) + \frac{1}{2} \varphi^+\varphi^- Y((\mathcal{D}^- \mathcal{D}^+ \\
& & \quad - \, \mathcal{D}^+ \mathcal{D}^-) \mu,x) \\
&=& Y(\mu,x) + \varphi^+ \left[ \mathcal{D}^+, Y(\mu,x) \right] + \varphi^-\left[ \mathcal{D}^-, Y(\mu,x) \right] \\
& & \quad + \, \frac{1}{2} \varphi^+\varphi^- \left( \left[ \mathcal{D}^-, \left[ \mathcal{D}^+ , Y(\mu,x) \right] \right] - \left[ \mathcal{D}^+, \left[ \mathcal{D}^-, Y(\mu,x) \right] \right] \right)  \\
&=& \sum_{n \in \mathbb{Z}} \biggl( J(n)  + \varphi^+ \left[ G^+(-1/2), J(n) \right] + \varphi^-\left[ G^-(-1/2 ), J(n) \right]   \\
& & \quad +\,  \frac{1}{2} \varphi^+\varphi^- \biggl( \bigl[ G^-(-1/2), \bigl[ G^+(-1/2) , J(n) \bigr] \bigr] - \bigl[ G^+(-1/2), \bigl[ G^-(-1/2), \\
& & \quad J(n) \bigr] \bigr] \biggr) \biggr)x^{-n-1}  \\
&=& \sum_{n \in \mathbb{Z}} \biggl( J(n)  - \varphi^+ G^+(n-1/2) + \varphi^-G^-(n-1/2 ) +  \frac{1}{2} \varphi^+\varphi^- \biggl( -\bigl[ G^-(-1/2), \\
& & \quad G^+(n-1/2) \bigr] - \bigl[ G^+(-1/2),  G^-(n-1/2) \bigr] \biggr) \biggr) x^{-n-1} \\
&=& \sum_{n \in \mathbb{Z}} \biggl( J(n)  - \varphi^+ G^+(n- 1/2) + \varphi^-G^-(n-1/2) +  \frac{1}{2} \varphi^+\varphi^- \Bigl( -(2L(n-1) \\
& & \quad + \,  n J(n-1)) - (2L(n-1) - n J(n-1)) \Bigr)\biggr) x^{-n-1}  \\
&=& \sum_{n \in \mathbb{Z}} \biggl( J(n) x^{-n-1} - \varphi^+ G^+(n+ 1/2) x^{-n-2} + \varphi^-G^-(n+1/2) x^{-n-2}  \\
& & \quad + \,  \frac{1}{2} \varphi^+\varphi^- ( -4L(n-1))x^{-n-1} \biggr)  \\
\end{eqnarray*}
which gives (\ref{homo-mu}).

Finally, the $G^\pm(-1/2)$-derivative properties (\ref{G(-1/2)-derivative}) follow {}from axiom (\ref{Dpm-compatibility}) for $(V, Y(\cdot,x), \mathbf{1}, \tau^{(+)}, \tau^{(-)})$, i.e., $\mathcal{D}^\pm = G^\pm (-1/2)$, and Proposition \ref{Dpm-derivative-prop} applied to the $N=2$ vertex superalgebra $(V,\tilde{Y}(\cdot,(x, \varphi^+, \varphi^-)), \mathbf{1})$, proving that $(V, \tilde{Y}(\cdot,(x,\varphi^+$, $\varphi^-)), \mathbf{1}, \mu)$ is an $N=2$ Neveu-Schwarz vertex operator superalgebra with two odd formal variables.  
\end{proof}

We note that we also have 
\begin{equation}
\mu  = - \frac{1}{2} G^-(1/2) \tau^{(+)}.
\end{equation}

\begin{rema}{\em
In the proof of Proposition \ref{get-a-vosa-with}, we could have proved the Jacobi identity for $(V, \tilde{Y}(\cdot, (x, \varphi^+, \varphi^-)), \mathbf{1})$ directly using (\ref{Y-from-VOSA-without}) and the commutation relations and bracket properties for $\mathcal{D}^\pm$.  However, this is a long, tedious, and not very illuminating calculation.  Thus we have instead chosen to prove the Jacobi identity using the equivalent properties of supercommutativity and the $\mathcal{D}^\pm$-bracket-derivative properties.   The ease of proving Proposition \ref{get-a-vosa-with}  in this manner is one of the reasons we chose first to develop the notion of $N=2$ vertex superalgebra with two odd formal variables and formulate weak supercommutativity directly {}from this notion, rather than first formulate the notion of $N=2$ vertex superalgebra without odd formal variables and then derive the notion with two odd formal variables and the property of weak supercommutativity with odd formal variables {}from the weak supercommutativity without odd formal variables.  }
\end{rema}

Let $\mathrm{VSA}_2(\varphi^+, \varphi^-)$ denote the category of $N=2$ vertex superalgebras with two odd formal variables, let $\mathrm{VSA}_2$ denote the category of $N=2$ vertex superalgebras without odd formal variables,  let $\mathrm{VOSA}_2(\varphi^+, \varphi^-,c)$ denote the category of $N=2$ Neveu-Schwarz vertex operator superalgebras with two odd formal variables and with central charge $c \in \mathbb{C}$, and let $\mathrm{VOSA}_2(c)$ denote the category of $N=2$ Neveu-Schwarz vertex operator superalgebras without odd formal variables and with central charge $c \in \mathbb{C}$.   Note that $\mathrm{VOSA}_2(\varphi^+, \varphi^-,c)$ is a subcategory of $\mathrm{VSA}_2 (\varphi^+, \varphi^-)$ and $\mathrm{VOSA}_2(c)$ is a subcategory of $\mathrm{VSA}_2$. 

Define $F_{\varphi^\pm =0}: \mathrm{VSA}_2(\varphi^+, \varphi^-) \longrightarrow \mathrm{VSA}_2$ by 
\begin{eqnarray*}
F_{\varphi^\pm =0}: (V,Y(\cdot,(x,\varphi^+, \varphi^-)), \mathbf{1}) &\mapsto& (V,Y(\cdot,(x,0,0)),
\mathbf{1}, \mathcal{D}^+, \mathcal{D}^-)\\
\gamma &\mapsto& \gamma 
\end{eqnarray*}
for $(V,Y(\cdot,(x,\varphi^+, \varphi^-)), \mathbf{1})$ an object in $\mathrm{VSA}_2(\varphi^+, \varphi^-)$ and $\gamma$ a morphism, where $\mathcal{D}^\pm \in (\mathrm{End} \, V)^1$ are defined by $\mathcal{D}^\pm v = v^\pm_{-3/2} \mathbf{1}$.   Proposition \ref{get-a-vosa-without} shows that $F_{\varphi^\pm =0}$ takes objects in $\mathrm{VSA}_2 (\varphi^+,\varphi^-)$ to objects in $\mathrm{VSA}_2$.  It is clear that $F_{\varphi^\pm =0}$ takes morphisms to morphisms and that $F_{\varphi^\pm =0}$ is a functor.   Moreover, Proposition \ref{get-a-vosa-without} implies that $F_{\varphi^\pm =0}$ restricted to $\mathrm{VOSA}(\varphi^+, \varphi^-,c)$ gives a functor {}from $\mathrm{VOSA}_2(\varphi^+$, $\varphi^-,c)$ to $\mathrm{VOSA}_2(c)$ via
\[F_{\varphi^\pm =0} \! : \! (V,Y(\cdot,(x,\varphi^+, \varphi^-)), \mathbf{1}, \mu) \mapsto (V,Y(\cdot,(x,0,0)),
\mathbf{1}, -G^+(-1/2)\mu, G^-(-1/2)\mu).\]

Define $F_{\varphi^+, \varphi^-} : \mathrm{VSA}_2 \longrightarrow \mathrm{VSA}_2(\varphi^+, \varphi^-)$ by 
\begin{eqnarray*}
F_{\varphi^+, \varphi^-} : (V,Y(\cdot,x), \mathbf{1}, \mathcal{D}^+, \mathcal{D}^-) &\mapsto& (V,\tilde{Y} (\cdot,(x,\varphi^+, \varphi^-)), \mathbf{1}), \\
\gamma &\mapsto& \gamma
\end{eqnarray*} 
where $\tilde{Y}(v,(x,\varphi^+, \varphi^-))$ is defined by (\ref{Y-from-VOSA-without}).  Proposition \ref{get-a-vosa-with} shows that $F_{\varphi^+, \varphi^-}$ takes  objects in $\mathrm{VSA}_2$ to objects in $\mathrm{VSA}_2 (\varphi^+,$ $ \varphi^-)$.   It is easy to show that $F_{\varphi^+, \varphi^-}$ takes morphisms to morphisms and that $F_{\varphi^+, \varphi^-}$ is a functor. Moreover, Proposition \ref{get-a-vosa-with} implies that $F_{\varphi^+, \varphi^-}$ restricted to $\mathrm{VOSA}(c)$ gives a functor {}from $\mathrm{VOSA}_2(c)$ to $\mathrm{VOSA}_2(\varphi^+, \varphi^-,c)$ via
\[F_{\varphi^+, \varphi^-} : (V,Y(\cdot,x), \mathbf{1}, \tau^{(+)}, \tau^{(-)}) \mapsto (V,\tilde{Y} (\cdot,
(x,\varphi^+, \varphi^-)), \mathbf{1}, \frac{1}{2} G^+ (1/2) \tau^{(-)}) . \]

For any category $\mathcal{C}$, let $1_{\mathcal{C}}$ denote the identity functor on the category $\mathcal{C}$.   It is easy to show that $F_{\varphi^\pm =0} \circ F_{\varphi^+, \varphi^-} = 1_{\mathrm{VSA}_2}$ and $F_{\varphi^+, \varphi^-} \circ F_{\varphi^\pm =0} = 1_{\mathrm{VSA}_2(\varphi^+, \varphi^-)}$, and that the corresponding restrictions of $F_{\varphi^\pm =0}$ and $F_{\varphi^+, \varphi^-}$ to the subcategories of $N=2$ Neveu-Schwarz vertex operator superalgebras with and without odd formal variables, respectively, satisfy $F_{\varphi^\pm =0} \circ F_{\varphi^+, \varphi^-} = 1_{\mathrm{VOSA}_2(c)}$ and $F_{\varphi^+, \varphi^-} \circ F_{\varphi^\pm =0} = 1_{\mathrm{VOSA}_2(\varphi^+, \varphi^-,c)}$.  Thus we have the following theorem:

\begin{thm}\label{superalgebras} 
The two categories $\mathrm{VSA}_2(\varphi^+, \varphi^-)$ and $\mathrm{VSA}_2$ are isomorphic.   In addition, for any $c \in \mathbb{C}$, the categories $\mathrm{VOSA}_2(\varphi^+, \varphi^-,c)$ and $\mathrm{VOSA}_2(c)$ are isomorphic. 
\end{thm}

\section{The nonhomogeneous coordinate system}\label{nonhomo-section}

So far our formulations have been in what we call the ``homogeneous" coordinate system, denoted by formal $N=2$ variables $(x, \varphi^+, \varphi^-)$ and the ``homogeneous" basis for the $N=2$ Neveu-Schwarz algebra given by $L_n, J_n, G^\pm_{n-1/2}$, for $n \in \mathbb{Z}$ and $d$ and relations (\ref{Virasoro-relation})--(\ref{Neveu-Schwarz-relation-last}).   In this section, we transfer some of our results to the ``nonhomogeneous" $N=2$ formal variables which we denote by $(x, \varphi^{(1)}, \varphi^{(2)})$ and the nonhomogeneous  basis for the $N=2$ Neveu-Schwarz algebra given by $L_n, J_n, G^{(j)}_{n-1/2}$, for $n \in \mathbb{Z}$ and $j = 1,2$, and $d$ and relations (\ref{Virasoro-relation2})--(\ref{transformed-Neveu-Schwarz-relation-last}).   This is a standard transformation in $N=2$ superconformal field theory \cite{C}, \cite{DRS}.  We chose to call these coordinate systems ``homogeneous" and ``nonhomogeneous", respectively, due to the transformation properties of a nonhomogeneous $N=2$ superconformal function on the corresponding nonhomogeneous superconformal operators $D^{(j)}$ for $j=1,2$ as described in \cite{B-moduli}, and due to the action of the $J_n$ terms, for $n \in \mathbb{Z}$,  in the algebra of infinitesimals as described in Remark \ref{homo-infinitesimals-remark}.   In addition, below in Remark \ref{last-nonhomo-terminology-remark}, we give another motivation for our terminology which is that the $J(0)$-weight (induced {}from the $J(0)$-grading) is no longer homogeneous on the components of the nonhomogeneous vertex operators.

\subsection{$N=2$ vertex superalgebras in the nonhomogeneous coordinates}\label{nonhomo-vertex-superalgebra-section}

The transformation {}from homogeneous to nonhomogeneous coordinates is given by
\begin{equation}\label{substitute-varphi}
\varphi^{(1)} =  \frac{1}{\sqrt{2}}(\varphi^+ +   \varphi^-) \qquad \mbox{and} \qquad \varphi^{(2)} =  - \frac{i}{\sqrt{2}} (\varphi^+ -  \varphi^-) , 
\end{equation}
or equivalently
\begin{equation}
\varphi^\pm = \frac{1}{\sqrt{2}}(\varphi^{(1)} \pm i \varphi^{(2)}).
\end{equation}
In addition, given an $N=2$ vertex superalgebra with two odd formal variables $(V, Y(\cdot, (x, \varphi^+, \varphi^-)), \mathbf{1}, \mu)$, for $v \in V$, we define $v_{n-1/2}^{(j)}, v_n^{(1,2)} \in \mathrm{End} \, V$, for $j = 1,2$, and $n \in \mathbb{Z}$, using the components of the vertex operator $Y(v, (x, \varphi^+, \varphi^-))$ as follows
\begin{equation}\label{substitute-v1}
v_{n - \frac{1}{2}}^{(1)} = \frac{1}{\sqrt{2}} \left( v_{n - \frac{1}{2}}^+  + v_{n - \frac{1}{2}}^- \right), \quad \mathrm{and} \quad  v_{n - \frac{1}{2}}^{(2)} = \frac{i}{\sqrt{2}} \left( v_{n - \frac{1}{2}}^+ - v_{n - \frac{1}{2}}^- \right)
\end{equation}
and
\begin{equation} \label{substitute-v3}
v_{n}^{(1,2)} = -i v_n^{+-},
\end{equation}
or equivalently
\begin{equation}
v^\pm_{n - \frac{1}{2}} = \frac{1}{\sqrt{2}} \left( v_{n - \frac{1}{2}}^{(1)}  \mp i v_{n - \frac{1}{2}}^{(2)} \right), \quad \mathrm{and} \quad v_{n}^{+-} = i v_n^{(1,2)}. 
\end{equation} 

Then letting
\begin{equation}
Y(v,(x,\varphi^{(1)},\varphi^{(2)})) = \sum_{n \in \mathbb{Z}} \Bigl( v_n + \varphi^{(1)}  v_{n - \frac{1}{2}}^{(1)} +  \varphi^{(2)} v_{n - \frac{1}{2}}^{(2)} + \varphi^{(1)} \varphi^{(2)} v_{n - 1}^{(1,2)} \Bigr) x^{-n-1}
\end{equation} 
we have that
\begin{equation}
Y(v, (x, \varphi^{(1)}, \varphi^{(2)})) = Y(v, (x,\varphi^+, \varphi^-)) |_{\varphi^\pm = \frac{1}{\sqrt{2}} (\varphi^{(1)}  \pm i \varphi^{(2)})}.
\end{equation}

Performing these substitutions we have the following notion of {\it $N=2$ vertex superalgebra with two odd formal variables in the nonhomogeneous coordinate system} which is equivalent to the notion of $N=2$ vertex superalgebra with two odd formal variables in the homogeneous coordinate system.

\begin{defn}\label{nonhomo-vertex-superalgebra-definition}
{\em An} $N=2$ vertex superalgebra over $\bigwedge_*$ and with two odd formal variables in the nonhomogeneous coordinate system {\em consists of a $\mathbb{Z}_2$-graded $\bigwedge_*$-module (graded by} sign $\eta${\em) 
\begin{equation}\label{nonhomo-first-n2-superalgebra-with}
V = V^0 \oplus V^1
\end{equation} 
equipped, first, with a linear map 
\begin{eqnarray}\label{nohomo-operator-with}
V &\longrightarrow&  (\mbox{End} \; V)[[x,x^{-1}]][\varphi^{(1)}, \varphi^{(2)}] \\
v  &\mapsto&  Y(v,(x,\varphi^{(1)}, \varphi^{(2)})) \nonumber
\end{eqnarray}
with 
\begin{equation}
Y(v,(x,\varphi^{(1)}, \varphi^{(2)})) = \sum_{n \in \mathbb{Z}} \Bigl( v_n  + \varphi^{(1)} v_{n - \frac{1}{2}}^{(1)} + \varphi^{(2)}  v_{n - \frac{1}{2}}^{(2)} + \varphi^{(1)} \varphi^{(2)} v_{n-1}^{(1,2)}  \Bigr)  x^{-n-1},
\end{equation}
where for the $\mathbb{Z}_2$-grading of $\mathrm{End} \; V$ induced {}from that of $V$, we have
\begin{equation}
v_n, v_n^{(1,2)} \in (\mbox{End} \; V)^{\eta(v)}, \quad \mathrm{and} \quad 
v_{n - \frac{1}{2}}^{(j)} \in  (\mbox{End} \; V)^{(\eta(v) + 1) \mbox{\begin{footnotesize} mod 
\end{footnotesize}} 2},
\end{equation}
for $j = 1,2$, and for $v$ of homogeneous sign in $V$, $x$ is an even formal variable, and $\varphi^{(j)}$, for $j=1,2$  are odd formal variables, and where $Y(v,(x,\varphi^{(1)}, \varphi^{(2)}))$ denotes the} nonhomogeneous vertex operator associated with $v$.  {\em We also have a distinguished element $\mathbf{1}$ in $V$ (the} vacuum vector{\em ).  The following conditions are assumed for $u,v \in V$:  the} truncation conditions: {\em for $j = 1,2$
\begin{equation}\label{nonhomo-truncation}
u_n v = u_{n}^{(1,2)} v =  u_{n-\frac{1}{2}}^{(j)} v = 0 \qquad \mbox{for $n \in \mathbb{Z}$ sufficiently large,} 
\end{equation}
that is
\begin{equation}
Y(u, (x, \varphi^{(1)}, \varphi^{(2)}))v \in V((x))[\varphi^{(1)}, \varphi^{(2)}];
\end{equation}
next, the following} vacuum property{\em :
\begin{equation}\label{nonhomo-vacuum-identity}
Y(\mathbf{1}, (x, \varphi^{(1)}, \varphi^{(2)})) = \mathrm{id}_V ;
\end{equation}
the} creation property {\em holds:
\begin{eqnarray}
Y(v,(x,\varphi^{(1)},\varphi^{(2)})) \mathbf{1} &\in& V[[x]][\varphi^{(1)},\varphi^{(2)}] \label{nonhomo-creation-property1} \\
\lim_{(x,\varphi^{(1)},\varphi^{(2)}) \rightarrow 0} Y(v,(x, \varphi^{(1)},\varphi^{(2)})) \mathbf{1} &=& v ;  \label{nonhomo-creation-property2}
\end{eqnarray} 
and finally the} Jacobi identity {\em holds:  
\begin{eqnarray}\label{nonhomo-Jacobi-identity}
\\
\lefteqn{\quad x_0^{-1} \delta \biggl( \frac{x_1 - x_2 - \varphi^{(1)}_1 \varphi^{(1)}_2 - \varphi^{(2)}_1 \varphi^{(2)}_2}{x_0} \biggr) Y(u,(x_1, \varphi^{(1)}_1, \varphi^{(2)}_1)) Y(v,(x_2, \varphi^{(1)}_2, \varphi^{(2)}_2)) }\nonumber\\
& & \quad -  (-1)^{\eta(u)\eta(v)} x_0^{-1} \delta \biggl( \frac{x_2 -  x_1 + \varphi^{(1)}_1 \varphi^{(1)}_2 + \varphi^{(2)}_1 \varphi^{(2)}_2}{-x_0} \biggr) Y(v,(x_2, \varphi^{(1)}_2,\varphi^{(2)}_2)) \nonumber \\
& & \quad \cdot Y(u,(x_1,\varphi^{(1)}_1,\varphi^{(2)}_1)) \nonumber\\
\quad &=&  \! \! \! \! x_2^{-1} \delta \biggl( \frac{x_1 - x_0 - \varphi^{(1)}_1 \varphi^{(1)}_2 - \varphi^{(2)}_1 \varphi^{(2)}_2}{x_2} \biggr) Y(Y(u,(x_0, \varphi^{(1)}_1 - \varphi^{(1)}_2,\varphi^{(2)}_1 - \varphi^{(2)}_2))v, \nonumber\\
& & \quad (x_2, \varphi^{(1)}_2,\varphi^{(2)}_2)) , \nonumber
\end{eqnarray}
for $u,v$ of homogeneous sign in $V$.
}
\end{defn}

The $N=2$ vertex superalgebra with two odd formal variables in the nonhomogeneous coordinate system just defined is denoted by \[(V,Y(\cdot,(x,\varphi^{(1)}, \varphi^{(2)})),\mathbf{1}).\]

\begin{rema} 
{\em
It is clear that the transformations  (\ref{substitute-varphi}), (\ref{substitute-v1}), (\ref{substitute-v3}), 
induce an isomorphism between the category of $N=2$ vertex superalgebras with two odd formal variables in the homogeneous coordinate system and the category of $N=2$ vertex superalgebras with two odd formal variables in the nonhomogeneous coordinate system.   }
\end{rema}

\begin{rema}\label{superconformal-remark-nonhomo}
{\em The substitution  $x_0 = x_1 - x_2 - \varphi_1^{(1)} \varphi_2^{(1)} - \varphi_1^{(2)} \varphi_2^{(2)}$ can be thought of as the even part of a nonhomogeneous ``superconformal shift of  $(x_1, \varphi_1^{(1)}, \varphi_1^{(2)})$ by $(x_2, \varphi_2^{(1)}, \varphi_2^{(2)})$".  That is, if we let $f(x, \varphi^{(1)}, \varphi^{(2)}) = (\tilde{x}, \tilde{\varphi}^{(1)}, \tilde{\varphi}^{(2)})$ with $\tilde{x} \in (\bigwedge_*[[x, x^{-1}]][\varphi^{(1)}, \varphi^{(2)}])^0$ and $\tilde{\varphi}^{(j)} \in (\bigwedge_*[[x, x^{-1}]][\varphi^{(1)}, \varphi^{(2)}])^1$, for $j = 1,2$, then $f$ is said to be {\it $N=2$ superconformal} in the nonhomogeneous coordinates $x, \varphi^{(1)}$ and $\varphi^{(2)}$ if and only if $f$ satisfies 
\begin{multline}\label{nice-superconformal-condition-nonhomo}
D^{(1)} \tilde{\varphi}^{(1)} - D^{(2)} \tilde{\varphi}^{(2)} = D^{(1)}\tilde{\varphi}^{(2)} + D^{(2)} \tilde{\varphi}^{(1)} \\
= D^{(j)}\tilde{x} - \tilde{\varphi}^{(1)} D^{(j)} \tilde{\varphi}^{(1)} - \tilde{\varphi}^{(2)} D^{(j)} \tilde{\varphi}^{(2)} =  0 
\end{multline}     
for 
\begin{equation}\label{define-D-derivative-nonhomo}
D^{(j)} = \frac{\partial}{\partial \varphi^{(j)}} + \varphi^{(j)} \frac{\partial}{\partial x} \quad \mbox{for $j = 1,2$}
\end{equation} 
(see \cite{B-moduli}).  For a formal superanalytic vector-valued function in the nonhomogeneous coordinates $f(x, \varphi^{(1)}, \varphi^{(2)}) = (\tilde{x}, \tilde{\varphi}^{(1)}, \tilde{\varphi}^{(2)})$, the conditions (\ref{nice-superconformal-condition-nonhomo}) are equivalent to requiring that $f$ transform the super-differential operators $D^{(j)}$ for $j = 1,2$, homogeneously of degree one in $D^{(1)}$ and in $D^{(2)}$, so that $D^{(1)} = (D^{(1)} \tilde{\varphi}^{(1)}) \tilde{D}^{(1)} - (D^{(2)} \tilde{\varphi}^{(1)}) \tilde{D}^{(2)}$ and $D^{(2)} = (D^{(1)} \tilde{\varphi}^{(2)}) \tilde{D}^{(1)} - (D^{(2)} \tilde{\varphi}^{(2)}) \tilde{D}^{(2)}$.  We observe that $f(x_1,\varphi_1^{(1)}, \varphi_1^{(2)}) = (x_1 - x_2 - \varphi_1^{(1)} \varphi_2^{(1)} - \varphi_1^{(2)} \varphi_2^{(2)}, \varphi_1^{(1)} - \varphi_2^{(1)}, \varphi_1^{(2)} - \varphi_2^{(2)})$ is formally $N=2$ superconformal in $x_1, \varphi_1^{(1)}$ and $\varphi_1^{(2)}$ since it satisfies (\ref{nice-superconformal-condition-nonhomo}). }
\end{rema}

We have the following consequences of the definition of $N=2$ vertex superalgebra with two odd formal variables in the nonhomogeneous coordinate system:
The {\it supercommutator formula} for $N=2$ vertex operators in the nonhomogeneous coordinates
\begin{multline}\label{nonhomo-bracket}
[ Y(u, (x_1,\varphi^{(1)}_1, \varphi^{(2)}_1)), Y(v,(x_2,\varphi^{(1)}_2, \varphi^{(2)}_2))]  \\
=  \mbox{Res}_{x_0} x_2^{-1} \delta \biggl( \frac{x_1 - x_0 - \varphi^{(1)}_1 \varphi^{(1)}_2 - \varphi^{(2)}_1 \varphi^{(2)}_2}{x_2} \biggr) \\
\cdot Y(Y(u,(x_0, \varphi^{(1)}_1 - \varphi^{(1)}_2, \varphi^{(2)}_1 - \varphi^{(2)}_2))v,(x_2, \varphi^{(1)}_2, \varphi^{(2)}_2)) . 
\end{multline}
The {\it iterate formula} for $N=2$ vertex operators in the nonhomogeneous coordinates
\begin{multline}\label{nonhomo-iterate}
Y(Y(u,(x_0, \varphi^{(1)}_1 - \varphi^{(1)}_2, \varphi^{(2)}_1 - \varphi^{(2)}_2))v,(x_2, \varphi_2^{(1)}, \varphi^{(2)}_2)) \\
= \mathrm{Res}_{x_1} \biggl( x_0^{-1} \delta \biggl( \frac{x_1 - x_2 - \varphi^{(1)}_1 \varphi^{(1)}_2 - \varphi^{(2)}_1 \varphi^{(2)}_2}{x_0} \biggr)Y(u,(x_1, \varphi^{(1)}_1, \varphi^{(2)}_1))  \\
\hspace{.5in} \cdot Y(v,(x_2, \varphi^{(1)}_2,\varphi^{(2)}_2)) 
-  (-1)^{\eta(u)\eta(v)} x_0^{-1} \delta \biggl( \frac{x_2 - x_1 + \varphi^{(1)}_1 \varphi^{(1)}_2 + \varphi^{(2)}_1 \varphi^{(2)}_2}{-x_0} \biggr) \\
\hspace{2.45in}\cdot Y(v,(x_2, \varphi^{(1)}_2,\varphi^{(2)}_2)) Y(u,(x_1,\varphi^{(1)}_1,\varphi^{(2)}_1)) \biggr).
\end{multline}

We have the following corollary to Proposition \ref{Dpm-derivative-prop}:
\begin{cor}\label{nonhomo-derivative-cor}
Let $(V, Y(\cdot, (x, \varphi^{(1)}, \varphi^{(2)})), \mathbf{1})$ be an $N=2$ vertex superalgebra in the nonhomogeneous coordinate system and let $\mathcal{D}^{(j)}$, for $j=1,2$, be the odd endomorphisms of $V$ defined by
\begin{equation}\label{nonhomo-define-Ds}
\mathcal{D}^{(j)} (v) = v^{(j)}_{-\frac{3}{2}} \mathbf{1} \qquad \mbox{for $v \in V$.}
\end{equation}
Then 
\begin{equation}
\mathcal{D}^{(1)} = \frac{1}{\sqrt{2}} (\mathcal{D}^+ + \mathcal{D}^-) \quad \mathrm{and} \quad \mathcal{D}^{(2)} =  \frac{i}{\sqrt{2}} (\mathcal{D}^+ - \mathcal{D}^-), 
\end{equation}
i.e.,
\begin{equation}
 \mathcal{D}^\pm = \frac{1}{\sqrt{2}}(\mathcal{D}^{(1)} \mp i \mathcal{D}^{(2)}),
\end{equation}
and for $j = 1,2$,
\begin{equation}\label{nonhomo-Ds-derivative-property}
Y(\mathcal{D}^{(j)} v,(x,\varphi^{(1)}, \varphi^{(2)})) = \biggl( \frac{\partial}{\partial \varphi^{(j)}} + \varphi^{(j)} \frac{\partial}{\partial x} \biggr) Y(v,(x,\varphi^{(1)}, \varphi^{(2)}))  . 
\end{equation} 
\end{cor}

\begin{rema}\label{geometry-correspondence-remark-nonhomo}
{\em The $\mathcal{D}^{(j)}$-derivative properties (\ref{nonhomo-Ds-derivative-property}), for $j = 1,2$, show that the operators $\mathcal{D}^{(j)} \in (\mathrm{End} \; V)^1$ correspond to the $N=2$ superconformal operators $D^{(j)}$ as defined in (\ref{define-D-derivative-nonhomo});  see Remarks \ref{superconformal-remark-nonhomo} and \ref{geometry-correspondence-remark}. }
\end{rema}

Note that for $j = 1,2$
\begin{equation}
\left( \frac{\partial}{\partial \varphi^{(j)}} + \varphi^{(j)} \frac{\partial}{\partial x} \right)^2 =  \frac{\partial }{\partial x},
\end{equation}
and for $j, k = 1,2$
\begin{equation}\label{Djk-bracket-get-D-nonhomo}
\left[ \mathcal{D}^{(j)}, \mathcal{D}^{(k)} \right] = \delta_{j,k} 2 \mathcal{D},
\end{equation}
where $\mathcal{D} = \frac{1}{2} \left[ \mathcal{D}^+, \mathcal{D}^- \right]$.  Of course the $\mathcal{D}$-derivative property still holds under the change of variables {}from homogeneous to nonhomogeneous; that is we still have
\begin{equation}\label{nonhomo-D-derivative}
Y(\mathcal{D}v, (x, \varphi^{(1)}, \varphi^{(2)})) = \frac{\partial}{\partial x} Y(v, (x, \varphi^{(1)}, \varphi^{(2)})) .
\end{equation}

\begin{rema}\label{N2-subalgebras-representation-remark-nonhomo}
{\em As noted in Remark \ref{N2-subalgebras-representation-remark}, an $N=2$ vertex superalgebra with two odd formal variables is a representation of $\mathfrak{osp}_{\bigwedge_*}(2|2)_{<0}$.  In the nonhomogeneous basis, this representation is given by $\mathcal{D}^{(j)} \mapsto G^{(j)} _{-1/2}$, for $j=1,2$, and $\mathcal{D} \mapsto L_{-1}$; see Remark \ref{N2-subalgebras-remark-nonhomo}.
}
\end{rema}

In addition, we have that
\begin{equation}
\mathcal{D}^{(j)} (\mathbf{1})  = 0 \quad \mbox{for $j= 1,2$},
\end{equation}
\begin{eqnarray}
\lefteqn{\ Y(e^{x_0 \mathcal{D} + \varphi_0^{(1)} \mathcal{D}^{(1)} + \varphi_0^{(2)} \mathcal{D}^{(2)}}v, (x,\varphi^{(1)}, \varphi^{(2)})) } \label{nonhomo-exponential-L(-1)-G(-1/2)} \\
&=& e^{x_0 \frac{\partial}{\partial x} + \varphi^{(1)}_0 \left(
\frac{\partial}{\partial \varphi^{(1)}} + \varphi^{(1)} \frac{\partial}{\partial x} \right)+ \varphi^{(2)}_0 \left(
\frac{\partial}{\partial \varphi^{(2)}} + \varphi^{(2)} \frac{\partial}{\partial x} \right)} Y(v,(x,\varphi^{(1)}, \varphi^{(2)})) \nonumber \\  
&=& Y(v,(x + x_0 + \varphi_0^{(1)} \varphi^{(1)} + \varphi_0^{(2)} \varphi^{(2)}, \varphi_0^{(1)} + \varphi^{(1)}, \varphi_0^{(2)} + \varphi^{(2)})) , \nonumber
\end{eqnarray}
and
\begin{equation}\label{for-skew-symmetry-nonhomo}
e^{x \mathcal{D} + \varphi^{(1)} \mathcal{D}^{(1)} +  \varphi^{(2)} \mathcal{D}^{(2)}} v \; = \; Y(v,(x,\varphi^{(1)}, \varphi^{(2)})) \mathbf{1} .
\end{equation}
The following {\it skew supersymmetry} holds for nonhomogeneous $N=2$ vertex operators 
\begin{multline}\label{nonhomo-skew-supersymmetry}
Y(u,(x,\varphi^{(1)}, \varphi^{(2)}))v \\
=(-1)^{\eta(u) \eta(v)}  e^{x \mathcal{D} + \varphi^{(1)} \mathcal{D}^{(1)} + \varphi^{(2)} \mathcal{D}^{(2)}} Y(v,(-x,-\varphi^{(1)}, - \varphi^{(2)}))u 
\end{multline}
for $u,v$ of homogeneous sign in $V$.

In analogy to Proposition \ref{bracket-prop}, we have the following {\em $\mathcal{D}^{(j)}$-} and {\em $\mathcal{D}$-bracket-derivative properties} for $j = 1,2$:
\begin{eqnarray}
\quad \left[ \mathcal{D}^{(j)}, Y(v, (x, \varphi^{(1)}, \varphi^{(2)})) \right] &=& \left( \frac{\partial}{\partial \varphi^{(j)}} - \varphi^{(j)} \frac{\partial}{\partial x} \right) Y(v, (x, \varphi^{(1)}, \varphi^{(2)})) \label{D*-bracket-derivative}\\
\left[ \mathcal{D}, Y(v, (x, \varphi^{(1)}, \varphi^{(2)})) \right] &=& \frac{\partial}{\partial x} Y(v, (x,\varphi^{(1)}, \varphi^{(2)})) , \label{nonhomo-D-bracket-derivative}
\end{eqnarray}
and the following {\em $\mathcal{D}^{(j)}$-} and {\em $\mathcal{D}$-bracket properties} for $j = 1,2$: 
\begin{eqnarray}
\left[ \mathcal{D}^{(j)}, Y(v, (x, \varphi^{(1)}, \varphi^{(2)})) \right] &=& Y(\mathcal{D}^{(j)} v, (x, \varphi^{(1)}, \varphi^{(2)})) \label{D*-bracket}\\
& & \quad - \, 2 \varphi^{(j)} Y(\mathcal{D}v, (x, \varphi^{(1)}, \varphi^{(2)})) \nonumber\\
\left[ \mathcal{D}, Y(v, (x, \varphi^{(1)}, \varphi^{(2)})) \right] &=& Y(\mathcal{D}v, (x, \varphi^{(1)}, \varphi^{(2)})) \label{nonhomo-D-bracket}.
\end{eqnarray}
The following conjugation formula holds
\begin{eqnarray}\label{conjugate-shift-nonhomo}
\lefteqn{ e^{x_0  \mathcal{D} + \varphi^{(1)}_0  \mathcal{D}^{(1)} +  \varphi^{(2)}_0  \mathcal{D}^{(2)}} Y(v,(x,\varphi^{(1)}, \varphi^{(2)})) e^{- x_0  \mathcal{D} - \varphi_0^{(1)}  \mathcal{D}^{(1)} -  \varphi^{(2)}_0  \mathcal{D}^{(2)}}}  \label{nonhomo-conjugation}\\
&=& Y(e^{(x_0 + 2 \varphi^{(1)} \varphi_0^{(1)} + 2 \varphi^{(2)} \varphi_0^{(2)})  \mathcal{D} + \varphi^{(1)}_0  \mathcal{D}^{(1)} + \varphi^{(2)}_0  \mathcal{D}^{(2)}}v,(x,\varphi^{(1)}, \varphi^{(2)})) \nonumber \\
&=& Y(v,(x +x_0  + \varphi^{(1)} \varphi_0^{(1)} + \varphi^{(2)} \varphi_0^{(2)}, \varphi^{(1)} + \varphi_0^{(1)}, \varphi^{(2)} + \varphi_0^{(2)})).  \nonumber
\end{eqnarray}

Again using the notation 
\[Y(v,x) = Y(v,(x, \varphi^+, \varphi^-)) |_{\varphi^+ = \varphi^- = 0} = Y(v,(x, \varphi^{(1)}, \varphi^{(2)})) |_{\varphi^{(1)} = \varphi^{(2)} = 0},\] 
we have for $j = 1,2$
\begin{equation}
\left[ \mathcal{D}^{(j)}, Y(v,x)\right] = Y(\mathcal{D}^{(j)}v, x),
\end{equation}
and 
\begin{multline}\label{nonhomo-odd-relationship1}
Y(v,(x,\varphi^{(1)}, \varphi^{(2)})) = \sum_{n \in \mathbb{Z}} \Bigl( v_n  +
\varphi^{(1)}  [\mathcal{D}^{(1)} , v_n ]  +  \varphi^{(2)} [\mathcal{D}^{(2)}, v_n ]  \\
- \varphi^{(1)} \varphi^{(2)} [\mathcal{D}^{(1)}, [ \mathcal{D}^{(2)}, v_n ] ] \Bigr) x^{-n-1},
\end{multline}
i.e., 
\begin{equation}\label{nonhomo-Dv-bracket-condition}
v^{(j)}_{n - 1/2} =  [\mathcal{D}^{(j)}, v_n] \quad \mbox{for $j=1,2$}, \quad \mathrm{and} \quad v^{(1,2)}_{n-1} = -  [\mathcal{D}^{(1)}, [ \mathcal{D}^{(2)}, v_n ] ],  
\end{equation}  
or equivalently  
\begin{equation}
v^{(j)}_{n - 1/2} =  (\mathcal{D}^{(j)} v)_n \quad \mbox{for $j=1,2$}, \quad \mathrm{and} \quad 
v^{(1,2)}_{n-1} = - (\mathcal{D}^{(1)} \mathcal{D}^{(2)} v)_n,
\end{equation}
i.e.,
\begin{multline}\label{constructing-nonhomo-operators}
Y(v,(x,\varphi^{(1)}, \varphi^{(2)})) = Y(v,x) + \varphi^{(1)} Y(\mathcal{D}^{(1)} v,x) + \varphi^{(2)} Y(\mathcal{D}^{(2)} v,x) \\
- \varphi^{(1)} \varphi^{(2)} Y(\mathcal{D}^{(1)} \mathcal{D}^{(2)} v,x) .
\end{multline}

For $u,v \in V$ and $n \in \mathbb{Z}$, we have for instance
\begin{eqnarray}
\left[u^{(j)}_{-\frac{1}{2}}, v^{(j)}_{n-\frac{1}{2}} \right] &=& (-1)^{\eta(u)+1}(u^{(j)}_{-\frac{1}{2}} v)^{(j)}_{n-\frac{1}{2}} \quad \mbox{for $j =1,2$,} \label{for-mistake-remark}\\
\left[u^{(1,2)}_{-1}, v_n^{(1,2)}\right] &=& (u^{(1,2)}_{-1}v)_n^{(1,2)}
\end{eqnarray}
and in particular,
\begin{eqnarray}
\left[u^{(j)}_{-\frac{1}{2}}, v^{(j)}_{-\frac{1}{2}} \right] &=& (-1)^{\eta(u)+1} (u^{(j)}_{-\frac{1}{2}} v)^{(j)}_{-\frac{1}{2}} \quad \mbox{for $j =1,2$,} \label{form-Lie4}\\
\left[u^{(1,2)}_{-1}, v_{-1}^{(1,2)}\right] &=& (u^{(1,2)}_{-1}v)_{-1}^{(1,2)}. \label{form-Lie5}
\end{eqnarray}
Thus the operators $u^{(1)}_{-1/2}$ form a Lie superalgebra, as do the the operators $u^{(2)}_{-1/2}$, and the operators $u^{(1,2)}_{-1}$, respectively.   

\begin{rema} {\em There is a mistake in the equations (34) and (35) in \cite{B-vosas} which give the corresponding formulas to (\ref{for-mistake-remark}) and (\ref{form-Lie4}) in the $N=1$ case.  In \cite{B-vosas}, the term $(-1)^{\eta(u) +1}$ was erroneously omitted in each of these formulas.}
\end{rema}

Of course weak supercommutativity and weak associativity can be formulated in the nonhomogeneous coordinate system as follows:  Let $(V,Y(\cdot,(x, \varphi^{(1)}, \varphi^{(2)})), \mathbf{1})$ be an $N=2$ vertex  superalgebra with two odd formal variables in the nonhomogeneous coordinate system.  There exists $k \in \Z$ such that 
\begin{equation}\label{weak-supercommutativity-formula-nonhomo}
(x_1 - x_2 - \varphi_1^{(1)} \varphi_2^{(1)} - \varphi_1^{(2)} \varphi_2^{(2)})^k \Bigl[Y(u,(x_1,\varphi_1^{(1)}, \varphi_1^{(2)}))
,Y(v,(x_2,\varphi_2^{(1)}, \varphi_2^{(2)})) \Bigr] = 0 ,
\end{equation}
and this weak supercommutativity follows {}from the truncation condition (\ref{nonhomo-truncation}) and the Jacobi identity (\ref{nonhomo-Jacobi-identity}).   For $u, w \in V$, there exists $k \in \Z$  such that for any $v \in V$ 
\begin{eqnarray}\label{weak-associativity-formula-nonhomo}
\lefteqn{ (x_0 + x_2 + \varphi_1^{(1)} \varphi_2^{(1)} + \varphi_1^{(2)} \varphi_2^{(2)})^k Y(Y(u,(x_0,\varphi_1^{(1)} - \varphi_2^{(1)}, \varphi_1^{(2)} - \varphi_2^{(2)})v,}\\
& & \quad  (x_2,\varphi_2^{(1)}, \varphi_2^{(2)}))w \nonumber\\
\qquad &=&  \! \!  (x_0 + x_2 + \varphi_1^{(1)} \varphi_2^{(1)} + \varphi_1^{(2)} \varphi_2^{(2)})^k Y(u,(x_0 + x_2 +\varphi_1^{(1)} \varphi_2^{(1)} + \varphi_1^{(2)} \varphi_2^{(2)} ,\nonumber\\
& & \quad  \varphi_1^{(1)}, \varphi_1^{(2)})) Y(v,(x_2,\varphi_2^{(1)}, \varphi_2^{(2)}))w  , \nonumber
\end{eqnarray}
and this weak associativity follows {}from the truncation condition (\ref{nonhomo-truncation}) and the Jacobi identity (\ref{nonhomo-Jacobi-identity}).   Thus we have the following corollary to Proposition \ref{equals-Jacobi1} and Theorem \ref{get-Jacobi2}:

\begin{cor}\label{equals-Jacobi1-nonhomo}  
In the presence of the other axioms in the definition of $N=2$ vertex superalgebra with two odd formal variables in the nonhomogeneous coordinate system, the Jacobi identity (\ref{nonhomo-Jacobi-identity}) is equivalent to weak supercommutativity (\ref{weak-supercommutativity-formula-nonhomo}) and weak associativity (\ref{weak-associativity-formula-nonhomo}).  Furthermore, the Jacobi identity for an $N=2$ vertex superalgebra with two odd formal variables in the nonhomogeneous coordinate system follows {}from weak supercommutativity (\ref{weak-supercommutativity-formula-nonhomo}) in the presence of the other axioms together with the $\mathcal{D}^{(j)}$-bracket-derivative properties (\ref{D*-bracket-derivative}), for $j=1,2$.  In particular, in the definition of the notion of $N=2$ vertex superalgebra with two odd formal variables in the nonhomogeneous coordinate system, the Jacobi identity can be replaced by these properties.
\end{cor}

In addition, we have the obvious statements of rationality of products, supercommutativity, rationality of iterates, and associativity obtained {}from performing the change of coordinates {}from homogeneous to nonhomogeneous in Section \ref{duality-section}, as well as the analogue of Theorem \ref{duality}.

\begin{rema}\label{nonhomo-equivalence-remark}
{\em The notion of $N=2$ vertex superalgebra is not equivalent to the notion of ``$N=2$ superconformal vertex algebra" given in \cite{Kac1997} since in \cite{Kac1997} the $\mathcal{D}^{(j)}$-bracket-derivative relation in is given (in our notation) by 
\begin{equation}
\left[ \mathcal{D}^{(j)}, Y(v, (x, \varphi^{(1)}, \varphi^{(2)})) \right] = \left( \frac{\partial}{\partial \varphi^{(j)}} + \varphi^{(j)} \frac{\partial}{\partial x} \right) Y(v, (x,  \varphi^{(1)}, \varphi^{(2)}))
\end{equation} 
(\cite{Kac1997} eq. (5.9.7) and p. 184) which differs {}from (\ref{D*-bracket-derivative}) by a minus sign.  However, in this case, the minus sign is very important since it determines the correspondence between the endomorphisms $\mathcal{D}^{(j)}$ and the $N=2$ superderivations $D^{(j)}$.  Our notion of $N=2$ vertex superalgebra is equivalent to the notion of ``$N_K=2$ SUSY vertex algebra" in \cite{HK} by Corollary \ref{equals-Jacobi1-nonhomo}. }
\end{rema}

\subsection{$N=2$ Neveu-Schwarz vertex operator superalgebras in the nonhomogeneous coordinates}

Now let
\begin{eqnarray}
G^{(1)}(n - 1/2) &=& \frac{1}{\sqrt{2}} \left( G^+(n - 1/2) + G^- (n - 1/2) \right) \label{substitute-G1}\\
G^{(2)}(n - 1/2) &=& \frac{i}{\sqrt{2}} \left( G^+(n - 1/2) - G^- (n - 1/2) \right), \label{substitute-G2}
\end{eqnarray}
or equivalently
\begin{equation}
G^\pm(n - 1/2) = \frac{1}{\sqrt{2}} \left( G^{(1)}(n -1/2) \mp i G^{(2)} (n -1/2) \right).
\end{equation}

Performing these substitutions along with (\ref{substitute-varphi}), (\ref{substitute-v1}), (\ref{substitute-v3}), we have the following notion of {\it $N=2$ Neveu-Schwarz vertex operator superalgebra with two odd formal variables in the nonhomogeneous coordinate system} which is equivalent to the notion of $N=2$ Neveu-Schwarz vertex operator superalgebra with two odd formal variables in the homogeneous coordinate system given in Definition \ref{VOSA-definition}.

\begin{defn}
{\em An} $N = 2$ Neveu-Schwarz vertex operator superalgebra over $\bigwedge_*$ and with two odd formal variables in the nonhomogeneous coordinate system {\em is a $\frac{1}{2} \mathbb{Z}$-graded $\bigwedge_*$-module (graded by} weights{\em)  
\begin{equation}
V = \coprod_{n \in \frac{1}{2}\mathbb{Z}} V_{(n)} 
\end{equation}
such that 
\begin{equation}
\dim V_{(n)} < \infty \qquad \mbox{for $n \in \frac{1}{2} \mathbb{Z}$,} 
\end{equation}
\begin{equation}
V_{(n)} = 0 \qquad \mbox{for $n$ sufficiently negative} , 
\end{equation}
equipped with an $N=2$ vertex superalgebra structure $(V, Y(\cdot, (x,$ $\varphi^{(1)}, \varphi^{(2)})), \mathbf{1})$ in nonhomogeneous coordinates, and a distinguished vector $\mu \in V_{(1)}^0$ (the {\em  $N=2$ Neveu-Schwarz element} or {\em  $N=2$ superconformal element}), satisfying the following conditions:  the $N=2$ Neveu-Schwarz algebra relations in the nonhomogeneous basis hold: 
\begin{eqnarray}
\left[L(m) ,L(n) \right] \! &=&  \! (m - n)L(m + n) \! + \!  \frac{1}{12} (m^3 - m) \delta_{m + n, 0} \; c_V  ,\label{nonhomo-n2-first}\\
\left[J(m) , J(n) \right]  \!  &=&  \! \frac{1}{3} m \delta_{m + n , 0} \;c_V, \\ 
\left[L(m) , J(n) \right]  \!  &=&  \! -nJ(m+n) ,\\
\qquad \qquad \left[L(m),G^{(j)}(n + 1/2)\right] \! &=&  \! \Bigl(\frac{m}{2} - n - \frac{1}{2} \Bigr) G^{(j)} (m + n + 1/2) , \\
\left[ J(m) , G^{(1)} (n + 1/2)\right]   \! &=&  \! - i G^{(2)} (m + n + 1/2),\\
\left[ J(m) , G^{(2)} (n + 1/2)\right]   \! &=& \!  i G^{(1)} (m + n + 1/2), \\
& & \hspace{-1.9in} \left[ G^{(j)}(m + 1/2) , G^{(j)} (n - 1/2) \right]  \ \, =  \ \,2L(m + n)  + \frac{1}{3}(m^2 + m) \delta_{m + n , 0}  \;c _V,  \\
& & \hspace{-1.9in} \left[ G^{(1)}(m + 1/2) , G^{(2)}(n - 1/2) \right]  \ \ = \ \, -i(m-n+1) J(m+n),  \label{nonhomo-n2-last}
\end{eqnarray}
for $m,n \in \mathbb{Z}$, and $j=1,2$, where 
\begin{eqnarray}
J(n) = \mu_n, &\quad&  iG^{(2)}(n - 1/2) = \mu^{(1)}_{n-\frac{1}{2}}, \\
- iG^{(1)}(n - 1/2) = \mu^{(2)}_{n-\frac{1}{2}}, &\quad& 2iL(n) = \mu_{n}^{(1,2)} , 
\end{eqnarray}
i.e., 
\begin{multline}\label{nonhomo-mu}
Y(\mu,(x,\varphi^{(1)},\varphi^{(2)})) = \sum_{n \in \mathbb{Z}} \Bigl( J(n) x^{- n - 1} + i\varphi^{(1)} G^{(2)} (n+1/2) x^{-n-2}  \\ 
- i\varphi^{(2)} G^{(1)}(n+1/2) x^{-n-2}    + 2 i \varphi^{(1)} \varphi^{(2)}  L(n) x^{- n - 2} \Bigr)
\end{multline}
and $c_V \in \mathbb{C}$ (the} central charge{\em);  for $n \in \frac{1}{2} \mathbb{Z}$ and $v \in V_{(n)}$
\begin{equation}
L(0)v = nv
\end{equation}
and in addition, $V_{(n)}$ is the direct sum of eigenspaces for $J(0)$ such that if $v \in V_{(n)}$ is also an eigenvector for $J(0)$ with eigenvalue $k$, i.e., if
\begin{equation}\label{J(0)-grading-nonhomo}
J(0)v = kv, \quad \mbox{then $k \equiv 2n \; \mathrm{mod} \; 2$};
\end{equation} 
and finally for $j=1,2$, the} $G^{(j)}(-1/2)$-derivative properties {\em hold:
\begin{equation}\label{G-derivative}
\biggl( \frac{\partial}{\partial \varphi^{(j)}} + \varphi^{(j)} \frac{\partial}{\partial x} \biggr) Y(v,(x,\varphi^{(1)},\varphi^{(2)})) = Y(G^{(j)}(- 1/2)v,(x,\varphi^{(1)},\varphi^{(2)})) .\\
\end{equation} }
\end{defn}

\medskip

An $N=2$ Neveu-Schwarz vertex operator superalgebra with two odd formal variables in the nonhomogeneous coordinate system is denoted by
\[(V,Y(\cdot,(x,\varphi^{(1)},\varphi^{(2)})),\mathbf{1},\mu) . \]

\begin{rema}\label{transform-isomorphism-remark} 
{\em
It is clear that the transformations  (\ref{substitute-varphi}), (\ref{substitute-v1}), (\ref{substitute-v3}), 
(\ref{substitute-G1}), and (\ref{substitute-G2}) induce an isomorphism between the category of $N=2$ Neveu-Schwarz vertex operator superalgebras with two odd formal variables in nonhomogeneous coordinates and the category of $N=2$ Neveu-Schwarz vertex operator superalgebras with two odd formal variables in homogeneous coordinates.   }
\end{rema}

Using the notation and results of Section \ref{nonhomo-vertex-superalgebra-section}, we observe that 
\begin{equation}\label{nonhomo-D-correspondence}
G^{(j)}(-1/2) = \mathcal{D}^{(j)} \quad \mbox{for $j = 1,2$}.
\end{equation}

\begin{rema}\label{HK-remark} {\em  Our notion of $N=2$ Neveu-Schwarz vertex operator superalgebra with two odd formal variables in the nonhomogeneous coordinate system coincides with the notion of ``strongly conformal $N_K=2$ SUSY vertex algebra" given in \cite{HK} and \cite{He}, with the exception of the additional grading condition (\ref{J(0)-grading-nonhomo}).  However, in \cite{HK} and \cite{He} there are several confusing discrepancies concerning the vertex operator corresponding to the $N=2$ superconformal element $\mu$.  To see the discrepancies, we note that in \cite{HK} and \cite{He} their $S^{(j)}$'s correspond to our $\mathcal{D}^{(j)}$'s for $j = 1,2$, and the $\mathcal{D}^{(j)}$-bracket-derivative properties in \cite{HK} and \cite{He} coincide with ours  (see p. 107 of \cite{HK}) (as opposed to in \cite{Kac1997} where they do not; see Remark \ref{nonhomo-equivalence-remark}).  Thus the formulation of vertex operators with two odd formal variables given on p.115 of \cite{HK} and on p. 8 of \cite{He} correspond to our equation (\ref{constructing-nonhomo-operators}).  However, in \cite{HK} and \cite{He}, Kac and Heluani conclude (for instance, in equation (2.6.4) in \cite{HK}), that 
\begin{equation}\label{HK-mu}
Y(\sqrt{-1} J_{(-1)} |0\rangle,z,\theta^1, \theta^2) =  \sqrt{-1} J(z) + \theta^1 G^{(2)}(z) - \theta^2 G^{(1)}(z) + 2 \theta^1 \theta^2 L(z),
\end{equation}
where $J(z) = \sum_{n \in \mathbb{Z}} J_n z^{-n-2}$, and $L(z) = \sum_{n \in \mathbb{Z}} L_n z^{-n-2}$.  Now, if one assumes $G^{(j)} (z) = \sum_{n \in \frac{1}{2} + \mathbb{Z}} G_n^{(j)} z^{-n- \frac{3}{2}}$, then  equation (\ref{HK-mu}) is equivalent to 
\begin{multline}
Y(\mu,(x,\varphi^{(1)},\varphi^{(2)})) = \sum_{n \in \mathbb{Z}} \Bigl( J(n) x^{- n - 1}-  i\varphi^{(1)} G^{(2)} (n+1/2) x^{-n-2} \\
+ i\varphi^{(2)} G^{(1)}(n+1/2) x^{-n-2}  - 2 i \varphi^{(1)} \varphi^{(2)}  L(n) x^{- n - 2} \Bigr)
\end{multline}
instead of (\ref{nonhomo-mu}) as we have.  On the other hand if one assumes that Kac and Heluani are using the rather bizarre notation that $G^{(j)} (z) = -\sum_{n \in \frac{1}{2} + \mathbb{Z}} G_n^{(j)} z^{-n- \frac{3}{2}}$ as the equations (2.6.2) of \cite{HK} and the equations on the bottom of p.7 of \cite{He} would seem to indicate, and in addition one assumes they have changed their notation midstream to have $L(z) = -\sum_{n \mathbb{Z}} L_n z^{-n-2}$ as their equation $L(z) = (1/2) Y(S^2 S^1 i \mu, z)$ on p. 114 of \cite{HK} and p.7 of \cite{He} would seem to indicate, then the vertex operator corresponding to $\mu$ would indeed be correct! Note that in this last statement we are using property (\ref{nonhomo-annihilation}) below.}
\end{rema}

\begin{rema}\label{huang-milas-remark2} {\em  Our notion of $N=2$ Neveu-Schwarz vertex operator superalgebra with two odd formal variables in the nonhomogeneous coordinate system is similar in spirit to the notion of ``$N=2$ superconformal vertex operator superalgebra with odd formal variables" given in \cite{HM}.  However their are a couple of typos in the definition and consequences given in \cite{HM}.  In \cite{HM} the vertex operators incorporate odd variables according to 
\begin{multline}\label{huang-milas}
Y(v, (x, \varphi_1, \varphi_2)) = Y(v,x) + \varphi_1 Y(G_1(-1/2)v,x) \\
+ \varphi_2 Y(G_2(-1/2)v,x) + \varphi_1 \varphi_2 Y(G_1(-1/2) G_2(-1/2)v,x);
\end{multline}
see p. 366 of \cite{HM}.  In terms of our notation, $G_1(-1/2) = G^{(2)}(-1/2)$ and $G_2(-1/2) = G^{(1)}(-1/2)$, $\varphi_1 = \varphi^{(2)}$ and $\varphi_2 = \varphi^{(1)}$.  Thus by comparing (\ref{huang-milas})  with (\ref{constructing-nonhomo-operators}) in light of (\ref{nonhomo-D-correspondence}), one can see that there is a typo in the sign of the last term of (\ref{huang-milas}).  Also, their is an additional typo in the vertex operator with odd formal variables corresponding to $\mu$.}
\end{rema}

We have the following consequences of the definition of $N=2$ Neveu-Schwarz vertex operator superalgebra with two odd formal variables in the nonhomogeneous coordinate system: 
\begin{equation}\label{nonhomo-annihilation}
G^{(j)}(n + 1/2) \mathbf{1}  = 0, \quad \mbox{for $n \geq -1$ and $j = 1,2$},
\end{equation}
(cf. (\ref{kill-vacuum})).  There exist two vectors $\tau^{1} = G^{(1)}(-3/2) \mathbf{1} = i G^{(2)}(-1/2) \mu$ and $\tau^{2} = G^{(2)}(-3/2)$ $\mathbf{1} = -i G^{(1)}(-1/2) \mu$ in $V_{(3/2)}^1$ such that
\begin{multline}\label{tau1}
Y(\tau^{1},(x,\varphi^{(1)},\varphi^{(2)})) =  \sum_{n \in \mathbb{Z}} \Bigl( G^{(1)}(n+1/2) x^{- n - 2} + 2\varphi^{(1)} L(n) x^{-n-2}  \\ 
- i(n+1)\varphi^{(2)} J(n) x^{-n-2}    - (n+2)  \varphi^{(1)} \varphi^{(2)} G^{(2)}(n+ 1/2) x^{- n - 3} \Bigr)
\end{multline}
\begin{multline}\label{tau2}
Y(\tau^{2},(x,\varphi^{(1)},\varphi^{(2)})) =   \sum_{n \in \mathbb{Z}} \Bigl( G^{(2)}(n+ 1/2) x^{- n - 2} + 2\varphi^{(2)}L(n) x^{-n-2}  \\ 
+ i(n+1) \varphi^{(1)}  J(n) x^{-n-2}   +  (n+2)  \varphi^{(1)} \varphi^{(2)} G^{(1)}(n+ 1/2) x^{- n - 3} \Bigr),
\end{multline}
as well as $\omega = L(-2) \mathbf{1} = -\frac{i}{2} G^{(2)} (-1/2) G^{(1)} (-1/2) \mu $  in $V_{(2)}$ with 
\begin{multline}
Y(\omega,(x,\varphi^{(1)},\varphi^{(2)})) =  \sum_{n \in \mathbb{Z}} \Bigl( L(n) x^{- n - 2} - \frac{1}{2}(n+2)  \varphi^{(1)} G^{(1)}(n+ 1/2) x^{-n-3}  \\ 
- \frac{1}{2} (n+2) \varphi^{(2)} G^{(2)}(n+ 1/2 ) x^{-n-3}   + \frac{i}{2}(n+1)(n+2) \varphi^{(1)} \varphi^{(2)} J(n) x^{- n - 3} \Bigr).
\end{multline}
Note that 
\begin{equation}
\tau^{1} = \frac{1}{\sqrt{2}} ( \tau^{(+)} + \tau^{(-)}) \quad \mbox{and} \quad \tau^{2} = \frac{i}{\sqrt{2}} ( \tau^{(+)} - \tau^{(-)}) ,
\end{equation}
i.e. 
\begin{equation}
\tau^{(\pm)} = \frac{1}{\sqrt{2}} ( \tau^{1} \mp i \tau^{2}).
\end{equation}
We have 
\begin{eqnarray}
\lefteqn{\quad Y(\mu, (x,\varphi^{(1)}, \varphi^{(2)}))\mu } \label{nonhomo-mu-singular}\\
&=& \! \! \! \! J(1) \mu x^{-2} + i\varphi^{(1)} G^{(2)}(- 1/2) \mu  x^{-1} - i \varphi^{(2)} G^{(1)} (- 1/2) \mu x^{-1}  \nonumber\\
& & \quad + \, 2i\varphi^{(1)} \varphi^{(2)} \mu x^{-2} + 2i \varphi^{(1)} \varphi^{(2)} L(-1) \mu x^{-1} + y (x, \varphi^{(1)}, \varphi^{(2)}) \nonumber\\
\quad &=& \! \! \! \! \frac{1}{3}c_V \mathbf{1} x^{-2}  + \varphi^{(1)} G^{(1)}(- 3/2) \mathbf{1} x^{-1} + \varphi^{(2)} G^{(2)}(-3/2) \mathbf{1} x^{-1}  \nonumber \\
& & \quad +\,  2i \varphi^{(1)} \varphi^{(2)} J(-1) \mathbf{1} x^{-2}  + 2i\varphi^{(1)} \varphi^{(2)} J(-2) \mathbf{1} x^{-1} + y(x, \varphi^{(1)}, \varphi^{(2)}) \nonumber \\
&=& \! \! \! \! \frac{1}{3}c_V \mathbf{1} x^{-2} + \varphi^{(1)} \tau^{1} x^{-1} + \varphi^{(2)} \tau^{2} x^{-1}  + 2i\varphi^{(1)} \varphi^{(2)} \mu x^{-2} \nonumber\\
& & \quad + \, 2i\varphi^{(1)} \varphi^{(2)} L(-1) \mu x^{-1} + y (x, \varphi^{(1)}, \varphi^{(2)}) \nonumber
\end{eqnarray}
where $y(x, \varphi^{(1)}, \varphi^{(2)}) \in V[[x]][\varphi^{(1)}, \varphi^{(2)}]$. 

Since $G^{(j)}(-1/2) = \mathcal{D}^{(j)}$ for $j = 1,2$, the identities for $\mathcal{D}^{(j)}$ apply using $G^{(j)}(-1/2)$ for $j=1,2$, respectively (and as before $\mathcal{D} = L(-1)$).  That is, we now have the $L(-1)$-derivative property (\ref{nonhomo-D-derivative}), skew supersymmetry (\ref{nonhomo-skew-supersymmetry}),  the  $G^{(j)}(-1/2)$- and $L(-1)$-bracket-derivative properties (\ref{D*-bracket-derivative}) and (\ref{nonhomo-D-bracket-derivative}), and the $G^{(j)}(-1/2)$- and $L(-1)$-bracket properties (\ref{D*-bracket}) and (\ref{nonhomo-D-bracket}) with $\mathcal{D}^{(j)}= G^{(j)}(-1/2)$, for $j = 1,2$, and $\mathcal{D} = L(-1)$.

In addition, 
\begin{multline} \label{L0-bracket-nonhomo} 
\left[ L(0), Y(v, (x, \varphi^{(1)}, \varphi^{(2)})) \right] = Y\biggl( \Bigl( L(0)+ \frac{1}{2} \varphi^{(1)} G^{(1)}(-1/2)  \\
+ \frac{1}{2} \varphi^{(2)} G^{(2)} (-1/2) + xL(-1)\Bigr) v, (x, \varphi^{(1)}, \varphi^{(2)})\biggr) 
\end{multline}
\begin{multline} \label{J0-bracket-nonhomo} 
\left[ J(0), Y(v, (x, \varphi^{(1)}, \varphi^{(2)})) \right] = Y\left(\left(J(0)- i \varphi^{(1)} G^{(2)}(-1/2) \right. \right. \\
\left. \left. + \, i  \varphi^{(2)} G^{(1)} (-1/2)  +  2 i \varphi^{(1)} \varphi^{(2)} L(-1)\right)v, (x, \varphi^{(1)}, \varphi^{(2)} )\right) 
\end{multline}
\begin{multline}\label{G-1/2-bracket-nonhomo} 
\left[ G^{(1)}(1/2), Y(v, (x, \varphi^{(1)}, \varphi^{(2)})) \right] = Y\left(\left(G^{(1)}(1/2)- 2\varphi^{(1)} L(0) + i \varphi^{(2)} J(0) \right. \right.  \\
\left. \left. -  \varphi^{(1)} \varphi^{(2)} G^{(2)} (-1/2)   + x G^{(1)} (-1/2)  -  2 x \varphi^{(1)} L(-1) \right) v, (x, \varphi^{(1)}, \varphi^{(2)})\right)
\end{multline}
\begin{multline}\label{G-1/2-bracket-nonhomo2} 
\left[ G^{(2)}(1/2), Y(v, (x, \varphi^{(1)}, \varphi^{(2)})) \right]  = Y\left(\left(G^{(2)}(1/2)- 2\varphi^{(2)} L(0) - i \varphi^{(1)} J(0)\right. \right.   \\
\left. \left. +  \varphi^{(1)} \varphi^{(2)} G^{(1)} (-1/2) + x G^{(2)} (-1/2) -  2 x \varphi^{(2)} L(-1) \right) v, (x, \varphi^{(1)}, \varphi^{(2)})\right) 
.
\end{multline}

\begin{rema}\label{Lie-algebras-remark-nonhomo}
{\em Recall that by (\ref{form-Lie1}), the operators $u_0$ form a Lie superalgebra denoted $\mathcal{L}_0$; see Remark \ref{Lie-algebras-remark}.  By (\ref{form-Lie4}), the operators $u^{(j)}_{-1/2}$, for $j = 1,2$ also form two Lie superalgebras; denote these by $\mathcal{L}^{(j)}_{-1/2}$, for $j=1,2$, respectively.  Finally, by (\ref{form-Lie5}), the operators $u^{(1,2)}_{-1}$ form a Lie superalgebra which by (\ref{substitute-v3}) is isomorphic to the Lie superalgebra generated by $u^{+-}_{-1}$ and denoted by $\mathcal{L}^{+-}_{-1}$; see Remark \ref{Lie-algebras-remark}.   Since $\tau^{j}_0 = G^{(j)}(-1/2) = \mathcal{D}^{(j)}$, by (\ref{nonhomo-Dv-bracket-condition}), we have that $\mathcal{L}^{(j)}_{-1/2}$, for $j=1,2$, are subalgebras of  $\mathcal{L}_0$.  }
\end{rema}

We continue to use the notation $\mathrm{wt} \ v = n$ for $v\in V$ satisfying $L(0)v = nv$, and $\mathrm{wt}^J \ v = n$ for $v \in V$ satisfying $J(0)v = nv$.  {}From (\ref{L0-bracket-nonhomo}), and using the $L(-1)$- and $G^{(j)}(-1/2)$-derivative properties, for $j=1,2$, we obtain the following:
\begin{eqnarray}
\mbox{wt} \; v_n  =  \mbox{wt} \; v_n^{(1,2)} &=& \mbox{wt} \; v - n - 1 ,\label{L(0)-grading-nonhomo} \\
\mbox{wt} \; v_{n- \frac{1}{2}}^{(j)}  &=&  \mbox{wt} \; v - (n - \frac{1}{2}) - 1 , \quad \mbox{for $j=1,2$} \label{L(0)-grading-nonhomo2}
\end{eqnarray}
for $n \in \mathbb{Z}$ and for $v \in V$ of homogeneous weight with respect to $L(0)$.
{}From (\ref{J0-bracket-nonhomo}), and using the $L(-1)$-, $G^{(j)}(-1/2)$-derivative properties, we obtain the following:
\begin{equation}
\mathrm{wt}^J \ v_n = \mathrm{wt}^J \ v^{(1,2)}_n = \mathrm{wt}^J \ v  \label{J(0)-grading-nonhomo1} 
\end{equation}
for $n \in \mathbb{Z}$ and for $v \in V$ of homogeneous weight with respect to the $J(0)$ grading.  However $v^{(j)}_{n-1/2}$ for $j=1,2$, are in general nonhomogeneous operators with respect to the $J(0)$ grading.  In fact, for $v,w \in V$ of homogeneous weight with respect to $J(0)$, we have
\begin{eqnarray}
J(0) v_{n - \frac{1}{2}}^{(1)} w &=& (\mathrm{wt}^J v + \mathrm{wt}^J  w) v_{n - \frac{1}{2}}^{(1)} w  - i v_{n-\frac{1}{2}}^{(2)} w  \label{J(0)-grading-nonhomo2}\\
J(0) v_{n - \frac{1}{2}}^{(2)} w &=& (\mathrm{wt}^J v + \mathrm{wt}^J  w) v_{n - \frac{1}{2}}^{(2)} w + i v_{n-\frac{1}{2}}^{(1)} w .\label{J(0)-grading-nonhomo3}
\end{eqnarray}

\begin{rema}\label{last-nonhomo-terminology-remark}
{\em  The action of $J(0)$ on the endomorphisms $v^{(j)}_{n -1/2}$ given in (\ref{J(0)-grading-nonhomo2}) and (\ref{J(0)-grading-nonhomo3}) versus the action of $J(0)$ on $v^\pm_{j - 1/2}$ given in 
 (\ref{J(0)-grading2}) provide another motivation for our terminology in calling the $(x, \varphi^+, \varphi^-)$ coordinate system homogeneous and the $(x, \varphi^{(1)}, \varphi^{(2)})$ coordinate system nonhomogeneous.  That is, for $v \in V$ homogeneous with respect to the $J(0)$-grading, the endomorphisms $v^\pm_{n-1/2}$ are also homogeneous with respect to this grading whereas the endomorphisms $v^{(j)}_{n -1/2}$, for $j = 1,2$ are not.}
\end{rema}

As a consequence of the above $L(0)$-grading properties (\ref{L(0)-grading-nonhomo}) and (\ref{L(0)-grading-nonhomo2}), we have the following property in analogy to the homogeneous case (\ref{conjugate-by-L(0)}) 
\begin{equation} \label{conjugate-by-L(0)-nonhomo} 
x_0^{2L(0)} Y(v, (x, \varphi^{(1)}, \varphi^{(2)})) x_0^{-2 L(0)} = Y(x_0^{2L(0)}v, (x_0^2 x,
x_0 \varphi^{(1)}, x_0 \varphi^{(2)})) .
\end{equation}
However {}from the $J(0)$-grading properties (\ref{J(0)-grading-nonhomo1})--(\ref{J(0)-grading-nonhomo3}), we have the following property which is very different {}from the homogeneous case given in (\ref{conjugate-by-J(0)}),
\begin{multline}\label{conjugate-by-J(0)-nonhomo}
x_0^{J(0)} Y(v, (x, \varphi^{(1)}, \varphi^{(2)})) x_0^{-J(0)} = Y\biggl(x_0^{J(0)}v, \biggl(x, \frac{\varphi^{(1)}}{2} (x_0 + x_0^{-1}) \\
+ \frac{i\varphi^{(2)}}{2}(x_0 - x_0^{-1}), -\frac{i\varphi^{(1)}}{2}(x_0 - x_0^{-1}) + \frac{\varphi^{(2)}}{2}(x_0 + x_0^{-1}) \biggr)\biggr).
\end{multline}
Equation (\ref{conjugate-by-J(0)-nonhomo}) can be expressed as
\begin{multline}\label{conjugate-by-J(0)-hyperbolic}
e^{x_0J(0)} Y(v, (x, \varphi^{(1)}, \varphi^{(2)})) e^{-x_0J(0)} \\
\qquad \qquad = Y(e^{x_0J(0)}v, (x, \varphi^{(1)} \cosh x_0 + i \varphi^{(2)} \sinh x_0 , -i\varphi^{(1)} \sinh x_0 + \varphi^{(2)} \cosh x_0)) .
\end{multline}

For $a,b \in (\bigwedge_*^0)^\times$ recall the map $\gamma_{(a,b)}: V \rightarrow V$ defined by (\ref{define-gamma}), that is $\gamma_{(a,b)} (v) = a^{2L(0)} b^{J(0)} v$ for $v \in V$.  Writing $b = e^\beta$ for $\beta \in \bigwedge_*^0$, and using (\ref{conjugate-by-L(0)-nonhomo}) and (\ref{conjugate-by-J(0)-nonhomo}), we have that 
\begin{multline}\label{gamma-iso-nonhomo}
\gamma_{(a,e^\beta)} \circ Y(v, (x, \varphi^{(1)}, \varphi^{(2)})) \circ \gamma_{(a,e^\beta)}^{-1} \\
\qquad = Y(\gamma_{(a,e^\beta)}(v), (a^2x, a \varphi^{(1)} \cosh \beta + i a \varphi^{(2)} \sinh \beta, - i a \varphi^{(1)} \sinh \beta + a \varphi^{(2)} \cosh \beta)),
\end{multline}
for $a \in (\bigwedge_*^0)^\times$ and $\beta \in \bigwedge_*^0$.
Thus we have the following lemma:

\begin{lem}\label{nonhomo-iso-lemma}
Let $(V, Y(\cdot, (x, \varphi^{(1)}, \varphi^{(2)})), \mathbf{1}, \mu)$ be an $N=2$ Neveu-Schwarz vertex operator superalgebra in the nonhomogeneous basis.  Then $(V, Y(\cdot, (a^2x, a \varphi^{(1)} \cosh \beta + i a \varphi^{(2)} \sinh \beta, - i a \varphi^{(1)} \sinh \beta + a \varphi^{(2)} \cosh \beta)), \mathbf{1}, a^2\mu)$ is isomorphic to  $(V, Y(\cdot, (x, \varphi^{(1)},$ $\varphi^{(2)})), \mathbf{1}, \mu)$, for $a \in (\bigwedge_*^0)^\times$ and $\beta \in \bigwedge_*^0$.
\end{lem}

\begin{rema}\label{nonhomo-auto-remark2} {\em Given an $N=2$ Neveu-Schwarz vertex operator superalgebra in nonhomogeneous coordinates $(V, Y(\cdot,$ $(x, \varphi^{(1)}, \varphi^{(2)})), \mathbf{1}, \mu)$, the automorphisms of the nonhomogeneous $N=2$ Neveu-Schwarz algebra (\ref{nonhomo-auto1})--(\ref{nonhomo-auto3}), give rise to the following $N=2$ Neveu-Schwarz vertex operator superalgebras which are isomorphic to $(V, Y(\cdot, (x, \varphi^{(1)}, \varphi^{(2)})), \mathbf{1}, \mu)$:  We have the following family of $N=2$ Neveu-Schwarz vertex operator superalgebras in nonhomogeneous coordinates which have the same $N=2$ Neveu-Schwarz algebra element $\mu$
\begin{multline}
(V, Y(\cdot, (x, \varphi^{(1)}, \varphi^{(2)})), \mathbf{1}, \mu) \cong (V,Y(\cdot, (x, \varphi^{(1)} \cosh \beta + i \varphi^{(2)} \sinh \beta, \\
- i \varphi^{(1)} \sinh \beta +  \varphi^{(2)} \cosh \beta)), \mathbf{1}, \mu)
\end{multline}
for $\beta \in \bigwedge_*^0$, which correspond to the isomorphisms $\gamma_{(1,e^\beta)}$ defined in (\ref{define-gamma}).  In addition, we have the isomorphism 
\begin{equation}
(V, Y(\cdot, (x, \varphi^{(1)},- \varphi^{(2)})), \mathbf{1}, \mu) \cong (V, Y(\cdot, (x, \varphi^{(1)}, \varphi^{(2)})), \mathbf{1}, -\mu) 
\end{equation}
which does not preserve the $N=2$ Neveu-Schwarz algebra element but does preserve the ``Virasoro element" $\omega = L(-2) \mathbf{1}$.  And, of course, the composition of these isomorphisms gives the following continuous family of isomorphic $N=2$ Neveu-Schwarz vertex operator superalgebras in nonhomogeneous coordinates which preserve $\omega$ but not $\mu$
\begin{multline}
(V, Y(\cdot, (x, \varphi^{(1)}, \varphi^{(2)})), \mathbf{1}, \mu) \cong (V, Y(\cdot, (x,  - i \varphi^{(1)} \sinh \beta +  \varphi^{(2)} \cosh \beta, \\
\varphi^{(1)} \cosh \beta + i  \varphi^{(2)} \sinh \beta)), \mathbf{1}, -\mu)
\end{multline}
for $\beta \in \bigwedge_*^0$. Cf. Remark \ref{auto-remark2}.  }
\end{rema}

One can see some aspects of the correspondence with the geometry of $N=2$ super-Riemann spheres with tubes in the nonhomogeneous coordinates by recalling {}from \cite{B-moduli} that such objects are topological superspaces whose transition functions are $N=2$ superconformal, meaning they transform the superdifferential operators $D^{(j)}$, for $j=1,2$,  defined by (\ref{define-D-derivative-nonhomo}), homogeneously of degree one in each of the $D^{(j)}$'s for $j=1,2$; see Remark \ref{superconformal-remark-nonhomo}.  As noted in Remark \ref{geometry-correspondence-remark-nonhomo} the $\mathcal{D}^{(j)} = G^{(j)} (-1/2)$-derivative properties in the notion of $N=2$ (Neveu-Schwarz) vertex (operator) superalgebra with two odd formal variables in the nonhomogeneous coordinate system show some of the correspondence with the notion of $N=2$ superconformality in the nonhomogeneous coordinate system.  In addition, here we note that in \cite{B-moduli} we prove the following proposition:

\begin{prop}\label{moduli-prop-nonhomo}(\cite{B-moduli})
There is a bijection between formal $N=2$ superconformal functions in the nonhomogeneous coordinate system vanishing at zero and invertible in a neighborhood of zero and expressions of the form
\begin{multline}\label{change-of-variables-nonhomo}
\exp \Biggl(-  \sum_{n \in \Z} \Bigl( A^{(1)}_n L_n(x,\varphi^{(1)}, \varphi^{(2)}) + A^{(2)}_n J_n(x,\varphi^{(1)}, \varphi^{(2)}) \\
+  M^{(1)}_{n - \frac{1}{2}} G^{(1)}_{n -\frac{1}{2}} (x, \varphi^{(1)}, \varphi^{(2)}) 
+  M^{(2)}_{n - \frac{1}{2}} G^{(2)}_{n -\frac{1}{2}} (x, \varphi^{(1)}, \varphi^{(2)}) \Bigr) \Biggr)  \cdot  \\
(a_0^{(1)})^{-2L_0(x,\varphi^{(1)}, \varphi^{(2)})} \cdot  (a_0^{(2)})^{-J_0(x,\varphi^{(1)}, \varphi^{(2)})} \cdot (x, \varphi^{(1)}, \varphi^{(2)}) 
\end{multline}
for $(a_0^{(1)}, a_0^{(2)}) \in ((\bigwedge_*^0)^\times)^2/\langle (-1,-1) \rangle$, and $A_n^{(j)} \in \bigwedge_*^0$ and $M^{(j)}_{n- 1/2} \in \bigwedge_*^1$, for $n \in \Z$, and $j = 1,2$
where 
\begin{eqnarray}
L_n(x,\varphi^{(1)},\varphi^{(2)}) \! \! \! \! &=&\! \! \! \!  - \biggl( x^{n+ 1} \frac{\partial}{\partial x} + (\frac{n + 1}{2})x^n \Bigl(  \varphi^{(1)} \frac{\partial }{\partial \varphi^{(1)}}  + \varphi^{(2)} \frac{\partial }{\partial \varphi^{(2)} }  \Bigr) \biggr) \label{L-notation-nonhomo} \\
\quad J_n(x,\varphi^{(1)},\varphi^{(2)}) \! \! \! \! &=& \! \! \! \! i x^n \biggl( \varphi^{(1)}  \frac{\partial }{\partial \varphi^{(2)}} -   \varphi^{(2)}  \frac{\partial }{\partial \varphi^{(1)}} \biggr)  \label{J-notation-nonhomo}  \\
\qquad \ \ \ \ G^{(1)}_{n-\frac{1}{2}} (x,\varphi^{(1)},\varphi^{(2)}) \! \! \! \! &=& \! \! \! \! - \biggl( x^n \Bigl(   \frac{\partial }{\partial \varphi^{(1)}}  -  \varphi^{(1)}  \frac{\partial}{\partial x} \Bigr)  - nx^{n-1} \varphi^{(1)} \varphi^{(2)}   \frac{\partial }{\partial \varphi^{(2)}}  \biggr)   \\
\qquad \ \ \ \ G^{(2)}_{n-\frac{1}{2}} (x,\varphi^{(1)},\varphi^{(2)}) \! \! \! \! &=& \! \! \! \! - \biggl( x^n \Bigl(   \frac{\partial }{\partial \varphi^{(2)}}  -  \varphi^{(2)}  \frac{\partial}{\partial x} \Bigr)  +  nx^{n-1} \varphi^{(1)} \varphi^{(2)}   \frac{\partial }{\partial \varphi^{(1)}}  \biggr)   \label{G-notation-nonhomo}
\end{eqnarray}
are superderivations in  $\mbox{Der} (\bigwedge_*[[x,x^{-1}]] [\varphi^{(1)},\varphi^{(2)}])$, for $n \in \mathbb{Z}$, which give a representation of the $N=2$ Neveu-Schwarz algebra with central charge zero in the nonhomogeneous basis; that is  (\ref{L-notation-nonhomo})--(\ref{G-notation-nonhomo}) satisfy (\ref{Virasoro-relation2})--(\ref{transformed-Neveu-Schwarz-relation-last}) with $d = 0$. Similarly, there is a bijection between formal $N=2$ superconformal functions in the nonhomogeneous coordinate system vanishing at $(\infty, 0,0)$ and invertible in a neighborhood of $(\infty, 0,0)$ and expressions of the form
\begin{multline}\label{change-of-variables-nonhomo-infty}
\exp \Biggl( \sum_{n \in \Z} \Bigl( A^{(1)}_n L_{-n}(x,\varphi^{(1)}, \varphi^{(2)}) - A^{(2)}_n J_{-n} (x,\varphi^{(1)}, \varphi^{(2)}) \\
+  iM^{(1)}_{n - \frac{1}{2}} G^{(1)}_{-n +\frac{1}{2}} (x, \varphi^{(1)}, \varphi^{(2)}) 
+  iM^{(2)}_{n - \frac{1}{2}} G^{(2)}_{-n +\frac{1}{2}} (x, \varphi^{(1)}, \varphi^{(2)}) \Bigr) \Biggr)  \cdot  \\
(a_0^{(1)})^{2L_0(x,\varphi^{(1)}, \varphi^{(2)})} \cdot  (a_0^{(2)})^{-J_0(x,\varphi^{(1)}, \varphi^{(2)})} \cdot \left(\frac{1}{x}, \frac{i\varphi^{(1)}}{x},\frac{i \varphi^{(1)}}{x} \right) .
\end{multline}
\end{prop}

Let $\mathcal{SC}(2,0)$ denote the set of formal $N=2$ superconformal functions vanishing at zero and invertible in a neighborhood of zero.   It is clear that $\mathcal{SC}(2,0)$ is a group under composition.  Let $A^{(j)} = \{A^{(j)}_n\}_{n \in \mathbb{Z}_+}$ and $M^{(j)}= \{M^{(j)}_{n-1/2}\}_{n \in \mathbb{Z}_+}$, for $j = 1,2$.  Let $(\bigwedge_*^k)^\infty$ denote the set of all infinite sequences in $\bigwedge_*^k$ for $k = 0,1$, and let $\bigwedge_*^\infty = (\bigwedge_*^0)^\infty \oplus (\bigwedge_*^1)^\infty$.  Let
\begin{equation}
\mathcal{G} = \mbox{$(\bigwedge_*^0)^\times)^2/\langle (-1,-1) \rangle \times (\bigwedge_*^\infty)^2$}.
\end{equation} 
Define the map 
\begin{equation}
\hat{E}_2: \mathcal{G} \longrightarrow \mathcal{SC}(2,0)
\end{equation}
by
\begin{multline}
\hat{E}_2(a_0^{(1)}, a_0^{(2)}, A^{(1)}, A^{(2)}, M^{(1)}, M^{(2)}) \\
= \exp \Biggl(-  \sum_{n \in \Z} \Bigl( A^{(1)}_n L_n(x,\varphi^{(1)}, \varphi^{(2)}) + A^{(2)}_n J_n(x,\varphi^{(1)}, \varphi^{(2)}) \\
+  M^{(1)}_{n - \frac{1}{2}} G^{(1)}_{n -\frac{1}{2}} (x, \varphi^{(1)}, \varphi^{(2)}) 
+  M^{(2)}_{n - \frac{1}{2}} G^{(2)}_{n -\frac{1}{2}} (x, \varphi^{(1)}, \varphi^{(2)}) \Bigr) \Biggr)  \cdot  \\
(a_0^{(1)})^{-2L_0(x,\varphi^{(1)}, \varphi^{(2)})} \cdot  (a_0^{(2)})^{-J_0(x,\varphi^{(1)}, \varphi^{(2)})} \cdot (x, \varphi^{(1)}, \varphi^{(2)}) .
\end{multline}
It is clear that $\hat{E}_2$ is a bijection.  Let $g, h \in \mathcal{G}$ with $g=(a_0^{(1)}, a_0^{(2)}, A^{(1)}, A^{(2)}, M^{(1)}, M^{(2)})$ and $h=(b_0^{(1)}, b_0^{(2)}$, $B^{(1)}, B^{(2)}, N^{(1)}, N^{(2)})$.  We define the {\it $N=2$ composition at zero} map
\begin{eqnarray}
\circ_0: \mathcal{G} \times \mathcal{G} &\longrightarrow& \mathcal{G} \\
(g,h) &\mapsto& g \circ_0 h \nonumber
\end{eqnarray}
by 
\begin{eqnarray}\label{define-circ}
\lefteqn{\quad g\circ_0 h }\\
\quad \quad &=& \! \! (a_0^{(1)}, a_0^{(2)}, A^{(1)}, A^{(2)}, M^{(1)}, M^{(2)}) \circ_0 (b_0^{(1)}, b_0^{(2)}, B^{(1)}, B^{(2)}, N^{(1)}, N^{(2)})  \nonumber \\
&=&\! \! \hat{E}_2^{-1} \left(  \hat{E}_2 (b_0^{(1)}, b_0^{(2)}, B^{(1)}, B^{(2)}, N^{(1)}, N^{(2)}) \circ \hat{E}_2(a_0^{(1)}, a_0^{(2)}, A^{(1)}, \right. \nonumber \\
& &  \hspace{3in} \left.  A^{(2)}, M^{(1)}, M^{(2)}) \right).\nonumber
\end{eqnarray}

For $(A^{(1)}, A^{(2)}, M^{(1)}, M^{(2)}) \in (\bigwedge_*^\infty)^2$, define
\begin{multline}\label{T-notation}
T((A^{(1)}, A^{(2)}, M^{(1)}, M^{(2)}), (x, \varphi^{(1)}, \varphi^{(2)}))\\
=-  \sum_{n \in \Z} \Bigl( A^{(1)}_n L_n(x,\varphi^{(1)}, \varphi^{(2)}) + A^{(2)}_n J_n(x,\varphi^{(1)}, \varphi^{(2)}) \\
+  M^{(1)}_{n - \frac{1}{2}} G^{(1)}_{n -\frac{1}{2}} (x, \varphi^{(1)}, \varphi^{(2)}) 
+  M^{(2)}_{n - \frac{1}{2}} G^{(2)}_{n -\frac{1}{2}} (x, \varphi^{(1)}, \varphi^{(2)}) \Bigr).
\end{multline}

The next proposition follows from Propositions 6.12 and 6.15  in \cite{B-moduli}, translated from homogeneous to nonhomogeneous coordinates.  

\begin{prop}(\cite{B-moduli})\label{N=2-group-zero-prop}
The set $\mathcal{G} = (\bigwedge_*^0)^\times)^2/\langle (-1,-1) \rangle \times (\bigwedge_*^\infty)^2$ is a group under the composition operation $\circ_0$ defined by (\ref{define-circ}).   In fact, $(\mathcal{G}, \circ_0)$ is isomorphic to the opposite group of  the group of formal $N=2$ superconformal functions vanishing at zero and invertible in a neighborhood of zero under composition, i.e., 
\begin{equation}
(\mathcal{G}, \circ_0) \cong (\mathcal{SC}(2,0), \circ)^{\mathrm{op}}.
\end{equation}
Furthermore for $g, h \in \mathcal{G}$ with $g = (a_0^{(1)}, a_0^{(2)}, A^{(1)}, A^{(2)}, M^{(1)}, M^{(2)})$ and 
$h = (b_0^{(1)}, b_0^{(2)}, B^{(1)}, B^{(2)}, N^{(1)}, N^{(2)})$, we have that if
\begin{equation}
g \circ_0 h = (c_0^{(1)}, c_0^{(2)}, C^{(1)}, C^{(2)},O^{(1)}, O^{(2)}) 
\end{equation}
then $(c_0^{(1)}, c_0^{(2)},C^{(1)}, C^{(2)},O^{(1)}, O^{(2)}) \in  \mathcal{G}$ is given by
\begin{eqnarray}
\lefteqn{e^{T((C^{(1)}, C^{(2)}, O^{(1)}, O^{(2)}), (x, \varphi^{(1)}, \varphi^{(2)}))}
 \cdot (c_0^{(1)})^{-2L_0(x,\varphi^{(1)}, \varphi^{(2)})} \cdot  (c_0^{(2)})^{-J_0(x,\varphi^{(1)}, \varphi^{(2)})} }\\
&=& e^{T((A^{(1)}, A^{(2)}, M^{(1)}, M^{(2)}), (x, \varphi^{(1)}, \varphi^{(2)}))} \cdot (a_0^{(1)})^{-2L_0(x,\varphi^{(1)}, \varphi^{(2)})}   \nonumber \\
& & \quad \cdot  (a_0^{(2)})^{-J_0(x,\varphi^{(1)}, \varphi^{(2)})} \cdot e^{T((B^{(1)}, B^{(2)}, N^{(1)}, N^{(2)}), (x, \varphi^{(1)}, \varphi^{(2)}))} \hspace{.2in} \nonumber\\
& & \quad \cdot (b_0^{(1)})^{-2L_0(x,\varphi^{(1)}, \varphi^{(2)})} \cdot  (b_0^{(2)})^{-J_0(x,\varphi^{(1)}, \varphi^{(2)})}. \nonumber
\end{eqnarray}
\end{prop}

Similarly, letting $\mathcal{SC}(2,\infty)$ denote the set of formal $N=2$ superconformal functions vanishing at infinity and invertible in a neighborhood of infinity.   Let $I_2(x,\varphi^{(1)}, \varphi^{(2)}) = (1/x, i\varphi^{(1)}/x, i \varphi^{(2)}/x)$.   Then $\mathcal{SC}(2,\infty)$ is a group under the operation 
\begin{equation}\label{infinity-composition}
f_1 \cdot f_2 = f_1 \circ I_2^{-1} \circ f_2 .
\end{equation}  
Again let $g, h \in \mathcal{G}$ with $g=(a_0^{(1)}, a_0^{(2)}, A^{(1)}$, $A^{(2)}, M^{(1)}, M^{(2)})$ and $h=(b_0^{(1)}, b_0^{(2)}$, $B^{(1)}, B^{(2)}, N^{(1)}, N^{(2)})$.  We define the {\it $N=2$ composition at infinity} map
\begin{eqnarray}
\circ_\infty: \mathcal{G} \times \mathcal{G} &\longrightarrow& \mathcal{G} \\
(g,h) &\mapsto& g \circ_\infty h \nonumber
\end{eqnarray}
as follows:  let
\begin{multline}
(c_0^{(1)}, c_0^{(2)},C^{(1)}, -C^{(2)},-iO^{(1)}, -iO^{(2)}) \\
= \hat{E}_2^{-1} \left(  \hat{E}_2 (b_0^{(1)}, b_0^{(2)}, B^{(1)}, B^{(2)}, N^{(1)}, N^{(2)}) \circ I_2^{-1} \circ \hat{E}_2(a_0^{(1)}, a_0^{(2)}, A^{(1)}, \right. \\
 \left.  A^{(2)}, M^{(1)}, M^{(2)}) \circ I_2^{-1}\right)
\end{multline} 
then define
\begin{eqnarray}\label{define-circ-infty}
\lefteqn{\quad g\circ_\infty h} \\
\qquad &=&(a_0^{(1)}, a_0^{(2)}, A^{(1)}, A^{(2)}, M^{(1)}, M^{(2)}) \circ_\infty (b_0^{(1)}, b_0^{(2)}, B^{(1)}, B^{(2)}, N^{(1)}, N^{(2)})  \nonumber \\
&=& (c_0^{(1)}, c_0^{(2)}, C^{(1)}, C^{(2)},O^{(1)}, O^{(2)}) . \nonumber
\end{eqnarray}

The next proposition follows from Proposition 6.20 and Corollary 6.21 in \cite{B-moduli}, translated from homogeneous to nonhomogeneous coordinates.  

\begin{prop}(\cite{B-moduli})\label{N=2-group-infinity-prop}
The set $\mathcal{G} = (\bigwedge_*^0)^\times)^2/\langle (-1,-1) \rangle \times (\bigwedge_*^\infty)^2$ is a group under the composition operation $\circ_\infty$ defined by (\ref{define-circ-infty}).  In fact, $(\mathcal{G}, \circ_\infty)$ is isomorphic to the opposite group of  the group of formal $N=2$ superconformal functions vanishing at infinity and invertible in a neighborhood of infinity under the composition (\ref{infinity-composition}), i.e.
\begin{equation}
(\mathcal{G}, \circ_\infty) \cong (\mathcal{SC}(2,\infty), \cdot)^{\mathrm{op}}.
\end{equation}
Furthermore for $g, h \in \mathcal{G}$ with $g = (a_0^{(1)}, a_0^{(2)}, A^{(1)}, A^{(2)}, M^{(1)}, M^{(2)})$ and $h = (b_0^{(1)}, b_0^{(2)}, B^{(1)}, B^{(2)}, N^{(1)}, N^{(2)})$, we have that if
\begin{equation}
g \circ_\infty h = (c_0^{(1)}, c_0^{(2)}, C^{(1)}, C^{(2)},O^{(1)}, O^{(2)}) 
\end{equation}
then $(c_0^{(1)}, c_0^{(2)},C^{(1)}, C^{(2)},O^{(1)}, O^{(2)}) \in  \mathcal{G}$ is given by
\begin{eqnarray}
\lefteqn{\quad e^{T((C^{(1)}, -C^{(2)}, -iO^{(1)},-i O^{(2)}), I_2(x, \varphi^{(1)}, \varphi^{(2)}))}
 \cdot (c_0^{(1)})^{2L_0(x,\varphi^{(1)}, \varphi^{(2)})} }\\
 & & \quad \cdot  (c_0^{(2)})^{-J_0(x,\varphi^{(1)}, \varphi^{(2)})} \nonumber\\
&=& e^{T((A^{(1)}, -A^{(2)}, -iM^{(1)}, -iM^{(2)}), I_2(x, \varphi^{(1)}, \varphi^{(2)}))}  \cdot (a_0^{(1)})^{2L_0(x,\varphi^{(1)}, \varphi^{(2)})} \nonumber \\
& & \quad \cdot  (a_0^{(2)})^{-J_0(x,\varphi^{(1)}, \varphi^{(2)})} \cdot e^{T((B^{(1)}, -B^{(2)}, -iN^{(1)}, -iN^{(2)}), I_2(x, \varphi^{(1)}, \varphi^{(2)}))} \nonumber\\
& & \quad \cdot (b_0^{(1)})^{2L_0(x,\varphi^{(1)}, \varphi^{(2)})} \cdot  (b_0^{(2)})^{-J_0(x,\varphi^{(1)}, \varphi^{(2)})}. \nonumber
\end{eqnarray}
\end{prop}

In light of Proposition  \ref{moduli-prop-nonhomo}, properties such as (\ref{gamma-iso-nonhomo}) can be rewritten as 
\begin{multline}\label{LJ-change-of-variables-nonhomo}
a^{2L(0)} b^{J(0)}  Y(v, (x, \varphi^{(1)}, \varphi^{(2)})) a^{-2L(0)} b^{-J(0)}  \\
= Y(a^{2L(0)} b^{J(0)} v, a^{-2L_0(x, \varphi^{(1)}, \varphi^{(2)})} b^{-J_0(x, \varphi^{(1)}, \varphi^{(2)})} (x,  \varphi^{(1)},  \varphi^{(2)})),
\end{multline}
where using (\ref{L-notation-nonhomo}) and (\ref{J-notation-nonhomo}), we have the following $N=2$ superconformal change of variables vanishing at zero and invertible in a neighborhood of zero in the nonhomogeneous coordinate system
\begin{eqnarray}
\lefteqn{\qquad a^{-2L_0(x, \varphi^{(1)}, \varphi^{(2)})} b^{-J_0(x, \varphi^{(1)}, \varphi^{(2)})} (x,  \varphi^{(1)},  \varphi^{(2)})}\\
&=&  a^{\left(2x\frac{\partial}{\partial x} + \varphi^{(1)} \frac{\partial}{\partial \varphi^{(1)}} + \frac{1}{2} \varphi^{(2)} \frac{\partial}{\partial \varphi^{(2)}}\right)} b^{-i\left( \varphi^{(1)} \frac{\partial}{\partial \varphi^{(2)}} - \varphi^{(2)} \frac{\partial}{\partial \varphi^{(1)}}\right)} (x,  \varphi^{(1)}, \varphi^{(2)}) \nonumber\\
&=& \biggl(a^2 x, \frac{a}{2}  \varphi^{(1)} (b + b^{-1}) + \frac{ia}{2} \varphi^{(2)} (b - b^{-1}) , \frac{-ia}{2}  \varphi^{(1)} (b - b^{-1}) \nonumber \\
& & \qquad + \, \frac{a}{2} \varphi^{(2)} (b + b^{-1}) \biggr). \nonumber
\end{eqnarray}
Another example of change of variables interpretation is that (\ref{conjugate-shift-nonhomo}) can be expressed as
\begin{eqnarray}\label{shift-change-of-variables-nonhomo}
\lefteqn{ \qquad \qquad e^{x_0 L(-1) + \varphi^{(1)}_0  G^{(1)}(-1/2) +  \varphi^{(2)}_0  G^{(2)}(-1/2)} Y(v,(x,\varphi^{(1)}, \varphi^{(2)})) }\\
& & \qquad \cdot e^{- x_0 L(-1)- \varphi_0^{(1)}  G^{(1)}(-1/2) -  \varphi^{(2)}_0  G^{(2)}(-1/2)}\nonumber  \\
&=& \! \! \! \! \! Y(v,e^{-x_0 L_{-1}(x, \varphi^{(1)}, \varphi^{(2)})  - \varphi^{(1)}_0  G^{(1)}_{-1/2} (x, \varphi^{(1)}, \varphi^{(2)})  -  \varphi^{(2)}_0  G^{(2)}_{-1/2}(x, \varphi^{(1)}, \varphi^{(2)})}  (x , \varphi^{(1)}, \varphi^{(2)})). \nonumber
\end{eqnarray} 

\begin{rema}\label{change-of-variables-remark-nonhomo}{\em 
We see {}from the discussion above that one can think of (\ref{LJ-change-of-variables-nonhomo}) and (\ref{shift-change-of-variables-nonhomo}) as a change of variables formula related to a certain change of variables of the form (\ref{change-of-variables-nonhomo}) or (\ref{change-of-variables-nonhomo-infty}), respectively.  See also Remark \ref{change-of-variables-remark}.
}
\end{rema}

\begin{rema}{\em
Proposition  \ref{moduli-prop-nonhomo} and Remark \ref{change-of-variables-remark-nonhomo} reflect the connection between $N=2$ Neveu-Schwarz vertex operator superalgebras with two odd formal variables and $N=2$ superconformal functions.  In a Section \ref{N=1-section}, we show that there is an equivalent notion of ``$N=2$ vertex superalgebra with one odd formal variable" that is related to a different but equivalent superfunctional setting, namely that of $N=1$ superanalytic (as opposed to superconformal) functions.   
 }
\end{rema}

\begin{rema}
{\em Given an $N=2$ Neveu-Schwarz vertex operator superalgebra with two odd formal variables in nonhomogeneous coordinates $(V, Y(\cdot, (x, \varphi^{(1)}, \varphi^{(2)})), \mathbf{1}$, $\mu)$, and setting $\varphi^{(1)} = \varphi^{(2)} = 0$, we have the notion of ``$N=2$ Neveu-Schwarz vertex operator superalgebra without odd formal variables in the nonhomogeneous basis" given by $(V, Y(\cdot, (x, 0, 0)), \mathbf{1}, \tau^{1} = G^{(1)}(-3/2)\mathbf{1}, \tau^{2} = G^{(2)}(-3/2) \mathbf{1})$.  Note that we require that the two elements $\tau^{j}$  for $j = 1,2$, in $V_{(3/2)}$ be specified in contrast to the corresponding notion with formal variables in which we specify just one vector giving rise to the $N=2$ Neveu-Schwarz algebra in the nonhomogeneous basis, namely $\mu \in V_{(1)}$.  This is because, in analogy to the homogeneous case (see Remark \ref{need-both-if-without}), if we do not include odd formal variables we need these two vectors in order to give a full generating set for the $N=2$ Neveu-Schwarz algebra in the nonhomogeneous basis.  However, unlike in the homogeneous case, in the nonhomogeneous case, we could also specify $\mu$ and one of $\tau^{1}$ or $\tau^{2}$.  Of course via Remark \ref{transform-isomorphism-remark}  and Theorem \ref{superalgebras}, we see that the category of $N=2$ (Neveu-Schwarz) vertex (operator) superalgebras without odd formal variables in the homogeneous and in the nonhomogeneous coordinate basis are isomorphic.
} 
\end{rema}

We have the following corollary to Proposition \ref{mu-prop}: 

\begin{cor}\label{nonhomo-mu-cor}
Let $(V, Y(\cdot, (x, \varphi^{(1)}, \varphi^{(2)})), \mathbf{1})$ be an $N=2$ vertex superalgebra with two odd formal variables in the nonhomogeneous coordinates system, and let $\mathcal{D}^{(j)}$, for $j=1,2$, be defined by (\ref{nonhomo-define-Ds}).  Let $\mathcal{D} = \frac{1}{2} [ \mathcal{D}^{(1)}, \mathcal{D}^{(2)}]$.  If $\mu \in V$, satisfies
\begin{multline}\label{nonhomo-for-mu-prop}
Y(\mu, (x,\varphi^{(1)}, \varphi^{(2)}))\mu =  \frac{1}{3}c \mathbf{1} x^{-2}  + i\varphi^{(1)} \mathcal{D}^{(2)}  \mu x^{-1} -i \varphi^{(2)} \mathcal{D}^{(1)}  \mu x^{-1} \\
 + 2i\varphi^{(1)} \varphi^{(2)} \mu x^{-2} + 2i\varphi^{(1)} \varphi^{(2)} \mathcal{D} \mu x^{-1} + y(x, \varphi^{(1)}, \varphi^{(2)}) 
\end{multline}
for some $c \in \mathbb{C}$ and $y \in V[[x]][\varphi^{(1)}, \varphi^{(2)}]$, that is if
\begin{eqnarray}
\mu_1 \mu &=& \frac{1}{3} c \mathbf{1} \\
\mu_{-\frac{1}{2}}^{(1)} \mu &=&i \mathcal{D}^{(2)}  \mu \ = \  i \mu_{-\frac{3}{2}}^{(2)} \mathbf{1} \\
\mu_{-\frac{1}{2}}^{(2)} \mu &=&-i \mathcal{D}^{(1)}  \mu \ = \ -i \mu_{-\frac{3}{2}}^{(1)} \mathbf{1} \\
\mu_0^{(1,2)} \mu &=& 2i \mu \\
\mu_{-1}^{(1,2)} \mu  &=& 2i \mathcal{D} \mu  \ = \ 2i\mu_{-2}\mathbf{1}
\end{eqnarray}
and
\begin{equation}
\mu_{n-1} \mu = \mu_{n-\frac{1}{2}}^{(j)} \mu = \mu_n^{(1,2)} \mu =  0 \qquad \mbox{for $n\geq 1$ and $j = 1,2$},
\end{equation}
then defining $L(n), J(n), G^{(j)}(n - 1/2) \in \mathrm{End} \, V$, for $n \in \mathbb{Z}$ and $j=1,2$, by $\mu_n = J(n)$, $\mu_{n-1/2}^{(1)} = iG^{(2)}(n - 1/2)$, $\mu_{n-1/2}^{(2)} = - iG^{(1)}(n - 1/2)$ and $\mu_n^{(1,2)} = 2iL(n)$, we have that $L(n), J(n)$, and $G^{(j)}(n-1/2)$ satisfy the relations for the $N=2$ Neveu-Schwarz algebra in the nonhomogeneous basis (\ref{nonhomo-n2-first})--(\ref{nonhomo-n2-last}) with central charge $c$. 
\end{cor}

\section{$N=1$ superanalytic fields and $N=2$ (Neveu-Schwarz) vertex (operator) superalgebras with one odd formal variable}\label{N=1-section}

In \cite{B-thesis} (see also \cite{B-vosas}), we introduced the notion of ``$N=1$ Neveu-Schwarz vertex operator superalgebra with odd formal variables".  If $(V, Y(\cdot, (x, \varphi^{(1)}, \varphi^{(2)})), \mathbf{1}, \mu)$ is an $N=2$ Neveu-Schwarz vertex operator superalgebra with two odd formal variables in the nonhomogeneous coordinate system, then $(V, Y(\cdot, (x, \varphi^{(1)}, 0)), \mathbf{1}, \tau^{1})$ and $(V, Y(\cdot, (x, 0, \varphi^{(2)})), \mathbf{1}, \tau^{2})$ are both $N=1$ Neveu-Schwarz vertex operator superalgebras with one odd formal variable where $\tau^j = G^{(j)} (-3/2) \mathbf{1}$ for $j = 1,2$, i.e. $\tau^1 = i G^{(2)}(-1/2) \mu$ and $\tau^2 = - iG^{(1)}(-1/2) \mu$; see (\ref{tau1}) and (\ref{tau2}).  In this section, we will see that there is in fact a deeper connection between $N=1$ Neveu-Schwarz vertex operator superalgebras and $N=2$ Neveu-Schwarz vertex operator superalgebras.   In particular, this will lead us to the notion of ``$N=2$ (Neveu-Schwarz) vertex (operator) superalgebras with one odd formal variable" motivated by $N=1$ superanalytic geometry or equivalently continuously deformed $N=1$ superconformal geometry. 

\subsection{The $N=1$ Neveu-Schwarz algebra and formal $N=1$ superconformal functions}

The $N=1$ Neveu-Schwarz algebra is the Lie superalgebra which has a basis consisting of the central element $d$, even elements $L_n$ and odd elements $G_{n + 1/2}$ for $n \in \mathbb{Z}$, and commutation relations  
\begin{eqnarray}
\left[L_m ,L_n \right] &=& (m - n)L_{m + n} + \frac{1}{12} (m^3 - m) \delta_{m + n 
, 0} \; d , \label{Virasoro-relation-N1} \\
\left[ L_m, G_{n + \frac{1}{2}} \right] &=& \left(\frac{m}{2} - n - \frac{1}{2} \right) G_{m+n+\frac{1}{2}} ,\\
\left[ G_{m + \frac{1}{2}} , G_{n - \frac{1}{2}} \right] &=&2L_{m + n}  + \frac{1}{3} (m^2 + m) \delta_{m + n , 0} \; d ,  \label{N1-Neveu-Schwarz-relation-last}
\end{eqnarray}
for $m, n \in \mathbb{Z}$.  Note that there are two copies of the $N=1$ Neveu-Schwarz algebra in the $N=2$ Neveu-Schwarz algebra, namely $\mathrm{span}_{\bigwedge_*}\{L_n, G^{(j)}_{n-1/2}, d \; | \; n \in \mathbb{Z} \}$ for $j =1$ and for $j = 2$.  

\begin{rema}\label{N1-subalgebras-remark}
{\em The subalgebra  $\mathrm{span}_{\bigwedge_*} \{ L_{-1}, G_{- 1/2}, L_0, G_{1/2}, L_1 \}$ of the $N=1$ Neveu-Schwarz algebra is the orthogonal-symplectic algebra 
$\mathfrak{osp}_{\bigwedge_*} (1|2)$; see \cite{Kac}, \cite{B-memoir}.  This is the Lie superalgebra of ``infinitesimal $N=1$ superconformal transformations of the $N=1$ super-Riemann sphere" as shown in \cite{CR}, see also \cite{B-memoir}.  We will denote the subalgebra of $\mathfrak{osp}_{\bigwedge_*} (1|2)$ given by $\mathrm{span}_{\bigwedge_*} \{ L_{-1}, G_{- 1/2} \}$ by $\mathfrak{osp}_{\bigwedge_*} (1|2)_{<0}$.
}
\end{rema}

Formally, we say that a power series $f(x,\varphi) = (\tilde{x}, \tilde{\varphi})$ with $\tilde{x} \in (\bigwedge_*[[x,x^{-1}]] [ \varphi])^0$ and $\tilde{\varphi} \in (\bigwedge_*[[x, x^{-1}]][\varphi])^1$ is {\it $N=1$ superanalytic}.  In addition, $f$ is {\it $N=1$ superconformal} in $x$ and $\varphi$ if and only if it satisfies 
\begin{equation}\label{nice-superconformal-condition-N1}
D_1\tilde{x} - \tilde{\varphi} D_1 \tilde{\varphi} =  0 
\end{equation}     
where
\begin{equation}\label{define-D-derivative-N1}
D_1 = \frac{\partial}{\partial \varphi} + \varphi \frac{\partial}{\partial x} 
\end{equation} 
(see \cite{Fd}, \cite{B-memoir} and cf. Remark \ref{superconformal-remark-nonhomo}).  For a superanalytic power series $f$, the condition (\ref{nice-superconformal-condition-N1}) is equivalent to requiring that $f$ transform the $N=1$ super-differential operator $D_1$ homogeneously of degree one \cite{B-memoir}.  (In the literature for the $N=1$ case $D_1$ is usually denoted $D$; however, we have already used $D$ to denote $\frac{\partial}{\partial x}$, and thus we will instead use $D_1$.)   

In \cite{B-memoir} (see also \cite{B-thesis}), the following proposition was proved:
\begin{prop}\label{moduli-prop-N1}(\cite{B-memoir})
There is a bijection between formal $N=1$ superconformal functions vanishing at zero and invertible in a neighborhood of zero and expressions of the form
\begin{equation}\label{change-of-variables-N1}
\exp \biggl(-  \sum_{n \in \Z} \Bigl( A_n L_n(x,\varphi) +  M_{n - \frac{1}{2}} G_{n -\frac{1}{2}} (x, \varphi) 
 \Bigr) \biggr)  \cdot  a_0^{-2L_0(x,\varphi)} \cdot   (x, \varphi) 
\end{equation}
for $a_0 \in (\bigwedge_*^0)^\times$, and $A_n \in \bigwedge_*^0$ and $M_{n- 1/2} \in \bigwedge_*^1$, for $n \in \Z$, where 
\begin{eqnarray}
L_n(x,\varphi) &=&  - \biggl( x^{n+ 1} \frac{\partial}{\partial x} + \Bigl(\frac{n + 1}{2} \Bigr)x^n \varphi  \frac{\partial }{\partial \varphi }   \biggr) \label{L-notation-N1} \\
G_{n-\frac{1}{2}} (x,\varphi) &=&  - x^n \Bigl(   \frac{\partial }{\partial \varphi}  -  \varphi  \frac{\partial}{\partial x} \Bigr)    \label{G-notation-N1}
\end{eqnarray}
are superderivations in  $\mbox{Der} (\bigwedge_*[[x,x^{-1}]] [\varphi])$, for $n \in \mathbb{Z}$, which give a representation of the $N=1$ Neveu-Schwarz algebra with central charge zero; that is  (\ref{L-notation-N1}) and (\ref{G-notation-N1}) satisfy (\ref{Virasoro-relation-N1})--(\ref{N1-Neveu-Schwarz-relation-last}) with $d = 0$.  Similarly, there is a bijection between formal $N=1$ superconformal functions vanishing at $(\infty, 0)$ and invertible in a neighborhood of $(\infty, 0)$ and expressions of the form
\begin{equation}\label{change-of-variables-N1-infty}
\exp \biggl(\sum_{n \in \Z} \Bigl( A_n L_{-n}(x,\varphi) +  iM_{n - \frac{1}{2}} G_{-n +\frac{1}{2}} (x, \varphi) 
 \Bigr) \biggr)  \cdot  a_0^{2L_0(x,\varphi)} \cdot  \left(\frac{1}{x}, \frac{i\varphi}{x} \right) .
\end{equation}
\end{prop}

\subsection{Alternate notions of $N=1$ superconformality and continuous deformations of $N=1$ (Neveu-Schwarz) vertex (operator) superalgebras with one odd formal variable}

The key property of the operator $D_1$ in $N=1$ superconformal field theory is simply that it satisfies $(D_1)^2 = \frac{\partial}{\partial x}$.  Let $s \in (\bigwedge_*^0)^\times$ and $\sigma \in \bigwedge_*^1$.   Define 
\begin{equation}\label{deformed-D}
D_{(s,\sigma)} = \frac{1}{s} \frac{\partial}{\partial \varphi} + (s \varphi + \sigma)\frac{\partial}{\partial x}. 
\end{equation} 
Then $D^{(s,\sigma)}$ is an odd superderivation in $\mathrm{Der} \bigwedge_* [[x, x^{-1}]] [\varphi]$ and 
\begin{equation}\label{deformed-condition}
(D_{(s,\sigma)})^2 = \frac{\partial}{\partial x}.
\end{equation}
In fact, it is not hard to verify that a superderivation $T \in \mathrm{Der} \in \bigwedge_* [[x, x^{-1}]] [\varphi]$ satisfies $T^2 = \frac{\partial}{\partial x}$ if and only if $T = D_{(s,\sigma)}$ for some $s \in (\bigwedge_*^0)^\times$ and $\sigma \in \bigwedge_*^1$.  This continuous family of derivations satisfying (\ref{deformed-condition}) was first given in \cite{B-announce} (see also \cite{B-thesis}). If we let $s = 1$ and $\sigma = 0$, then $D_{(1,0)} = D_1$ is the standard superderivation used to define $N=1$ superconformality and $D_{(s,\sigma)}$ is a deformation of $D_{(1,0)}$ by the $N=1$ superanalytic (but not $N=1$ superconformal) transformation 
\begin{equation}
\varphi \mapsto s\varphi + \sigma .
\end{equation}

Since all $D_{(s, \sigma)}$ satisfy (\ref{deformed-condition}) it makes sense to consider $N=1$ superconformal field theory -- the geometric and algebraic aspects of the theory -- where instead of $D_{(1,0)}$, we use $D_{(s, \sigma)}$.  We will call such a theory {\it $N=1$ $D_{(s,\sigma)}$-superconformal}.  For the purposes of this paper, we restrict ourselves to $\sigma = 0$, and {}from now on will write $D_{(s,0)} = D_s$.

In \cite{B-announce}--\cite{B-iso} the author gave a rigorous study of the moduli space of genus-zero $N=1$ super-Riemann surfaces with tubes, modulo $N=1$ superconformal equivalence.  These ``$N=1$ superspheres with tubes" model certain worldsheets, and a tube (representing a propagating $N=1$ superstring) is $N=1$ superconformally equivalent to a puncture on the $N=1$ supersphere with a local $N=1$ superconformal coordinate chart vanishing at the puncture \cite{B-memoir}.   In \cite{B-memoir} (see also \cite{B-thesis}), the author proved Proposition \ref{moduli-prop-N1} which shows 
that these local $N=1$ superconformal coordinates can be expressed in terms of exponentials of certain $N=1$ superderivations, and that these $N=1$ superderivations give a representation of the $N = 1$  Neveu-Schwarz algebra with zero central charge \cite{NS}.  In \cite{B-memoir} a sewing operation on this moduli space was defined using local $N=1$ superconformal coordinates and an interpretation of sewing in terms of the exponentials of representatives of $N=1$ Neveu-Schwarz algebra elements was formulated.  Using this supergeometric structure, in \cite{B-iso} the author introduced the notion of {\it $N=1$ supergeometric vertex operator superalgebra with central charge $c \in \mathbb{C}$} and proved the following theorem.

\begin{thm}\label{iso-thm} (\cite{B-iso})
Assuming the convergence of certain projective factors, the category of $N=1$ supergeometric vertex operator superalgebras with central charge $c \in \mathbb{C}$ is isomorphic to the category  of (superalgebraic) $N=1$ Neveu-Schwarz vertex operator superalgebras with central charge $c \in\mathbb{C}$ (\cite{B-vosas}).
\end{thm}

The proof of this theorem involves algebraic, differential geometric, and analytic techniques.  As part of this rigorization of the correspondence between $N=1$ Neveu-Schwarz vertex operator superalgebras and the worldsheet supergeometry of genus-zero superconformal field theory, the author was able to formulate the odd variable components for $N=1$ Neveu-Schwarz vertex operator superalgebras \cite{B-vosas} so that the formal algebraic notions reflected the differential supergeometry. 

With the motivation above, we can repeat the work of the author in \cite{B-announce}--\cite{B-iso} using the operator $D_s$ for $s \in (\bigwedge_*^0)^\times$, or in other words deforming the odd variable by $\varphi \mapsto s\varphi$.   In particular, we will say that a formal power series $f(x,\varphi) = (\tilde{x}, \tilde{\varphi}) \in \bigwedge_*[[x, x^{-1}]][\varphi]$ is $N=1$ $D_s$-superconformal in $x$ and $\varphi$ if and only if $f$ satisfies 
\begin{equation}\label{nice-superconformal-condition-Ds}
D_s\tilde{x} - s\tilde{\varphi} D_s s\tilde{\varphi} =  0 
\end{equation}     
for 
\begin{equation}\label{define-D-derivative-Ds}
D_s = \frac{1}{s} \frac{\partial}{\partial \varphi} + s\varphi \frac{\partial}{\partial x}. 
\end{equation} 
For a formal superanalytic function $f$ in one even variable and one odd variables the condition (\ref{nice-superconformal-condition-Ds}) is equivalent to requiring that $f$ transform the $N=1$ super-differential operator $D_s$ homogeneously of degree one.

We have the following corollary to Proposition \ref{moduli-prop-N1}:

\begin{cor}\label{moduli-prop-Ds}
There is a bijection between formal $N=1$ $D_s$-superconformal functions vanishing at zero and invertible in a neighborhood of zero and expressions of the form
\begin{equation}\label{change-of-variables-Ds}
\exp \Biggl(-  \sum_{n \in \Z} \Bigl( A_n L_{s,n} (x,\varphi) +  M_{n - \frac{1}{2}} G_{s,n -\frac{1}{2}} (x, \varphi) 
 \Bigr) \Biggr)  \cdot  a_0^{-2L_{s,0}(x,\varphi)} \cdot   (x, \varphi) 
\end{equation}
for $a_0 \in (\bigwedge_*^0)^\times$, and $A_n \in \bigwedge_*^0$ and $M_{n- 1/2} \in \bigwedge_*^1$, for $n \in \Z$, where 
\begin{eqnarray}
L_{s,n}(x,\varphi) &=&  L_n(x, s\varphi) = - \biggl( x^{n+ 1} \frac{\partial}{\partial x} + \Bigl(\frac{n + 1}{2} \Bigr)x^n \varphi  \frac{\partial }{\partial \varphi }   \biggr)  \label{L-notation-Ds} \\
&=& L_n(x, \varphi) \nonumber \\
G_{s,n-\frac{1}{2}} (x,\varphi) &=& G_{n-\frac{1}{2}} (x,s\varphi)=  - x^n \Bigl( \frac{1}{s}  \frac{\partial }{\partial \varphi}  - s \varphi  \frac{\partial}{\partial x} \Bigr)    \label{G-notation-Ds}
\end{eqnarray}
are superderivations in  $\mbox{Der} (\bigwedge_*[[x,x^{-1}]] [\varphi])$, for $n \in \mathbb{Z}$, which give a representation of the $N=1$ Neveu-Schwarz algebra with central charge zero; that is  (\ref{L-notation-Ds}) and (\ref{G-notation-Ds}) satisfy (\ref{Virasoro-relation-N1})--(\ref{N1-Neveu-Schwarz-relation-last}) with $d = 0$.  Similarly, there is a bijection between formal $N=1$ $D_s$-superconformal functions vanishing at $(\infty, 0)$ and invertible in a neighborhood of $(\infty, 0)$ and expressions of the form
\begin{equation}\label{change-of-variables-Ds-infty}
\exp \Biggl(\sum_{n \in \Z} \Bigl( A_n L_{s,-n} (x,\varphi) +  iM_{n - \frac{1}{2}} G_{s,-n +\frac{1}{2}} (x, \varphi) 
 \Bigr) \Biggr)  \cdot  a_0^{2L_{s,0}(x,\varphi)}  \cdot  \left(\frac{1}{x}, \frac{i\varphi}{x} \right) .
\end{equation}
\end{cor}

Following the work of the author in \cite{B-announce}--\cite{B-iso} in this setting, this then gives rise to the algebraic notions of ``$N=1$ $D_s$-vertex superalgebra" and ``$N=1$ $D_s$-Neveu-Schwarz vertex operator superalgebra" as we present below.

\begin{defn}\label{N1-vertex-superalgebra-definition}
{\em For $s \in (\bigwedge_*^0)^\times$, an} $N=1$ $D_s$-vertex superalgebra over $\bigwedge_*$ and with odd formal variables {\em consists of a $\mathbb{Z}_2$-graded $\bigwedge_*$-module (graded by} sign $\eta${\em) 
\begin{equation}\label{first-n1-superalgebra-with}
V = V^0 \oplus V^1
\end{equation} 
equipped, first, with a linear map 
\begin{eqnarray}\label{N1-operator-with}
V &\longrightarrow&  (\mbox{End} \; V)[[x,x^{-1}]][\varphi] \\
v  &\mapsto&  Y_s(v,(x,\varphi)) = \sum_{n \in \mathbb{Z}} \Bigl( v_n  + \varphi v_{n - \frac{1}{2}}  \Bigr)  x^{-n-1},\nonumber
\end{eqnarray}
where for the $\mathbb{Z}_2$-grading of $\mathrm{End} \; V$ induced {}from that of $V$, we have
\begin{equation}
v_n \in (\mbox{End} \; V)^{\eta(v)}, \quad \mbox{and} \quad v_{n - \frac{1}{2}} \in  (\mbox{End} \; V)^{(\eta(v) + 1) \mbox{\begin{footnotesize} mod 
\end{footnotesize}} 2},
\end{equation}
for $v$ of homogeneous sign in $V$, $x$ is an even formal variable, $\varphi$ is an odd formal variables, and where $Y_s(v,(x,\varphi))$ denotes the} $D_s$-vertex operator associated with $v$.  {\em We also have a distinguished element $\mathbf{1}$ in $V$ (the} vacuum vector{\em ).  The following conditions are assumed for $u,v \in V$:  the} truncation conditions: {\em
\begin{equation}\label{N1-truncation}
u_n v =   u_{n-\frac{1}{2}} v = 0 \qquad \mbox{for $n \in \mathbb{Z}$ sufficiently large,} 
\end{equation}
that is
\begin{equation}
Y_s(u, (x, \varphi))v \in V((x))[\varphi];
\end{equation}
next, the following} vacuum property{\em :
\begin{equation}\label{N1-vacuum-identity}
Y_s(\mathbf{1}, (x, \varphi)) = \mathrm{id}_V ;
\end{equation}
the} creation property {\em holds:
\begin{equation}\label{N1-creation-property} 
Y_s(v,(x,\varphi)) \mathbf{1} \in V[[x]][\varphi] \quad \mathrm{and} \quad 
\lim_{(x,\varphi) \rightarrow 0} Y_s(v,(x, \varphi)) \mathbf{1} = v ;  
\end{equation} 
and finally the} Jacobi identity {\em holds:  
\begin{eqnarray}\label{N1-Jacobi-identity}
\lefteqn{\ \ x_0^{-1} \delta \biggl( \frac{x_1 - x_2 - s^2\varphi_1 \varphi_2 }{x_0} \biggr) Y_s(u,(x_1, \varphi_1)) Y_s(v,(x_2, \varphi_2)) }\\
& & \ \  -  (-1)^{\eta(u)\eta(v)} x_0^{-1} \delta \biggl( \frac{x_2 -  x_1 + s^2\varphi_1 \varphi_2}{-x_0} \biggr) Y_s(v,(x_2, \varphi_2))Y_s(u,(x_1,\varphi_1)) \nonumber\\
\qquad &=& \! \! \!  x_2^{-1} \delta \biggl( \frac{x_1 - x_0 - s^2\varphi_1 \varphi_2}{x_2} \biggr) Y_s(Y_s(u,(x_0, \varphi_1 - \varphi_2))v,(x_2, \varphi_2)) , \nonumber
\end{eqnarray}
for $u,v$ of homogeneous sign in $V$.
}
\end{defn}

The $N=1$ $D_s$-vertex superalgebra just defined is denoted by 
\[(V,Y_s(\cdot,(x,\varphi)),\mathbf{1}).\]

\begin{rema}
{\em  The cases $s = \pm1$ coincide with the notion of an $N=1$ vertex superalgebra in the usual sense (cf. \cite{B-vosas} and \cite{HK}).   In the next section we will see that the case when $s$ is allowed to vary over all of $(\bigwedge_*^0)^\times$ simultaneously is also significant for its relationship to the $N=2$ Neveu-Schwarz algebra.  
}
\end{rema}

As with $N=2$ vertex superalgebras, there are many useful consequences of the definition.  We will give just a few of the analogous consequences below which will need later. 

Taking $\mathrm{Res}_{x_1}$ of the Jacobi identity and using the $\delta$-function identity (\ref{delta-2-terms-with-phis}) which still holds with the substitutions $\varphi_1^+= s\varphi_1 $, $\varphi_2^- = s\varphi_2$, $\varphi_1^- = 0$, and $\varphi_2^+ = 0$, we obtain the following {\it iterate formula} for $N=1$ $D_s$-vertex operators with odd formal variables
\begin{multline}\label{iterate-N1}
Y_s(Y_s(u,(x_0, \varphi_1 - \varphi_2))v,(x_2, \varphi_2)\\
= \mathrm{Res}_{x_1} \biggl( x_0^{-1} \delta \biggl( \frac{x_1 - x_2 - s^2\varphi_1 \varphi_2}{x_0} \biggr)Y_s(u,(x_1, \varphi_1)) Y_s(v,(x_2, \varphi_2)) \\
-  (-1)^{\eta(u)\eta(v)} x_0^{-1} \delta \biggl( \frac{x_2 - x_1 + s^2 \varphi_1 \varphi_2 }{-x_0} \biggr) Y_s(v,(x_2, \varphi_2))Y_s(u,(x_1,\varphi_1)) \biggr).
\end{multline}

\begin{prop}\label{N1-derivative-prop}
Let $(V, Y_s(\cdot, (x, \varphi)), \mathbf{1})$ be an $N=1$ $D_s$-vertex superalgebra and let $\mathcal{D}_s$ be the odd endomorphism of $V$ defined by
\begin{equation}\label{define-Ds}
\mathcal{D}_s (v) = \frac{1}{s} v_{-\frac{3}{2}} \mathbf{1} \qquad \mbox{for $v \in V$.}
\end{equation}
Then 
\begin{equation}\label{Ds-derivative-property}
Y_s(\mathcal{D}_s v,(x,\varphi)) = \biggl(\frac{1}{s} \frac{\partial}{\partial \varphi} + s \varphi \frac{\partial}{\partial x} \biggr) Y_s(v,(x,\varphi))  . 
\end{equation} 
Furthermore, letting 
\begin{equation}\label{define-D-N1}
\mathcal{D} = \frac{1}{2} \left[ \mathcal{D}_s, \mathcal{D}_s \right] = \mathcal{D}_s^2,
\end{equation}
we have the $\mathcal{D}$-derivative property 
\begin{equation}\label{D-derivative-N1}
Y_s(\mathcal{D}v, (x, \varphi)) = \frac{\partial}{\partial x} Y_s(v, (x, \varphi)) ,
\end{equation}
and
\begin{equation}\label{D-is-what-it-should-be2}
\mathcal{D}(v) = v_{-2}\mathbf{1}.
\end{equation}
In particular, $\mathrm{span}_{\bigwedge_*} \{ \mathcal{D}, \mathcal{D}_s \}$ satisfies the commutation relations for the Lie superalgebra $\mathfrak{osp}_{\bigwedge_*}(1|2)_{<0}$, and thus $V$ is a representation for this Lie superalgebra.
\end{prop}

\begin{proof}
Using the iterate formula (\ref{iterate-N1}), the vacuum property (\ref{N1-vacuum-identity}), Proposition \ref{expansion-prop} in the case $\varphi_1^+= s\varphi_1 $, $\varphi_2^- = s\varphi_2$, $\varphi_1^- = 0$, and $\varphi_2^+ = 0$, (\ref{delta-2-terms-with-phis}) with $x_0=0$, $\varphi_1^+= s\varphi_1 $, $\varphi_2^- = s\varphi_2$, $\varphi_1^- = 0$, and $\varphi_2^+ = 0$ (which is well-defined), and  (\ref{delta-substitute}), again with $x_0=0$, $\varphi_1^+= s\varphi_1 $, $\varphi_2^- = s\varphi_2$, $\varphi_1^- = 0$, and $\varphi_2^+ = 0$, we have
\begin{eqnarray*}
\lefteqn{Y_s(\mathcal{D}_s v,(x_2,\varphi_2)) }\\
&=&  \frac{1}{s}Y(v_{-\frac{3}{2}} \mathbf{1},(x_2,\varphi_2)) \\
&=&  \frac{1}{s} \frac{\partial}{\partial \varphi_1}  \mathrm{Res}_{x_0} x_0^{-1} Y_s(Y_s(v,(x_0, \varphi_1 - \varphi_2)) \mathbf{1},  (x_2, \varphi_2)) \biggr|_{\varphi_1= 0}\\
&=&  \frac{1}{s} \frac{\partial}{\partial \varphi_1}  \mathrm{Res}_{x_0} x_0^{-1}  \mathrm{Res}_{x_1} \biggl( x_0^{-1} \delta \biggl( \frac{x_1 - x_2 - s^2\varphi_1 \varphi_2}{x_0} \biggr)Y_s(v,(x_1, \varphi_1))  \\
& &\quad  - \, x_0^{-1} \delta \biggl( \frac{x_2 - x_1 + s^2 \varphi_1 \varphi_2 }{-x_0} \biggr) Y_s(v,(x_1,\varphi_1)) \biggr) \biggr|_{\varphi_1 = 0}\\
&=&  \frac{1}{s} \frac{\partial}{\partial \varphi_1} \mathrm{Res}_{x_1} \biggl( \Bigl( (x_1 - x_2 - s^2\varphi_1 \varphi_2)^{-1}\\
& & \quad  -  \, (-x_2 + x_1 - s^2\varphi_1 \varphi_2)^{-1} \Bigr)  Y_s(v,(x_1,\varphi_1)) \biggr) \biggr|_{\varphi_1 = 0}\\
&=&  \frac{1}{s} \frac{\partial}{\partial \varphi_1}  \mathrm{Res}_{x_1} \biggl( x_2^{-1} \delta \left(\frac{x_1 - s^2 \varphi_1 \varphi_2}{x_2} \right)  Y_s(v,(x_1,\varphi_1)) \biggr) \biggr|_{\varphi_1 = 0}\\
&=& \frac{1}{s}  \frac{\partial}{\partial \varphi_1}  \mathrm{Res}_{x_1} \biggl( x_1^{-1} \delta \left(\frac{x_2 +  s^2 \varphi_1 \varphi_2}{x_1} \right)  Y_s(v,(x_1,\varphi_1)) \biggr) \biggr|_{\varphi_1 = 0}\\
&=&  \frac{1}{s} \frac{\partial}{\partial \varphi_1}  \mathrm{Res}_{x_1} \biggl( x_1^{-1} \delta \left(\frac{x_2 +  s^2 \varphi_1 \varphi_2}{x_1} \right)  Y_s(v,(x_2 + s^2 \varphi_1 \varphi_2,\varphi_1)) \biggr) \biggr|_{\varphi_1 = 0}\\
&=& \frac{1}{s} \frac{\partial}{\partial \varphi_1}  Y_s(v,(x_2 + s^2 \varphi_1 \varphi_2,\varphi_1)) \biggr|_{\varphi_1 = 0}\\
&=& \frac{1}{s} \frac{\partial}{\partial \varphi_1}  \left( Y_s(v,(x_2,\varphi_1))  + s^2 \varphi_1 \varphi_2 \frac{\partial}{\partial x_2} Y_s(v,(x_2 ,\varphi_1))  \right)\biggr|_{\varphi_1 = 0}\\
&=& \biggl(  \frac{1}{s} \frac{\partial}{\partial \varphi_1}  Y_s(v,(x_2 ,\varphi_1)) 
+   s \varphi_2 \frac{\partial}{\partial x_2} Y_s(v,(x_2 ,\varphi_1))    \biggr) \biggr|_{\varphi_1 = 0}\\
&=& \sum_{n \in \mathbb{Z}} \left( \frac{1}{s} v_{n-\frac{1}{2}} x_2^{-n-1}  +  s \varphi_2 \frac{\partial}{\partial x_2}v_n x_2^{-n-1} \right) \\
&=& \left(  \frac{1}{s}\frac{\partial}{\partial \varphi_2} + s\varphi_2 \frac{\partial}{\partial x_2} \right)  Y_s(v,(x_2,\varphi_2 )),
\end{eqnarray*}
proving (\ref{Ds-derivative-property}).

To prove the $\mathcal{D}$-derivative property (\ref{D-derivative-N1}) we note that
\begin{eqnarray*}
Y_s(\mathcal{D} v, (x, \varphi)) &=&Y_s ((\mathcal{D}_s)^2v, (x, \varphi)) \ = \ \left( \frac{1}{s} \frac{\partial}{\partial \varphi} + s \varphi \frac{\partial}{\partial x} \right)^2 Y_s(v,(x, \varphi)) \\
&=& \frac{\partial}{\partial x} Y(v, (x, \varphi)),
\end{eqnarray*}
proving (\ref{D-derivative-N1}).

{}From the $\mathcal{D}_s$-derivative property (\ref{Ds-derivative-property}), we have
\begin{equation}
(\mathcal{D}_s v)_{n-\frac{1}{2}} = -s n v_{n-1},
\end{equation}
and thus
\[\mathcal{D}(v) \ = \  \mathcal{D}_s^2 v  \  = \  \frac{1}{s} (\mathcal{D}_s v)_{-\frac{3}{2}} \mathbf{1}  \ = \ \frac{1}{s} (s v_{-2}) \mathbf{1}  \ = \ v_{-2} \mathbf{1}, \]
proving (\ref{D-is-what-it-should-be2}).  The fact that $\mathrm{span}_{\bigwedge_*} \{\mathcal{D}, \mathcal{D}_s\}$ satisfy the commutation relations for $\mathfrak{osp}_{\bigwedge_*}(1|2)_{<0}$ is obvious.
\end{proof}

\begin{rema}\label{geometry-correspondence-remark-N1}
{\em The $\mathcal{D}_s$-derivative property (\ref{Ds-derivative-property}) implies that the operator $\mathcal{D}_s \in (\mathrm{End} \; V)^1$ corresponds to the $N=1$ $s$-deformed superconformal operators $D_s$ as defined in (\ref{define-D-derivative-Ds}) whereas $\mathcal{D}$ corresponds to the usual conformal operator $\frac{\partial}{\partial x}$. }
\end{rema}

{}From Proposition \ref{Taylor} with $\varphi_0^+ =  \varphi_0$, $\varphi_0^- = s^2 \varphi_0$, $\varphi^+ = \varphi$ and $\varphi^- = 0$, and using the $\mathcal{D}_s$- and $\mathcal{D}$-derivative properties (\ref{Ds-derivative-property}) and (\ref{D-derivative-N1}), we have 
\begin{eqnarray}
Y_s(e^{x_0 \mathcal{D} + s\varphi_0 \mathcal{D}_s }v, (x,\varphi)) &=& e^{x_0 \frac{\partial}{\partial x} + s\varphi_0  \left( \frac{1}{s}
\frac{\partial}{\partial \varphi} + s\varphi \frac{\partial}{\partial x} \right)} Y_s (v,(x,\varphi)) \label{exponential-L(-1)-G(-1/2)-N1} \\ 
&=& Y_s(v,(x + x_0 + s^2\varphi_0 \varphi, \varphi_0 + \varphi)) . \nonumber
\end{eqnarray}

{}From the creation property and (\ref{exponential-L(-1)-G(-1/2)-N1}), we have
\begin{equation}\label{for-skew-symmetry-N1}
e^{x \mathcal{D} + s\varphi \mathcal{D}_s } v \; = \; Y_s(v,(x,\varphi))
\mathbf{1} .
\end{equation}

\begin{rema}\label{Jacobi-symmetry-remark-N1}
{\em The left-hand side of the Jacobi identity (\ref{N1-Jacobi-identity}) is invariant under  the transformation 
\begin{equation}
(u,v,x_0,x_1,x_2,\varphi_1,\varphi_2) \longleftrightarrow ((-1)^{\eta(u)
\eta(v)} v,u,- x_0,x_2,x_1,\varphi_2,\varphi_1) .
\end{equation}
Thus the right-hand side of the Jacobi identity must be symmetric with respect to this also.  }
\end{rema}

\begin{prop} {\bf (skew supersymmetry)}
Let $V$ be an $N=1$ $D_s$-vertex superalgebra.  Recall the operators $\mathcal{D}_s$ and $\mathcal{D}$ defined in Proposition \ref{N1-derivative-prop}.  We have
\begin{equation}\label{skew-supersymmetry-N1}
Y(u,(x,\varphi))v =(-1)^{\eta(u) \eta(v)}  e^{x \mathcal{D} + s\varphi \mathcal{D}_s } Y_s(v,(-x,-\varphi))u 
\end{equation}
for $u,v$ of homogeneous sign in $V$.
\end{prop}

\begin{proof}
Using the symmetry of the left-hand side of the Jacobi identity as noted in Remark \ref{Jacobi-symmetry-remark-N1}, property (\ref{delta-substitute}) with $\varphi_1^+ = s\varphi_1$, $\varphi_2^- = s \varphi_s$, $\varphi_1^- = 0$, and $\varphi_2^+ = 0$, and (\ref{exponential-L(-1)-G(-1/2)-N1}), we obtain the following:
\begin{eqnarray*}
\lefteqn{x_2^{-1} \delta \biggl( \frac{x_1 - x_0 - s^2\varphi_1 \varphi_2}{x_2} \biggr) Y_s(Y_s(u,(x_0,  \varphi_1 - \varphi_2))v,(x_2, \varphi_2))}\\
&=& (-1)^{\eta(u) \eta(v)}  x_1^{-1} \delta \biggl( \frac{x_2 + x_0 - s^2\varphi_2 \varphi_1 }{x_1} \biggr) Y_s(Y_s(v,(-x_0,  \varphi_2 - \varphi_1))u,(x_1, \varphi_1)) \\
&=& (-1)^{\eta(u) \eta(v)}  x_1^{-1} \delta \biggl( \frac{x_2 + x_0 - s^2 \varphi_2 \varphi_1 }{x_1} \biggr) Y_s(Y_s(v,(-x_0,  \varphi_2 - \varphi_1))u,(x_2 + x_0 \\
& & \quad - \, s^2 \varphi_2 \varphi_1 , \varphi_1)) \\
&=& (-1)^{\eta(u) \eta(v)}  x_1^{-1} \delta \biggl( \frac{x_2 + x_0 -s^2 \varphi_2 \varphi_1}{x_1} \biggr) Y_s(e^{x_0 \mathcal{D} + s(\varphi_1 - \varphi_2) \mathcal{D}_s } Y_s(v,(-x_0,  \varphi_2 \\
& & \quad  - \,\varphi_1))u,(x_2 , \varphi_2)) .
\end{eqnarray*}  
Taking $\mathrm{Res}_{x_1}$ and using  (\ref{delta-2-terms-with-phis}) in this setting (i.e., with $\varphi_1^+ = s\varphi_1$, $\varphi_2^- = s \varphi_s$, $\varphi_1^- = 0$, and $\varphi_2^+ = 0$), we have
\begin{multline}
Y_s(Y_s(u,(x_0,  \varphi_1 - \varphi_2))v,(x_2, \varphi_2)) \\
= 
(-1)^{\eta(u) \eta(v)} Y_s(e^{x_0 \mathcal{D} + s(\varphi_1 - \varphi_2) \mathcal{D}_s } Y_s(v,(-x_0,  \varphi_2 - \varphi_1))u,(x_2 , \varphi_2)).
\end{multline}
Then acting on $\mathbf{1}$, taking the limit as $(x_2, \varphi_2)$ goes to zero, and using the creation property (\ref{N1-creation-property}), the result follows.
\end{proof}

\begin{prop}\label{bracket-prop-N1}
Let $V$ be an $N=1$ $D_s$-vertex superalgebra and let $\mathcal{D}_s$ and $\mathcal{D}$ be defined by (\ref{define-Ds}) and (\ref{define-D-N1}), respectively.   Then for $v \in V$ we have the following $\mathcal{D}_s$- and $\mathcal{D}$-bracket-derivative properties
\begin{eqnarray}
\left[ \mathcal{D}_s, Y_s(v, (x, \varphi)) \right] &=& \left(\frac{1}{s} \frac{\partial}{\partial \varphi} - s \varphi \frac{\partial}{\partial x} \right) Y_s(v, (x, \varphi)) \label{Ds-bracket-derivative-N1}\\
\left[ \mathcal{D}, Y_s(v, (x, \varphi)) \right] &=& \frac{\partial}{\partial x} Y_s(v, (x,\varphi)) ,
\end{eqnarray}
and the following $\mathcal{D}_s$- and $\mathcal{D}$-bracket properties
\begin{eqnarray}
\left[ \mathcal{D}_s, Y_s(v, (x, \varphi)) \right] &=& Y(\mathcal{D}_s v, (x, \varphi)) - 2s \varphi Y(\mathcal{D}v, (x, \varphi)) \label{Ds-bracket-N1} \\
\left[ \mathcal{D}, Y_s(v, (x, \varphi)) \right] &=& Y_s(\mathcal{D}v, (x, \varphi)) \label{D-bracket-N1}.
\end{eqnarray}
\end{prop}

\begin{proof}
Let  $T = x \mathcal{D} + s\varphi \mathcal{D}_s$.  Since $(\partial/\partial x) T = \mathcal{D}$, we have
\begin{equation}\label{for-proof-1-N1}
\frac{\partial}{\partial x} e^T = \mathcal{D} e^T + e^T \frac{\partial}{\partial x} .
\end{equation}
However since $(\partial/\partial \varphi)T = s\mathcal{D}_s$, we have that for $n \in \mathbb{N}$,
\begin{eqnarray*}
\frac{\partial}{\partial \varphi} T^n &=& s\mathcal{D}_s T^{n-1} + \sum_{k=1}^{n-1} \left( [T^k, s\mathcal{D}_s] T^{n-1-k} + s\mathcal{D}_s  T^{n-1}\right) + T^n \frac{\partial}{\partial \varphi} \\
&=&  \sum_{k=0}^{n-1} \left( 2ks^2 \varphi \mathcal{D} T^{k-1} T^{n-1-k} + s\mathcal{D}_s  T^{n-1}\right) + T^n \frac{\partial}{\partial \varphi} \\
&=& n(n-1) s^2 \varphi \mathcal{D} T^{n-2} + ns\mathcal{D}_s  T^{n-1} + T^n \frac{\partial}{\partial \varphi} ,
\end{eqnarray*}
and thus
\begin{equation}\label{for-proof-2-N1}
\frac{\partial}{\partial \varphi} e^T = (s^2\varphi \mathcal{D} + s\mathcal{D}_s) e^T + e^T \frac{\partial}{\partial \varphi} .
\end{equation}
Therefore for $u,v \in V$ of homogeneous sign in $V$, using skew supersymmetry (\ref{skew-supersymmetry-N1}), equations (\ref{for-proof-1-N1}) and (\ref{for-proof-2-N1}), and the $\mathcal{D}_s$-derivative property (\ref{Ds-derivative-property}), we obtain
\begin{eqnarray*}
\lefteqn{\left(\frac{1}{s} \frac{\partial}{\partial \varphi} - s\varphi \frac{\partial}{\partial x} \right) Y_s(u, (x, \varphi))v}\\
&=& (-1)^{\eta(v) \eta(u)}   \left( \frac{1}{s} \frac{\partial}{\partial \varphi} -s \varphi \frac{\partial}{\partial x} \right) e^{x \mathcal{D} + s\varphi \mathcal{D}_s} Y_s(v,(-x,-\varphi))u \\
&=& (-1)^{\eta(v) \eta(u)}   \biggl(\Bigl(s\varphi \mathcal{D} + \mathcal{D}_s - s\varphi \mathcal{D}  \Bigr) e^{x \mathcal{D} + s\varphi \mathcal{D}_s } Y_s(v,(-x,-\varphi))u \\
& & \quad + \, e^{x \mathcal{D} +s \varphi \mathcal{D}_s } \left( \frac{1}{s} \frac{\partial}{\partial \varphi}  - s \varphi \frac{\partial}{\partial x}  \right) Y_s(v,(-x,-\varphi))u  \biggr) \\
&=& (-1)^{\eta(v) \eta(u)}   \biggl( \mathcal{D}_s  e^{x \mathcal{D} + s\varphi \mathcal{D}_s} Y_s(v,(-x,-\varphi))u \\
& & \quad - \,  e^{x \mathcal{D} + s\varphi \mathcal{D}_s} \left( \frac{1}{s} \frac{\partial}{\partial (-\varphi)}  + s(-\varphi) \frac{\partial}{\partial (-x)}  \right) Y_s(v,(-x,-\varphi))u  \biggr) \\
\end{eqnarray*}
\begin{eqnarray*}
&=& (-1)^{\eta(v) \eta(u)}   \biggl(\mathcal{D}_s e^{x \mathcal{D} + s\varphi \mathcal{D}_s } Y_s(v,(-x,-\varphi))u - e^{x \mathcal{D} + s\varphi \mathcal{D}_s }Y_s(\mathcal{D}_s v,(-x,-\varphi))u  \biggr) \\
&=&  \mathcal{D}_s Y_s(u,(x,\varphi))v  -  (-1)^{\eta(u)}   Y_s(u,(x,\varphi)) \mathcal{D}_s v \\
&=& \left[ \mathcal{D}_s , Y_s(u,(x,\varphi)) \right]v .
\end{eqnarray*}
Using this $\mathcal{D}_s$-bracket property and the $\mathcal{D}$-derivative property (\ref{D-derivative-N1}), we have
\begin{eqnarray*}
\frac{\partial}{\partial x} Y(u, (x, \varphi))
&=& - \left(\frac{1}{s} \frac{\partial}{\partial \varphi} -s \varphi \frac{\partial}{\partial x} \right)^2 Y_s(u, (x, \varphi))\\
&=& [ \mathcal{D}_s, [\mathcal{D}_s, Y_s(u, (x, \varphi)) ]]\\
&=& \frac{1}{2} [ [\mathcal{D}_s, \mathcal{D}_s], Y_s(u, (x, \varphi)) ]\\
&=& [\mathcal{D}, Y(u, (x, \varphi)) ].
\end{eqnarray*}

Finally (\ref{Ds-bracket-N1}) and (\ref{D-bracket-N1}) follow {}from the $\mathcal{D}_s$- and $\mathcal{D}$-derivative properties (\ref{Ds-derivative-property}) and (\ref{D-derivative-N1}), respectively. 
\end{proof}

\begin{rema}\label{Ds-construction-remark} 
{\em Let $(V, Y_s(\cdot, (x, \varphi)), \mathbf{1})$ be an $N=1$ $\mathcal{D}_s$-vertex superalgebra and let $Y_s(\cdot,x) = Y_s(\cdot,(x, 0))$.  Proposition \ref{N1-derivative-prop} implies that
\begin{equation}\label{constructing-s-operators}
Y_s (v, (x, \varphi)) = Y_s(v,x) + s \varphi Y_s(\mathcal{D}_sv,x).  
\end{equation}
This shows that if $(V,Y_1(\cdot,(x,\varphi)),\mathbf{1})$ is an $N=1$ vertex superalgebra in the usual sense (i.e. satisfies Definition \ref{N1-vertex-superalgebra-definition} with $s = 1$) , then $(V,Y_1(\cdot,(x,s\varphi)),\mathbf{1})$ is an $N=1$ $D_s$-vertex superalgebra.   In fact this implies that in making a choice as to how to add the odd variable components to a vertex superalgebra which is a representation of $\mathfrak{osp}_{\bigwedge_*} (1|2)_{<0}$ one is choosing a superconformal structure $D_s$.  Thus, on the one hand, the choice of how one constructs the odd formal variable component of an $N=1$ vertex superalgebra determines whether it is $D_1$-superconformal or $D_s$-superconformal for $s \neq \pm1$.  On the other hand, for any given $N=1$ vertex superalgebra, one can construct a continuous one-parameter family of $N=1$ $D_s$-vertex superalgebras parameterized by $s \in (\bigwedge_*^0)^\times$.  }
\end{rema}

\begin{defn}
{\em An} $N = 1$ $D_s$-Neveu-Schwarz vertex operator superalgebra over $\bigwedge_*$ and with odd formal variables {\em is a $\frac{1}{2} \mathbb{Z}$-graded $\bigwedge_*$-module (graded by} weights{\em)  
\begin{equation}
V = \coprod_{n \in \frac{1}{2}\mathbb{Z}} V_{(n)} 
\end{equation}
such that 
\begin{equation}
\dim V_{(n)} < \infty \qquad \mbox{for $n \in \frac{1}{2} \mathbb{Z}$,} 
\end{equation}
\begin{equation}
V_{(n)} = 0 \qquad \mbox{for $n$ sufficiently negative} , 
\end{equation}
equipped with an $N=1$ $D_s$-vertex superalgebra structure $(V, Y_s(\cdot, (x,\varphi)), \mathbf{1})$, and a distinguished vector $\tau \in V_{(3/2)}^1$ (the {\em  $N=1$ Neveu-Schwarz element} or {\em  $N=1$ $D_s$-superconformal element}), satisfying the following conditions:  the $N=1$ Neveu-Schwarz algebra relations hold: 
\begin{eqnarray}
\left[L(m) ,L(n) \right] \! \! \! &=&  \! \! \! (m - n)L(m + n) \! + \!  \frac{1}{12} (m^3 - m) \delta_{m + n, 0} \; c_V  ,\label{n1-first}\\
\bigl[L(m),G(n + 1/2)\bigr]  \! \! \! &=&  \! \! \! \Bigl(\frac{m}{2} - n - \frac{1}{2} \Bigr) G(m + n + 1/2) , \\
\qquad \ \ \ \bigl[ G(m + 1/2) , G(n - 1/2) \bigr]  \! \! \! &=&  \! \! \!  2L(m + n)  + \frac{1}{3}(m^2 + m) \delta_{m + n , 0}  \;c _V,  \label{n1-last}
\end{eqnarray}
for $m,n \in \mathbb{Z}$,  where 
\begin{equation}
G(n - 1/2) = \tau_n,\quad \mbox{and} \quad   2sL(n) = \tau_{n+ \frac{1}{2}} , \quad \mbox{for $n \in \mathbb{Z}$}
\end{equation}
i.e., 
\begin{equation}\label{N1-tau}
Y_s(\tau,(x,\varphi)) = \sum_{n \in \mathbb{Z}} \left(G (n+1/2) x^{-n-2}   + 2 s \varphi  L(n) x^{- n - 2} \right)
\end{equation}
and $c_V \in \mathbb{C}$ (the} central charge{\em);  for $n \in \frac{1}{2} \mathbb{Z}$ and $v \in V_{(n)}$
\begin{equation}
L(0)v = nv ;
\end{equation}
and finally the} $G(-1/2)$-derivative property {\em holds:
\begin{equation}\label{G-derivative-N1}
\biggl( \frac{1}{s} \frac{\partial}{\partial \varphi} + s\varphi  \frac{\partial}{\partial x} \biggr) Y_s(v,(x,\varphi)) = Y_s(G(- 1/2)v,(x,\varphi)) .\\
\end{equation} }
\end{defn}

\medskip

The $N=1$ $D_s$-Neveu-Schwarz vertex operator superalgebra with odd formal variables is denoted by 
\[(V,Y_s(\cdot,(x,\varphi)),\mathbf{1},\tau).\] 

\begin{rema}  
{\em 
If $s = \pm1$, then $(V,Y_s(\cdot,(x,\varphi)),\mathbf{1},s\tau)$ is an $N=1$ Neveu-Schwarz vertex operator superalgebra in the usual sense as defined in \cite{B-announce} (see also \cite{B-thesis} and \cite{B-vosas}).  Thus given an $N=1$ Neveu-Schwarz vertex operator superalgebra $(V, Y(\cdot,(x,\varphi)), \mathbf{1}, \tau)$ with odd formal variables, the one parameter family of deformed $N=1$ vertex algebras given by the $N=1$ $\mathcal{D}_s$-Neveu-Schwarz vertex operator superalgebras $(V, Y(\cdot,(x,s\varphi)), \mathbf{1}, \tau)$ (see Remark \ref{Ds-construction-remark}) continuously deforms between the two $N=1$ Neveu-Schwarz vertex operator superalgebras for $s = \pm 1$.  }
\end{rema}

If we set $\varphi = 0$, then for all $s \in(\bigwedge_*^0)^\times$, we have that $(V, Y_s(\cdot,(x,0)), \mathbf{1}, \tau)$ is an $N=1$ Neveu-Schwarz vertex operator superalgebra without odd variables in the usual sense as defined in \cite{B-vosas} (see also \cite{B-thesis}).   Thus it makes sense to denote $Y_s (\cdot, (x,0))$ simply by $Y(v,x)$.

\begin{rema}\label{N1-from-N2-remark}
{\em As mentioned as motivation in the beginning of the section, if $(V, Y(\cdot, (x, \varphi^{(1)}, \varphi^{(2)})), \mathbf{1}, \mu)$ is an $N=2$ Neveu-Schwarz vertex operator superalgebra with two odd formal variables in nonhomogeneous coordinates, then $(V, Y(\cdot, (x$, $\varphi^{(1)}, 0)), \mathbf{1}, \tau^{1})$ and $(V, Y(\cdot, (x, 0, \varphi^{(2)})), \mathbf{1}, \tau^{2})$ are both $N=1$ Neveu-Schwarz vertex operator superalgebras with one odd formal variable where $\tau^j = G^{(j)} (-3/2) \mathbf{1}$ for $j = 1,2$, i.e. $\tau^1 = i G^{(2)}(-1/2) \mu$ and $\tau^2 = - iG^{(1)}(-1/2) \mu$.  In fact  $(V, Y(\cdot, (x, -\varphi^{(1)}$, $0)), \mathbf{1}, -\tau^{1})$ and $(V, Y(\cdot, (x, 0, -\varphi^{(2)})), \mathbf{1}, -\tau^{2})$ are also $N=1$ Neveu-Schwarz vertex operator superalgebras.  But this is simply due to the fact that if $(V, Y(\cdot, (x, \varphi)), \mathbf{1}, \tau)$ is an $N=1$ Neveu-Schwarz vertex operator superalgebra then so is $(V, Y(\cdot, (x, -\varphi))$, $\mathbf{1},- \tau)$.)  That is, it is due to the action of the group of automorphisms of the $N=1$ Neveu-Schwarz algebra given by $\mathbb{Z}_2$.   And in fact $(V,Y(\cdot, (x, \varphi)), \mathbf{1}, \tau)$ and  $(V, Y(\cdot, (x, -\varphi))$, $\mathbf{1},- \tau)$ are isomorphic as $N=1$ Neveu-Schwarz vertex operator superalgebras.  Moreover, similar statement can be made about $N=1$ vertex superalgebra structures on $N=2$ vertex superalgebras. }
\end{rema} 

The statements in the preceding remark are certainly not new or surprising, however as we shall see below and in the next section, there is a much deeper connection between $N=1$ $D_s$-Neveu-Schwarz vertex operator superalgebras and $N=2$ Neveu-Schwarz vertex operator superalgebras.  

Now for $j = 1,2$, define 
\begin{equation}
\pi^{(j)} Y(v, (x, \varphi^{(1)}, \varphi^{(2)})) = \sum_{n \in \mathbb{Z}} \Bigl( v_n + \varphi^{(j)} v^{(j)}_{n - \frac{1}{2}}   \Bigr) x^{-n-1} = Y(v,(x,\varphi^{(j)})) ,
\end{equation}
and for $\beta  \in \bigwedge_*^0$, define
\begin{eqnarray}
s_{(1,1)}(\beta) = i \sinh \beta & & s_{(1,2)}(\beta) = \cosh \beta\\
s_{(2,1)}(\beta) = \cosh \beta & & s_{(2,2)}(\beta) = -i \sinh \beta.
\end{eqnarray}
For $j, k = 1,2$, let
\begin{multline}
Y_{s_{(j,k)}(\beta)}(\cdot, (x, \varphi))\\
=\pi^{(j)} Y(\cdot, (x, \varphi^{(1)} \cosh \beta + i \varphi^{(2)} \sinh \beta, \varphi^{(2)} \cosh \beta - i \varphi^{(1)} \sinh \beta)) |_{\varphi^{(k)} = 0} .
\end{multline}

The following proposition follows {}from Lemma \ref{nonhomo-iso-lemma} and Remark \ref{N1-from-N2-remark}:

\begin{prop} Let $(V, Y(\cdot, (x, \varphi^{(1)}, \varphi^{(2)})), \mathbf{1}, \mu)$ be an $N=2$ Neveu-Schwarz vertex operator superalgebra.  For $\beta \in \bigwedge_*^0$ and $j,k = 1,2$, let $Y_{s_{(j,k)} (\beta)} (\cdot, (x,\varphi))$ be defined as above.  Then $(V, Y_{s_{(j,k)}(\beta)} (\cdot, (x, \varphi)), \mathbf{1}, \tau =- iG^{(k)}(-1/2) \mu)$ is an $N=1$ $D_{s_{(j,k)}(\beta)}$-Neveu-Schwarz vertex operator superalgebra.
\end{prop}

Thus we see that the $J(0)$ operator gives rise to continuous families of $D_s$-deformed $N=1$ Neveu-Schwarz vertex operator superalgebras.  In the next section we will investigate this phenomenon {}from another point of view, and see that it is a reflection of the $N=1$ superanalytic geometry underlying the notion of $N=2$ Neveu-Schwarz vertex operator superalgebra.



\begin{rema}{\em Just as we have deformed the odd variable in the notion of $N=1$ vertex superalgebra while maintaining the $N=1$ superconformality property of the resulting $\mathcal{D}_s$ operator (or equivalently maintaining the property that the resulting algebraic structure of $N=1$ $D_s$-vertex superalgebra is still naturally a representation of $\mathfrak{osp}_{\bigwedge_*} (1|2)_{<0}$), we could also deform the two odd variables in the notion of $N=2$ vertex superalgebra while maintaining the $N=2$ superconformality property of the resulting operators which would act as deformed versions of $\mathcal{D}^{(j)}$ for $j = 1,2$ (or equivalently maintaining the property that the resulting algebraic structure which we could call an ``$N=2$ $D^{(j)}_s$-vertex superalgebra" is still naturally a representation of $\mathfrak{osp}_{\bigwedge_*} (2|2)_{<0}$).  And we could extend this to include a representation of the $N=2$ Neveu-Schwarz algebra just as we did in the $N=1$ case with the notion of $N=2$ $D_s$-Neveu-Schwarz vertex operator superalgebra.  However, we do not follow this obvious extension in this paper.}
\end{rema}

\subsection{The $N=2$ Neveu-Schwarz algebra arising {}from continuous deformations of the $N=1$ Neveu-Schwarz algebra}

Regarding $N=1$ $D_s$-Neveu-Schwarz vertex operator superalgebras as continuously deformed
$N=1$ Neveu-Schwarz vertex operator superalgebras, a natural question arises: Is there a (necessarily)
non-superconformal infinitesimal change of coordinates operator associated to this deformation?  That is, just as $L_0(x,\varphi) = - ( x \frac{\partial}{\partial x} + \frac{\varphi}{2} \frac{\partial}{\partial \varphi})$ is the infinitesimal local superconformal coordinate corresponding to a deformation 
\[ a^{-2 L_0(x,\varphi)} (x,\varphi) = (a^2 x, a \varphi)\]
for $a \in (\bigwedge_*^0)^\times$, in what setting can we expect to find an infinitesimal superanalytic coordinate $\varphi \frac{\partial}{\partial \varphi}$ giving the deformation
\[s^{ \varphi \frac{\partial}{\partial \varphi}} (x,\varphi) = (x,s\varphi)\]
for $s \in (\bigwedge_*^0)^\times$?  We will find that such an operator naturally occurs in the $N=2$ superconformal setting.

To deform {}from an $N=1$ Neveu-Schwarz vertex operator superalgebra into an $N=1$ $D_s$-Neveu-Schwarz vertex operator superalgebra, we need to extend the representation of the $N=1$ Neveu-Schwarz algebra in terms of infinitesimal $N=1$ superconformal local coordinates  to include a term that behaves as $\varphi \frac{\partial}{\partial\varphi}$.  By Proposition \ref{moduli-prop-N1}, the superderivations  $L_n (x,\varphi)$ and $G_{n-\frac{1}{2}} (x,\varphi)$, for $n \in \mathbb{Z}$, given by (\ref{L-notation-N1}) and (\ref{G-notation-N1}) give a representation of the $N=1$ Neveu-Schwarz Lie superalgebra with central charge zero in terms of infinitesimal superconformal transformations in $\mathrm{Der}(\bigwedge_* [[x, x^{-1}]][\varphi])$.   If we let
\begin{equation}
J_0 (x,\varphi)  = \varphi \frac{\partial}{\partial \varphi} ,  
\end{equation}
then
\begin{equation}
[ L_m(x,\varphi),J_0(x,\varphi)] = 0,
\end{equation}
and
\begin{eqnarray}
\left[ J_0(x,\varphi) , G_{n+\frac{1}{2}} (x,\varphi)\right] &=& x^{n+1} \left( \varphi \frac{\partial}{\partial x} + \frac{\partial}{\partial \varphi} \right) \\
&=& -i G^*_{n+\frac{1}{2}} (x,\varphi), \nonumber
\end{eqnarray}
where we define
\begin{equation}
G^*_{n-\frac{1}{2}} (x,\varphi) = ix^n \left( \frac{\partial}{\partial \varphi} + \varphi \frac{\partial}{\partial x}\right)  \quad \mbox{for $n \in \mathbb{Z}$}.
\end{equation}
We then also have that 
\begin{eqnarray}
\left[ L_m(x,\varphi), G^*_{n + \frac{1}{2}}(x,\varphi) \right]
&=& (- n +  \frac{m-1}{2} ) G^*_{m+n+\frac{1}{2}}(x,\varphi) ,\\
\left[ J_0(x,\varphi), G^*_{n + \frac{1}{2}} (x,\varphi)\right] 
&=&  iG_{n+ \frac{1}{2}} (x,\varphi),\\
\left[G^*_{m+\frac{1}{2}}(x,\varphi), G^*_{n - \frac{1}{2}} (x,\varphi)\right]
&=& 2L_{m+n} (x,\varphi),
\end{eqnarray}
\begin{eqnarray}
\left[G_{m+\frac{1}{2}}(x,\varphi), G^*_{n - \frac{1}{2}}(x,\varphi) \right]
&=& -i(m-n +1) x^{m +n} \varphi \frac{\partial}{\partial \varphi} \\
&=& -i(m-n+1) J_{m +n}(x,\varphi), \nonumber
\end{eqnarray}
where we define
\begin{equation}
J_n (x,\varphi)  = x^n \varphi \frac{\partial}{\partial \varphi} \quad \mbox{for $n \in \mathbb{Z}$}. 
\end{equation}

Finally, we have that 
\begin{eqnarray}
\left[J_m(x,\varphi), J_n(x,\varphi) \right] 
&=& 0 , \\
\left[L_m(x,\varphi),J_n(x,\varphi) \right]
&=& -nJ_{m+n}(x,\varphi) , \\
\left[J_m(x,\varphi),G_{n + \frac{1}{2}} (x,\varphi)\right] 
&=& -iG^*_{m+n+\frac{1}{2}}(x,\varphi), \\
\left[J_m(x,\varphi),G^*_{n + \frac{1}{2}}(x,\varphi) \right] 
&=& iG_{m+n+\frac{1}{2}}(x,\varphi). 
\end{eqnarray}

In conclusion, the operators $L_n (x, \varphi)$ and $G_{n-\frac{1}{2}} (x, \varphi)$, for $n \in \mathbb{Z}$, given by (\ref{L-notation-N1}) and (\ref{G-notation-N1}), along with $J_0(x, \varphi) = \varphi \frac{\partial}{\partial \varphi}$  generate the basis of operators
\begin{eqnarray}
L_n (x, \varphi)  &=& - \left( x^{n + 1} \frac{\partial}{\partial x} + (\frac{n + 1}{2})\varphi x^n \frac{\partial}{\partial \varphi} \right) \label{recap-L}  \\ 
J_n (x, \varphi)   &=& x^n \varphi \frac{\partial}{\partial \varphi}  \\
G_{n-\frac{1}{2}} (x, \varphi)  &=& - x^n \left( \frac{\partial}{\partial \varphi} - \varphi \frac{\partial}{\partial x}\right) \\
G^*_{n-\frac{1}{2}} (x, \varphi)  &=&  i x^n \left( \frac{\partial}{\partial \varphi} + \varphi \frac{\partial}{\partial x}\right) \label{recap-G*}
\end{eqnarray}
which satisfy the commutation relations for the $N=2$ Neveu-Schwarz algebra in the nonhomogeneous basis (\ref{Virasoro-relation2})--(\ref{transformed-Neveu-Schwarz-relation-last}) under the identification $L_n (x, \varphi) \mapsto L_n$, $J_n (x, \varphi) \mapsto J_n$, $G_{n + \frac{1}{2}} (x, \varphi) \mapsto G^{(1)}_{n+\frac{1}{2}}$, $G^*_{n + \frac{1}{2}} (x, \varphi) \mapsto G^{(2)}_{n+ \frac{1}{2}}$, for $n \in \mathbb{Z}$, and $d \mapsto 0$. 

Note that 
\begin{eqnarray}
G^*_{n-\frac{1}{2}} (x, \varphi)  &=&  i x^n \left( \frac{\partial}{\partial \varphi} + \varphi \frac{\partial}{\partial x}\right) \ = \  - x^n \left( \frac{\partial}{\partial i \varphi} - i \varphi \frac{\partial}{\partial x}\right) \\
&=& G_{n-\frac{1}{2}} (x, i \varphi). \nonumber
\end{eqnarray}

In fact, it is straightforward to show that there are exactly two extensions of the set of operators $S = \{L_n(x, \varphi), G_{n - \frac{1}{2}} (x, \varphi) \; | \; n \in \mathbb{Z} \}$ with $L_n(x, \varphi)$, and $G_{n - \frac{1}{2}} (x, \varphi)$ given by (\ref{L-notation-N1}) and (\ref{G-notation-N1}), respectively,  by an even superderivation $J(0) \in \mathrm{Der} (\bigwedge_* [[x, x^{-1}]][\varphi])$ such that $S \cup \{J(0) \}$ generate a representation of the $N=2$ Neveu-Schwarz Lie superalgebra with central charge zero.    These two are given by $J(0) = \pm \varphi \frac{\partial}{\partial \varphi}$, reflecting the Virasoro preserving automorphism of the $N=2$ Neveu-Schwarz algebra given by (\ref{nonhomo-auto2}).

Note that for $n \in \mathbb{Z}$, we have 
\begin{eqnarray}
x^n \frac{\partial}{\partial x} &=& - L_{n-1}(x, \varphi) - \frac{n}{2} J_{n-1}(x, \varphi) \\
x^n \varphi \frac{\partial}{\partial x} &=& \frac{1}{2} \left(G_{n- \frac{1}{2}} (x, \varphi) - i G^*_{n-\frac{1}{2}} (x, \varphi) \right) \\
x^n  \frac{\partial}{\partial \varphi} &=& \frac{1}{2} \left(-G_{n- \frac{1}{2}} (x, \varphi) - i G^*_{n-\frac{1}{2}} (x, \varphi) \right)\\
x^n \varphi \frac{\partial}{\partial \varphi} &=& J_n(x, \varphi);
\end{eqnarray}
that is the superderivations given by (\ref{recap-L})--(\ref{recap-G*}) generate $\mathrm{Der} (\bigwedge_* [[x, x^{-1}]] [\varphi])$.  And writing $G^\pm_{n-\frac{1}{2}} (x, \varphi) = \frac{1}{\sqrt{2}} ( G_{n-\frac{1}{2}} (x, \varphi) \mp i G^*_{n- \frac{1}{2}} (x, \varphi))$, we have that
\begin{equation}
x^n \varphi \frac{\partial}{\partial x} = \frac{1}{\sqrt{2}} G^+_{n- \frac{1}{2}} (x, \varphi), \quad \mathrm{and} \quad 
x^n  \frac{\partial}{\partial \varphi} = - \frac{1}{\sqrt{2}} G^-_{n- \frac{1}{2}} (x, \varphi).
\end{equation}

We recall {}from for instance \cite{B-memoir} that a superanalytic function in one even variable $z \in \bigwedge_*^0$ and one odd variable $\theta \in \bigwedge_*^1$ is one that has well-defined derivatives with respect to $z$ and $\theta$.  This is a much less strict condition than that for superconformality which requires that the function, in addition to being superanalytic, satisfy the much more stringent condition that it transform the operator $D_1 = \partial / \partial \theta + \theta \partial / \partial z$ homogeneously of degree one.  

\begin{prop}\label{moduli-prop-N2-one-variable}
There is a bijection between formal $N=1$ superanalytic (not necessarily superconformal) functions 
\begin{equation}\label{where-H-lives}
f(x, \varphi) = (\tilde{x}, \tilde{\varphi}) \in (x \mbox{$\bigwedge_*^0$}[[x]] \oplus \varphi x \mbox{$\bigwedge_*^1$}[[x]], \ x\mbox{$\bigwedge_*^1$}[[x]] \oplus \varphi \mbox{$\bigwedge_*^0$}[[x]]),
\end{equation} 
(that is vanishing at $(x, \varphi) = (0,0)$ and with even part vanishing at $x = 0$) which are invertible in a neighborhood of zero, and expressions of the form
\begin{multline}\label{change-of-variables-superanalytic}
\exp \Biggl(-  \sum_{n \in \Z} \Bigl( A^{(1)}_n L_n(x,\varphi) + A^{(2)}_n J_n(x,\varphi) +  M^{(1)}_{n - \frac{1}{2}} G^{(1)}_{n -\frac{1}{2}} (x, \varphi) \\
+  M^{(2)}_{n - \frac{1}{2}} G^{(2)}_{n -\frac{1}{2}} (x, \varphi) \Bigr) \Biggr)  \cdot  (a_0^{(1)})^{-2L_0(x,\varphi)} \cdot  (a_0^{(2)})^{-J_0(x,\varphi)} \cdot (x, \varphi) 
\end{multline}
for $(a_0^{(1)}, a_0^{(2)}) \in ((\bigwedge_*^0)^\times)^2/\langle (-1,-1) \rangle$, and $A_n^{(j)} \in \bigwedge_*^0$ and $M^{(j)}_{n- 1/2} \in \bigwedge_*^1$, for $n \in \Z$, and $j = 1,2$
where 
\begin{eqnarray}
L_n(x,\varphi) &=&  - \biggl( x^{n+ 1} \frac{\partial}{\partial x} + (\frac{n + 1}{2})x^n \varphi \frac{\partial }{\partial \varphi}  \biggr) \label{L-notation-N1-N2} \\
J_n(x,\varphi) &=& x^n \varphi \frac{\partial }{\partial \varphi}  \label{J-notation-N1-N2}  \\
G^{(j)}_{n-\frac{1}{2}} (x,\varphi) &=& (-i)^{j+1} x^n \Bigl(   \frac{\partial }{\partial \varphi}  +  (-1)^{j} \varphi  \frac{\partial}{\partial x} \Bigr) \label{G-notation-N1-N2}  \end{eqnarray}
are superderivations in  $\mbox{Der} (\bigwedge_*[[x,x^{-1}]] [\varphi])$, for $n \in \mathbb{Z}$, which give a representation of the $N=2$ Neveu-Schwarz algebra with central charge zero in the nonhomogeneous basis; that is  (\ref{L-notation-N1-N2})--(\ref{G-notation-N1-N2}) satisfy (\ref{Virasoro-relation2})--(\ref{transformed-Neveu-Schwarz-relation-last}) with $d = 0$.  Similarly, there is a bijection between formal $N=1$ superanalytic functions 
\begin{equation}\label{where-H-lives2}
f(x, \varphi) = (\tilde{x}, \tilde{\varphi}) \in (x^{-1} \mbox{$\bigwedge_*^0$}[[x^{-1}]] \oplus \varphi x^{-2} \mbox{$\bigwedge_*^1$}[[x^{-1}]], \ x^{-1}\mbox{$\bigwedge_*^1$}[[x^{-1}]] \oplus \varphi x^{-1}\mbox{$\bigwedge_*^0$}[[x^{-1}]]),
\end{equation} 
(that is vanishing at $(\infty, 0)$ and such that $x\frac{\partial}{\partial \varphi}\tilde{x}$ also vanishes as $x \rightarrow \infty$) which are invertible in a neighborhood of $(\infty, 0)$ and expressions of the form
\begin{multline}\label{change-of-variables-superanalytic-infty}
\exp \Biggl(\sum_{n \in \Z} \Bigl( A^{(1)}_n L_{-n}(x,\varphi) - A^{(2)}_n J_{-n}(x,\varphi) +  iM^{(1)}_{n - \frac{1}{2}} G^{(1)}_{-n +\frac{1}{2}} (x, \varphi) \\
+  iM^{(2)}_{n - \frac{1}{2}} G^{(2)}_{-n +\frac{1}{2}} (x, \varphi) \Bigr) \Biggr)  \cdot  (a_0^{(1)})^{2L_0(x,\varphi)} \cdot  (a_0^{(2)})^{-J_0(x,\varphi)} \cdot  \left(\frac{1}{x}, \frac{i\varphi}{x} \right) .
\end{multline}

\end{prop}

\begin{proof}
An $N=1$ superanalytic function vanishing at zero, invertible in a neighborhood of zero satisfying $f(x, \varphi)|_{x = 0} = (0, b_0\varphi)$ for some $b_0 \in (\bigwedge_*^0)^\times$, can be written as $f(x, \varphi) = (\tilde{x}, \tilde{\varphi})$ where
\begin{eqnarray}
\tilde{x} &=& a_0 \biggl( x + \sum_{n  \in \mathbb{Z}_+} a_n x^{n+1} + \varphi \sum_{n \in \mathbb{Z}_+} k_{n-\frac{1}{2}} x^n \biggr) \label{superanalytic1}\\
\tilde{\varphi} &=& b_0 \biggl( \sum_{n  \in \mathbb{Z}_+} m_{n-\frac{1}{2}} x^{n} + \varphi \Bigl( 1+ \sum_{n \in \mathbb{Z}_+} b_n x^n\Bigr) \biggr), \label{superanalytic2}
\end{eqnarray}
for $a_0, b_0 \in (\bigwedge_*^0)^\times$ and $a_n, b_n \in \bigwedge_*^0$ and $\alpha_n, \beta_n \in \bigwedge_*^1$ for $n \in \mathbb{Z}_+$.

It is clear that (\ref{change-of-variables-superanalytic}) is superanalytic and vanishing at zero with invertible even coefficient of $x$ and invertible odd coefficient of $\varphi$.  It remains to show that any superanalytic function $f(x, \varphi) = (\tilde{x}, \tilde{\varphi})$ given by (\ref{superanalytic1}) and (\ref{superanalytic2}) can be expressed as (\ref{change-of-variables-superanalytic}).

We first rewrite (\ref{change-of-variables-superanalytic}) as
\begin{multline}\label{change-of-variables-superanalytic2}
\exp \Biggl(\sum_{n \in \Z} \Bigl( A_n x^{n+1} \frac{\partial}{\partial x} + B_n \varphi x^n \frac{\partial}{\partial \varphi} +  K_{n - \frac{1}{2}} \varphi x^n \frac{\partial}{\partial x}  \\
+  M_{n - \frac{1}{2}} x^n \frac{\partial}{\partial \varphi}\Bigr) \Biggr)  \cdot  (a_0^{(1)})^{-2L_0(x,\varphi)} \cdot  (a_0^{(2)})^{-J_0(x,\varphi)} \cdot (x, \varphi),
\end{multline}
where 
\begin{eqnarray}
A_n &=& A_n^{(1)} \label{change-basis-first}\\
B_n &=& \left( \frac{n+1}{2} \right) A_n^{(1)} - A_n^{(2)}\\
K_{n-\frac{1}{2}} &=&  - (M^{(1)}_{n-\frac{1}{2}} + i M^{(2)}_{n-\frac{1}{2}})\\
M_{n- \frac{1}{2}} &=& M^{(1)}_{n-\frac{1}{2}} - i M^{(2)}_{n-\frac{1}{2}}. \label{change-basis-last}
\end{eqnarray}
Then given $f(x, \varphi) = (\tilde{x}, \tilde{\varphi})$ of the form (\ref{superanalytic1}) and (\ref{superanalytic2}) (i.e., of the form (\ref{where-H-lives})), equating $f$ with (\ref{change-of-variables-superanalytic2}), we have the following system of equations for the unknown quantities $a_0^{(1)}, a_0^{(2)} \in  ((\bigwedge_*^0)^\times)^2$, and $A_n, B_n \in \bigwedge_*^0$, $K_{n - \frac{1}{2}}, M_{n - \frac{1}{2}} \in \bigwedge_*^1$, for $n \in \mathbb{Z}_+$,  in terms of the known quantities $b_0, a_0 \in  ((\bigwedge_*^0)^\times)^2$, and $a_n, b_n \in \bigwedge_*^0$, $k_{n - \frac{1}{2}}, m_{n - \frac{1}{2}} \in \bigwedge_*^1$, for $n \in \mathbb{Z}_+$, obtained by expanding (\ref{change-of-variables-superanalytic2})
\begin{eqnarray}
a_0 &=& (a_0^{(1)})^2 \label{find-a}\\
b_0 &=& a_0^{(1)} a_0^{(2)}  \label{find-b}\\
a_n &=& A_n + p^{(1)}_n(A_1,\dots, A_{n-1}, B_1, \dots, B_{n-1}, M_{\frac{1}{2}}, \dots, M_{n - \frac{1}{2}}, \label{find-A} \\
& & \quad  K_{\frac{1}{2}}, \dots, K_{n - \frac{1}{2}}) \nonumber \\
b_n &=& B_n + p^{(2)}_n(A_1,\dots, A_{n-1}, B_1, \dots, B_{n-1}, M_{\frac{1}{2}}, \dots, M_{n - \frac{1}{2}}, \\
& & \quad K_{\frac{1}{2}}, \dots, K_{n - \frac{1}{2}}) \nonumber \\
\qquad k_{n-\frac{1}{2}} &=& K_{n-\frac{1}{2}} + q^{(1)}_n(A_1,\dots, A_{n-1}, B_1, \dots, B_{n-1}, M_{\frac{1}{2}}, \dots, M_{n - \frac{3}{2}}, \\
& & \quad K_{\frac{1}{2}}, \dots, K_{n - \frac{3}{2}}) \nonumber 
\end{eqnarray}
\begin{eqnarray}
\qquad m_{n-\frac{1}{2}} &=& M_{n-\frac{1}{2}} + q^{(2)}_n(A_1,\dots, A_{n-1}, B_1, \dots, B_{n-1}, M_{\frac{1}{2}}, \dots, M_{n - \frac{3}{2}}, \label{find-M} \\
& & \quad K_{\frac{1}{2}}, \dots, K_{n - \frac{3}{2}}),\nonumber
\end{eqnarray}
where
\begin{equation}
p^{(j)}_n \in \mbox{$\bigwedge_*$}[A_1,\dots, A_{n-1}, B_1, \dots, B_{n-1}, M_{\frac{1}{2}}, \dots, M_{n - \frac{1}{2}}, K_{\frac{1}{2}}, \dots, K_{n - \frac{1}{2}}]
\end{equation}
and 
\begin{equation}
q^{(j)}_n  \in \mbox{$\bigwedge_*$}[A_1,\dots, A_{n-1}, B_1, \dots, B_{n-1}, M_{\frac{1}{2}}, \dots, M_{n - \frac{3}{2}}, K_{\frac{1}{2}}, \dots, K_{n - \frac{3}{2}}]
\end{equation}
for $j = 1,2$.  Equations (\ref{find-A})--(\ref{find-M}) have a unique solution found recursively by solving for $K_\frac{1}{2}$, $M_\frac{1}{2}$, $A_1$, $B_1$, $K_\frac{3}{2}$, $M_\frac{3}{2}$, $A_2$, $B_2, \dots$.  It is now clear that given $a_n$, $b_n$, $k_{n - \frac{1}{2}}$, $m_{n- \frac{1}{2}}$ for $n \in \mathbb{Z}_+$, there is a unique solution for $A_n$, $B_n$, $K_{n-\frac{1}{2}}$, $M_{n-\frac{1}{2}}$ for $n \in \mathbb{Z}_+$ and vice versa.  Thus using (\ref{change-basis-first})--(\ref{change-basis-last}), there is a unique solution for $A_n^{(1)}$, $A_n^{(2)}$, $M^{(1)}_{n-\frac{1}{2}}$, $M^{(2)}_{n-\frac{1}{2}}$ for $n \in \mathbb{Z}_+$.  Next, we note that the map 
\begin{eqnarray*}
\mbox{$((\bigwedge_*^0)^\times)^2$} & \longrightarrow & \mbox{$((\bigwedge_*^0)^\times)^2/\langle(-1,-1)\rangle$}\\
(a_0, b_0) & \mapsto & (a_0^{(1)}, a_0^{(2)}) 
\end{eqnarray*}
given by (\ref{find-a}) and (\ref{find-b}) is a bijection.   This proves the first statement of the proposition.  The fact that (\ref{L-notation-N1-N2})--(\ref{G-notation-N1-N2}) give a representation of the $N=2$ Neveu-Schwarz Lie superalgebra with central charge zero is easy to check (as was done at the beginning of this section).

Finally we prove the last statement of the proposition by noting that $f(x,\varphi)$ is superanalytic, invertible in a neighborhood of $(\infty, 0)$ and of the form (\ref{where-H-lives2}) if and only if $f (1/x, -i\varphi/x)$ is superanalytic invertible in a neighborhood of zero and of the form (\ref{where-H-lives}).  Thus by the first part of the proposition, $f(1/x, -i\varphi/x)$ can be expressed uniquely as (\ref{change-of-variables-superanalytic}).  Composing $f(1/x, -i\varphi/x)$ with $(1/x, i\varphi/x)$, one obtains (\ref{change-of-variables-superanalytic}).
\end{proof}

\begin{rema}{\em
We note that the case of $N=1$ superanalytic functions written as exponentials of infinitesimals giving a representation of the $N=2$ Neveu-Schwarz algebra is slightly more complicated than the case of $N=2$ superconformal functions written as exponentials of infinitesimals giving a representation of the $N=2$ Neveu-Schwarz algebra.  The reason for the added condition that the first bijection given in Proposition \ref{moduli-prop-N2-one-variable} is with superanalytic functions that vanish not only at $(x,\varphi) = (0,0)$ but that also have even part vanishing at $x = 0$ is that general superanalytic functions vanishing at $(x,\varphi) = (0,0)$ are given by $f_1 \circ f_2(x, \varphi)$ where $f_1(x, \varphi) = ( x + k_{-1/2} \varphi, \varphi)$ for $k_{-1/2} \in \bigwedge_*^1$ and $f_2$ is vanishing at $(x, \varphi) = (0,0)$ but also has even part vanishing at $x = 0$.   Since 
\begin{eqnarray}
f_1(x, \varphi) &=&  (x+ k_{-\frac{1}{2}} \varphi, \varphi) \ = \ e^{k_{-\frac{1}{2}} \varphi \frac{\partial}{\partial x}} \cdot (x, \varphi) \\
&=& e^{\frac{1}{2}k_{-\frac{1}{2}}  \left( G^{(1)}_{-\frac{1}{2}}(x, \varphi) - i G^{(2)}_{-\frac{1}{2}} (x, \varphi) \right)} \cdot (x, \varphi).\nonumber
\end{eqnarray}
And in the spirit of \cite{B-memoir} we would normally think of $G^{(j)}_{n - 1/2}(x, \varphi)$ as an infinitesimal acting at infinity for $j = 1,2$ and $n \leq 0$.  However, by taking out the superconformality condition we need these infinitesimals acting at zero in order to get all superanalytic functions vanishing at zero and invertible in a neighborhood of zero.  Similarly, to get all superanalytic functions vanishing at infinity and invertible in a neighborhood of zero, we need to take $f_1 \circ f_2(x, \varphi)$ with $f_1(x, \varphi) = ( x - i k_{1/2} \varphi x,  \varphi)$ for $k_{1/2} \in \bigwedge_*^1$  and $f_2(x, \varphi) = (\tilde{x}, \tilde{\varphi})$ vanishing at infinity, invertible in a neighborhood of infinity and with $x\frac{\partial}{\partial \varphi}\tilde{x}$ also vanishing as $x \rightarrow \infty$.  Note then that 
\begin{eqnarray}
f_1(x, \varphi) &=&  (x- ik_{\frac{1}{2}} \varphi x , \varphi) \ = \ e^{-ik_{\frac{1}{2}} \varphi x \frac{\partial}{\partial x}} \cdot (x, \varphi) \\
&=& e^{-i\frac{1}{2}k_{\frac{1}{2}}  \left( G^{(1)}_{\frac{1}{2}}(x, \varphi) - i G^{(2)}_{\frac{1}{2}} (x, \varphi) \right)} \cdot (x, \varphi)\nonumber
\end{eqnarray}
and 
\begin{eqnarray}
e^{-i\frac{1}{2}k_{\frac{1}{2}}  \left( G^{(1)}_{\frac{1}{2}}(x, \varphi) - i G^{(2)}_{\frac{1}{2}} (x, \varphi) \right)} \cdot \Bigl(\frac{1}{x}, \frac{i\varphi}{x} \Bigr) = \left( \frac{1}{x} + k_{\frac{1}{2}} \frac{i \varphi}{x}, \frac{i \varphi}{x} \right).
\end{eqnarray}
Thus in symmetry with the case of functions vanishing at zero, the infinitesimals $G^{(j)}_{ 1/2}(x, \varphi)$ for $j = 1,2$ act at infinity as well as zero.}
\end{rema}

Let $\mathcal{SA}_{\mathrm{res}}(1,0)$ denote the set of formal $N=1$ superanalytic functions vanishing at zero and invertible in a neighborhood of zero restricted to only include those that also have even part vanishing at $x=0$.   It is clear that $\mathcal{SA}_{\mathrm{res}}(1,0)$ is a group under composition. Again let $\mathcal{G} = (\bigwedge_*^0)^\times)^2/\langle (-1,-1) \rangle \times (\bigwedge_*^\infty)^2$.
Define the map 
\begin{equation}
\hat{E}_1: \mathcal{G} \longrightarrow \mathcal{SA}_\mathrm{res}(1,0)
\end{equation}
by
\begin{multline}
\hat{E}_1(a_0^{(1)}, a_0^{(2)}, A^{(1)}, A^{(2)}, M^{(1)}, M^{(2)}) \\
= \exp \Biggl(-  \sum_{n \in \Z} \Bigl( A^{(1)}_n L_n(x,\varphi) + A^{(2)}_n J_n(x,\varphi) +  M^{(1)}_{n - \frac{1}{2}} G^{(1)}_{n -\frac{1}{2}} (x, \varphi) \\
+  M^{(2)}_{n - \frac{1}{2}} G^{(2)}_{n -\frac{1}{2}} (x, \varphi) \Bigr) \Biggr)  \cdot  (a_0^{(1)})^{-2L_0(x,\varphi)} \cdot  (a_0^{(2)})^{-J_0(x,\varphi)} \cdot (x, \varphi) .
\end{multline}
It is clear that $\hat{E}_1$ is a bijection.  Let $g, h \in \mathcal{G}$ with $g=(a_0^{(1)}, a_0^{(2)}, A^{(1)}, A^{(2)}, M^{(1)}, M^{(2)})$ and $h=(b_0^{(1)}, b_0^{(2)}$, $B^{(1)}, B^{(2)}, N^{(1)}, N^{(2)})$.  We define the {\it $N=1$ composition at zero} map
\begin{eqnarray}
\circ_0' : \mathcal{G} \times \mathcal{G} &\longrightarrow& \mathcal{G} \\
(g,h) &\mapsto& g \circ_0' h \nonumber
\end{eqnarray}
by 
\begin{eqnarray}\label{define-circ-prime}
\lefteqn{\quad g\circ_0' h} \\
\qquad \quad &=&(a_0^{(1)}, a_0^{(2)}, A^{(1)}, A^{(2)}, M^{(1)}, M^{(2)}) \circ_0 (b_0^{(1)}, b_0^{(2)}, B^{(1)}, B^{(2)}, N^{(1)}, N^{(2)})  \nonumber \\
&=& \hat{E}_1^{-1} \left(  \hat{E}_1 (b_0^{(1)}, b_0^{(2)}, B^{(1)}, B^{(2)}, N^{(1)}, N^{(2)}) \circ \hat{E}_1(a_0^{(1)}, a_0^{(2)}, A^{(1)}, \right. \nonumber \\
& & \hspace{3in} \left.  A^{(2)}, M^{(1)}, M^{(2)}) \right).\nonumber
\end{eqnarray}

For $(A^{(1)}, A^{(2)}, M^{(1)}, M^{(2)}) \in (\bigwedge_*^\infty)^2$, we will use the notation
\begin{multline}
T((A^{(1)}, A^{(2)}, M^{(1)}, M^{(2)}), (x, \varphi)) \\
=-  \sum_{n \in \Z} \Bigl( A^{(1)}_n L_n(x,\varphi) + A^{(2)}_n J_n(x,\varphi) +  M^{(1)}_{n - \frac{1}{2}} G^{(1)}_{n -\frac{1}{2}} (x, \varphi) 
+  M^{(2)}_{n - \frac{1}{2}} G^{(2)}_{n -\frac{1}{2}} (x, \varphi) \Bigr)
\end{multline}
in analogy (but not to be confused) with (\ref{T-notation}).

\begin{prop}\label{composition-switch-prop}
Let $f_1(x, \varphi) = (\tilde{x}, \tilde{\varphi}) \in (\bigwedge_*((x))[\varphi] )^0 \oplus (\bigwedge_*((x))[\varphi] )^1$, and let $f_2 \in \mathcal{SA}_\mathrm{res}(1,0)$ be given by
\begin{equation}
f_2(x, \varphi) = e^{T(A^{(1)}, A^{(2)}, M^{(1)}, M^{(2)}), (x, \varphi))} \cdot  (a_0^{(1)})^{-2L_0(x,\varphi)} \cdot  (a_0^{(2)})^{-J_0(x,\varphi)} \cdot (x, \varphi)
\end{equation}
for $(a_0^{(1)}, a_0^{(2)}, A^{(1)}, A^{(2)}, M^{(1)}, M^{(2)})  \in \mathcal{G}$.  Then
\begin{eqnarray}\label{composition-switch}
\lefteqn{\quad f_1 \circ f_2 (x, \varphi) }\\
\qquad \ &=& f_1\left(e^{T(A^{(1)}, A^{(2)}, M^{(1)}, M^{(2)}), (x, \varphi))} \cdot  (a_0^{(1)})^{-2L_0(x,\varphi)} \cdot  (a_0^{(2)})^{-J_0(x,\varphi)} \cdot(x, \varphi)\right) \nonumber \\
&=& e^{T(A^{(1)}, A^{(2)}, M^{(1)}, M^{(2)}), (x, \varphi))} \cdot  (a_0^{(1)})^{-2L_0(x,\varphi)} \cdot  (a_0^{(2)})^{-J_0(x,\varphi)} \cdot f_1(x, \varphi) .\nonumber
\end{eqnarray}
More generally, if $f_1(x, \varphi) = (\tilde{x}, \tilde{\varphi}) \in (\bigwedge_*[[x, x^{-1}]][\varphi] )^0 \oplus (\bigwedge_*[[x, x^{-1}]][\varphi] )^1$ and the composition $f_1 \circ f_2$ is well-defined, then (\ref{composition-switch}) holds.
\end{prop}

\begin{proof}
Since $T((A^{(1)}, A^{(2)}, M^{(1)}, M^{(2)}), (x, \varphi)) \in (\mathrm{Der}(\bigwedge_*((x))[\varphi] ))^0$, the automorphism property, Proposition 3.1 of \cite{B-memoir} holds in this setting.  Thus, in the proof of Proposition 3.12 in \cite{B-memoir},  we can replace $T_H$ with the more general even derivation $T((A^{(1)}, A^{(2)}, M^{(1)}, M^{(2)}), (x, \varphi))$, and the result still follows.  The final statement follows by linearity. 
\end{proof}

The next proposition follows from Proposition \ref{composition-switch-prop}.

\begin{prop}\label{N=1-group-zero-prop}
The set $\mathcal{G} = (\bigwedge_*^0)^\times)^2/\langle (-1,-1) \rangle \times (\bigwedge_*^\infty)^2$ is a group under the composition operation $\circ_0'$ defined by (\ref{define-circ-prime}).   In fact, $(\mathcal{G}, \circ_0')$ is isomorphic to the opposite group of  the group of formal $N=1$ superanalytic functions vanishing at zero, invertible in a neighborhood of zero, and with even part vanishing at $x=0$ under composition, i.e., 
\begin{equation}
(\mathcal{G}, \circ_0') \cong (\mathcal{SA}_\mathrm{res}(1,0), \circ)^{\mathrm{op}}.
\end{equation}
Furthermore for $g, h \in \mathcal{G}$ with $g = (a_0^{(1)}, a_0^{(2)}, A^{(1)}, A^{(2)}, M^{(1)}, M^{(2)})$ and 
$h = (b_0^{(1)}, b_0^{(2)}, B^{(1)}, B^{(2)}, N^{(1)}, N^{(2)})$, we have that if
\begin{equation}
g \circ_0' h = (c_0^{(1)}, c_0^{(2)}, C^{(1)}, C^{(2)},O^{(1)}, O^{(2)}) 
\end{equation}
then $(c_0^{(1)}, c_0^{(2)},C^{(1)}, C^{(2)},O^{(1)}, O^{(2)}) \in  \mathcal{G}$ is given by
\begin{eqnarray}
\lefteqn{e^{T((C^{(1)}, C^{(2)}, O^{(1)}, O^{(2)}), (x, \varphi))} \cdot (c_0^{(1)})^{-2L_0(x,\varphi)} \cdot  (c_0^{(2)})^{-J_0(x,\varphi)} }\\
&=& e^{T((A^{(1)}, A^{(2)}, M^{(1)}, M^{(2)}), (x, \varphi))} \cdot (a_0^{(1)})^{-2L_0(x,\varphi)} \cdot  (a_0^{(2)})^{-J_0(x,\varphi)} \nonumber \\
 & & \quad \cdot e^{T((B^{(1)}, B^{(2)}, N^{(1)}, N^{(2)}), (x, \varphi))} \cdot (b_0^{(1)})^{-2L_0(x,\varphi)} \cdot  (b_0^{(2)})^{-J_0(x,\varphi)}. \nonumber
\end{eqnarray}
\end{prop}

Similarly, letting $\mathcal{SA}_\mathrm{res}(1,\infty)$ denote the set of formal $N=1$ superanalytic functions vanishing at infinity and invertible in a neighborhood of infinity and such that $x (\partial/\partial \varphi)$ of the even component vanishes as the even variable $x$ goes to infinity.   Let $I_1(x,\varphi) = (1/x, i\varphi/x)$.   Then $\mathcal{SA}_\mathrm{res}(1,\infty)$ is a group under the operation 
\begin{equation}\label{infinity-composition-prime}
f_1 \cdot' f_2 = f_1 \circ I_1^{-1} \circ f_2 .
\end{equation}  
Again let $g, h \in \mathcal{G}$ with $g=(a_0^{(1)}, a_0^{(2)}, A^{(1)}$, $A^{(2)}, M^{(1)}, M^{(2)})$ and $h=(b_0^{(1)}, b_0^{(2)}$, $B^{(1)}, B^{(2)}, N^{(1)}, N^{(2)})$.  We define the {\it $N=1$ composition at infinity} map
\begin{eqnarray}
\circ_\infty': \mathcal{G} \times \mathcal{G} &\longrightarrow& \mathcal{G} \\
(g,h) &\mapsto& g \circ_\infty' h \nonumber
\end{eqnarray}
as follows:  let
\begin{multline}
(c_0^{(1)}, c_0^{(2)},C^{(1)}, -C^{(2)},-iO^{(1)}, -iO^{(2)}) \\
= \hat{E}_1^{-1} \left(  \hat{E}_1 (b_0^{(1)}, b_0^{(2)}, B^{(1)}, B^{(2)}, N^{(1)}, N^{(2)}) \circ I_1^{-1} \circ \hat{E}_1(a_0^{(1)}, a_0^{(2)}, A^{(1)}, \right. \\
 \left.  A^{(2)}, M^{(1)}, M^{(2)}) \circ I_1^{-1}\right)
\end{multline} 
then define
\begin{eqnarray}\label{define-circ-infty-prime}
\lefteqn{\quad g\circ_\infty' h} \\
\qquad &=&(a_0^{(1)}, a_0^{(2)}, A^{(1)}, A^{(2)}, M^{(1)}, M^{(2)}) \circ_\infty' (b_0^{(1)}, b_0^{(2)}, B^{(1)}, B^{(2)}, N^{(1)}, N^{(2)})  \nonumber \\
&=& (c_0^{(1)}, c_0^{(2)}, C^{(1)}, C^{(2)},O^{(1)}, O^{(2)}) . \nonumber
\end{eqnarray}

\begin{prop}\label{composition-switch-prop-infinity}
Let $f_1(x, \varphi) = (\tilde{x}, \tilde{\varphi}) \in (\bigwedge_*((x^{-1}))[\varphi] )^0 \oplus (\bigwedge_*((x))[\varphi] )^1$, and let $f_2 \in \mathcal{SA}_\mathrm{res}(1,\infty)$ be given by
\begin{equation}
f_2(x, \varphi) = e^{T(A^{(1)}, -A^{(2)}, -iM^{(1)}, -iM^{(2)}), I_1(x, \varphi))} \cdot  (a_0^{(1)})^{2L_0(x,\varphi)} \cdot  (a_0^{(2)})^{-J_0(x,\varphi)} \cdot (x, \varphi)
\end{equation}
for $(a_0^{(1)}, a_0^{(2)}, A^{(1)}, A^{(2)}, M^{(1)}, M^{(2)})  \in \mathcal{G}$.  Then
\begin{eqnarray}\label{composition-switch-infinity}
\lefteqn{\qquad \ \ f_1 \circ f_2 (x, \varphi) }\\
\ \  &=& \! \! \! \! f_1\left(e^{T(A^{(1)}, -A^{(2)}, -iM^{(1)}, -iM^{(2)}), I_1(x, \varphi))} \cdot  (a_0^{(1)})^{2L_0(x,\varphi)} \cdot  (a_0^{(2)})^{-J_0(x,\varphi)} \cdot(x, \varphi)\right) \nonumber \\
&=& \! \! \! \! e^{T(A^{(1)}, -A^{(2)}, -iM^{(1)}, -iM^{(2)}), I_1 (x, \varphi))} \cdot  (a_0^{(1)})^{2L_0(x,\varphi)} \cdot  (a_0^{(2)})^{-J_0(x,\varphi)} \cdot f_1(x, \varphi) .\nonumber
\end{eqnarray}
More generally, if $f_1(x, \varphi) = (\tilde{x}, \tilde{\varphi}) \in (\bigwedge_*[[x, x^{-1}]][\varphi] )^0 \oplus (\bigwedge_*[[x, x^{-1}]][\varphi] )^1$ and the composition $f_1 \circ f_2$ is well-defined, then (\ref{composition-switch-infinity}) holds.
\end{prop}

\begin{proof}
Since $T((A^{(1)}, -A^{(2)}, -iM^{(1)}, -iM^{(2)}), I_1(x, \varphi)) \in (\mathrm{Der}(\bigwedge_*((x^{-1}))[\varphi] ))^0$, the automorphism property, Proposition 3.1 of \cite{B-memoir} holds in this setting.  Thus, in the proof of Proposition 3.20 in \cite{B-memoir},  we can replace $\bar{T}$ with the more general $T((A^{(1)}, -A^{(2)}, -iM^{(1)}, -iM^{(2)}), I_1 (x, \varphi))$, and the result still follows.  The final, statement follows by linearity. 
\end{proof}

The next proposition follows from Proposition \ref{composition-switch-prop-infinity}.

\begin{prop}\label{N=1-group-infinity-prop}
The set $\mathcal{G} = (\bigwedge_*^0)^\times)^2/\langle (-1,-1) \rangle \times (\bigwedge_*^\infty)^2$ is a group under the composition operation $\circ_\infty'$ defined by (\ref{define-circ-infty-prime}).  In fact, $(\mathcal{G}, \circ_\infty')$ is isomorphic to the opposite group of  the group of formal $N=1$ superanalytic functions vanishing at infinity and invertible in a neighborhood of infinity and such that $x (\partial/\partial \varphi)$ of the even component vanishes as the even variable $x$ goes to infinity under the composition (\ref{infinity-composition-prime}), i.e.
\begin{equation}
(\mathcal{G}, \circ_\infty') \cong (\mathcal{SA}_\mathrm{res} (1,\infty), \cdot')^{\mathrm{op}}.
\end{equation}
Furthermore for $g, h \in \mathcal{G}$ with $g = (a_0^{(1)}, a_0^{(2)}, A^{(1)}, A^{(2)}, M^{(1)}, M^{(2)})$ and $h = (b_0^{(1)}, b_0^{(2)}, B^{(1)}, B^{(2)}, N^{(1)}, N^{(2)})$, we have that if
\begin{equation}
g \circ_\infty' h = (c_0^{(1)}, c_0^{(2)}, C^{(1)}, C^{(2)},O^{(1)}, O^{(2)}) 
\end{equation}
then $(c_0^{(1)}, c_0^{(2)},C^{(1)}, C^{(2)},O^{(1)}, O^{(2)}) \in  \mathcal{G}$ is given by
\begin{eqnarray}
\lefteqn{e^{T((C^{(1)}, -C^{(2)}, -iO^{(1)},-i O^{(2)}), I_1(x, \varphi))}
 \cdot (c_0^{(1)})^{2L_0(x,\varphi)} \cdot  (c_0^{(2)})^{-J_0(x,\varphi)} }\\
\qquad \ \ &=& e^{T((A^{(1)}, -A^{(2)}, -iM^{(1)}, -iM^{(2)}), I_1(x, \varphi))} \cdot (a_0^{(1)})^{2L_0(x,\varphi)}  \cdot  (a_0^{(2)})^{-J_0(x,\varphi)} \nonumber\\
& & \quad \cdot e^{T((B^{(1)}, -B^{(2)}, -iN^{(1)}, -iN^{(2)}), I_1(x, \varphi))} \cdot (b_0^{(1)})^{2L_0(x,\varphi)} \cdot  (b_0^{(2)})^{-J_0(x,\varphi)}. \nonumber
\end{eqnarray}
\end{prop}

Comparing Propositions \ref{N=2-group-zero-prop} and \ref{N=1-group-zero-prop}, and comparing Propositions \ref{N=2-group-infinity-prop} and \ref{N=1-group-infinity-prop}, we have the following theorem:

\begin{thm} The groups $(\mathcal{G}, \circ_0)$ and $(\mathcal{G}, \circ_0')$ are isomorphic and thus the groups $(\mathcal{SC}(2,0), \circ)$ and $(\mathcal{SA}_\mathrm{res}(1,0), \circ)$ are isomorphic.  Similarly, the groups $(\mathcal{G}, \circ_\infty)$ and $(\mathcal{G}, \circ_\infty')$ are isomorphic and thus the groups $(\mathcal{SC}(2,\infty), \cdot)$ and $(\mathcal{SA}_\mathrm{res} (1,\infty), \cdot')$ are isomorphic.   

That is to say, the group of $N=2$ superconformal functions vanishing at zero and invertible in a neighborhood of zero is isomorphic to the group of $N=1$ superanalytic functions vanishing at zero and invertible in a neighborhood of zero and with even component vanishing when the even variable is set equal to zero.

Similarly, the group of $N=2$ superconformal functions vanishing at infinity and invertible in a neighborhood of infinity is isomorphic to the group of $N=1$ superanalytic functions vanishing at infinity and invertible in a neighborhood of infinity and with  $x(\partial/\partial \varphi)$ of the even component vanishing when the even variable approaches infinity.
\end{thm}

\begin{rema}\label{compare-to-DRS-remark}  {\em It is important to note that we do {\it not} use the transformation used in for instance \cite{DRS} to achieve a reduction {}from $N=2$ superconformal functions to $N=1$ superanalytic functions.  The transformation used in \cite{DRS} and \cite{He} is the nonsuperconformal change of variables
\begin{equation}
\tilde{x} = x - \frac{i}{2}\varphi^{(1)} \varphi^{(2)} = x + \frac{1}{2}\varphi^+ \varphi^-, \quad \mathrm{and} \quad 
\tilde{\varphi}^{(j)} = \varphi^{(j)}, \quad \mbox{for $j=1,2$},
\end{equation}
which results in the following infinitesimals
\begin{eqnarray}
L_n(\tilde{x}, \tilde{\varphi}^{(1)}, \tilde{\varphi}^{(2)}) &=& L_n(x, \varphi^{(1)}, \varphi^{(2)}) \\
J_n(\tilde{x}, \tilde{\varphi}^{(1)}, \tilde{\varphi}^{(2)}) &=& J_n(x, \varphi^{(1)}, \varphi^{(2)}) \\
G^{(1)}_{n-\frac{1}{2}} (\tilde{x}, \tilde{\varphi}^{(1)}, \tilde{\varphi}^{(2)}) &=& G^{(1)}_{n-\frac{1}{2}} (x, \varphi^{(1)}, \varphi^{(2)}) - \frac{i}{2} \biggl( x^n \varphi^{(2)} \frac{\partial}{\partial x} \\
& & \quad - n x^{n-1} \varphi^{(1)} \varphi^{(2)} \frac{\partial}{\partial \varphi^{(1)}} \biggr) \nonumber\\
G^{(2)}_{n-\frac{1}{2}} (\tilde{x}, \tilde{\varphi}^{(1)}, \tilde{\varphi}^{(2)}) &=& G^{(2)}_{n-\frac{1}{2}} (x, \varphi^{(1)}, \varphi^{(2)}) + \frac{i}{2} \biggl( x^n \varphi^{(1)} \frac{\partial}{\partial x} \\
& & \quad + n x^{n-1} \varphi^{(1)} \varphi^{(2)} \frac{\partial}{\partial \varphi^{(2)}} \biggr), \nonumber
\end{eqnarray} 
using the notation (\ref{L-notation-nonhomo})--(\ref{G-notation-nonhomo}).  Instead, we achieve a correspondence between the $N=2$ superconformal and $N=1$ superanalytic settings by continuously deforming the odd variable in the $N=1$ superconformal setting, i.e using the transformation $\varphi \mapsto s\varphi$ for $s \in (\bigwedge^0)^\times$.   Thus we give another way of looking at the supergeometric phenomenon studied in \cite{DRS}.  That is Proposition \ref{moduli-prop-N2-one-variable} shows that the $N=2$ superconformal geometry is equivalent to $N=1$ superanalytic geometry, but we use a different means to prove this through the notion of continuously deformed $N=1$ superconformality, namely $D_s$-superconformality. }
\end{rema}

\subsection{$N=2$ (Neveu-Schwarz) vertex (operator) superalgebras with one odd formal variable}
  
{}From the previous section, we see that the infinitesimal $N=1$ superanalytic transformations give a representation of the $N=2$ Neveu-Schwarz Lie algebra just as the infinitesimal $N=2$ superconformal transformations do.  This motivates the following notions of $N=2$ vertex superalgebra with one odd formal variable and $N=2$ Neveu-Schwarz vertex operator superalgebra with one odd formal variable.

\begin{defn}\label{vertex-superalgebra-definition-one-variable}
{\em An} $N=2$ vertex superalgebra over $\bigwedge_*$ and with one odd formal variable {\em consists of an $N=1$ vertex superalgebra $(V, Y(\cdot, \varphi), \mathbf{1})$ such that in addition to $\mathcal{D}_1 \in (\mathrm{End} \; V)^1$ given by Proposition \ref{N1-derivative-prop} (which we will denote by $\mathcal{D}^{(1)}$, that is $\mathcal{D}^{(1)} v = v_{-3/2} \mathbf{1}$), there exists $\mathcal{D}^{(2)} \in (\mathrm{End} \; V)^1$ such that 
\begin{equation}\label{first-D-axiom-one-var}
\bigl[\mathcal{D}^{(j)}, \mathcal{D}^{(k)} \bigr] = 2\delta_{j,k} \mathcal{D},
\end{equation}
and such that in addition to the $\mathcal{D}^{(1)}$-bracket relation
\begin{equation}
\bigl[ \mathcal{D}^{(1)}, Y(v, (x,\varphi)) \bigr] = Y(\mathcal{D}^{(1)}v, (x, \varphi)) - 2 \varphi Y(\mathcal{D}v, (x, \varphi))
\end{equation}
arising naturally {}from the $N=1$ vertex superalgebra structure (see Proposition \ref{bracket-prop-N1}), we have that the} $\mathcal{D}^{(2)}$-bracket relation {\em holds:
\begin{equation}\label{second-D-axiom-one-var}
\bigl[ \mathcal{D}^{(2)}, Y(v, (x, \varphi)) \bigr] = Y(\mathcal{D}^{(2)}v,(x, \varphi)).
\end{equation}
 }
\end{defn}

The $N=2$ vertex superalgebra with one odd formal variable just defined is denoted by \[(V,Y(\cdot,(x,\varphi)),\mathbf{1},\mathcal{D}^{(1)}, \mathcal{D}^{(2)}).\]

\begin{rema}\label{N2-subalgebras-representation-remark-one-var}
{\em An $N=2$ vertex superalgebra with one odd formal variables is a representation of $\mathfrak{osp}_{\bigwedge_*}(2|2)_{<0}$.  This representation is given by $\mathcal{D}^{(j)} \mapsto G^{(j)} _{-1/2}$, for $j=1,2$, and $\mathcal{D} \mapsto L_{-1}$; see Remark \ref{N2-subalgebras-remark-nonhomo}.
}
\end{rema}

As a consequence of the definition, if $V$ is an $N=2$ vertex superalgebra with one odd formal variable, {}from (\ref{constructing-s-operators}) and the fact that $V$ is an $N=2$ vertex superalgebra, letting $Y(v, x) = Y(v, (x, 0))$, we have that 
\begin{equation}\label{constructing-one-var-operators}
Y(v, (x, \varphi))  = Y(v,x) + \varphi Y(\mathcal{D}^{(1)}v, x).
\end{equation}

\begin{defn}
{\em An} $N = 2$ Neveu-Schwarz vertex operator superalgebra over $\bigwedge_*$ and with one odd formal variable {\em is a $\frac{1}{2} \mathbb{Z}$-graded $\bigwedge_*$-module (graded by} weights{\em)  
\begin{equation}\label{graded-one-var}
V = \coprod_{n \in \frac{1}{2}\mathbb{Z}} V_{(n)} 
\end{equation}
such that 
\begin{equation}
\dim V_{(n)} < \infty \qquad \mbox{for $n \in \frac{1}{2} \mathbb{Z}$,} 
\end{equation}
\begin{equation}\label{third-n2-with-one-var}
V_{(n)} = 0 \qquad \mbox{for $n$ sufficiently negative} , 
\end{equation}
equipped with an $N=2$ vertex superalgebra structure $(V, Y(\cdot, (x,\varphi)),\mathbf{1}, \mathcal{D}^{(1)}, \mathcal{D}^{(2)})$, and two distinguished vectors $\tau^j \in V_{(3/2)}^1$ (the {\em  $\tau^j $ $N=2$ Neveu-Schwarz elements} or {\em  $\tau^j$ $N=2$ superconformal elements}), satisfying the following conditions:  the $N=2$ Neveu-Schwarz algebra relations in the nonhomogeneous basis (\ref{Virasoro-relation2})--(\ref{transformed-Neveu-Schwarz-relation-last}) hold for $L(n), J(n) \in  (\mathrm{End} \; V)^0$ and $G^{(j)}(n + 1/2) \in (\mathrm{End} \; V)^1$, for $n \in \mathbb{Z}$, $j=1,2$, and $c_V \in \mathbb{C}$ (the} central charge{\em), where
\begin{equation}\label{first-n2-one-var}
G^{(j)}(n - 1/2) = \tau^j_n,  \qquad \mbox{for $j=1,2$} 
\end{equation}
\begin{equation}
2L(n) = \tau^1_{n+ \frac{1}{2}}, \quad \mbox{and} \quad  i(n+1)J(n) = \tau^2_{n+ \frac{1}{2}},
\end{equation}
for $n \in \mathbb{Z}$, i.e., 
\begin{eqnarray}
Y(\tau^1,(x,\varphi)) &=& \sum_{n \in \mathbb{Z}} \left( G^{(1)}(n+1/2) + 2\varphi L(n) \right) x^{-n-2}\\
Y(\tau^2,(x,\varphi)) &=& \sum_{n \in \mathbb{Z}} \left( G^{(2)}(n+1/2) + i(n+1)\varphi J(n) \right) x^{-n-2},
\end{eqnarray}
and
\begin{equation}\label{D-G-compatibility-one-var}
G^{(j)}(-1/2) = \mathcal{D}^{(j)}, \qquad \mbox{for $j = 1,2$;}
\end{equation}
for $n \in \frac{1}{2} \mathbb{Z}$ and $v \in V_{(n)}$
\begin{equation}\label{L0-grading-one-var}
L(0)v = nv
\end{equation}
and in addition, $V_{(n)}$ is the direct sum of eigenspaces for $J(0)$ such that if $v \in V_{(n)}$ is also an eigenvector for $J(0)$ with eigenvalue $k$, i.e., if
\begin{equation}\label{J0-grading-one-var}
J(0)v = kv, \quad \mbox{then $k \equiv 2n \; \mathrm{mod} \; 2$}.
\end{equation} 
}
\end{defn}

\medskip

An $N=2$ Neveu-Schwarz vertex operator superalgebra with one odd formal variable is denoted by
\[(V,Y(\cdot,(x,\varphi)),\mathbf{1},\tau^1, \tau^2) . \]

Note that in particular, (\ref{D-G-compatibility-one-var}) implies that the $G(-1/2)$-derivative property holds:
\begin{equation}
\left( \frac{\partial}{\partial \varphi} + \varphi \frac{\partial}{\partial x} \right) Y(v, (x, \varphi)) = Y(G^{(1)}(-1/2)v, (x, \varphi)),
\end{equation}
and thus the $L(-1)$-derivative property also holds:
\begin{equation}
\frac{\partial}{\partial x}  Y(v, (x, \varphi)) = Y(L(-1)v, (x, \varphi)).
\end{equation}

\subsection{Isomorphism between the categories of $N=2$ (Neveu-Schwarz) vertex (operator) superalgebras with one odd formal variable and without odd formal variables}\label{iso-section-one-var}

The following proposition is immediate:

\begin{prop}\label{get-a-vosa-without-one-var}
Let $(V,Y(\cdot,(x,\varphi)), \mathbf{1}, \mathcal{D}^{(1)}, \mathcal{D}^{(2)})$ be an $N=2$ vertex superalgebra with one odd formal variable.  Then $(V,Y(\cdot,(x,0)), \mathbf{1}, \mathcal{D}^{(1)}, \mathcal{D}^{(2)})$ is an $N=2$ vertex superalgebra without odd formal variables in the nonhomogeneous coordinate system.  That is  setting $\mathcal{D}^\pm = \frac{1}{\sqrt{2}} (\mathcal{D}^{(1)} \mp i \mathcal{D}^{(2)})$, then $(V,Y(\cdot,(x,0)), \mathbf{1}, \mathcal{D}^{+}, \mathcal{D}^{-})$ is an $N=2$ vertex superalgebra without odd formal variables in the homogeneous coordinate system.

Let $(V,Y(\cdot,(x,\varphi)), \mbox{\bf 1},\tau^1, \tau^2)$ be an $N=2$ Neveu-Schwarz vertex
operator superalgebra with one odd formal variable.  Then $(V,Y(\cdot,(x,0)), \mathbf{1}, \tau^1, \tau^2)$ is an $N=2$ Neveu-Schwarz vertex operator superalgebra without odd formal variables in the nonhomogeneous coordinate system.  That is setting $\tau^{(\pm)} = \frac{1}{\sqrt{2}} (\tau^1 \mp i \tau^2)$, then $(V,Y(\cdot,(x,0)), \mathbf{1}, \tau^{(+)}, \tau^{(-)})$ is an $N=2$ Neveu-Schwarz vertex operator superalgebra without odd formal variables in the homogeneous coordinate system.
\end{prop}

On the other hand we have:

\begin{prop}\label{get-a-vosa-with-one-var}
Let $(V,Y(\cdot,(x,\varphi^{(1)}, \varphi^{(2)})), \mathbf{1})$ be an $N=2$ vertex superalgebra with two odd formal variables in the nonhomogeneous coordinate system.  Then letting $\mathcal{D}^{(j)}$ be defined by (\ref{nonhomo-define-Ds}), we have that $(V,Y(\cdot,(x,\varphi^{(1)},0)), \mathbf{1}, \mathcal{D}^{(1)}, \mathcal{D}^{(2)})$ is an $N=2$ vertex superalgebra with one odd formal variable.  

Let $(V,Y(\cdot,(x,\varphi^{(1)}, \varphi^{(2)})), \mbox{\bf 1}, \mu)$ be an $N=2$ Neveu-Schwarz vertex
operator superalgebra with two odd formal variables in the nonhomogeneous coordinate system.  Then letting $\tau^{1} = - iG^{(2)}(-1/2) \mu$ and  $\tau^2 = -i G^{(1)}(-1/2)\mu$, we have that $(V,Y(\cdot,(x,\varphi^{(1)},0)), \mathbf{1}, \tau^1, \tau^2)$ is an $N=2$ Neveu-Schwarz vertex operator superalgebra with one odd formal variable. 
\end{prop}

\begin{proof}
Setting the odd variable $\varphi^{(2)}$ equal to zero in the axioms for an $N=2$ vertex superalgebra $(V,Y(\cdot,(x,\varphi^{(1)}, \varphi^{(2)})),  \mathbf{1})$, it is clear that $(V,Y(\cdot,(x,\varphi^{(1)},0)), \mathbf{1}, \mathcal{D}^{(1)}$, $\mathcal{D}^{(2)})$ satisfies all of the axioms for an $N=2$ vertex superalgebra with one odd formal variable except for the axioms involving $\mathcal{D}^{(j)}$, for $j = 1,2$.  The remaining axioms (\ref{first-D-axiom-one-var}) and (\ref{second-D-axiom-one-var}) follow {}from consequences (\ref{Djk-bracket-get-D-nonhomo}) and (\ref{D*-bracket}) of the definition of $N=2$ nonhomogeneous vertex superalgebra with two odd formal variables.

Similarly, it is trivial that $(V,Y(\cdot,(x,\varphi^{(1)},0)), \mathbf{1}, \tau^1, \tau^2)$ satisfies axioms (\ref{graded-one-var})--(\ref{third-n2-with-one-var}), (\ref{L0-grading-one-var}) and (\ref{J0-grading-one-var}).  By the proof above, $V$ is an $N=2$ vertex superalgebra with one odd formal variable.   Axioms (\ref{first-n2-one-var})--(\ref{D-G-compatibility-one-var})  follow {}from consequences of the definition of an $N=2$ Neveu-Schwarz vertex operator superalgebra with two odd formal variables which show that the vertex operators for $\tau^1$ and $\tau^2$ are given by (\ref{tau1}) and (\ref{tau2}), respectively, and that by (\ref{nonhomo-D-correspondence}) $G^{(j)}(-1/2) = \mathcal{D}^{(j)}$. 
\end{proof}

Let $\mathrm{VSA}_2(\varphi^{(1)}, \varphi^{(2)})$ denote the category of $N=2$ vertex  superalgebras with two odd formal variables in the nonhomogeneous coordinate system, and let $\mathrm{VSA}_2(\varphi)$ denote the category of $N=2$ vertex  superalgebras with one odd formal variable.    Let $\mathrm{VOSA}_2(\varphi^{(1)}, \varphi^{(2)},c)$ denote the category of $N=2$ Neveu-Schwarz vertex operator superalgebras with two odd formal variables in the nonhomogeneous coordinate system and with central charge $c \in \mathbb{C}$,  and let $\mathrm{VOSA}_2(\varphi,c)$ denote the category of $N=2$ Neveu-Schwarz vertex operator superalgebras with one odd formal variable and with central charge $c \in \mathbb{C}$.   Also recall the notation {}from the end of Section \ref{iso-section}.

We have the following theorem:

\begin{thm}\label{superalgebras-one-var} 
The categories $\mathrm{VSA}_2(\varphi^+, \varphi^-)$, $\mathrm{VSA}_2(\varphi^{(1)}, \varphi^{(2)})$, $\mathrm{VSA}_2(\varphi)$, and $\mathrm{VSA}_2$ are all isomorphic to each other.   In addition, for any $c \in \mathbb{C}$, the subcategories $\mathrm{VOSA}_2(\varphi^+, \varphi^-,c)$, $\mathrm{VOSA}_2(\varphi^{(1)}, \varphi^{(2)}, c)$, $\mathrm{VOSA}_2(\varphi,c)$ and $\mathrm{VOSA}_2(c)$ are all isomorphic to each other. 
\end{thm}

\begin{proof}
Define $F_{\varphi =0}: \mathrm{VSA}_2(\varphi) \longrightarrow \mathrm{VSA}_2$ by 
\begin{eqnarray*}
F_{\varphi =0}: (V,Y(\cdot,(x,\varphi)), \mathbf{1}, \mathcal{D}^{(1)}, \mathcal{D}^{(2)}) &\mapsto& (V,Y(\cdot,(x,0)),\mathbf{1}, \mathcal{D}^+, \mathcal{D}^-)\\
\gamma &\mapsto& \gamma 
\end{eqnarray*}
for $(V,Y(\cdot,(x,\varphi)), \mathbf{1}, \mathcal{D}^{(1)}, \mathcal{D}^{(2)})$ an object in $\mathrm{VSA}_2(\varphi)$ and $\gamma$ a morphism, where $\mathcal{D}^\pm = \frac{1}{\sqrt{2}} (\mathcal{D}^{(1)} \mp i \mathcal{D}^{(2)})$.   Proposition \ref{get-a-vosa-without-one-var} shows that $F_{\varphi=0}$ takes objects in $\mathrm{VSA}_2 (\varphi)$ to objects in $\mathrm{VSA}_2$.  It is clear that $F_{\varphi =0}$ takes morphisms to morphisms and that $F_{\varphi =0}$ is a functor.   Moreover, Proposition \ref{get-a-vosa-without-one-var} implies that $F_{\varphi =0}$ restricted to $\mathrm{VOSA}_2(\varphi,c)$ gives a functor {}from $\mathrm{VOSA}_2(\varphi,c)$ to $\mathrm{VOSA}_2(c)$.

Define $F_{\varphi^{(2)} = 0} : \mathrm{VSA}_2(\varphi^{(1)}, \varphi^{(2)}) \longrightarrow \mathrm{VSA}_2(\varphi)$ by 
\begin{eqnarray*}
F_{\varphi^{(2)} = 0} : (V,Y(\cdot,(x, \varphi^{(1)}, \varphi^{(2)})), \mathbf{1}) &\mapsto& (V,Y (\cdot,(x,\varphi^{(1)}, 0)), \mathbf{1}, \mathcal{D}^{(1)}, \mathcal {D}^{(2)}), \\
\gamma &\mapsto& \gamma
\end{eqnarray*} 
where $\mathcal{D}^{(j)}$, for $j = 1,2$ are defined by (\ref{nonhomo-define-Ds}).  Proposition \ref{get-a-vosa-with-one-var} shows that $F_{\varphi^{(2)}=0}$ takes  objects in $\mathrm{VSA}_2(\varphi^{(1)}, \varphi^{(2)})$ to objects in $\mathrm{VSA}_2 (\varphi)$.   It is easy to show that $F_{\varphi^{(2)} = 0}$ takes morphisms to morphisms and is a functor.  Moreover, Proposition \ref{get-a-vosa-with-one-var} implies that $F_{\varphi^{(2)} = 0}$ restricted to $\mathrm{VOSA}_2(\varphi^{(1)}, \varphi^{(2)}, c)$ gives a functor {}from $\mathrm{VOSA}_2(\varphi^{(1)}, \varphi^{(2)}, c)$ to $\mathrm{VOSA}_2(\varphi,c)$ via
\[F_{\varphi^{(2)} = 0} : (V,Y(\cdot,(x,\varphi^{(1)}, \varphi^{(2)})), \mathbf{1}, \mu) \mapsto (V,Y(\cdot,
(x,\varphi)), \mathbf{1},\tau^1, \tau^2), \]
where $\tau^1 = -iG^{(2)} (-1/2) \mu$ and $\tau^2 = -iG^{(1)} (-1/2) \mu$.

Let $F_{\pm \mapsto (1,2)} : \mathrm{VSA}_2(\varphi^+, \varphi^-) \longrightarrow \mathrm{VSA}_2(\varphi^{(1)}, \varphi^{(2)})$ be the functor given by the change of coordinates {}from homogeneous to nonhomogeneous.   

Finally, let $F_\varphi: \mathrm{VSA}_2 \rightarrow \mathrm{VSA}_2 (\varphi)$ be given by
$F_\varphi = F_{\varphi^{(2)} = 0} \circ F_{\pm \mapsto (1,2)} \circ F_{\varphi^+, \varphi^-}$.  It is easy to show that $F_{\varphi =0} \circ F_\varphi = 1_{\mathrm{VSA}_2}$.  By (\ref{Y-from-VOSA-without}), (\ref{constructing-nonhomo-operators}) and (\ref{constructing-one-var-operators}), we have that $F_\varphi \circ F_{\varphi =0} = 1_{\mathrm{VSA}_2(\varphi)}$.  It is then obvious that the corresponding restrictions of $F_{\varphi =0}$ and $F_\varphi$ to the subcategories of $N=2$ Neveu-Schwarz vertex operator superalgebras with one odd formal variable and without odd formal variables, respectively, satisfy $F_{\varphi =0} \circ F_\varphi  = 1_{\mathrm{VOSA}_2(c)}$ and $F_\varphi \circ F_{\varphi =0} = 1_{\mathrm{VOSA}_2(\varphi,c)}$.   By Theorem \ref{superalgebras} the result follows.
\end{proof}

\end{document}